\documentclass[a4paper,11pt, leqno]{article}
\usepackage[latin1]{inputenc}
\usepackage{articlemacro}
\usepackage{comment}

\excludecomment{shortversion}
\includecomment{longversion}
\newcommand{\Warning}{\noindent\emph{Warning}:\ }

\newcommand{\affoid}[1]{(#1,#1^+)}

\newcommand{\an}{{}_{\textrm{a}}}
\newcommand{\na}{{}_{\textrm{na}}}
\newcommand{\Ltr}{{(\textup{Ltr})}}

\newcommand{\ang}[1]{{\langle #1 \rangle}}

\newcommand{\czero}{\cup \{0\}}

\newcommand{\cl}[1]{{\overline{\{#1\}}}}

\newcommand{\Cvee}{{\check{C}}}
\newcommand{\Hvee}{{\check{H}}}

\begin{document}

\title{Adic Spaces}

\author{Torsten Wedhorn}

\maketitle

%------------------------------------------------------------------

\textbf{This script is highly preliminary and unfinished. It is online only to give the audience of our lecture easy access to it. Therefore \textsc{usage is at your own risk}.}

\bigskip

\noindent{\scshape Abstract.\ }
Adic spaces are objects in the realm of non-archimedean analysis and have been developed by Roland Huber. The goal of these lecture notes is to give an introduction to adic spaces.

\bigskip

\noindent{\scshape Acknowledgement.\ } I do not claim any originality here. The main sources are R.~Huber's works (\cite{Hu_Habil}, \cite{Hu_Cont}, \cite{Hu_Gen},  \cite{Hu_Etale}, \cite{HuKne}) and M.~Strauch's introduction to adic spaces \cite{Str_Adic}. Occasionally I follow these sources almost verbatim. The main example in Section~\ref{ExampleAdicSpec} follows Section~2 of P.~Scholze's work on perfectoid spaces \cite{Scholze_Perfect}.

I am very grateful to R.~Huber for sending me a copy of~\cite{Hu_Habil} and to M.~Strauch for sending me a copy 
of~\cite{Str_Adic}. I am also very grateful for helpful remarks from K.~Buzzard, J.~Hilgert, E.~Lau, J.~Sch\"utt, and D.~Wortmann. J.~Sch\"utt also provided the proofs in the sections 2.1--2.3. Special thanks go to B.~Conrad who read the mansucript very thoroughly and made many helpful remarks.

%===========================================================
\section*{Introduction}

\subsection*{Notation}
All rings are commutative and with $1$. All complete topological groups/rings are Hausdorff.

%==========================================================

\tableofcontents

%=============================================================

\section{Valuations}

\subsection{Totally ordered groups}
\begin{definition}\label{DefTotOrdered}
A \emph{totally ordered group} (resp.~\emph{totally ordered monoid}) is a commutative group (resp.~commutative monoid) $\Gamma$ (whose composition law is written multiplicatively) together with a total order $\leq$ on $\Gamma$ such that
\[
\gamma \leq \gamma' \implies \gamma\delta \leq \gamma'\delta, \qquad\text{for all $\gamma,\gamma',\delta \in 
\Gamma$}
\]
Let $\Gamma, \Gamma'$ be totally ordered groups (resp.~monoids). A homomorphism of totally ordered groups 
(resp.~totally ordered monoids) is a homomorphism $f\colon \Gamma \to \Gamma'$ of groups (resp.~monoids) such 
that for all $\gamma_1,\gamma_2 \in \Gamma$ one has $\gamma_1 \leq \gamma_2 \implies f(\gamma_1) \leq 
f(\gamma_2)$.
\end{definition}

We obtain the category of totally ordered groups (resp.~totally ordered monoids). It is easy to see that a 
homomorphism is an isomorphism if and only if it is bijective.

If $\Gamma$ is a totally ordered group, $\gamma \in \Gamma$, and $H \subseteq \Gamma$ a subset, we write $\gamma < H$ if $\gamma < \delta$ for all $\delta \in H$. Similarly we define ``$\gamma \leq H$'', ``$\gamma > H$'', ``$\gamma \geq H$''.

We also define $\Gamma_{<\gamma} := \set{\delta \in \Gamma}{\delta < \gamma}$ and have analogous definitions for $\Gamma_{\leq\gamma}$, $\Gamma_{>\gamma}$, and $\Gamma_{\geq\gamma}$.

\begin{remark}
Let $\Gamma, \Gamma'$ be totally ordered groups. A group homomorphism $f\colon \Gamma \to \Gamma'$ is a homomorphism of totally ordered groups if and only if for all $\gamma \in \Gamma$ with $\gamma < 1$ one has $f(\gamma) < 1$.
\end{remark}

\begin{example}
\begin{ali}
\item
$\RR^{>0}$ is a totally ordered group with respect to multiplication and the standard order. $(\RR,+)$ is a totally 
ordered group w.r.t. the standard order. The logarithm $\RR^{>0} \to \RR$ is an isomorphism of totally ordered 
groups.
\item
If $\Gamma$ is a totally ordered group, then every subgroup with the induced order is again totally ordered. It is 
called a \emph{totally ordered subgroup}.

In particular if $M \subseteq \Gamma$ is a subset, then the intersection of all totally ordered subgroups of $\Gamma
$ containing $M$ is again a totally ordered subgroup, called the totally ordered subgroup \emph{generated by $M$}.
\item
Let $I$ be a well ordered set (e.g. $I = \{1,\dots,n\}$ with the standard order) and let $(\Gamma_i)_{i \in I}$ be a 
family of totally ordered groups. We endow $\prod_{i \in I} \Gamma_i$ with the lexicographic order (i.e. $(\gamma_i) 
< (\gamma'_i)$ if and only if $\gamma_j < \gamma'_j$, where $j$ is the smallest element of $I$ such that $
\gamma_j \ne \gamma'_j$). Then $\prod_{i \in I} \Gamma_i$ is a totally ordered group.

The product order on $\prod_{i \in I} \Gamma_i$ (i.e., $(\gamma_i)_i \leq (\gamma'_i)_i$ if and only if $\gamma_i \leq 
\gamma'_i$ for all $i \in I$) is not a total order (except if there exists only one index $i$ such that $\Gamma_i \ne 
\{1\}$).
\begin{longversion}
\item
More generally, let $I$ be a totally ordered index set and let $(\Gamma_i)_{i \in I}$ be a family of totally ordered 
groups. We define the \emph{Hahn product} $\Hbf_i \Gamma_i$ as the subgroup of those $\gamma = (\gamma_i)_i \in 
\prod_i \Gamma_i$ such that $\supp\gamma := \set{i \in I}{\gamma_i \ne 1}$ is well-ordered. Then $\Hbf_i \Gamma_i
$ becaomes a totally ordered group if one defines $\gamma = (\gamma_i)_i < 1$ whenever for the first element $j$ 
in the well-ordering of $\supp \gamma$ one has $\gamma_j < 1$.
\end{longversion}
\item
If $\Gamma$ is a totally ordered group, then $\Gamma_{\leq 1}$ and $\Gamma_{\geq 1}$ are totally ordered 
monoids.
\end{ali}
\end{example}

\begin{prop}
Every totally ordered group is order-isomorphic to a subgroup of a Hahn product of copies of the totally ordered 
group $\RR^{>0}$ over a suitable totally ordered set $I$.
\end{prop}

Ordered monoids will mainly come up in the following construction.

\begin{rem}
Let $\Gamma$ be a totally ordered group. We add an element $0$ to $\Gamma$ and define on the disjount union $\Gamma \cup \{0\}$ the structure of a totally ordered monoid as follows. Restricted to $\Gamma$ it is the given structure. We extend the multiplication by $0\cdot \gamma := \gamma\cdot 0 := 0$ for all $\gamma \in \Gamma \cup \{0\}$. The total order is extended by defining $0 \leq \gamma$ for all $\gamma \in \Gamma \cup \{0\}$. Then $\Gamma \cup \{0\}$ is a totally ordered monoid.

Let $\Gamma'$ be a second totally ordered group. Every isomorphism of totally ordered monoids $f\colon \Gamma \cup \{0\} \iso \Gamma' \cup \{0\}$ sends $0$ to $0$ and its restriction to $\Gamma$ yields an isomorphism $\Gamma \iso \Gamma'$ of totally ordered groups. Conversely, every isomorphism $\Gamma \to \Gamma'$ of totally ordered groups can be extended to an isomorphism $\Gamma \cup \{0\} \to \Gamma' \cup \{0\}$ of totally ordered monoids by sending $0$ to $0$.

In particular $\Gamma$ and $\Gamma'$ are isomorphic if and only if $\Gamma \cup \{0\}$ and $\Gamma' \cup \{0\}$ 
are isomorphic.
\end{rem}

\begin{rem}\label{PropTotallyOrdered}
Let $\Gamma$ be a totally ordered group, $\gamma,\delta \in \Gamma$.
\begin{ali}
\item
$\gamma < 1 \iff \gamma^{-1} > 1$.
\item
\begin{equation}
\begin{aligned}
\gamma,\delta \leq 1 \implies \gamma\delta \leq 1, &\qquad \gamma < 1, \delta \leq 1 \implies \gamma\delta < 1,\\
\gamma,\delta \geq 1 \implies \gamma\delta \geq 1, &\qquad \gamma > 1, \delta \geq 1 \implies \gamma\delta > 1.
\end{aligned}
\end{equation}
\item
$\Gamma$ is torsion free.
\end{ali}
\end{rem}

\begin{proof}
If $\delta \leq 1$, then $\gamma\delta \leq \gamma\cdot 1 = \gamma$. Similarly if we replace $\leq$ by $\geq$. This 
shows~(2). (1) and~(3) are immediate consequences of~(2).
\end{proof}

% \begin{longversion}
% \begin{rem}\label{TopOrderedGroup}
% Let $\Gamma$ be a totally ordered group. We endow $\Gamma \cup \{0\}$ with the following topology. $U \subset 
%\Gamma \cup \{0\}$ is open if any only if $0 \notin U$ or there exists $\gamma \in \Gamma$ such that $\Gamma_{<
%\gamma} := \set{\delta \in \Gamma}{\delta < \gamma} \subseteq U$.
% 
% Thus the induced topology on $\Gamma$ is the discrete topology. A basis for the topology are the sets $\
%{\gamma\}$ and $\Gamma_{<\gamma}$ for $\gamma \in \Gamma$.
% 
% The topology determines the total order: For two element $\gamma, \delta \in \Gamma$ one has $\gamma < 
%\delta$ if and only if there exists an open neighborhood of $0$ containing $\gamma$ but not $\delta$.
% 
% In particular an isomorphism $\Gamma \cup \{0\} \to \Gamma' \cup \{0\}$ of monoids is an isomorphism of totally ordered monoids if and only if it is a homeomorphism.
% \end{rem}
% \end{longversion}

\begin{remdef}
A subgroup $\Delta$ of a totally ordered group $\Gamma$ is called \emph{isolated} or \emph{convex} if the following 
equivalent conditions are satisfied for all $\delta, \delta', \gamma \in \Gamma$.
\begin{eli}
\item
$\delta \leq \gamma \leq 1$ and $\delta \in \Delta$ imply $\gamma \in \Delta$.
\item
$\delta, \gamma \leq 1$ and $\delta\gamma \in \Delta$ imply $\delta,\gamma \in \Delta$.
\item
$\delta \leq \gamma \leq \delta'$ and $\delta,\delta' \in \Delta$ imply $\gamma \in \Delta$.
\end{eli}
\end{remdef}

\begin{proof}
The equivalence of (i) and (iii) is clear.
\proofstep{(i) $\implies$ (ii)}
As $\delta, \gamma \leq 1$ one has $\delta\gamma \leq \gamma \leq 1$ and hence $\gamma \in \Delta$ by~(i). Then 
$\delta \in \Delta$ because $\Delta$ is a group.
\proofstep{(ii) $\implies$ (i)}
Let $\delta \leq \gamma \leq 1$ with $\delta \in \Delta$. Then $\delta\gamma^{-1} \leq 1$ and $(\delta\gamma^{-1})
\gamma = \delta \in \Delta$. This implies $\gamma \in \Delta$ by~(ii).
\end{proof}

\begin{example}
\begin{ali}
\item
If $\Gamma$ is a totally ordered group, then $1$ and $\Gamma$ are convex subgroups. If $\Gamma = \RR_{>0}$ 
(or any non-trivial subgroup thereof), these are the only convex subgroups.
\item
Let $\Gamma_1,\dots,\Gamma_n$ be totally ordered groups and endow $\Gamma = \prod_i \Gamma_i$ with the 
lexicographic order (where $\{1,\dots,n\}$ is well ordered in the standard way). Then for all $r = 1,\dots,n+1$ the 
subgroup $\prod_{i=r}^n \Gamma_i$ is a convex subgroup of $\Gamma$. If some $\Gamma_i$ for $i > 1$ is non-
trivial, then $\Gamma_1$ is not a convex subgroup of $\Gamma$.
\end{ali}
\end{example}

\begin{rem}\label{GeneratedConvex}
Let $\Gamma$ be a totally ordered group, let $H$ be a subgroup. Then the convex subgroup of $\Gamma$ generated by $H$ consists of those $\gamma \in \Gamma$ such that there exist $h,h' \in H$ with $h \leq \gamma \leq h'$ (it is easy to check that this is a subgroup of $\Gamma$).
\end{rem}

\begin{rem}
Let $\Gamma$ be a totally ordered group and let $\Delta$, $\Delta'$ be convex subgroups. Then $\Delta \subseteq 
\Delta'$ or $\Delta' \subseteq \Delta$.
\end{rem}

\begin{proof}
Assume there exist $\delta \in \Delta \setminus \Delta'$ and $\delta' \in \Delta' \setminus \Delta$. After possibly 
replacing these elements by their inverse we may assume that $\delta, \delta' < 1$. And after possibly swapping $
\Delta$ with $\Delta'$ we may assume that $\delta < \delta'$. But then $\delta' \in \Delta$ because $\Delta$ is 
convex.
\end{proof}

\begin{rem}
Let $\Gamma$ be a totally ordered group.
\begin{ali}
\item
If $f\colon \Gamma \to \Gamma'$ is a homomorphism of totally ordered groups, then $\ker(f)$ is a convex subgroup of $\Gamma$.
\item
Let $\Delta \subseteq \Gamma$ be a convex subgroup and let $f\colon \Gamma \to \Gamma/\Delta$ be the canonical homomorphism. Then there exists a unique total order on $\Gamma/\Delta$ such that $f(\Gamma_{\leq 1}) = (\Gamma/\Delta)_{\leq 1}$. Then $f$ is a homomorphism of totally ordered groups.
\end{ali}
\end{rem}

\begin{remark}\label{InverseQuotTotallyOrdered}
Let $\Delta \subseteq \Gamma$ be a convex subgroup and let $f\colon \Gamma \to \Gamma/\Delta$ be the canonical homomorphism. For all $\gamma \in \Gamma$ one then has
\[
f(\Gamma_{\geq \gamma}) = (\Gamma/\Delta)_{\geq \gamma\Delta}
\]
and hence
\[
f^{-1}((\Gamma/\Delta)_{<\gamma\Delta}) = \bigcap_{\delta \in \Delta}\Gamma_{<\gamma\delta}.
\]
\end{remark}

\begin{defi}
Let $\Gamma$ be a totally ordered group. The number of convex subgroups $\ne 1$ of $\Gamma$ is called the 
\emph{height of $\Gamma$}: $\height\Gamma \in \NN_0 \cup \{\infty\}$.
\end{defi}

Clearly $\height\Gamma = 0$ if and only if $\Gamma = 1$. The height of $\RR_{>0}$ or (of any non-trivial subgroup) 
is $1$.

\begin{longversion}
In fact there is the following converse.

\begin{prop}\label{Height1Groups}
For a totally ordered group $\Gamma \ne 1$ the following assertions are equivalent.
\begin{eli}
\item
$\Gamma$ has height $1$.
\item
There exists an injective homomorphism $\Gamma \mono \RR_{>0}$ of totally ordered groups.
\item
$\Gamma$ is archimedean, i.e. for all $\gamma,\delta \in \Gamma_{<1}$ there exists an integer $m > 0$ such that $
\delta^m < \gamma$.
\end{eli}
\end{prop}

\begin{proof}
\cite{Bou_AC}~VI, 4.5, Prop.~8.
\end{proof}
\end{longversion}

\begin{rem}
Let $\Delta$ be a convex subgroup of a totally ordered group $\Gamma$. Then
\[
\height \Gamma = \height\Delta + \height\Gamma/\Delta
\]
In particular if $\Gamma$ is the lexicographically ordered product of totally ordered groups $\Delta$ and $\Delta'$, 
then
\[
\height \Gamma = \height \Delta + \height \Delta'
\]
\end{rem}

\begin{defi}\label{DefCofinal}
Let $\Gamma$ be a totally ordered group and let $H$ be a subgroup of $\Gamma$. We say that $\gamma \in \Gamma \czero$ is \emph{cofinal for $H$}, if for all $h \in H$ there exists $n \in \NN$ such that $\gamma^n < h$.
\end{defi}

Clearly, $0$ is cofinal for $H$ for every subgroup of $\Gamma$, and no $\gamma \geq 1$ is cofinal for any subgroup of $\Gamma$.

\begin{remark}\label{TopTotOrdGroup}
Sometimes it is convenient to consider for a totally ordered group $\Gamma$ the following topology on $\Gamma \czero$. A subset $U$ of $\Gamma \czero$ is open if $0 \notin U$ or if there exists $\gamma \in \Gamma$ such that $\Gamma_{<\gamma} \subseteq U$.

Then an element $\gamma \in \Gamma$ is cofinal for $\Gamma$ if and only if $0$ is in the closure of $\set{\gamma^n}{n \in \NN}$.
\end{remark}

\begin{example}\label{Height1Cofinal}
Let $\Gamma$ be a totally ordered group of height $1$. As $\Gamma$ is archimedean (Proposition~\ref{Height1Groups}), every element $\gamma \in \Gamma$ with $\gamma < 1$ is confinal for $\Gamma$.
\end{example}

\begin{rem}\label{RemCofinal}
Let $\Gamma$ be a totally ordered group, $\gamma \in \Gamma$ with $\gamma < 1$, and let $H$ be a subgroup of $\Gamma$. Then $\gamma$ is cofinal in $H$ if and only if the convex subgroup generated by $\gamma$ (Remark~\ref{GeneratedConvex}) contains $\set{h \in H}{h \leq 1}$ (or, equivalently, $H$).

Indeed, let $\gamma$ be cofinal in $H$ and let $h \in H$ with $h < 1$. Then there exists $n \in \NN$ such that $\gamma^n < h < 1$ and thus $h$ is contained in the convex subgroup generated by $\gamma$. Conversely, assume that the convex subgroup generated by $\gamma$ contains $\set{h \in H}{h \leq 1}$. It suffices to show that for all $h \in H$ with $h < 1$ there exists $n \in \NN$ with $\gamma^n < h$. By hypothesis there exists $m \in \ZZ$ such that $\gamma^m \leq h$. One necessarily has $m \geq 1$ because $\gamma < 1$. Then $\gamma^{m+1}< h$.
\end{rem}

\begin{prop}\label{CofinalConvex}
Let $\Gamma$ be a totally ordered group, $\Gamma' \subseteq \Gamma$ a convex subgroup, let $\gamma \in \Gamma$ be cofinal for $\Gamma'$, and let $\Delta \subsetneq \Gamma'$ be a proper convex subgroup. Then $\delta\gamma$ is cofinal for $\Gamma'$ for all $\delta \in \Delta$.
\end{prop}

\begin{proof}
Let $\gamma_0 \in \Gamma'$ be an element with $\gamma_0 < \Delta$ and hence $\Delta < \gamma_0^{-1}$. Then there exists $n \in \NN$ with $\gamma^n < \gamma_0$ and hence $(\delta\gamma)^{2n} < \gamma_0^{-1}\gamma^{2n} < \gamma^n$ which implies the claim.
\end{proof}

\begin{cor}\label{CofinalQuot}
Let $\Gamma$ be a totally ordered group, $\Delta \subsetneq \Gamma$ a proper convex subgroup. If $\gamma \in \Gamma$ is cofinal for $\Gamma$, then the image of $\gamma$ in $\Gamma/\Delta$ is cofinal for $\Gamma/\Delta$.
\end{cor}

%--------------------------------------------------------------------

\subsection{Valuations}
\begin{defi}
Let $A$ be a ring. A \emph{valuation of $A$} is a map $\abs\colon A \to \Gamma \cup \{0\}$, where $\Gamma$ is a 
totally ordered group, such that
\begin{dli}
\item
$|a+b| \leq \max(|a|,|b|)$ for all $a,b \in A$.
\item
$|ab| = |a||b|$ for all $a,b \in A$.
\item
$|0| = 0$ and $|1| = 1$.
\end{dli}
The subgroup of $\Gamma$ generated by $\im(\abs) \setminus \{0\}$ is called the \emph{value group of $\abs$}. It is 
denoted by $\Gamma_\abs$.

The set $\supp(\abs) := \abs^{-1}(0)$ is called the \emph{support of $\abs$}.
\end{defi}

\begin{longversion}
The second condition shows that $|u| \in \Gamma$ for every unit $u \in A^{\times}$ and that $\abs\rstr{A^{\times}}
\colon A^{\times} \to \Gamma$ is a homomorphism of groups. Moreover $|-1||-1| = 1$ shows that $|-1| = 1$ by 
\ref{PropTotallyOrdered}. Thus we have for all $a \in A$
\[
|-a| = |a|.
\]
\end{longversion}

\begin{example}
Let $A$ be a ring and $\pfr$ be a prime ideal of $A$. Then
\[
a \sends \begin{cases} 1,&a \notin \pfr;\\ 0,&a \in \pfr \end{cases}
\]
is a valuation with value group $1$. Every valuation on $A$ of this form is called a \emph{trivial valuation}.

If $F$ is a finite field, $F^\times$ is a torsion group. Thus on $F$ there exists only the trivial valuation.
\end{example}

\begin{rem}\label{PropVal}
Let $\abs$ be a valuation on $A$.
\begin{ali}
% \begin{longversion}
% \item
% $\Atilde(\abs) := \set{a \in A}{|a| \leq 1}$ is a subring of $A$.
% \item
% For every $\gamma \in \Gamma$ the set $\set{a \in A}{|a| \leq \gamma}$ is an $\Atilde(\abs)$-submodule of $A$.
% \end{longversion}
\item
Let $a,b \in A$ with $|a| \ne |b|$. Then $|a+b| = \max(|a|,|b|)$. Indeed let $|a| < |b|$ and assume that $|a+b| < |b|$. 
Then $|b| = |a+b-a| \leq \max(|a+b|,|a|) < |b|$. Contradiction.
\item
Let $\varphi\colon B \to A$ be a homomorphism of rings. Then $\abs \circ \varphi$ is a valuation on $B$ and $
\supp(\abs \circ \varphi) = \varphi^{-1}(\supp(\abs))$.
\end{ali}
\end{rem}

\begin{rem}\label{ResFieldVal}
Let $\abs\colon A \to \Gamma \cup \{0\}$ be a valuation of a ring $A$. Then $\supp(\abs)$ is a prime ideal in $A$ 
(and of $\set{a \in A}{|a| \leq 1}$). Denote by $K$ the field of fractions of $A/\supp(\abs)$. Define a valuation
\[
\abs'\colon K \to \Gamma \cup \{0\}, \qquad \frac{\abar}{\bbar} \sends |a||b|^{-1},
\]
where $\abar, \bbar \in A/\supp(\abs)$ are the images of elements $a,b \in A$. [Note that for $a \in A$, $a' \in 
supp(\abs)$ one has $|a+a'| = |a|$ by \ref{PropVal}. This and the multiplicativity of $\abs$ shows that $\abs'$ is well 
defined.]

The support of $\abs'$ is the zero ideal and $\Gamma_\abs = \Gamma_{\abs'}$.

Conversely let $\pfr \subset A$ be a prime ideal and let $\abs'$ be a valuation on $\kappa(\pfr)$. Then its support is 
the zero ideal. Composing $\abs'$ with the canonical ring homomorphism $A \to \kappa(\pfr)$ we obtain a valuation $
\abs$ on $A$.
\end{rem}

\begin{defi}\label{ValRing}
Let $\abs$ be a valuation on $A$. Then
\begin{align*}
\text{$K(\abs)$} &:= \Frac (A/\supp(\abs)),\\
\text{(resp.\ $A(\abs)$} &:= \set{x \in K(\abs)}{|x|' \leq 1},\\
\text{resp.\ $\mfr(\abs)$} &:= \set{x \in K(\abs)}{|x|' < 1},\\
\text{resp.\ $\kappa(\abs)$} &:= A(\abs)/\mfr(\abs)\ )
\end{align*}
is called the \emph{valued field} (resp.~\emph{valuation ring}, resp.~the \emph{maximal ideal}, resp.~the 
\emph{residue field}) \emph{of $\abs$}.

The \emph{height} (or \emph{rank}) \emph{of $\abs$} is defined as the height of the value group of $\abs$.
\end{defi}

\begin{propdef}\label{ValEquiv}
Two valuations $\abs_1$ and $\abs_2$ on a ring $A$ are called \emph{equivalent} if the following equivalent 
conditions are satisfied.
\begin{eli}
\item
There exists an isomorphism of totally ordered monoid $f\colon \Gamma_{\abs_1} \cup \{0\} \iso \Gamma_{\abs_2} 
\cup \{0\}$ such that $f \circ \abs_1 = \abs_2$.
\item
$\supp(\abs_1) = \supp(\abs_2)$ and $A(\abs_1) = A(\abs_2)$.
\item
For all $a,b \in A$ one has $|a|_1 \leq |b|_1$ if and only if $|a|_2 \leq |b|_2$.
\end{eli}
Note that $f$ in~(i) restricts to an isomorphism $\Gamma_{\abs_1} \iso \Gamma_{\abs_2}$ of totally ordered groups.
\end{propdef}

\begin{proof}
All three conditions imply that $\supp(\abs_1) = \supp(\abs_2) =: \pfr$. For $a \in A$ the valuations $|a|_1$ and $|a|
_2$ depend only on the image of $a$ in $A/\pfr$. Thus we may replace $A$ by $A/\pfr$. Then $A$ is an integral 
domain, $\supp(\abs_i) = 0$ for $i=1,2$, and we can extend $\abs_1$ and $\abs_2$ to the field of fractions $K$ of 
$A$. Thus we may assume that $\abs_1$ and $\abs_2$ are valuations on a field.

Clearly, (i) implies~(iii) and~(iii) implies~(ii). If~(ii) holds, then also $\set{a \in K}{|a|_1 \geq 1} = \set{a \in K}{|a|_2 
\geq 1}$ and hence $\abs_1\colon K^{\times} \to \Gamma_{\abs_1}$ and $\abs_2\colon K^{\times} \to 
\Gamma_{\abs_2}$ are surjective group homomorphisms with the same kernel. Thus there exists a unique group 
homomorphism $f\colon \Gamma_{\abs_1} \iso \Gamma_{\abs_2}$ such that $f \circ \abs_1 = \abs_2$. It maps 
elements $\leq 1$ to elements $\leq 1$ and thus is a homomorphism of totally ordered groups, and we can extend $f
$ to $\Gamma_{\abs_1} \cup \{0\}$ by setting $f(0) = 0$.
\end{proof}

\begin{longversion}
\begin{rem}
Let $A$ be a ring. For a valuation $\abs$ denote by $[\abs]$ its equivalence class. It follows from \ref{ResFieldVal} 
that $\abs \sends (\supp(\abs), \abs')$ yields a bijection
\begin{equation}\label{EquivSupp}
\begin{aligned}
&\set{[\abs]}{\text{$\abs$ valuation on $A$}}\\
\bijective &\set{(\pfr,[\abs'])}{\pfr \in \Spec A, \text{$\abs'$ valuation on $\Frac A/\pfr$}}.
\end{aligned}
\end{equation}
\end{rem}
\end{longversion}

\begin{example}
Let $\abs$ be a valuation on a field $k$ and let $\Gamma$ be its value group. Let $A = k\dlbrack T_1,\dots,T_r
\drbrack$ be the ring of formal power series in $r$ variables over $k$, and endow the group $\RR_{>0}^r \times \Gamma$ with the lexicographic order. Fix some real 
numbers $\rho_1,\dots,\rho_r$ with $0 < \rho_i < 1$. Then
\begin{align*}
\abs'\colon A &\to (\RR_{>0}^r \times \Gamma) \czero,\\
\sum_{\nbf = (n_1,\dots,n_r)}a_{\nbf}T_1^{n_1}\cdots T_r^{n_r} &\sends \rho_1^{m_1}\cdots \rho_r^{m_r}|a_{\mbf}|, 
\quad\text{where $\mbf = \inf\set{\mbf \in \ZZ^r}{a_{\mbf} \ne 0}$}
\end{align*}
is a valuation on $A$.
% Its restriction to $k(T_1,\dots,T_r)$ is the unique valuation $\abs'$ on $k(T_1,\dots,T_r)$ that extends $\abs$ on 
%$k$ and satisfies $|T_i|' = (\rho_i$, where $e_i \in \ZZ^n$ is the $i$-th standard basis vector. Note that $v$ itself is 
%not characterized by these properties.
\end{example}

\begin{longversion}
\begin{example}
Every totally ordered group $\Gamma$ is the value group of a valuation which can be constructed as follows. Let $k
$ be any field and set $C := k[\Gamma_{\leq 1}]$ (the monoid algebra of the monoid of elements $\leq 1$ in $
\Gamma$). By definition $C$ has a $k$-basis $(x_\gamma)_{\gamma \in \Gamma_{\leq 1}}$ and its multiplication is 
given by $x_\gamma x_\delta = x_{\gamma\delta}$. Define a valuation on $C$ by $|0| = 0$ and
\[
\bigl|\sum_{\gamma \in \Gamma_{\leq 1}}a_\gamma x_\gamma\bigr| := \sup\set{\gamma}{a_\gamma \ne 0}.
\]
Then $\abs$ is a valuation on $C$ with value group $\Gamma$ and support $0$.
\end{example}
\end{longversion}

%================================================================

\section{Valuation rings}

\subsection{Definition of valuation rings}

\begin{rem}
Let $A$ and $B$ local rings with $A \subseteq B$. Then we say that $B$ \emph{dominates} $A$ if $\mfr_B \cap A = 
\mfr_A$. Given a field $K$, the set of local subrings of $K$ is inductively ordered with respect to domination order.
\end{rem}

\begin{prop}\label{CharValRing}
Let $K$ be a field and let $A$ be a subring of $K$. The following conditions are equivalent.
\begin{equivlist}
\item
For every $a \in K^{\times}$ one has $a \in A$ or $a^{-1} \in A$.
\item
$\Frac A = K$ and the set of ideals of $A$ is totally ordered by inclusion.
\item
$A$ is local and a maximal element in the set of local subrings of $K$ with respect to the domination order. 
\item
There exists a valuation $\abs$ on $K$ such that $A = A(\abs)$.
\item
There exists an algebraically closed field $L$ and a ring homomorphism $h\colon A \to L$ which is maximal in the set of homomorphisms from subrings of $K$ to $L$ (here we define $h \leq h'$ for homomorphisms $h\colon A \to L$ and $h'\colon A' \to L$, where $A$, $A'$ are subrings of $K$, if $A \subseteq A'$ and $h = h'\rstr{A}$).
\end{equivlist}
\end{prop}

\begin{proof}
\proofstep{``(i) $\iff$ (iv)''}
Clearly (iv) implies (i). Conversely if (i) holds, it is easy to check that 
\[
|\cdot|\colon K\rightarrow K^\times/A^\times\cup\{0\},\qquad
 |x| =
  \begin{cases}
   0, & \text{if $x=0$;} \\
   [x] := x \mod A^\times, & \text{if $x\in K^\times$,}
  \end{cases}
\]
defines a valuation, where $K^\times/A^\times$ is ordered by $[x] \leq [y] :\Leftrightarrow xy^{-1}\in A$.

\proofstep{``(i) $\Rightarrow$ (ii)''}
$\Frac A = K$ is obvious. Let $\mathfrak{a},\mathfrak{b}\subseteq A$ be ideals. If $\mathfrak{a}$ is not a subset of $\mathfrak{b}$ and $a\in\mathfrak{a}\backslash \mathfrak{b}$, then $b^{-1}a \notin A$ holds for every $b\in\mathfrak{b}$ and thus $a^{-1}b\in A$. This implies $Ab\subseteq Aa$ for all $b \in \bfr$ and in particular $\mathfrak{b}\subseteq\mathfrak{a}$.

\proofstep{``(ii) $\Rightarrow$ (iii)''}
As the ideals of $A$ are totally ordered by inclusion, $A$ has to be local. Denote by $\mfr_A$ its maximal ideal. Let $B\subseteq K$ be a subring dominating $A$. For $x\in B\backslash\{0\}$ we can write $x=ab^{-1}$ with $a,b\in A$. We have to show that $x \in A$. In case $Aa\subseteq Ab$ it directly follows $x\in A$. On the other hand if $Ab\subseteq Aa$ we have $x^{-1}\in A$. But $x^{-1}\notin\mathfrak{m}_A$ since $B$ dominates $A$. This implies $x\in A^\times$.

\proofstep{``(iii) $\Rightarrow$ (v)''}
Let $\kappa(A)$ be the residue field of $A$, let $L$ be an algebraic closure of $\kappa(A)$, and let $h\colon A\rightarrow L$ the canonical homomorphism. Further let $A'\subseteq K$ be a subring containing $A$ and let $h'\colon A'\rightarrow L$ be an extension of $h$. For $\mathfrak{p}= \ker h'$ we get $\mathfrak{p}\cap A=\mathfrak{m}_A$. Hence $A'_\mathfrak{p}$ dominates $A_{\mathfrak{m}_A}=A$, which shows $A=A'_\mathfrak{p}$ and thus $A'=A$.

\proofstep{``(v) $\Rightarrow$ (i)''}
For this implication we will use the following easy remark: Let $R$ be a ring, $B$ an $R$-algebra, and $b\in B$. Then $b$ is integral over $A$ if and only if $b^{-1} \in A[b^{-1}]^\times$.

Now we are able to prove that (v) implies (i). First we show that $A$ is local with $\mathfrak{m}_A = \ker h$. If $h(x)\neq 0$ for $x\in A$ we get a natural extension of $h$ to $A_x$. By maximality of $h$ it then follows $A_x=A$ and thus $x\in A^\times$.

Now let $x\in K$ be integral over $A$. Then $A':=A[x]$ is a finite extension of $A$ and there exists a maximal ideal $\mathfrak{m}'$ of $A'$ so that $\mathfrak{m}'\cap A=\mathfrak{m}_A$. Hence $A'/\mathfrak{m}'$ is a finite extension of $A/\mathfrak{m}$ and because $L$ is algebraically closed there exists an extension of $h$ to $A'/\mathfrak{m}'$ and thereby to $A'$. Thus $A=A'$ and $x\in A$.

Now, if $x\in K$ is not integral over $A$ we have $x^{-1}\notin A[x^{-1}]^\times$ by the above remark. Hence there exists a maximal ideal $\mathfrak{m}'$ of $A[x^{-1}]$ containing $x^{-1}$. Since $A\rightarrow A[x^{-1}]\rightarrow A[x^{-1}]/\mathfrak{m}'$ is surjective, its kernel must be $\mathfrak{m}_A$ and we obtain a natural extension of $h$ to $A[x^{-1}]$. Therefore $A[x^{-1}]=A$ and $x^{-1}\in A$.
\end{proof}

\begin{defi}
An integral domain $A$ with field of fractions $K$ is called \emph{valuation ring} (\emph{of $K$}) if $K$ and $A$ 
satisfy the equivalent conditions in Prop.~\ref{CharValRing}.
\end{defi}

\begin{rem}
If $A$ is a valuation ring, then up to equivalence there exists a unique valuation $\abs$ on $K := \Frac A$ such that 
$A = A(\abs)$ (Proposition~\ref{ValEquiv}). Its value group $\Gamma_\abs$ is isomorphic $\Gamma_A := K^\times/
A^\times$ with the total order defined by $xA^{\times} \leq yA^{\times} :\iff xy^{-1} \in A$. The valuation is trivial if 
and only if $A = K$.
\end{rem}

\begin{rem}\label{ValResidueField}
Let $A$ be a valuation ring, $K := \Frac A$, $k := A/\mfr_A$, and let $\pi\colon A \to k$ be the canonical map. Then 
the map
\begin{align*}
\set{B}{\text{$B$ valuation ring of $k$}} &\to \set{\Btilde}{\text{$\Btilde$ valuation ring of $K$ contained in $A$}},\\
B &\sends \pi^{-1}(B)
\end{align*}
is (well defined and) bijective.
\end{rem}

%----------------------------------------------------------------

\subsection{Examples of valuation rings}

\begin{lem}\label{Lemdvr}
Let $A$ be a local ring whose maximal ideal is generated by one element $p$. Assume that $\bigcap_{n \geq 1}Ap^n 
= 0$ (which is automatic if $A$ is noetherian). Define
\[
v\colon A\rightarrow\mathbb{Z}\cup\{\infty\}, \qquad v(x):=
\begin{cases}
n, & \text{if $x\in Ap^n\backslash Ap^{n+1}$;}\\
\infty, & \text{otherwise}.
\end{cases}
\]
Then the only ideals of $A$ are $0$ and $Ap^n$ for $n \geq 0$. Moreover $p$ is either nilpotent or $A$ is a valuation ring and its valuation is given by $x \sends \gamma^{v(x)}$, where $\gamma$ is any real number with $0 < \gamma < 1$.
\end{lem}

\begin{proof}
Since $\bigcap_{n\geq 1}Ap^n$=0 we have $v(x)=\infty$ if and only if $x=0$. Let $\mathfrak{a}$ be a non-zero ideal of $A$ and $a\in\mathfrak{a}$ such that $v(a)$ is minimal. Clearly $\mathfrak{a}\subseteq Ap^{v(a)}$. Furthermore, there exists an $u\in A$ with $a=up^{v(a)}$ and clearly $v(u) = 0$, i.e., $u \in A^{\times}$. This shows $p^{v(a)}\in Aa$ and thus $\mathfrak{a}=Ap^{v(a)}$. 
As we have seen that every element $a \in A$ can be written as $a=up^{v(a)}$ with $u \in A^\times$, it also follows that if $p$ is not nilpotent, $A$ will be an integral domain with valuation $v$.
\end{proof}

\begin{prop}\label{chardvr}
Let $A$ be an integral domain which is not a field. Then the following assertions are equivalent.
\begin{eli}
\item
$A$ is a noetherian valuation ring.
\item
$A$ is a local principal domain.
\item
$A$ is a valuation ring and $\Gamma_A$ is isomorphic to the totally ordered group $\ZZ$.
\item
$A$ is local, its maximal ideal $\mfr$ is a principal ideal, and $\bigcap_{n>0}\mfr^n = 0$.
\end{eli}
\end{prop}

\begin{proof}
\proofstep{``(iii) $\Rightarrow$ (ii),(i)''}
Without loss of generality we may assume that $\Gamma_{A}$ is the (additive) group $\mathbb{Z}$. Let $\pi \in A$ with $|\pi|=1$. For $x\in A$ we have $|x|=n=|\pi^n|$ which implies $x=u\pi^n$ for a unit $u\in A^\times$. Now (ii) and (i) follow as in the proof of Lemma~\ref{Lemdvr}.

\proofstep{``(ii) $\Rightarrow$ (iv)''}
Let $\mathfrak{m}_A=A\pi$. Then $\pi$ is the only prime element of $A$ up to multiplication with units of $A$ and every non-zero ideal has the form $A\pi^n$. As $A$ is factorial, the only element of $A$ dividing all powers of the prime element $\pi$ is $0$, i.e., $\bigcap_{n\geq 1}\mathfrak{m}_A^n=0$.

\proofstep{``(iv) $\Rightarrow$ (iii)''}
This follows directly from Lemma~\ref{Lemdvr} with $v$ as valuation.

\proofstep{``(i) $\Rightarrow$ (ii)''}
Let $\mathfrak{a}\subseteq A$ be an ideal generated by elements $x_1,\dots,x_n\in A$. Since the ideals $x_i A$ are totally ordered by inclusion, $\afr$ is generated by one of the elements $x_i$. Thus $A$ is a principal domain.
\end{proof}

\begin{defi}
Let $A$ be an integral domain which is not a field. If the equivalent properties of Prop.~\ref{chardvr} are satisfied, $A
$ is called a \emph{discrete valuation ring} and its valuation is called a \emph{discrete valuation}. It is called normed 
if its value group is $\ZZ$.
\end{defi}

\begin{example}\label{PAdicVal}
Let $A$ be a factorial ring, $K = \Frac A$. Let $p \in A$ be a prime element. Let $v_p\colon K \to \ZZ \cup \{\infty\}$ 
be the $p$-adic valuation. Fix a real number $0 < \rho < 1$ and define $|x| := \rho^{v_p(x)}$ for $x \in K$. Then $
\abs_p$ is a valuation (also called the $p$-adic valuation). Its valuation ring is $A_{(p)}$. Its equivalence class does 
not depend on $\rho$.
\end{example}

\begin{example}\label{ValPID}
Let $A$ be a principal ideal domain and $K$ its field of fractions. The valuation rings of $K$ containing $A$ and 
distinct from $K$ are the rings $A_{(p)}$, where $p$ is a prime element of $A$.

In particular, every valuation ring of the field $\QQ$ and distinct from $\QQ$ is of the form $\ZZ_{(p)}$ (because 
every subring of $\QQ$ contains $\ZZ$).
\end{example}

\begin{example}
Let $k$ be a field. The only valuations on $k(T)$ whose restriction to $k$ is trivial are (up to equivalence) the $p$-
adic valuations, where $p \in k[T]$ is an irreducible polynomial, and the valuation $f/g \sends \rho^{\deg(g) - \deg(f)}$ 
for $f,g \in k[T]$, $g \ne 0$ and $0 < \rho < 1$ a real number.
\end{example}

%------------------------------------------------------------

\subsection{Ideals in valuation rings}

\begin{prop}\label{IdealVal}
Let $K$ be a field, let $\abs$ be a valuation on $K$, $A$ its valuation ring, and let $\Gamma$ be the value group of 
$\abs$. Consider the following condition on subsets $M$ of $\Gamma$
\[
\gamma \in M, \delta \leq \gamma \implies \delta \in M.\tag{$\ast$}
\]
Then the map $M \sends \afr(M) := \set{x \in K}{|x| \in M \cup \{0\}}$ yields a bijection between the subsets of $
\Gamma$ satisfying $(\ast)$ and the set of $A$-submodules of $K$. Under this bijection one has the following 
correspondences.
\begin{ali}
\item
The ideals of $A$ correspond to those subsets $M$ satisfying $(\ast)$ contained in $\Gamma_{\leq 1}$.
\item
The monogenic (or, equivalently, the finitely generated) $A$-submodules of $K$ correspond to the subsets of the 
form $\Gamma_{\leq \gamma}$ for some $\gamma \in \Gamma$. In this case $\afr(\Gamma_{\leq \gamma})$ is 
generated by any element $x \in K$ with $|x| = \gamma$.
\end{ali}
\end{prop}

For instance the maximal ideal of $A$ corresponds to $\Gamma_{<1}$.

\begin{proof}
For every $A$-subomdule $\mathfrak{b}$ of $K$ we define the set $M(\mathfrak{b}):=\{|x|;x\in\mathfrak{b}\backslash\{0\}\}$. We note that these sets satisfy ($\ast$). We claim that this defines an inverse map to $M \sends \afr(M)$. Thus we have to show:
\begin{dli}
\item
$M(\mathfrak{a}(N))=N$ for every subset $N$ of $\Gamma$ satisfying ($\ast$).
\item
$\mathfrak{a}(M(\mathfrak{b}))=\mathfrak{b}$ for every $A$-submodule $\mathfrak{b}$ of $K$.
\end{dli}
Now~(a) is obvious, as is $\mathfrak{b}\subseteq\mathfrak{a}(M(\mathfrak{b}))$. On the other hand for $x\in\mathfrak{a}(M(\mathfrak{b}))$ we have $|x|\in M(\mathfrak{b})$. Hence there exists $y\in\mathfrak{b}$ with $|x|=|y|$ which implies $x=uy$ where $|u|=1$. This finally shows $x\in Ay\subseteq\mathfrak{b}$.

Furthermore (1) is obvious and concerning (2) we observe that for $x\in K$ with $|x|=\gamma$ the relation $x^{-1}\mathfrak{a}(\Gamma_{\leq\gamma})\subseteq A$ holds.
\end{proof}

\begin{rem}
Let $K$ be a field and let $A$ be a valuation ring of $K$. Then every ring $B$ with $A \subseteq B \subseteq K$ is 
a valuation ring.
\end{rem}

\begin{prop}\label{PrimeIdealVal}
Let $K$ be a field, let $\abs$ be a valuation on $K$, $A$ its valuation ring, and let $\Gamma$ be the value group of 
$\abs$. Define
\begin{align*}
\Rcal &:= \set{\text{$B$ (valuation) ring}}{A \subseteq B \subseteq K};\\
\Ical &:= \set{\Delta}{\text{$\Delta$ convex subgroup of $\Gamma$}}.
\end{align*}.
\begin{ali}
\item
The maps
\begin{align*}
\Spec A \to \Rcal, &\qquad \pfr \sends A_{\pfr},\\
\Rcal \to \Spec A, &\qquad B \sends \mfr_B
\end{align*}
are well defined mutually inverse inclusion reversing bijections.
\item
The maps
\begin{align*}
\Ical \to \Rcal, &\qquad \Delta \sends A(K \ltoover{\abs} \Gamma \to \Gamma/\Delta),\\
\Rcal \to \Ical, &\qquad B \sends |B^{\times}|
\end{align*}
are well defined mutually inverse inclusion preserving bijections.
\end{ali}
\end{prop}

\begin{proof}
\proofstep{(1)}
For $B\in\mathcal{R}$ let $x\in \mfr_B$. Then $x^{-1} \in K\setminus B \subseteq K \setminus A$ and hence $x \in\mathfrak{m}_A$. This shows $\mathfrak{m}_B\subseteq\mathfrak{m}_A$. Furthermore $\mathfrak{m}_B = \mfr_B \cap A$ is a prime ideal in $A$. The map $\mathcal{R}\rightarrow\mathrm{Spec}\ A$ is inclusion reversing. Now clearly $A_{\mathfrak{m}_B}\subseteq B$ holds. Conversely for $x\in B\backslash A$ we have $x^{-1}\in A$ and $x^{-1}\notin\mathfrak{m}_B$, so that $x\in A_{\mathfrak{m}_B}$. Thus $A_{\mathfrak{m}_B}=B$.

On the other hand let $\mathfrak{p}\in\mathrm{Spec}\ A$ and $B:=A_\mathfrak{p}$, then $\pfr = \mathfrak{m}_B \cap A = \mfr_B$ which finishes the proof of (1).

\proofstep{(2)}
Let $B\in\mathcal{R}$. We may assume $\Gamma_A=K^\times / A^\times$ and $\Gamma_B=K^\times/ B^\times$, so that we get a surjective homomorphism of ordered groups $\lambda\colon\Gamma_A\rightarrow\Gamma_B$ with kernel $B^\times/A^\times=:H_B$. Furthermore we have $|\cdot|_B=\lambda\circ|\cdot|_A$ and we can identify $\Gamma_A /H_B$ and $\Gamma_B$. Therefore $H_B$ defines $|\cdot|_B$ up to equivalence, hence determines $B$ uniquely. Conversely for every $\Delta\in\mathcal{I}$ the mapping $K\xrightarrow[]{|\cdot|}\Gamma\rightarrow\Gamma/\Delta$ is a valuation on $K$ whose ring contains $A$ since the map $\Gamma\rightarrow\Gamma/\Delta$ is a homomorphism of ordered groups.
\end{proof}

%-----------------------------------------------------------------

\begin{longversion}
\subsection{Extension of valuation}

\begin{prop}\label{ExistExtVal}
Let $K$ be a field, let $\abs$ be a valuation on $K$, and let $K'$ be a field extension of $K$. Then there exists a 
valuation $\abs'$ on $K'$ whose restriction to $K$ is equivalent to $\abs$. In other words there exists a commutative 
diagram
\begin{equation}\label{EqExtVal}
\begin{aligned}\xymatrix{
K' \ar[r]^{\abs'} & \Gamma_{\abs'}\czero \\
K \ar[u] \ar[r]^{\abs} & \Gamma_{\abs}\czero \ar[u]_{\lambda}\ ,
}\end{aligned}
\end{equation}
where $\lambda$ is an injective homomorphism of totally ordered groups.

Moreover, if $x_1,\dots,x_n$ are elements in $K'$ that are algebraic independent over $K$, and $\gamma_1,\dots,
\gamma_n$ are any elements of $\Gamma$. Then $\abs'$ may be chosen such that $|x_i|' = \gamma_i$ for all $i = 
1,\dots,n$.
\end{prop}

\begin{proof}
\cite{Bou_AC}~VI, \S2.4, Prop.~4.
\end{proof}

\begin{prop}\label{ExtValAlgebraic}
Let $K$ be a field, let $K' \supset K$ be an algebraic extension, let $\abs'$ be a valuation on $K'$ and let $\abs$ be 
its restriction to $K$. Let $k$ (resp.~$k'$) be the residue field of $A(\abs)$ (resp.~$A(\abs')$).
\begin{ali}
\item
$\Gamma_{\abs'}/\Gamma_\abs$ is a torsion group and $\height(\Gamma_{\abs'}) = \height(\Gamma_{\abs})$.
\item
$k'$ is an algebraic extension of $k$.
\item
$\abs$ is trivial (resp.~of height $1$) if and only if $\abs'$ is trivial (resp.~of height $1$).
\item
Let $K'$ be a finite extension of $K$.  Then
\[
[k':k](\Gamma_{\abs'}: \Gamma_\abs) \leq [K' : K];
\]
in particular $\Gamma_{\abs'}/\Gamma_\abs$ is a finite group and $k'$ is a finite extension of $k$. Moreover, $\abs$ 
is discrete if and only if $\abs'$ is discrete.
\end{ali}
\end{prop}

\begin{proof}
\cite{Bou_AC} VI, \S8.
\end{proof}

In general, an extension of a valuation on a field $K$ to a field extension $K'$ is not unique. In fact, this is never the 
case if $K'$ is transcendental over $K$ as the second assertion of Prop.~\ref{ExistExtVal} shows. But even if $K'$ 
is an algebraic extension we need some additional hypothesis. In fact this is closely connected to Hensel's lemma. 
Thus there is the following definition.

\begin{propdef}
Let $A$ be a local ring, $k$ its residue field. Then $A$ is called \emph{henselian} if the following equivalent 
conditions are satisfied.
\begin{eli}
\item
Every finite $A$-algebra is a (necessarily finite) product of local rings.
\item
For every monic polynomial $F \in A[T]$ with image $F_0 \in k[T]$ and for every decomposition $F_0 = G_0H_0$ 
into polynomials $G_0,H_0 \in k[T]$ that are prime to each other, there exists a unique pair $(G,H)$ of monic 
polynomials in $A[T]$ such that: $G_0$ and $H_0$ are the images of $G$ and $H$, respectively, $F = GH$, and the 
ideal of $A[T]$ generated by $G$ and $H$ is $A[T]$.
\item
For every smooth morphism of schemes $f\colon X \to \Spec A$ the canonical map $X(A) \to X(k)$ is surjective 
(where, with $S = \Spec A$, $X(A) = \Hom_S(S,X)$ and $X(k) = \Hom_S(\Spec k, X) = \Hom_k(\Spec k, X
\otimes_A k)$).
\end{eli}
\end{propdef}

\begin{proof}
\cite{EGA}~IV~(18.5.11), (18.5.13), (18.5.17).
\end{proof}

\begin{prop}\label{ExtValHensel}
Let $\abs$ be a valuation on a field $K$. Then $A(\abs)$ is henselian if and only if for every algebraic extension $K'$ 
of $K$ any two extensions of $\abs$ to $K'$ are equivalent.

If these equivalent conditions are satisfied, $A(\abs')$ is also henselian.
\end{prop}

\begin{proof}
See e.g. \cite{FuchsSalce}~II, Theorem~7.3
\end{proof}

\begin{rem}\label{CompleteHensel}
Let $\abs$ be a valuation on a field $K$. By Hensel's lemma, if $A(\abs)$ is complete and $\abs$ is of height $1$, 
then $A(\abs)$ is henselian.

But there exist valuations of height $2$ such that $A(\abs)$ is complete but there exist two non-equivalent 
extensions to a quadratic extension of $K$ (\cite{Bou_AC} VI, \S8, Exercise~1). In particular $A(\abs)$ is not 
henselian.
\end{rem}

For transcendental extensions we extend valuations almost ``arbitrarily'':

\begin{prop}
Let $K$ be a field, $\abs$ a valuation on $K$, $\Gamma$ its value group. Let $\Gamma'$ be a totally ordered group 
containing $\Gamma$ and let $\gamma \in \Gamma'$. Then there exists a valuation $\abs'$ on $K(T)$ with values in 
$\Gamma'$ such that $|T|' = \gamma$.
\end{prop}

\begin{proof}
One easily checks that
\begin{equation}\label{EqExtValTrans}
|\sum_i a_iT^i|' = \max_i\{|a_i|\gamma^i\}
\end{equation}
defines a valuation on $K(T)$ that extends $\abs$.
\end{proof}

It is not difficult to construct other extensions than~\eqref{EqExtValTrans} of $\abs$ to $K(T)$ (e.g.~\cite{Bou_AC} VI, \S3, 
Exercise~1).

The next two corollaries show that as special cases one can extend a valuation to $K(T)$ without changing the 
residue field or without changing the value group.

\begin{cor}
Let $K$ be a field, $\abs$ a valuation on $K$, $\Gamma$ its value group. Let $\Gamma'$ be a totally ordered group 
containing $\Gamma$ and let $\gamma \in \Gamma'$ be an element whose image in $\Gamma'/\Gamma$ has 
infinite order. Then there exists a unique valuation $\abs'$ on $K(T)$ extending $\abs$ such that $|T|' = \gamma$. 
One has $\kappa(\abs') = \kappa(\abs)$ and $\Gamma_{\abs'} = \Gamma_{\abs}\cdot \gamma^{\ZZ}$.
\end{cor}

\begin{proof}
We have only to show the uniqueness assertion. But using the hypothesis on $\gamma$ this follows easily from $|a
+b|' = \max\{|a|',|b|'\}$ if $|a|' \ne |b|'$.
\end{proof}

We may also extend $\abs$ by choosing $\gamma = 1$ in~\eqref{EqExtValTrans}.

\begin{cor}
Let $K$ be a field, $\abs$ a valuation on $K$, $\Gamma$ its value group. Then there exists a unique valuation $
\abs'$ on $K(T)$ extending $\abs$ such that $|T|' = 1$ and such that the image $t$ of $T$ in $\kappa(\abs')$ is 
transcendental over $\kappa(\abs)$. One has $\kappa(\abs') = \kappa(\abs)(t)$ and $\Gamma_{\abs'} = 
\Gamma_{\abs}$.
\end{cor}

\begin{proof}
\cite{Bou_AC} VI, \S10.1, Prop.~2
\end{proof}

\begin{rem}
Let $K$ be a field, $\abs$ a valuation on $K$, $\Gamma$ its value group. Let $K' = K(T_1,\dots,T_n)$ be a finitely 
generated purely transcendental extension of $K$. Arguing by induction one can therefore find for every all integers 
$r,s \geq 0$ with $r+s = n$ an extension $\abs'$ of $\abs$ to $K'$ such that
\begin{align*}
\Gamma_{\abs'} &= \Gamma \times \prod_{i=1}^r\ZZ,\qquad\text{lexicographically ordered},\\
\kappa(\abs') &\cong \kappa(\abs)(t_1,\dots,t_s)
\end{align*}
\end{rem}

Conversely there are the following restrictions.

\begin{prop}
Let $K$ be a field, $\abs$ a valuation on $K$, $\Gamma$ its value group, $k$ its residue field. Let $K'$ be an 
extension of $K$, $\abs'$ a valuation on $K'$ extending $\abs$, $\Gamma'$ its value group, $k'$ its residue field. 
Then
\begin{align*}
\trdeg(k'/k) + \dim_{\QQ}((\Gamma'/\Gamma) \otimes_{\ZZ} \QQ) &\leq \trdeg(K'/K),\\
\height(\Gamma') - \height(\Gamma) &\leq \dim_{\QQ}((\Gamma'/\Gamma) \otimes_{\ZZ} \QQ)
\end{align*}
\end{prop}

\begin{proof}
\cite{Bou_AC} VI, \S10.3.
\end{proof}

\begin{prop}\label{ExtendSpecializations}
Let $K_1 \mono K_2$ be a field extension, let $A_2$ be a valuation ring of $K_2$ and set $A_1 := K_1 \cap A_2$ 
(which is a valuation ring of $K_1$). For $i = 1,2$ let $\Gcal_i$ (resp.~$\Scal_i$) be the set of valuation rings of $K_i
$ which contain $A_i$ (resp.~which are contained in $A_i$). Then $B \sends B \cap K_1$ yields surjective maps $
\Gcal_2 \epi \Gcal_1$ and $\Scal_2 \epi \Scal_1$.

If $K_2$ is algebraic over $K_1$, then the map $\Gcal_2 \to \Gcal_1$ is bijective.
\end{prop}

\begin{proof}
Any valuation on a field can be extended to any field extension (Prop.~\ref{ExistExtVal}). If we apply this to the 
extension of residue field of $A_1$ and $A_2$, the surjectivity of $\Scal_2 \to \Scal_1$ follows from Remark~
\ref{ValResidueField}. Applying this to $K_1 \mono K_2$ itself, the surjectivity of $\Gcal_2 \to \Gcal_1$ follows from 
Proposition~\ref{PrimeIdealVal}.

Let $\Gamma_i$ be the value group (of the equivalence class) of the valuation on $K_i$ given by $A_i$. If $K_2 
\supset K_1$ is algebraic, then the cokernel of $\Gamma_1 \mono \Gamma_2$ is torsion and $\Delta \sends 
\Gamma_1 \cap \Delta$ is a bijection from the set of convex subgroups of $\Gamma_2$ onto the set of convex 
subgroups of $\Gamma_1$ (Proposition~\ref{ExtValAlgebraic}). Thus the bijectivity of $\Gcal_2 \to \Gcal_1$ follows 
from Proposition~\ref{PrimeIdealVal}.
\end{proof}

\end{longversion}

%====================================================================

\section{Spectral spaces}

Similarly as general schemes are obtained by gluing affine schemes, adic spaces will be obtained by gluing affinoid 
adic spaces. In both cases the underlying topological spaces are spectral spaces, which first have been singled out 
by Hochster (\cite{Ho_PrimeSpec}).

\subsection{Sober spaces}
We first recall the definition of irreducible spaces.

\begin{remdef}\label{DefIrred}
Let $X$ be a non-empty topological space. Then $X$ is called \emph{irreducible} if $X$ satisfies one of the following 
equivalent conditions.
\begin{eli}
\item
$X$ cannot be expressed as the union of two proper closed subsets.
\item
Every non-empty open subset of $X$ is dense.
\item
Any two non-empty open subsets of $X$ have non-empty intersection.
\end{eli}
A non-empty subset $Z$ of $X$ is called \emph{irreducible} if $Z$ is irreducible when we endow it with the induced 
topology.
\end{remdef}

\begin{lem}\label{IrredClosure}
Let $X$ be a topological space. A subspace $Y \subseteq X$ is irreducible if and only if its closure $\overline{Y}$ is 
irreducible.
\end{lem}

\begin{proof}
A subset $Z$ of $X$ is irreducible if and only if for any two open subsets $U$ and $V$ of $X$ with $Z \cap U \ne 
\emptyset$ and $Z \cap V \ne \emptyset$ we have $Z \cap (U \cap V) \ne \emptyset$. This implies the lemma 
because an open subset meets $Y$ if and only if it meets $\overline{Y}$.
\end{proof}

\begin{rem}
In particular, a subset of the form $\overline{\{x\}}$ for $x \in X$ is always irreducible.
\end{rem}

\begin{longversion}
\begin{rem}\label{IrredOpen}
Let $X$ be a topological space. If $U \subseteq X$ is an open subset and $Z \subseteq X$ is irreducible and 
closed, $Z \cap U$ is open in $Z$. Hence if $Z \cap U \ne \emptyset$, then $Z \cap U$ is an irreducible closed 
subset of $U$ whose closure in $X$ is $Z$. Together with Lemma~\ref{IrredClosure} this shows that there are 
mutually inverse bijective maps
\begin{equation}\label{BijIrredOpen}
\begin{aligned}
\{\text{$Y \subseteq U$ irreducible closed}\} &\bijective \{\text{$Z \subseteq X$ irreducible closed with $Z \cap U \ne 
\emptyset$}\} \\
Y & \mapsto \overline{Y} \quad\text{(closure in $X$)} \\
Z \cap U & \leftmapsto Z
\end{aligned}
\end{equation}
\end{rem}

\begin{remdef}\label{RemIrredComp}
An irreducible subset that is maximal w.r.t. inclusion is called an \emph{irreducible component}.

The set of irreducibly subsets is inductively ordered. Hence every irreducible subset is contained in an irreducible 
component. In particular $X$ is the union of its irreducible components. By Lemma~\ref{IrredClosure} irreducible 
components are always closed. 

A point $x \in X$ is called a \emph{maximal point} if its closure $\overline{\{x\}}$ is an irreducible component of $X$.
\end{remdef}
\end{longversion}

\begin{defi}\label{DefGeneric}
Let $X$ be an arbitrary topological space.
\begin{assertionlist}
\item
A point $x \in X$ is called \emph{closed} if the set $\{x\}$ is closed,
\item
We say that a point $\eta \in X$ is a \emph{generic point} if $\overline{\{\eta\}} = X$.
\item
Let $x$ and $x'$ be two points of $X$. We say that $x$ is a \emph{generization of $x'$} or that $x'$ is a 
\emph{specialization of $x$} if $x' \in \overline{\{x\}}$. We write $x \succeq x'$ or $x' \preceq x$.
\end{assertionlist}
\end{defi}

\begin{defi}
A topological space $X$ is called \emph{sober} if every irreducible closed subset of $X$ has a unique generic point.
\end{defi}

\begin{rem}\label{SoberKolm}
Recall that a topological space $X$ is called \emph{Kolmogorov} (or $T_0$) if for any two distinct points there exists 
an open set containing one of the points but not the other. 

A topological space is Kolmogorov if and only if every closed irreducible subset admits at most one generic point. 
Indeed, for $x,y \in X$ with $x \ne y$ one has $\overline{\{x\}} \ne \overline{\{y\}}$ if and only if there exists a closed 
subset containing one of the points but not the other one. By taking complements this is equivalent to the existence 
of an open subset containing one point but not the other.

In particular we see that every sober space is Kolmogorov.
\end{rem}

\begin{rem}\label{SpecialPartialOrder}
If $X$ is a sober space (or, more generally, a topological space where every closed irreducible subset has at most 
one generic point), then $\preceq$ is a partial order on $X$.
\end{rem}

% \begin{prop}\label{SoberPerm}
% \begin{ali}
% \item
% Any locally closed subspace of a sober space $X$ is again sober.
% \item
% Let $X$ be a topological space. Then finite unions and arbitrary intersections of sober subspaces of $X$ are again sober.
% \end{ali}
% \end{prop}
% 
% \begin{proof}
% Any subspace of a Kolmogorov space is again Kolmogorov. And in a Kolmogorov space the generic point of an irreducible subset is unique if it exists. 
% 
% \proofstep{(2)} Let $Y$ and $Y'$ be two sober subspaces of $X$ and let $Z$ be a closed irreducible subset of $Y \cup Y'$. Let $M$ (resp.~$M'$) be the set of maximal points of $Z \cap Y$ (resp.~of $Z \cap Y'$). 
% 
% Let $(X_i)_i$ be a family of sober subspaces. Let $Z$ be a closed irreducible subset of $\bigcap_i X_i$ and let $Z_i$ be the closure of $Z$ in $X_i$. By hypothesis $Z_i$ has a unique generic point $z_i$.
% 
% \proofstep{(1)} It is clear that any closed subspace of $X$ is sober. Let $U$ be an open subspace, let $Z \subseteq U$ an irreducible subset which is closed in $U$, and let $\overline{Z}$ be its closure in $X$. Then $\overline{Z}$ is irreducible and has a unique generic point $z$. Moreover $Z$ is open and dense in $\overline{Z}$ and hence $z \in Z$ and $z$ is a generic point of $Z$. It is unique because $U$ is Kolmogorov. Thus by~(2) any locally closed subspace is sober.
% \end{proof}

\begin{longversion}
\begin{example}
\begin{ali}
\item
If $X$ is a Hausdorff space, then the only irreducible subsets of $X$ are the sets $\{x\}$ for $x \in X$.
\item
Let $Z$ be an infinite set endowed with the topology such that the closed sets $\ne Z$ are all finite subsets. Then 
$Z$ is irreducible but has no generic point. In particular $Z$ is not sober. If we add a single point $\eta$ without 
changing the closed subsets $\ne Z \cup \{\eta\}$, then $Z \cup \{\eta\}$ is sober.

In particular we see that there are subspaces of sober spaces which are not sober.
\item
Let $X$ be a topological space in which $X$ and $\emptyset$ are the only open subsets. Then every point of $X$ is 
a generic point.
\end{ali}
\end{example}
\end{longversion}

%---------------------------------------------------------------------

\subsection{Spectral spaces}

\begin{defi}
Let $X$ be a topological space and let $f\colon X \to Y$ be a continuous map of topological spaces.
\begin{ali}
\item
$X$ is called \emph{quasi-separated} if the intersection of any two quasi-compact open subsets is again quasi-compact.
\begin{longversion}
\item
A subspace $Z$ of $X$ is called \emph{retro-compact} if the inclusion $Z \to X$ is quasi-compact.
\end{longversion}
\item
$f$ is called \emph{quasi-compact} if for every open quasi-compact subset $V$ of $Y$ the inverse image $f^{-1}(V)$ 
is quasi-compact.
\begin{longversion}
\item
$f$ is called \emph{quasi-separated} if for every open quasi-separated subset $V$ of $Y$ the inverse image $f^{-1}
(V)$ is quasi-separated.
\end{longversion}
\end{ali}
We abbreviate the property ``quasi-compact and quasi-separated'' to ``qcqs''. 
\end{defi}

\begin{defi}\label{DefSpectralSpace}
A topological space $X$ is called \emph{spectral} if it satisfies the following conditions.
\begin{dli}
\item
$X$ is quasi-compact.
\item
$X$ has a basis consisting of quasi-compact open subsets which is closed under finite intersections.
\item
$X$ is sober.
\end{dli}
$X$ is called \emph{locally spectral} if $X$ has an open covering by spectral spaces.
\end{defi}

\begin{rem}\label{RemSpectral}
\begin{ali}
\item
Any locally spectral space $X$ is sober and in particular Kolmogorov (Remark~\ref{SoberKolm}).
\item\label{SpectralQuasiSep}
A locally spectral space is spectral if and only if it is qcqs.
\begin{longversion}
\item
A finite sum of spectral spaces is again spectral.
\item
Every quasi-compact open subspace $U$ of a spectral space is again spectral: By~\eqref{SpectralQuasiSep}, $U$ 
inherits from $X$ a basis of quasi-compact open subsets which is closed under intersection. By~
\eqref{BijIrredOpen} any open subspace of a sober space is again sober.
\item
Every closed subspace of a spectral space is again spectral: This is clear as any closed subspace of a quasi-
compact is again quasi-compact and every closed subspace of sober space is sober.
\item
A locally spectral space $X$ has a basis consisting of spectral open subspaces and $X$ is sober. Note that 
Example~\ref{PatholExample} below shows that it may happen that such a basis is not stable under finite 
intersections (in oher words, $X$ is not quasi-separated).
\item
Every open subspace of a locally spectral space is locally spectral.
\end{longversion}
\end{ali}
\end{rem}

\begin{example}
Let $A$ be a ring and endow $\Spec A := \set{\pfr}{\text{$\pfr$ prime ideal}}$ with the usual topology (i.e. the closed sets of $\Spec A$ are of the form $V(\afr) := \set{\pfr}{\afr \subseteq \pfr}$ for ideal $\afr$ of $A$). Then $\Spec A$ is spectral: For $f \in A$ set $D(f) = \set{\pfr}{f \notin \pfr}$. Then the $(D(f))_{f \in A}$ form a basis of open quasi-compact subsets of $\Spec A$ stable under finite intersections.

For two prime ideals $\pfr$ and $\qfr$ of $A$, $\pfr$ is a specialization of $\qfr$ (or, equivalently, $\qfr$ is a 
generization of $\pfr$) if and only if $\qfr \subseteq \pfr$. For every irreducible closed subset $Z$ of $\Spec A$ there 
exists a unique prime ideal $\pfr \in Z$ such that $\pfr \subseteq \qfr$ for all $\qfr \in Z$.

\begin{longversion}
Let $X$ be a scheme. Then the underlying topological space of $X$ is spectral if and only if $X$ is qcqs.
\end{longversion}
\end{example}

In fact, Hochster has shown in~\cite{Ho_PrimeSpec} the converse.

\begin{theorem}\label{ThmSpectralSpace}
\begin{ali}
\item
For a topological space $X$ the following assertions are equivalent.
\begin{eli}
\item
$X$ is spectral.
\item
$X$ is homeomorphic to $\Spec A$ for some ring $A$.
\begin{longversion}
\item
$X$ is homeomorphic to the underlying topological space of a quasi-compact quasi-separated scheme.
\end{longversion}
\item
$X$ is the projective limit of finite Kolmogorov spaces.
\end{eli}
\item
A topological space is locally spectral if and only it is homeomorphic to the underlying topological space of a 
scheme.
\item
For every continuous quasi-compact map $f\colon X' \to X$ between spectral spaces there exists a ring 
homomorphism $\varphi\colon A \to A'$ such that the associated continuous map is $f$.
\end{ali}
\end{theorem}

\begin{example}\label{PatholExample}
Set $X := \NN \cup \{\infty\}$ endowed with the usual (total) order and define its set of open subsets by $\set{X_{< 
x}}{x \in X} \cup \{X\}$.

The only open subset containing $\infty$ is $X$ itself, in particular $X$ is quasi-compact. For all $x \in X$ the set 
$X_{\leq x}$ is open quasi-compact and these sets form a basis of the topoology stable under finite intersection. 
The non-empty closed subsets are the sets $X_{\geq x}$ for $x \in X$. All of them are irreducible with unique 
generic point $x$. Thus $X$ is a spectral space.

In fact, $X$ is the spectrum of a valuation ring with value group $\prod_{\NN} \ZZ$ (ordered lexicographically).

Note that $X_{< \infty}$ is open and not quasi-compact. If we glue two copies of $X$ along $X_{< \infty}$, then the 
resulting space is locally spectral and quasi-compact but it is not quasi-separated and in particular not spectral.
\end{example}

% \begin{longversion}
% \begin{rem}
% Every quasi-compact open subset of a spectral space is retro-compact (Remark~\ref{RemSpectral}~{SpectralQuasiSep}).
% \end{rem}
% 
% 
% \begin{defi}
% A continuous map $f\colon X' \to X$ between locally spectral spaces is called \emph{spectral} if for every open spectral subspace $U$ of $X$ and for every open spectral subspace $U'$ of $f^{-1}(U)$ the induced map $U' \to U$ is quasi-compact.
% \end{defi}
% 
% \begin{example}
% If $f$ is the underlying continuous map of a morphism of schemes, then $f$ is spectral.
% \end{example}
% 
% \begin{lem}
% A continuous map $f\colon X' \to X$ between spectral spaces is spectral if and only if $f$ is quasi-compact.
% \end{lem}
% 
% \begin{proof}
% Clearly a spectral map is quasi-compact. Conversely, assume that $f$ is quasi-compact and let $U$ and $U'$ be as above. By Hochster's therem $f$ is the underlying continuous map of a morphism $F$ of affine schemes. Its restriction $U' \to U$ is a morphism of open subschemes which are by definition quasi-compact and quasi-separated. Hence $U' \to U$ is quasi-compact by \cite{EGA}~I~(6.1.10)~(iii).
% \end{proof}
% \end{longversion}

%--------------------------------------------------------------

\subsection{Constructible topology}

\begin{defi}\label{DefConstructible}
Let $X$ be a spectral space.
\begin{ali}
\item
A subset $Z$ of $X$ is called \emph{constructible} if it is in the Boolean algebra $\Ccal$ of subsets of $X$ that is 
generated by all open quasi-compact subsets (i.e. $\Ccal$ is the smallest set of subsets of $X$ that contains all 
open quasi-compact subsets, is stable under finite intersections and under taking complements (this implies that $
\Ccal$ is also stable under finite unions)).
\item
A subset $Z$ of $X$ is called \emph{pro-constructible} (resp.~\emph{ind-constructible}) if it is an intersection 
(resp.~a union) of constructible subsets.
\end{ali}
If $X$ is a locally spectral space, a subset $Z$ of $X$ is called \emph{constructible} (resp.~\emph{pro-
constructible}, resp.~\emph{ind-constructible}) if there exists an open covering $(U_i)_i$ by spectral spaces such 
that $Z \cap U_i$ is a constructible (resp.~pro-constructible, resp.~ind-constructible) subset of $U_i$ for all $i \in I$.
\end{defi}

Using the quasi-compactness of spectral spaces, it is immediate that the second definition of constructible subsets in locally spectral spaces $X$ coincides with the first one if $X$ is spectral. 

%\begin{shortversion}
%\begin{rem}
%It can be shown that if a subset $Z$ of a locally spectral space is constructible (resp.~pro-constructible, resp.~ind-
%constructible), then $Z \cap U$ is constructible (resp.~pro-constructible, resp.~ind-constructible) for every open 
%subspace $U$ of $X$.
%\end{rem}
%\end{shortversion}

\begin{lem}\label{DescribeConstr}
Let $X$ be a spectral space. Then a subset $Z$ is constructible if and only if it is a finite union of subsets of the 
form $U \setminus V$, where $U$ and $V$ are quasi-compact open subsets of $X$.
\end{lem}

\begin{proof}
This is a special case of the following general principle: Let $X$ be any set and let $\Bcal$ be a Boolean algebra 
generated by a set $\Scal$ of subsets of $X$ such that $\Scal$ is stable under finite intersections and finite unions. 
Then $\Bcal$ is the set of finite unions of subsets of the form $U \setminus V$ with $U,V \in \Scal$.

Indeed, clearly every such set $U \setminus V$ is in $\Bcal$. So it suffices to show that the complement of any 
finite union of subsets $U \setminus V$ as above is again such a finite union. But if $C = \bigcup_{i\in I} U_i
\setminus V_i$ is a finite union, $U_i, V_i \in \Scal$, then, denoting complements by $-^c$,
\[
C^c = \bigcap_i U_i^c \cup V_i
\]
is the union of all sets of the form
\[
\bigl(\bigcap_{i\in I\setminus J} V_i\bigr) \cap \bigl(\bigcap_{i\in J} U_i^c\bigr) = \bigl(\bigcap_{i\in I\setminus J} V_i
\bigr) \setminus \bigl(\bigcup_{i\in J} U_i\bigr), \qquad J\subseteq I.
\]
As finite unions and finite intersections of subsets in $\Scal$ are again in $\Scal$, $C^c$ is again of the desired 
form.
\end{proof}

\begin{prop}\label{FiniteSpectral}
Let $X$ be a finite Kolmogorov space.
\begin{ali}
\item
$X$ is spectral.
\item
Every subset of $X$ is constructible.
\end{ali}
\end{prop}

\begin{proof}
\proofstep{(1)} Clearly, any open subset of $X$ is quasi-compact. Thus it remains to show that any closed 
irreducible subset $Z$ has a generic point. Assume that this is not the case, i.e. for all $z \in Z$, $U_z := Z 
\setminus \overline{\{z\}}$ is a non-empty open subset of $Z$ not containing $z$. As $Z$ is finite and irreducible, 
one has $\emptyset \ne \bigcap_z U_z$. But by definition, one clearly has $\bigcap_z U_z = \emptyset$. 
Contradiction.

\proofstep{(2)} Every subset of a finite Kolmogorov space is locally closed.
\end{proof}

\begin{cor}\label{OpenConstr}
Let $X$ be a spectral space.
\begin{ali}
\item
Every constructible subspace is quasi-compact.
\item
An open (resp.~a closed) subset $Z$ of $X$ is constructible if and only if $Z$ is quasi-compact (resp.~$X \setminus 
Z$ is quasi-compact).
\end{ali}
\end{cor}

\begin{proof}
Clearly~(1) implies~(2). To show~(1) it suffices by Lemma~\ref{DescribeConstr} to show that if $U$, $V$ are open 
quasi-compact, then $U \setminus V$ is quasi-compact. But $U \setminus V$ is closed in the quasi-compact subset 
$U$. Hence it is quasi-compact.
\end{proof}

% \begin{longversion}
% \begin{lem}\label{SpectralChevalley}
% Let $f\colon X' \to X$ be a quasi-compact map of spectral spaces. Then there exists a filtered projective system $(X'_i,p_{ij})_{i \in I}$ of spectral spaces and quasi-compact maps $f_i \colon X'_i \to X$ with $f_i \circ p_{ij} = f_j$ such that
% \begin{ali}
% \item
% $X' = \limproj_i X'_i$ and $f = f_i \circ p_i$ for all $i \in I$, where $p_i\colon X' \to X'_i$ is the canonical projection.
% \item
% For every constructible set $Z_i$ of $X'_i$ the image $f_i(Z_i)$ is a constructible subset of $X$.
% \item
% $f(X') = \bigcap_i f_i(X'_i)$.
% \end{ali}
% \end{lem}
% 
% \begin{proof}
% By Hochster's Theorem~\ref{ThmSpectralSpace} we may assume that $f$ is the map associated to a ring homomorphism $\varphi\colon A \to A'$. Moreover $A'$ is the filtered inductive of $A$-algebra $A'_i$ of finite presentation. Set $X'_i = \Spec A'_i$. Then $X'$ is the filtered projective limit of the $X'_i$. The second assertion then follows from Chevalley's theorem (\cite{EGA}~I~(7.1.4)) and the third is shown in~\cite{EGA}~I~(3.4.10).
% \end{proof}
% \end{longversion}
% 

\begin{rem}\label{RemConstTop}
Let $X$ be a locally spectral space.
\begin{ali}
\item
A subset $Z$ of $X$ is pro-constructible if and only if its complement is ind-constructible.
\item
Finite unions and arbitrary intersections of pro-constructible subsets are pro-constructible. Finite intersections and 
arbitrary unions of ind-constructible subsets are ind-constructible.
\begin{longversion}
\item
Every pro-constructible subset of $X$ is retro-compact in $X$. Conversely, a locally closed subset of $X$ is retro-compact if and only if it is pro-constructible.
\item
Let $(U_i)_i$ be an open covering of $X$. Then a subset $Z$ of $X$ is pro-constructible (resp.~ind-constructible) if 
and only if $Z \cap U_i$ is pro-constructible (resp.~ind-constructible) for all $i \in I$.
\end{longversion}
\end{ali}
\end{rem}

\begin{longversion}
\begin{proof}
The first two assertions are clear. For the other assertions we may assume by Hochster's Theorem~
\ref{ThmSpectralSpace} that $X$ is the underlying topological space of a scheme. Then all assertions follow from 
\cite{EGA}~I~(7.2.3)
\end{proof}
\end{longversion}

In particular we see that the ind-constructible (resp.~pro-constructible) subsets are the open (resp.~closed) subsets 
of a topology on $X$:

\begin{defi}\label{defConstTop}
Let $X$ be a locally spectral space. The topology on $X$ where the open subsets are the ind-constructible subsets 
is called the \emph{constructible topology}. $X$ endowed with the constructible topology is denoted by $X_{\rm 
cons}$.
\end{defi}

\begin{prop}\label{PropConstTop}
Let $X$ be a locally spectral space, $\Tcal$ its topology and $\Tcal_{\rm cons}$ its constructible topology.
\begin{ali}
\item
The topology $\Tcal_{\rm cons}$ is finer than $\Tcal$ (i.e., every closed subset of $X$ is pro-constructible).
\item
The constructible subsets are those subsets which are open and closed in the constructible topology (i.e., $Z 
\subseteq X$ is constructible if and only if $Z$ is pro-constructible and ind-constructible).
\item
For every open subspace $U$ of $X$, the topology of $\Tcal_{\rm cons}$ induced on $U$ is the same as the 
topology of $U_{\rm cons}$ (i.e. the inclusion is an open continuous map $U_{\rm cons} \to X_{\rm cons}$.
\item
$X_{\rm cons}$ is quasi-compact if and only if $X$ is quasi-compact.
\item
$X_{\rm cons}$ is Hausdorff if and only if $X$ is quasi-separated. In this case $X_{\rm cons}$ is locally compact 
and totally disconnected.
\end{ali}
\end{prop}

In particular if $X$ is a spectral space, then $X_{\rm cons}$ is a compact totally disconnected space.

\begin{shortversion}
\begin{proof}
Use Algebraic Geometry.
\end{proof}
\end{shortversion}

\begin{longversion}
\begin{proof}
By Hochster's Theorem~\ref{ThmSpectralSpace} we may assume that $X$ is the underlying topological space of a 
scheme. Then all assertions follow from \cite{EGA}~I~(7.2.12) and~(7.2.13).
\end{proof}
\end{longversion}

\begin{longversion}
\begin{cor}
Let $X$ be a quasi-compact locally spectral space, $(Z_i)_{i \in I}$ a family of pro-constructible subsets of $X$ and 
let $W$ be an ind-constructible subset of $X$ such that $\bigcap_{i \in I}Z_i \subseteq W$. Then there exists a finite 
subset $J$ of $I$ such that $\bigcap_{j \in J}Z_j \subseteq W$.
\end{cor}

\begin{proof}
The complement of $W$ is closed in the quasi-compact space $X_{\rm cons}$ and hence $W$ is quasi-compact in 
$X_{\rm cons}$. Thus already a finite number of complements of the $Z_i$ (which are open in $X_{\rm cons}$) 
suffice to cover the complement of $W$.
\end{proof}
\end{longversion}

% \begin{rem}
% Le $f\colon X' \to X$ be a quasi-compact map of locally spectral spaces. Then there exists an open covering $(U_i)_i$ of $X$ by spectral spaces and for all $i \in I$ a finite open covering $(U'_{ij})_j$ of $f^{-1}(U_i)$ by spectral spaces such that $f\rstr{U'_{ij}}\colon U'_{ij} to U_i$ is quasi-compact for all $i,j$.
% \end{rem}
% 
% \begin{proof}
% We may assume that $X$ is spectral. Then $X'$ is quasi-compact and a finite union of spectral spaces $U'_j$. If $U$ is an open quasi-compact subset, 
% \end{proof}
% 

\begin{longversion}
\begin{lem}
Let $X$ be a locally spectral space.
\begin{ali}
\item
Let $f\colon X' \to X$ be a qcqs map of locally spectral spaces. If $Z$ is a pro-constructible (resp.~ind-constructible) 
subset of $X$, then $f^{-1}(Z)$ is a pro-constructible (resp.~ind-constructible) subset of $X'$.

If $Z'$ is a pro-constructible subset of $X'$, then $f(Z')$ is pro-constructible.
\item
Assume that $X$ is spectral. Then conversely a subset $Z$ of $X$ is pro-constructible if and only if there exists a 
spectral space $X'$ and a quasi-compact map $f\colon X' \to X$ such that $Z = f(X')$.
\end{ali}
\end{lem}

\begin{proof}
To prove~(1) we may assume that $X$ is spectral. Then $X'$ is qcqs and hence spectral.

Thus for the proof of all assertions we may assume that $f\colon X' \to X$ is associated to a homomorphism of rings 
and the assertions are \cite{EGA}~I~(7.2.1) and \cite{EGA}~I~(7.2.3)~(vii).
\end{proof}
\end{longversion}

\begin{defi}\label{DefSpectralMap}
A continuous map $f\colon X \to Y$ of locally spectral spaces is called \emph{spectral} if for all open spectral 
subspaces $V \subseteq Y$ and $U \subseteq f^{-1}(V)$ the restriction $f\colon U \to V$ is quasi-compact.
\end{defi}

% \begin{longversion}
% \begin{rem}
% A quasi-compact continuous map $f\colon X \to Y$ of locally spectral spaces is spectral if and only if it is quasi-separated.
% \end{rem}
% 
% \begin{proof}
% Assume that $f$ is qcqs. Let $V \subseteq Y$ be a spectral open subspace. Then $f^{-1}(V)$ is an open spectral subspace of $X$, the restriction $f^{-1}(V) \to V$ is quasi-compact, and any open quasi-compact subset $U\subseteq f^{-1}(V)$ is retrocompact. Thus the restriction $U \to V$ is quasi-compact.
% 
% Conversely, assume that $f$ is spectral and quasi-compact. Let $V \subseteq Y$ be an open quasi-separated subset and let $(V_j)$ be an open cover of $V$ by spectral subspaces. To show that $f^{-1}(V)$ is quasi-separated, it suffices to show that the intersection of any two spectral subspaces $U_1$ and $U_2$ of $f^{-1}(V_j)$ is quasi-compact.
% \end{proof}
% \end{longversion}

\begin{prop}\label{CharSpectralMap}
Let $f\colon X \to Y$ be a continuous map of spectral spaces. Then the following assertions are equivalent.
\begin{eli}
\item
$f$ is spectral.
\item
$f\colon X_{\rm cons} \to Y_{\rm cons}$ is continuous.
\item
$f$ is quasi-compact.
\item
The inverse image of every constructible subset is constructible.
\item
The inverse image of every pro-constructible subset is pro-constructible.
\end{eli}
If these equivalent conditions are satisfied, then $f\colon X_{\rm cons} \to Y_{\rm cons}$ is proper. In particular the 
images of pro-constructible subsets are again pro-constructible.
\end{prop}

\begin{proof}
The implications ``(i) $\implies$ (iii) $\implies$ (iv) $\implies$ (v) $\iff$ (ii)'' are clear. Moreover,~(ii) implies~(iii) 
because the open quasi-compact subsets of $X$ are those subsets that are open in $X$ and open and closed in 
$X_{\rm cons}$.

Thus to show that all assertions are equivalent it suffices to show that (iii) implies (i). But if $f$ is quasi-compact, by 
Hochster's theorem we may assume that $f$ is associated to a ring homomorphism. Then the restriction of $f$ to 
open spectral subspaces $U$ and $V$ is the underlying topological map of a morphism of schemes between qcqs 
schemes and hence it is quasi-compact.

The last assertion is clear as any continuous map between compact spaces is proper (see below).
\end{proof}

Recall that a continuous map $f\colon X \to Y$ of topological spaces is called \emph{proper} (\cite{Bou_TG}
~Chap.~I~\S10) if it satisfies the following equivalent properties.
\begin{eli}
\item
For every topological space $Z$ the map $f \times \id_Z\colon X \times Z \to Y \times Z$ is closed.
\item
$f$ is closed and for all $y \in Y$ the fiber $f^{-1}(y)$ is quasi-compact.
\end{eli}

\begin{example}
Let $\varphi\colon A \to B$ be a homomorphism of rings. Then the corresponding scheme morphism $\Spec B \to 
\Spec A$ is quasi-compact and hence spectral.
\end{example}

%-------------------------------------------------------------

\subsection{Construction of spectral spaces}

\begin{lem}\label{CharSpectral}
Let $X$ be a quasi-compact Kolmogorov space which has a basis consisting of open quasi-compact subsets which 
is stable under finite intersections. Let $X'$ be the topological space with the same underlying set as $X$ and 
whose topology is generated by the open quasi-compact subsets of $X$ and their complements. Then the following 
assertions are equivalent.
\begin{eli}
\item
$X$ is spectral.
\item
$X'$ is compact with a basis of open and closed subspaces.
\item
$X'$ is quasi-compact.
\end{eli}
\end{lem}

\begin{proof}
Assertion~(i) implies~(ii) by Proposition~\ref{PropConstTop} and ``(ii) $\implies$ (iii)'' is clear. Let us show that (iii) implies (i). Let $Z$ be a closed irreducible subspace of $X$ and denote by $Z'$ the subspace of $X'$ with the same underlying set. As the topology of $X'$ is finer than the topology of $X$, $Z'$ is closed in $X'$ and hence quasi-compact.

Assume that $Z$ has no generic point. Then $\overline{\{z\}}$ is a proper subspace of $Z$ for all $z \in Z$. By hypothesis there exists a non-empty open quasi-compact subset $U_z$ of $Z$ that does not meet $\overline{\{z\}}$. In particular $Z = \bigcup_z (Z \setminus U_z)$. By hypothesis $Z' \setminus U_z$ is open in $Z'$ and as $Z'$ is quasi-compact, there exist finitely many of the sets $Z \setminus U_z$ that cover $Z$. This contradicts the irreducibility of $Z$.
\end{proof}

\begin{prop}\label{ProConstrSpectral}
Let $X$ be a spectral space and let $Z$ be a subspace. 
\begin{ali}
\item
$Z$ is pro-constructible if and only if $Z$ is spectral and the inclusion $Z \to X$ is spectral.
\item
Let $Z$ be a pro-constructible subspace of $X$. The closure of $Z$ is the set of specializations of points of $Z$.
\end{ali}
\end{prop}

\begin{proof}
\proofstep{(1)} The condition is sufficient by the last assertion of Proposition~\ref{CharSpectralMap}.
% By Remark~\ref{SoberPerm}, $Z$ is sober.
The inclusion $Z \to X$ is quasi-compact by Remark~\ref{RemConstTop}~(3). Thus $Z$ is a quasi-compact Kolmogorov 
space and admits a basis by open quasi-compact subsets stable under finite intersections. With the induced 
topology of $X_{\rm cons}$ the subspace $Z$ is compact. Therefore $Z$ is spectral by Lemma~\ref{CharSpectral}, 
and the inclusion $Z \to X$ is spectral by Proposition~\ref{CharSpectralMap}.

\proofstep{(2)} We have to show that if $x \in \overline{Z}$ then there exists $z \in Z$ with $x \in \overline{\{z\}}$. Set
\[
\Fcal := \set{U}{\text{$U$ quasi-compact open neigborhood of $x$}} \cup \{Z\}.
\]
Then any finite intersection of sets in $\Fcal$ is non-empty. Moreover all sets in $\Fcal$ are closed subsets of the 
compact space $X_{\rm cons}$. Thus the intersection of all quasi-compact open neigborhoods of $x$ in $X$ with $Z
$ contains a point $z$. Thus $x \in \overline{\{z\}}$.
\end{proof}

The following proposition will allow us to construct spectral spaces.

\begin{prop}\label{ConstructSpectral}
Let $X' = (X_0,\Tcal')$ be a quasi-compact space, let $\Ucal$ be a set of open and closed subspaces, let $\Tcal$ be the 
topology generated by $\Ucal$, and set $X = (X_0,\Tcal)$. Assume that $X$ is Kolmogorov. Then $X$ is a spectral 
space, $\Ucal$ is a basis of open quasi-compact subsets of $X$, and $X_{\rm cons} = X'$.
\end{prop}

\begin{proof}
We may replace $\Ucal$ by the set of finite intersections of sets in $\Ucal$. As the topology on $X$ is coarser than the topology of $X'$, every quasi-compact subset of $X'$ is also quasi-compact in $X$. Every set in $\Ucal$ is closed in $X'$ and thus quasi-compact in $X$. Thus $X$ is quasi-compact and $\Ucal$ is by construction a basis of open quasi-compact subsets of $X$ stable under finite intersection. Moreover $X$ is Kolmogorov by hypothesis.

By Lemma~\ref{CharSpectral} it remains to show that the topology $\Tcal''$ generated by the open quasi-compact subsets of $X$ and their complements is equal to $\Tcal'$. But $\Tcal''$ is Hausdorff (as $X$ is Kolmogorov there exists for $x \ne y$ an $U \in \Ucal$ with, say, $x \in U$ and $y \notin U$; then $U$ and $X' \setminus U$ are open sets of $\Tcal''$ separating $x$ and $y$) and coarser then $\Tcal'$. Hence $\id\colon X' \to X'' := (X_0,\Tcal'')$ is continuous and hence a homeomorphism because $X'$ is quasi-compact and $X''$ is Hausdorff.
\end{proof}

%============================================================

\section{Valuation spectra}

\subsection{The functor $A \sends \Spv A$}
\begin{definition}\label{DefValSpectrum}
Let $A$ be a ring. The \emph{valuation spectrum} $\Spv A$ is the set of all equivalence classes of valuations on $A
$ equipped with the topology generated by the subsets
\[
\Spv(A)(\frac{f}{s}) := \set{\abs \in \Spv A}{|f| \leq |s| \ne 0},
\]
where $f,s \in A$.
\end{definition}

Note for $s,t,f \in A$ we have $\Spv(A)(\frac{tf}{ts}) \subseteq \Spv(A)(\frac{f}{s})$, but in general these sets are not 
equal. For instance $\Spv(A)(\frac{1}{1}) = \Spv A$ but $\Spv(A)(\frac{s}{s}) = \set{v \in \Spv A}{v(s) \ne 0}$.

If $v$ is a point of $\Spv A$ we will often write $\abs_v$ instead of $v$ if we think of $v$ as an (equivalence class 
of an) absolute value on $A$.

%\begin{rem}
%Let $A$ be a ring and let $f \in A$. Then one has
%\begin{equation}\label{EqSpv}
%\begin{aligned}
%\set{v \in \Spv A}{v(f) \leq 1} &= \Spv(A)(\frac{f}{1}), \\
%\set{v \in \Spv A}{v(f) < 1} &= \bigcup_{, \\
%\set{v \in \Spv A}{v(f) > \gamma}, \qquad \set{v \in \Spv A}{v(f) \geq \gamma}
%\end{aligned}
%\end{equation}
%Indeed, for the first subset this holds by definition
%\end{rem}
%

\begin{example}
\begin{ali}
\item
Let $A = \QQ$ be the field of rational numbers. Then the only valuations on $\QQ$ are the $p$-adic valuations $
\abs_p$ for prime numbers $p$ and the trivial valuation $\abs_0$ (\ref{ValPID}). A set $\ne \emptyset$ is open if and 
only if it is the complement of a finite set of $p$-adic valuations. Hence $\Spv\QQ = \Spec\ZZ$.
\item
Let $A = \ZZ$ be the ring of integers. Then
\[
\Spv \ZZ = \Spv \QQ \cup \set{\abs_{0,p}}{\text{$p$ prime number}},
\]
where $\abs_{0,p}$ is induced by the trivial valuation on $\FF_p$. For all $p$ is $\abs_{0,p}$ a closed point ($
\Spv(\ZZ)(\frac{p}{p})$ is an open complement). Every open set which contains $\abs_{0,p}$ also contains $\abs_p$, 
i.e. $\abs_p$ is a generization of $\abs_{0,p}$. In fact $\cl{\abs_p} = \{\abs_p,\abs_{0,p}\}$.
\end{ali}
\end{example}

\begin{rem}\label{FunctorialitySpv}
Let $\varphi\colon A \to B$ be a ring homomorphism and denote by $\Spv(\varphi)\colon \Spv B \to \Spv A$ the map 
$\abs \sends \abs \circ \varphi$. For $f,s \in A$ one has
\begin{equation}\label{EqInvImageBasis}
\Spv(\varphi)^{-1}(\Spv(A)(\frac{f}{s})) = \Spv(B)(\frac{\varphi(f)}{\varphi(s)}).
\end{equation}
This shows that $\Spv(\varphi)$ is continuous. We will see in Proposition~\ref{SpvFunctor} that $\Spv(\varphi)$ is always quasi-compact.
\end{rem}

\begin{rem}
Let $A$ be a ring. Then~\eqref{EquivSupp} shows:
\begin{ali}
\item
Let $S \subset A$ be a multiplicative set and $\varphi$ be the canonical homomorphism $A \to S^{-1}A$. Then $
\Spv(\varphi)$ is a homeomorphism of $\Spv S^{-1}A$ onto the subspace $\set{v \in \Spv A}{\supp(v) \cap S = 
\emptyset}$.
\item
Let $\afr \subseteq A$ be an ideal and let $\varphi$ be the canonical homomorphism $A \to A/\afr$. Then $
\Spv(\varphi)$ is a homeomorphism of $\Spv A/\afr$ onto the subspace $\set{v \in \Spv A}{\supp(v) \supseteq \afr}$
\end{ali}
In other words, if $B = S^{-1}A$ or $B = A/\afr$, then
\[\xymatrix{
\Spv B \ar[r] \ar[d] &  \Spv A \ar[d] \\
\Spec B \ar[r] & \Spec A
}\]
is cartesian.
\end{rem}

\begin{rem}\label{SpvFieldIrred}
Let $K$ be a field. Then $\Spv K$ is irreducible and the trivial valuation is a generic point.

Indeed, let $v$ be a valuation on $K$. If for $f,s \in K$ one has $|f|_v \leq |s|_v \ne 0$, then $|f|_{\rm triv} \leq |s|
_{\rm triv} = 1$.
\end{rem}

\begin{rem}\label{SpvSpec}
Let $A$ be a ring. Consider $\Spv A \to \Spec A$, $x \sends \supp(x)$.
\begin{ali}
\item
For $s \in A$ let $D(s) = \set{\pfr \in \Spec A}{s \notin \pfr}$. The $D(s)$ form a basis of the topology of $\Spec A$. 
One has
\begin{equation}\label{EqSuppBasis}
\supp^{-1}(D(s)) = \set{v \in \Spv A}{|s|_v \ne 0} = (\Spv A)(\frac{0}{s}).
\end{equation}
In particular, $\supp$ is continuous.
\item
Its restriction to the subspace $T$ of trivial valuations on $A$ is a homeomorphism. Indeed, it is clearly bijective 
and $\supp(\Spv(A)(\frac{f}{s}) \cap T) = D(s)$.
\item
The canonical map $\iota\colon \Spv K(v) \to \Spv A$ is by Remark~\ref{FunctorialitySpv} a homeomorphism of $
\Spv K(v)$ onto $\set{w \in \Spv A}{\supp(w) = \supp(v)}$. This is an irreducible subspace by Remark~
\ref{SpvFieldIrred}.
\end{ali}
\end{rem}

\begin{prop}\label{SpvFunctor}
Let $A$ be a ring.
\begin{ali}
\item
The valuation spectrum $\Spv A$ is a spectral space. The sets of the form $\Spv(A)(\frac{f}{s})$ for $f,s \in A$ are open and quasi-compact. The Boolean algebra generated by them is the set of constructible sets, and this is also the Boolean algebra generated by the  sets $\set{v}{|f|_v \leq |g|_v}$ for $f,g \in A$.
\item
Any ring homomorphism $\varphi\colon A \to B$ induces a spectral map $\Spv(\varphi)\colon \Spv B \to \Spv A$ by 
$\abs \sends \abs \circ \varphi$. We obtain a contravariant functor from the category of rings to the category of 
spectral spaces and spectral maps.
\end{ali}
\end{prop}

\begin{proof}
We will use Proposition~\ref{ConstructSpectral} to show Assertion~(1).

\proofstep{(i)}
We first note that the Boolean algebra generated by the sets $\Spv(A)(\frac{f}{s})$ for $f,s \in A$ is the same as the 
Boolean algebra generated by the sets $\set{v}{|f|_v \leq |g|_v}$ for $f,g \in A$. Indeed, if $(\ )^c$ denotes the 
complement in $\Spv A$, then
\begin{align*}
\set{v}{|f|_v \leq |g|_v} &= \Spv(A)(\frac{f}{g}) \cup (\Spv(A)(\frac{0}{g}) \cup \Spv(A)(\frac{0}{f}))^c,\\
\Spv(A)(\frac{f}{s}) &= \set{v}{|f|_v \leq |s|_v} \cap (\set{v}{|s|_v \leq |0|_v})^c.\\
\end{align*}

\proofstep{(ii)}
Let $\Tcal$ be the topology defined above on $\Spv A$. For any two points in $\Spv A$ there exists a subset of the 
form $\Spv(A)(\frac{f}{s})$ for $f,s \in A$ containing one point but not the other (Proposition~\ref{ValEquiv}). Thus 
$X := (\Spv A, \Tcal)$ is Kolmogorov.

\proofstep{(iii)}
To the equivalence class $v$ of a valuation an $A$ we attach a binary relation $|_v$ on $A$ by $f |_v g$ iff $|f|_v 
\geq |g|_v$. We obtain a map $\rho\colon v \sends |_v$ of $\Spv A$ into the power set $\Pcal(A \times A)$ of $A 
\times A$. It is injective by Proposition~\ref{ValEquiv}. We identify $\Pcal(A \times A)$ with $\{0,1\}^{A \times A}$ 
and endow it with the product topology, where $\{0,1\}$ carries the discrete topology. Then $\Pcal(A \times A)$ is a 
compact space.

The image of $\rho$ consists of those binary relations $|$ such that for all $f,g,h \in A$ one has
\begin{dli}
\item
$f | g$ or $g | f$.
\item
If $f|g$ and $g|h$, then $f|h$.
\item
If $f|g$ and $f|h$, then $f|g+h$.
\item
If $f|g$, then $fh|gh$; and if $fh|gh$ and $0 \nmid h$, then $f | g$.
\item
$0 \nmid 1$.
\end{dli}
Each of these conditions (for fixed $f,g,h$) defines a closed subset of $\Pcal(A \times A)$: Let $\pi_{f,g}\colon 
\{0,1\}^{A \times A} \to \{0,1\}$ the projection onto the $(f,g)$-th component. Then for instance the elements of $
\{0,1\}^{A \times A}$ satisfying condition~(c) is the union of $\pi_{f,g+h}^{-1}(1)$ and of the complement of $\pi_{f,g}
^{-1}(1) \cap \pi_{f,h}^{-1}(1)$. Hence the image of $\rho$ is closed. Let $\Tcal'$ be the topology induced by $\Pcal(A 
\times A)$ on $\Spv A$. Then $X' := (\Spv A, \Tcal')$ is compact.

For $f,s \in A$ the subset $\Spv(A)(\frac{f}{s})$ is the intersection of $\Spv A$ with the open and closed subset of 
binary relations $|$ such that $s | f$ and $0 \nmid s$. Thus $\Spv(A)(\frac{f}{s})$ is open and closed in the compact 
space $X'$.

Using~(i) and~(ii), Assertion~(1) now follows from Proposition~\ref{ConstructSpectral}.

\proofstep{(iv)}
Assertion~(2) follows from Remark~\ref{FunctorialitySpv} and Proposition~\ref{CharSpectralMap}.
\end{proof}

\begin{corollary}\label{SuppQuasicompact}
The continuous map $\supp \colon \Spv A \to \Spec A$ is quasi-compact for every ring $A$.
\end{corollary}

\begin{proposition}\label{ExistValTensor}
Let $\varphi\colon R \to A$, $\psi\colon R \to B$ be ring homomorphisms and let $v \in \Spv A$, $w \in \Spv B$ with $v \circ \varphi = w \circ \psi =: u$. Then there exists $z \in \Spv(A \otimes_R B)$ such that its image in $\Spv A$ is $v$ and its image in $\Spv B$ is $w$.
\end{proposition}

%\begin{proof}
%We may replace $R$ by $K(u)$, $A$ by $K(v)$, and $B$ by $K(w)$ and assume that $R$, $A$, and $B$ are fields. It suffices a valuation $z$ on $A \otimes_R B/\mfr$, where $\mfr$ is a maximal ideal of $A \otimes_R B$, lying over $v$ and over $w$. 
%\end{proof}

%-----------------------------------------------------------------

\subsection{Specializations in valuation spectra}

\begin{rem}
Let $v,w \in \Spv A$ such that $v$ is a specialization of $w$. Then $\supp v \supseteq \supp w$.
\end{rem}

We start by considering specializations within the fiber of $\supp\colon \Spv A \to \Spec A$. Recall that $\supp^{-1}(\supp(v)) = \Spv(K(v))$. Hence we first look at specializations of valuations on a field.

\begin{prop}\label{SpecField}
Let $K$ be a field and let $v$ and $w$ be valuations on $K$. Then $v$ is a specialization of $w$ if and only if $A(v) \subseteq A(w)$.
\end{prop}

\begin{proof}
One has
\begin{align*}
A(v) \subseteq A(w) &\iff (\forall\ g \in A\colon v(g) \leq 1 \implies w(g) \leq 1)\\
&\iff (\forall\ f, 0 \ne s \in A\colon v(f) \leq v(s) \implies w(f) \leq w(s)) \\
&\iff (\forall\ f,s \in A\colon v \in \Spv(A)(\frac fs) \implies w \in \Spv(A)(\frac fs)).\qedhere
\end{align*}
\end{proof}

\begin{remdef}\label{DefVertical}
Let $v$ be a valuation on a ring $A$. A \emph{vertical} or \emph{secondary specialization} (resp.~\emph{vertical} or \emph{secondary generization}) of $v$ is a specialization (resp.~generization) $w$ of $v$ with $\supp(v) = \supp(w)$. Let $K(v) = \Frac A/(\supp v)$ and let $A(v) \subseteq K(v)$ the valuation ring of $v$.
\begin{ali}
\item
By Proposition~\ref{SpecField} one has homeomorphisms
\begin{equation}\label{VertGen}
\begin{aligned}
&\set{w \in \Spv A}{\text{$w$ vertical generization of $v$}} \\
\cong & \set{w \in \Spv K(v)}{\text{$w$ generization of $v$}}\\
\cong & \set{A(v) \subseteq B \subseteq K(v)}{\text{$B$ (valuation) ring}}\\
\cong & \set{H \subseteq \Gamma_v}{\text{$H$ convex subgroup of $\Gamma_v$}}.
\end{aligned}
\end{equation}
Here to $H$ corresponds the vertical generization $v/H$ of $v$ defined as follows
\[
v/H\colon A \to \Gamma_v/H \czero, \qquad
f \sends \begin{cases}
v(f) \mod{H},&\text{if $v(f) \ne 0$};\\
0,&\text{if $v(f) = 0$};
\end{cases}
\]
\item
By Remark~\ref{ValResidueField}, for the set of vertical specializations one has homeomorphisms
\begin{equation}\label{VertSpec}
\begin{aligned}
&\set{w \in \Spv A}{\text{$w$ vertical specialization of $v$}} \\
\cong & \set{w \in \Spv K(v)}{\text{$w$ specialization of $v$}}\\
\cong & \set{B \subseteq A(v) \subseteq K(v)}{\text{$B$ valuation ring of $K(v)$}}\\
\cong & \Spv(A(v)/\mfr(v)).
\end{aligned}
\end{equation}
\end{ali}
\end{remdef}

To consider specialization with changing support we need the following notion.

\begin{defi}
Let $v$ be a valuation on a ring $A$ and let $\Gamma_v$ its value group. Then the convex subgroup of $\Gamma_v$ 
generated by the elements $\Gamma_{v,\geq 1} \cap \im(v)$ is called the \emph{characteristic group of $v$}. It is 
denoted by $c\Gamma_v$.
\end{defi}

Hence $c\Gamma_v = 1$ if and only if $v(a) \leq 1$ for all $a \in A$.

\begin{example}
\begin{ali}
\item
Let $v$ be a valuation on a field. Then $\im(v) \setminus \{0\} = \Gamma_v$ and $c\Gamma_v = \Gamma_v$.
\item
Let $A$ be a valuation ring and let $v\colon A \to \Gamma_v \czero$ its valuation. Then $\Gamma_{v,\geq 1} \cap 
\im(v) = \{1\}$ and $c\Gamma_v = 1$.
\end{ali}
\end{example}

\begin{rem}
Let $v$ be a valuation on a ring $A$. Let $H \subseteq \Gamma_v$ be a subgroup. Define a map
\[
\abs_{v\rstr{H}}\colon A \to H \czero, \qquad  f \sends
\begin{cases}|f|_v,&\text{if $|f|_v \in H$};\\
0,&\text{if $|f|_v \notin H$}.\end{cases}
\]
Then $v\rstr{H}$ is a valuation if and only if $H$ is a convex subgroup and contains $c\Gamma_v$.

Indeed, one always has $v\rstr{H}(0) = 0$, $v\rstr{H}(1) = 1$. Assume that $v\rstr{H}$ is a valuation. Then $|fg|_{v
\rstr{H}} = |f|_{v\rstr{H}}|g|_{v\rstr{H}}$ is equivalent to $|f|_v|g|_v \in H \implies |f|_v,|g|_v \in H$, in particular $H$ is a 
convex subgroup. Moreover if there existed $f \in A$ with $|f|_v > 1$ and $|f|_v \notin H$, then $|f+1|_v = \max\{|f|
_v,1\} = |f|_v \notin H$ and hence $0 = |f+1|_{v\rstr{H}} = \max\{|f|_{v\rstr{H}},1\} = 1$, contradiction.

Conversely it is easy to see that if $H$ is convex and contains $c\Gamma_v$, then $v\rstr{H}$ is a valuation.
\end{rem}

\begin{remdef}\label{DefHorizontal}
Let $v$ be a valuation on a ring $A$ with value group $\Gamma_v$, let $\varphi\colon A \to K(v) = \Frac(A/\supp(v))
$ be the canonical homomorphism, and let $A(v) \subseteq K(v)$ be the valuation ring of $K(v)$. Let $H$ be a 
convex subgroup of $\Gamma_v$ and let $\pfr$ be the corresponding prime ideal of $A(v)$, i.e., $\pfr = \set{x \in 
A(v)}{\forall\,\delta \in H:|x|_v < \delta}$. Then $c\Gamma_v \subseteq H$ if and only if $\varphi(A) \subseteq 
A(v)_{\pfr}$. Let this be the case now.
\begin{ali}
\item
The valuation $v\rstr{H}$ is a specialization of $v$.

Indeed, if $|f|_{v\rstr{H}} \leq |s|_{v\rstr{H}} \ne 0$, then $|s|_v \in H$ and in particular $|s|_v \ne 0$. Assume $|f|_{v} > 
|s|_{v}$. Then $|f|_v \notin H$ and hence $|f|_v < 1$ because $H$ contains $c\Gamma_v$. Thus $|s|_v < |f|_v < 1$ 
and hence $|f|_v \in H$ because $H$ is convex. Contradiction.

We call $v\rstr{H}$ a \emph{primary} or a \emph{horizontal specialization of $v$} or $v$ a \emph{horizontal 
generalization of $v\rstr{H}$}.
\item
One has $\supp(v) \subseteq \supp(v\rstr{H}) = \varphi^{-1}(\pfr A(v)_{\pfr})$ with equality if and only if $\pfr = 0 \iff H 
= \Gamma_v$.
\item
$A(v)/\pfr$ is a valuation ring of $A(v)_{\pfr}/\pfr$, and $v\rstr{H}$ is the composition of the canonical ring 
homomorphism $\pi\colon A \to A(v)_{\pfr}/\pfr$ followed by the valuation given by $A(v)/\pfr$. In other words, one 
has $\ker(\pi) = \supp (v\rstr{H})$ and $\pi$ induces an extension of valued fields
\begin{equation}\label{ExtensionHorizontal}
(K(v\rstr{H}), A(v\rstr{H})) \mono (A(v)_\pfr/\pfr, A(v)/\pfr).
\end{equation}
\end{ali}
\end{remdef}

\begin{defi}
Let $A$ be a ring and let $v$ be a valuation on $A$. A subset $T$ of $A$ is called \emph{$v$-convex} if $t_1,t_2 \in 
T$, $s \in A$, and $v(t_1) \leq v(s) \leq v(t_2)$ imply $s \in T$.
\end{defi}

A subset $T$ containing $0$ is $v$-convex if and only if for $t \in T$ and $s \in A$ with $v(s) \leq v(t)$ one has $s 
\in T$.

If $v$ is a trivial valuation, the only $v$-convex sets $T$ with $0 \in T$ are $T = A$ and $T = \supp v$.

With this notion we can now describe the horizontal specializations of a given valuation.

\begin{prop}\label{SuppHorizontal}
Let $A$ be a ring and let $v$ be a valuation on $A$. Let $H$ be the set of $v$-convex prime ideals, ordered by 
reversed inclusion, and let $S$ be the set of horizontal specializations of $v$, ordered by being a specialization. Then $H$ and 
$S$ are totally ordered and the map
\[
\supp\colon S = \set{w}{\text{$w$ horizontal specialization of $v$}} \to H = \{\text{$v$-convex prime ideal of $A$}\}
\]
is a well-defined isomorphism of totally ordered sets.
\end{prop}

\begin{proof}
It is immediate that the support of a horizontal specialization of $v$ is $v$-convex and that $\supp\colon S \to H$ 
preserves the order. Moreover, every $v$-convex ideal clearly contains $\supp v$.

Let $\varphi\colon A \to K := \Frac(A/\supp(v))$ be the canonical homomorphism. By definition, the ordered set $S$ 
is isomorphic to the set $P$ of prime ideals $\pfr$ of $A(v)$ such that $\varphi(A) \subseteq A(v)_\pfr$, endowed 
with the reversed inclusion order. As the set of ideals in a valuation ring is totally ordered, $P$ and hence $S$ is 
totally ordered. Identifying $S$ with $P$, the map $\supp\colon S \to H$ becomes $\pfr \sends \varphi^{-1}(\pfr A_
\pfr)$. It is easy to check that an inverse map is given by sending a $v$-convex prime ideal $\qfr$ to the ideal of 
$A(v)$ generated by $\varphi(\qfr) \cup A(v)$.
\end{proof}

Our next goal is to show that every specialization is the horizontal specialization of a vertical specialization.

\begin{lem}\label{HorizontalVertical}
Let $w$ be a horizontal specialization of $v$.
\begin{ali}
\item\label{HorizontalVertical1}
Let $v'$ a vertical specialization of $v$. Then there exists a unique vertical specialization $w'$ of $w$ such that $w'$ 
is a horizontal specialization of $v'$.
\item\label{HorizontalVertical2}
Let $w'$ be a vertical specialization of $w$. Then there exists a vertical specialization $v'$ of $v$ such that $w'$ is a 
horizontal specialization of $v'$.
\end{ali}
\end{lem}

Indicating by $w \sends v$ if $v$ is a specialization, a visualization of this lemma is given by the following diagrams.
\[
\xymatrix{
v \ar@{|->}[r] \ar@{|->}[d] & w \ar@{|.>}[d] \\
v' \ar@{|.>}[r] & \exists! w'
}
\qquad
\xymatrix{
v \ar@{|->}[r] \ar@{|.>}[d] & w \ar@{|->}[d] \\
\exists\,v' \ar@{|.>}[r] & w'
}
\]

\begin{proof}
\proofstep{(1)}
As $\supp(w') = \supp(w)$ for every vertical specialization $w'$ of $w$, the uniqueness of $w'$ follows from Proposition~
\ref{SuppHorizontal}. Let $H'$ be the convex subgroup of $\Gamma_{v'}$ with $v = v'/H'$ and let $L$ be the convex 
subgroup of $\Gamma_v = \Gamma_{v'}/H'$ such that $w = v\rstr{L}$. Let $L'$ be the convex subgroup of $
\Gamma_{v'}$ containing $H'$ such that $L'/H' = L$. Then $c\Gamma_{v'} \subseteq L'$ and we set $w' = v'\rstr{L'}$.

\proofstep{(2)}
Consider the extension of valued fields $(K(w), A(w)) \mono (A(v)_\pfr/\pfr, A(v)/\pfr)$~\eqref{ExtensionHorizontal}. 
By Proposition~\ref{ExtendSpecializations} there exists a valuation ring $B$ of the field $A(v)_\pfr/\pfr$ with $B 
\subseteq A(v)/\pfr$ such that $B \cap K(w) = A(w')$. Let $v'$ be the valuation of $A$ with $\supp v' = \supp v$ such 
that $A(v')$ is the inverse image of $B$ in $A(v)$.
\end{proof}

\begin{cor}\label{ExistHorizontalGen}
Let $v$ be a valuation on a ring $A$ and let $\pfr$ be a prime ideal of $A$ with $\pfr \subseteq \supp v$. Then there 
exists a horizontal generization $w$ of $v$ such that $\supp w = \pfr$.
\end{cor}

\[\xymatrix{
\exists\,w \ar@{|.>}[r] \ar[d]^{\supp} & v \ar[d]^{\supp} \\
\pfr \ar@{|->}[r] & \supp(v) 
}\]

\begin{proof}
Let $R$ be the localization of $A/\pfr$ by the image of $\supp v$ in $A/\pfr$. This is a local ring whose field of 
fractions is $\Frac A/\pfr$. Then there exists a valuation ring $B$ of $\Frac A/\pfr$ dominating $R$ (Prop.~
\ref{CharValRing}). Let $w'$ be the correspnding valuation on $A$. Then $\supp w' = \pfr$ and we have a diagram
\[\xymatrix{
w' \ar@{|->}[r] & t_{\pfr} \ar@{|->}[d] \\
  & v
}\]
where $t_{\pfr}$ is the trivial valuation with support $\pfr$. By Lemma~\ref{HorizontalVertical}~\eqref{HorizontalVertical2} we now can fill in the lower left corner with a valuation $w$ as claimed.
%
%
%Let $\mfr$ be the maximal ideal of $B$. Then the canonical map $A \to R$ yields a field 
%extension $K(v) \mono B/\mfr$. Let $C$ be a valuation ring of $B/\mfr$ such that $C \cap K(v) = A(v)$ and let $w$ 
%be the valuation of $A$ with $\supp(w) = \pfr$ and $A(w) = \lambda^{-1}(C)$, where $\lambda\colon B \to B/\mfr$ is 
%the canonical homomorphism.
\end{proof}

\begin{prop}\label{Specialization}
Let $A$ be a ring, let $v$ be a valuation on $A$ and let $w$ be a specialization of $v$.
\begin{ali}
\item
Then $w$ is a horizontal specialization of a vertical specialization of $v$, i.e., there exist $v' \in \Spv A$ and convex 
subgroups $H, L \subseteq \Gamma_{v'}$ with $L \supseteq c\Gamma_{v'}$ such that $w = v'\rstr{L}$ and $v = v'/
H'$.
\item
Conversely, $w$ is a vertical specialization of a valuation $w'$ such that $w'$ is a horizontal specialization of $v$ or 
$c\Gamma_v = 1$ and $w'$ is a trivial valuation whose support contains $\supp(v\rstr{1})$.
\end{ali}
\end{prop}

Again we may visualize (1) by a diagram
\[\xymatrix{
v \ar@{|->}[dr] \ar@{|.>}[d]  \\
\exists\,v' \ar@{|.>}[r] & w
}\]

\begin{proof}
\proofstep{(2)}
We first consider the case $c\Gamma_v = 1$ and $v(a) \geq 1$ for all $a \in A \setminus \supp(w)$. Then $\supp(v
\rstr{1}) = \set{a \in A}{v(a) < 1} \subseteq \supp(w)$ and we may take for $w'$ the trivial valuation of $A$ with $
\supp(w') = \supp(w)$.

Thus we may assume that $c\Gamma_v \ne 1$ or $v(a) < 1$ for some $a \in A \setminus \supp(w)$. We first claim 
that for all $x,y \in A$ we have:
\begin{equation}\label{EqSpecProof1}
\text{$v(x) \leq v(y)$ and $w(x) \ne 0 = w(y)$} \implies v(x) = v(y) \ne 0.
\end{equation}
Indeed, $w(x) \ne 0 = w(y)$ implies $w \in \Spf(A)(\frac{y}{x})$ and hence $v \in \Spf(A)(\frac{y}{x})$ because $v$ is 
a generization of $w$. Thus $v(y) \leq v(x) \ne 0$. This shows the claim.

Next we claim that $\supp(w)$ is $v$-convex. Let $x,y \in A$ with $y \in \supp(w)$ and $v(x) \leq v(y)$. Assume that 
$x \notin \supp w$. Then we have
\begin{equation}\label{EqSpecProof2}
v(x) \leq v(y), \quad w(x) \ne 0 = w(y)
\end{equation}
and hence by~\eqref{EqSpecProof1} $v(x) = v(y) \ne 0$.

We now distinguish two cases. First we assume $c\Gamma_v \ne 1$. Then there exist $a \in A$ with $v(a) > 1$. 
Thus by~\eqref{EqSpecProof2} we obtain $v(x) \leq v(ay)$, $w(x) \ne 0 = w(ay)$ and hence $v(x) = v(ay)$ by~
\eqref{EqSpecProof1}. But this contradicts $v(x) = v(y)$. The second case is the existence of $a \in A \setminus 
\supp(w)$ with $v(a) > 0$. Then again applying~\eqref{EqSpecProof1} we obtain $v(ax) = v(y)$ which contradicts 
$v(x) = v(y)$. Hence we have shown that $\supp(w)$ is $v$-convex.

By Proposition~\ref{SuppHorizontal} there exists a horizontal specialization $u$ of $v$ such that $\supp u = \supp w$. 
By the definition of vertical specializations it suffices to prove that $w$ is a specialization of $u$. Let $f,s \in A$ with $w \in 
\Spv(A)(\frac{f}{s})$. Then $v \in \Spv(A)(\frac{f}{s})$ and hence $v(f) \leq v(s)$. As $u$ is a horizontal specialization 
of $v$, we find $u(f) \leq u(s)$. As $w(s) \ne 0$ and $\supp w = \supp u$, one has $u(s) \ne 0$. Hence $u \in \Spv(A)
(\frac{f}{s})$.

\proofstep{(1)}
If the specialization $w$ is a vertical specialization of a horizontal specialization of $v$, then we are done by 
Lemma~\ref{HorizontalVertical}~(2). Thus by~(2) and again using Lemma~\ref{HorizontalVertical}~(2) we may 
assume that $c\Gamma_v = 1$, that $w$ is trivial, and that $\supp(v\rstr{1}) \subseteq \supp(w)$. By Corollary~
\ref{ExistHorizontalGen} there exists a horizontal generization $u$ of $w$ such that $\supp u = \supp(v\rstr{1})$.

As $v\rstr{1}$ is trivial, $u$ is a vertical specialization of $v\rstr{1}$, which in turn is a horizontal specialization of $v
$. Thus again by Lemma~\ref{HorizontalVertical}~(2) there exists a vertical specialization $v'$ of $v$ such that $u$ 
is a horizontal specialization of $v'$. But then also $w$ is a horizontal specialization of $v'$.
\end{proof}

\subsection{The valuation spectrum of a scheme}

Let $X$ be a scheme. A \emph{valuation of $X$} is a pair $(x,v)$, where $x \in X$ and where $v$ is the equivalence class of a valuation on the residue field of $x$. The value group of $v$ is again denoted by $\Gamma_v$. The set of valuation on $X$ is denoted by $\Spv X$. We endow it with the topology generated by the sets of the form
\[
\set{(x,v) \in \Spv X}{x \in U, v(a(x)) \leq v(s(x)) \ne 0},
\]
where $U \subseteq X$ is open and $a,s \in \Oscr_X(U)$.

\begin{definition}\label{DefValSpecScheme}
The topological space $\Spv X$ is called \emph{valuation spectrum of $X$}.
\end{definition}

%----------------------------------------------------------------------------------------

% \begin{longversion}
% \subsection{The valuation spectrum of a finitely generated algebra over a field}
% \end{longversion}

%===========================================================

\section{Topological algebra}

\subsection{Topological groups}
\begin{definition}
A \emph{topological group} is a set $G$ endowed with the structure of a topological space and of a group such that 
the maps
\[
G \times G \to G, \quad (g,g') \sends gg', \qquad G \to G, \quad g \sends g^{-1}
\]
are continuous (where we endow $G \times G$ with the product topology).
\end{definition}

\begin{remark}
$G$ topological group, $a \in G$. Then left translation $g \sends ag$ and right tranlation $g \sends ga$ are homeomorphisms $G \to G$. In particular the topology of $G$ is uniquely determined by a fundamental system of neighborhoods of one element of $G$.
\end{remark}

The following trivial remark will be used very often.

\begin{remark}
Let $G$ be a topological group, let $H' \subseteq H \subseteq G$ be subgroups. If $H'$ is open in $G$, then $H$ is open in $G$.

Indeed, $H = \bigcup_{h\in H}hH'$.
\end{remark}

\begin{proposition}\label{SubGroupTop}
$G$ topological group, $H \subset G$ subgroup. We endow $G/H$ with the quotient topology (i.e. if $\pi\colon G \to G/H$ is the projection, $U \subseteq G/H$ is defined to be open if and only if $\pi^{-1}(U)$ is open).
\begin{assertionlist}
\item
The closure $\overline{H}$ is a subgroup of $G$.
\item
If $H$ is locally closed in $G$, then $H$ is closed in $G$.
\item
$H$ is open in $G$ (resp.~closed in $G$) if and only if $G/H$ is discrete (resp.~Hausdorff).
\item
If $H$ is a normal subgroup, $G/H$ is a topological group.
\end{assertionlist}
\end{proposition}

\begin{proof}
Let $a\colon G \times G \to G$ be the continuous map $(g,h) \sends gh^{-1}$. Then
\[
a(\overline{H} \times \overline{H}) = a(\overline{H \times H}) \subset \overline{a(H \times H)} = \overline{H}.
\]
This shows~(1).

To prove~(2) we may assume that $H$ is open and dense in $G$ (by replacing $G$ by the subgroup $\overline{H}$). 
Then for $g \in G$ the two cosets $gH$ and $H$ have nonempty intersection hence they are equal, i.e. $g \in H$.

(3): $G/H$ discrete $\iff$ $gH$ open in $G$ for all $g \in G$ $\iff$ $H$ open in $G$. If $G/H$ is Hausdorff, then $eH 
\in G/H$ is a closed point and its inverse image $H$ in $G$ is closed. Conversely, if $H$ is closed, then $H = HeH$ 
is a closed point in the quotient space $H\backslash G/H$. Hence its inverse image under the continuous map $G/H 
\times G/H \to H\backslash G/H$, $(g_1H,g_2H) \sends Hg_2^{-1}g_1H$ is closed. But this is the diagonal of $G/H 
\times G/H$.

(4) is clear.
\end{proof}

\begin{cor}\label{HausdorffTopGroup}
$G$ topological group, $e \in G$ neutral element. Equivalent:
\begin{equivlist}
\item
$G$ Hausdorff.
\item
$\{e\}$ is closed in $G$.
\item
$\{e\}$ is the intersection of all its neighborhoods.
\end{equivlist}
\end{cor}

\begin{proof}
The equivalence of~(i) and~(ii) is a special case of Proposition~\ref{SubGroupTop}~(3). ``(i) $\implies$ (iii)'' is 
obvious. Assume (iii) holds and let $g \ne e$. Then there exists a neighborhood $V$ of $e$ such that $g^{-1} \notin 
V$, i.e. $e \notin gV$. Hence $g$ is not in the closure of $\{e\}$.
\end{proof}

The following example of a topological group will be essential.

\begin{example}\label{TopBySubGroup}
Let $\Gamma$ be a non-empty partially ordered set which is right directed (i.e., for all $\gamma,\gamma' \in 
\Gamma$ there exists $\delta \in \Gamma$ with $\delta \geq \gamma$ and $\delta \geq \gamma'$). Let $G$ be a 
group and let $(G_\gamma)_{\gamma \in \Gamma}$ be a family of normal subgroups of $G$ such that $G_\delta 
\subseteq G_\gamma$ if $\delta \geq \gamma$. Equivalently, let $\Gamma$ be a set of normal subgroups such that 
for all $G_\gamma, G_{\gamma'} \in \Gamma$ there exists $G_{\delta} \in \Gamma$ such that $G_\delta \subseteq 
G_\gamma \cap G_{\gamma'}$.

Then there exists on $G$ a unique topology making $G$ into a topological group such that $(G_\gamma)_{\gamma 
\in \Gamma}$ is a fundamental system of neighborhoods of $e$ (this uses that all $G_\gamma$ are normal). 
Moreover:
\begin{assertionlist}
\item
By \ref{SubGroupTop} the $G_\gamma$ form a fundamental system of open and closed neighborhoods of $e$.
\item
$G$ is Hausdorff if and only if $\bigcap_\gamma G_\gamma = \{e\}$. In this case $G$ is totally disconnected 
by~(1).
\end{assertionlist}
\end{example}

Very often we will encounter the following special case.

\begin{example}\label{TopFilteredGroup}
Take in Example~\ref{TopBySubGroup} $\Gamma = \ZZ$, endowed with the standard order. Thus $G$ is given a 
descending filtration $\dots \supset G_n \supset G_{n+1} \supset \dots$ of normal subgroups $G_n$ for $n \in \ZZ$. 
As explained this defines on $G$ the structure of a topological group. Moreover:
\begin{assertionlist}
\item
We can endow $G$ with the structure of a pseudometric space as follows. For $g \in G$ set
\[
v(g) :=
\begin{cases}
-\infty,&\text{if $g \notin \bigcup_n G_n$};\\
n,&\text{if $g \in G_n \setminus G_{n+1}$};\\
\infty,&\text{if $g \in \bigcap_n G_n$}.
\end{cases}
\]
Fix a real number $\rho$ with $0 < \rho < 1$ and set
\[
d(g,h) := \rho^{v(gh^{-1})} \in \RR_{\geq 0} \cup \{\infty\}, \qquad g,h \in G
\]
Then clearly $d(g,g) = 0$, $d(g,h) = d(h,g)$ and
\begin{equation}\label{NonArch}
d(g,h) \leq \sup(d(g,k),d(h,k)), \qquad \text{for all $g,h,k \in G$}
\end{equation}
This a metric if and only if $\bigcap_n G_n = \{e\}$ and $\bigcup_n G_n = G$.
\item
Assume that $G$ is Hausdorff. Then the underlying topological space of $G$ is metrisable by a metric satisfying~\eqref{NonArch}. [Define the metric by composing $d$ with the map $\RR_{\geq 0} \cup \{\infty\} \to \RR_{\geq 0}$, 
$u \sends \inf(u,1)$.]
\end{assertionlist}
In the sequel we simply say that $G$ is a \emph{filtered group}, we call $(G_n)$ its \emph{filtration}, and we endow 
$G$ always with the topology defined above. A metric on $G$ satisfying~\eqref{NonArch} is called a \emph{non-
archimedean metric}.
\end{example}

\begin{rem}\label{MetrisableGroup}
In general is a topological group $G$ with unit $e$ metrisable if and only if $G$ is Hausdorff and $e$ has a countable fundamental system of neighborhoods (\cite{Bou_TG}~IX, \S3.1, Prop.~1). In this case there exist left-invariant metrics (i.e. $d(gh,gk) = d(h,k)$ for all $g,h,k \in G$) and right-invarint metrics (\cite{Bou_TG}~IX, \S3.1, Prop.~2).
\end{rem}

\begin{remdef}
Let $f\colon G \to H$ be a continuous homomorphism of topological groups. We endow $G/\ker(f)$ with the quotient topology induced by $G$ and $f(G)$ with the subspace topology induced by $H$. Then $f$ is called \emph{strict}, if the following equivalent conditions are satisfied.
\begin{equivlist}
\item
The bijective homomorphism $G/\ker(f) \to f(G)$ is a homeomorphism.
\item
The induced homomorphism $G \to f(G)$ is open.
\item
For all neighborhoods $U$ of $e$ in $G$ its image $f(U)$ is a nieghborhood of $e$ in $f(G)$.
\end{equivlist}
\end{remdef}

\begin{rem}
Let $f\colon G \to H$, $g\colon H \to K$ be strict homomorphisms of topological groups. If $\ker(g) \subseteq f(G)$, then it is not difficult to show that the composition $g \circ f$ is again strict.

But in general, $g \circ f$ is not strict, even if $f$ is injective and $g$ is surjective (e.g., let $x \in \RR$ be an irrational number, let $f$ be the inclusion $\ZZ x \mono \RR$ and let $g\colon \RR \to \RR/\ZZ$ be the canonical projection; use that $\ZZ + \ZZ x$ is dense in $\RR$). 
\end{rem}

%--------------------------------------------------------------------

\subsection{Topological rings}

If $A$ is a ring and $S,T \subseteq$ any subsets we denote by $S\cdot T$ the subgroup of $(A,+)$ generated by 
the elements $st$ for $s \in S$ and $t \in T$.

\begin{definition}
A \emph{topological ring} is a set $A$ endowed with the structure of a topological space and of a ring such that $(A,
+)$ is a topological group and the map
\[
A \times A \to A, \quad (a,a') \sends aa'
\]
is continuous.

Let $A$ be a topological ring. A topological $A$-module is an $A$-module $E$ endowed with a topology such that $
(E,+)$ is a topological group and such that the scalar multiplication $A \times E \to E$ is continuous.
\end{definition}

\begin{rem}
Let $A$ be a topological ring $A$.
\begin{ali}
\item
The homethety $A \to A$, $x \sends ax$ is continuous for all $a \in A$.
\item
The group of units $A^{\times}$ with the subspace topology is in general not a topological group (the map $x \sends 
x^{-1}$ is not necessarily continuous). Instead we endow $A^{\times}$ usually with the finest topology such that the 
two maps $A^{\times} \to A$, $x \sends x$ and $x \sends x^{-1}$ are continuous. Then $A^{\times}$ is a topological 
group.
\end{ali}
\end{rem}

\begin{defi}\label{DefTopField}
A \emph{topological field} is a field $K$ together with a topology that makes $K$ into a topological ring and such that 
$K^{\times} \to K^{\times}$, $x \sends x^{-1}$ is continuous (i.e. $K^{\times}$ with the subspace topology is a 
topological group).
\end{defi}

\begin{rem}
Let $A$ be a topological ring and let $I$ be an ideal. Similarly as in Proposition~\ref{SubGroupTop}~(1) one shows 
that its closure is also an ideal of $A$. In particular $\overline{\{0\}}$ is an ideal of $A$.
\end{rem}

\begin{rem}
Let $K$ be a topological field. As the closure of $\{0\}$ is an ideal, it must be either $\{0\}$ or $K$. Therefore if $K$ 
does not carry the chaotic topology (where the only open subsets are $\emptyset$ and $K$), then $K$ is Hausdorff.
\end{rem}

\begin{example}\label{FilteredRing}
Let $\Gamma$ be a partially ordered abelian group (written additively) whose order is right directed (e.g. if $\Gamma
$ is a totally ordered abelian group). Let $A$ be a ring and let $(A_\gamma)_{\gamma \in \Gamma}$ be a descending 
family of subgroups $A_\gamma$ of $(A,+)$ (descending means that $A_\delta \subseteq A_\gamma$ for $\delta 
\geq \gamma$) such that $A_\gamma A_\delta \subseteq A_{\gamma +\delta}$ for all $\gamma,\delta \in \Gamma$ 
and such that $1 \in A_0$.

Then $A_0$ is a subring of $A$ and $A_\gamma$ is an $A_0$-submodule of $A$. The set $B = \bigcup_{\gamma\in
\Gamma}A_\gamma$ is a subring of $A$. The set $\nfr = \bigcap_{\gamma} A_\gamma$ is an ideal of $B$ (for if $a 
\in A_\gamma$ and $x \in \nfr$ one has $x \in A_{\delta-\gamma}$ for all $\delta \in \Gamma$ and hence $ax \in A_
\gamma A_{\delta-\gamma} \subseteq A_\delta$ for all $\delta \in \Gamma$).

The ring $B$ together with the topology defined by the family of subgroups $(A_\gamma)_\gamma$ is a topological 
ring. Indeed, one has to show that the multiplication $B \times B \to B$ is continuous. Let $b_0,b'_0 \in B$ and 
hence $b_0 \in A_{\gamma}$ and $b'_0 \in A_{\gamma'}$ for some $\gamma,\gamma' \in \Gamma$. Let $\eps \in 
\Gamma$ be arbitrary and choose $\delta,\delta' \in \Gamma$ such that $\delta \geq \eps - \gamma'$, $\delta' \geq 
\eps - \gamma$, and $\delta' \geq \eps - \delta$. Then if $b,b' \in B$ with $b - b_0 \in A_\delta$ and $b' - b'_0 \in 
A_{\delta'}$ one has
\[
bb' - b_0b'_0 = (b-b_0)b'_0 + b_0(b' - b'_0) + (b - b_0)(b'-b'_0) \in A_{\delta+\gamma'} + A_{\gamma + \delta'} + 
A_{\delta + \delta'} \subseteq A_{\eps}. 
\]

Let $E$ be an $A$-module and endow $(E,+)$ with the structure of filtered topological group defined by the filtration 
$(A_\gamma E)_\gamma$. Then $\bigcup_\gamma A_\gamma E$ is a topological module over the topological ring 
$B$ (use the same proof).
\end{example}

We will also often encounter the following variant.

\begin{example}\label{LinTopRing}
Let $\Gamma$ be a set of ideals of $A$ such that for all $A_\gamma, A_{\gamma'} \in \Gamma$ there exists 
$A_{\delta} \in \Gamma$ such that $A_\delta \subseteq A_\gamma \cap A_{\gamma'}$. Then there exists on $A$ a 
unique topology that makes $A$ into a topological ring and such that the $A_\gamma$ form a fundamental system 
of neighborhoods of $0$ in $A$ (\cite{Bou_TG} III, \S6.3, Exemple~3). Such a topological ring $A$ is called 
\emph{linearly topologized} and $\Gamma$ is called a \emph{fundamental system of ideals}.

An ideal $I$ of $A$ is then open if and only if it contains some ideal $A_\gamma \in \Gamma$. 

If $(A, (A_\gamma))$ is a linearly topologized ring, then an ideal $I$ of $A$ is called \emph{ideal of definition} if it is 
open and for all $A_{\gamma}$ there exists an $n > 0$ such that $I^n \subseteq A_{\gamma}$. Note that an ideal of 
definition does not necessarily exist. Even if it exists, $I^n$ is not necessarily open.
\end{example}

A very important special case is the following.

\begin{definition}\label{AdicRing}
A topological ring $A$ is called \emph{adic} if there exists an ideal $I$ of $A$ such that $\set{I^n}{n \geq 0}$ is a 
fundamental system of neighborhoods of $0$ in $A$. The ideal $I$ is then an ideal of definition and this topology is 
called the \emph{$I$-adic topology}.

If $E$ is an $A$-module, the topology defined by $(I^nE)_n$ is called the \emph{$I$-adic topology on $E$}.
\end{definition}

\Warning Here we deviate from the terminology in EGA, where adic topological rings are always complete by 
definition.

The $I$-adic topology is also a special case of Example~\ref{FilteredRing} with $\Gamma = \ZZ$, $A_n = A$ for $n 
\leq 0$ and $A_n = I^n$ for $n > 0$.

\begin{rem}
Let $A$ be a ring and let $I$ and $J$ be ideals. Then the $J$-adic topology is finer than the $I$-adic topology if and 
only if there exists an integer $n > 0$ such that $J^n \subseteq I$.
\end{rem}

\begin{rem}
Let $A$ be a ring endowed with the $I$-adic topology for some ideal $I$ of $A$. Then $A$ is Hausdorff if and only if 
$\bigcap_{n \geq 0}I^n = 0$. If $A$ is noetherian, then by Krull's theorem
\[
\bigcap_{n \geq 0}I^n = \set{x \in A}{\exists\,a\in I\colon (1+a)x = 0}
\]
Thus the $I$-adic topology on a noetherian ring $A$ is Hausdorff if $I$ is contained in the Jacobson radical of $A$ or 
if $A$ has no zero divisors and $I \ne A$.
\end{rem}

The following notion yields also a topological ring.

\begin{defi}\label{DefNormedAlgebra}
Let $K$ be a field endowed with a valuation $\abs$ of height $1$. After replacing $\abs$ by an equivalent valuation 
we may assume that $\abs$ takes values in $\RR_{\geq0}$ (Proposition~\ref{Height1Groups}). A \emph{(non-
archimedean) normed $K$-algebra} is a unital $K$-algebra $A \ne 0$ together with a map $\norm\colon A \to 
\RR^{\geq 0}$, called \emph{norm}, such that
\begin{dli}
\item For $x \in A$ one has $\dvert x \dvert = 0$ if and only if $x = 0$.
\item $\dvert x + y \dvert \leq \max\{\dvert x \dvert,\dvert y \dvert\}$ for all $x,y \in A$. 
\item $\dvert \lambda x \dvert = |\lambda|\dvert x\dvert$ for all $x \in A$, $\lambda \in K$.
\item $\dvert xy\dvert \leq \dvert x\dvert\dvert y\dvert$ for all $x,y \in A$.
\item $\dvert 1 \dvert \leq 1$.
\end{dli}
\end{defi}

As usual the norm yields on $A$ the structure of a metric space and in particular that of a topological space. The 
properties of $\norm$ ensure that $A$ together with this topology is a topological ring.

Of course, the notion of a normed algebra is also used over fields with an archimedean absolute value. In this case 
one replaces condition~(b) by the usual triangle inequation.

\begin{rem}
Let $(A,\norm)$ be a normed algebra over a real-valued field $(K,\abs)$ with $A \ne 0$.
\begin{ali}
\item
Condition~(a) implies that $\dvert 1\dvert \ne 0$ and condition~(d) implies $\dvert 1 \dvert \leq \dvert 1 \dvert^2$ and 
hence $\dvert 1 \dvert = 1$ because of condition~(e). Thus $K \to A$, $\lambda \sends \lambda 1$ is an isometry.
\item
$A_0 := \set{a \in A}{\dvert a \dvert \leq 1}$ is an open subring of $A$. It is an algebra over the open subring 
$O_K := \set{\lambda \in K}{|\lambda|\leq 1}$ of $K$.
\end{ali}
\end{rem}

All of the above examples are a special case of the following type of ring.

\begin{defi}\label{NonArchRing}
A topological ring $A$ is called \emph{non-archimedean} if $A$ has a fundamental system of neighborhoods of $0$ 
consisting of subgroups of $(A,+)$.
\end{defi}

\begin{rem}\label{DefineNATop}
Let $A$ be a ring and let $\Gcal$ be a set of subgroups of $(A,+)$. Then $\Gcal$ is a fundamental system of neighborhoods of $0$ of a (unique) topology that makes $A$ into a topological ring if and only if $\Gcal$ satisfies the following properties.
\begin{dli}
\item
For all $G,G' \in \Gcal$ there exists $H \in \Gcal$ with $H \subseteq G \cap G'$.
\item
For all $x \in A$ and all $G \in \Gcal$ there exists $H \in \Gcal$ with $xH \subseteq G$.
\item
For all $G \in \Gcal$ there exists $H \in \Gcal$ with $H\cdot H \subseteq G$.
\end{dli}
\end{rem}

%-------------------------------------------------------------------

\subsection{Bounded sets and topologically nilpotent elements}

For every subset $T$ of a ring and for $n \in \NN$ we set
\begin{equation}\label{DefTn}
T(n) := \set{t_1t_2\cdots t_n}{t_i \in T}.
\end{equation}

\begin{defi}
Let $A$ be a topological ring. An element $x \in A$ is called \emph{topologically nilpotent} if $0$ is a limit of the 
sequence $(x^n)_{n \in \NN}$. We denote by $A^{oo}$ the set of topologically nilpotent elements of $A$.

More generally, we call a subset $T$ of $A$ \emph{topologically nilpotent} if there exists for every neighborhood $U$ 
of $0$ an $N \in \NN$ such that $T(n) \subseteq U$ for all $n \geq N$.
\end{defi}

Clearly $x \in A$ is topologically nilpotent if and only if $\{x\}$ is topologically nilpotent.

\begin{rem}\label{TopNilpotentLinTop}
Let $A$ be a linearly topologized ring. Then $x \in A$ is topologically nilpotent if and only if for every open ideal $J$ 
of $A$ the image of $x$ in $A/J$ is nilpotent. Therefore the set $A^{oo}$ of topologically nilpotent elements is an 
ideal of $A$.

Assume now that $A$ has an ideal of definition $I$. Then the following assertions are equivalent for $x \in A$.
\begin{eli}
\item
$x$ is topologically nilpotent.
\item
The image in $A/I$ is nilpotent.
\item
$x$ is contained in an ideal of definition.
\end{eli}
Indeed, clearly one has ``(iii) $\implies$ (i) $\implies$~(ii)''. Let $x \in A$ with $x^m \in I$. Then $I + Ax$ is open and 
$(I + Ax)^m \subseteq I$. Thus $I + Ax$ is an ideal of definition.

In particular we see that $A^{oo}$ is the inverse image of the nilradical of $A/I$ and hence it is an open ideal of $A$. 
Moreover, $A^{oo}$ is itself an ideal of definition (and then clearly the largest ideal of definition) if and only if the 
nilradical of $A/I$ is nilpotent. This is automatic if the nilradical of $A/I$ is finitely generated (e.g., if $A/I$ is 
noetherian).
\end{rem}

\begin{defi}\label{DefBounded}
Let $A$ be a topological ring. A subset $B$ of $A$ is called \emph{bounded} if for every neighborhood $U$ of $0$ in 
$A$ there exists an open neighborhood $V$ of $0$ in $A$ such that $vb \in U$ for all $v \in V$ and $b \in B$.

An element $x \in A$ is called \emph{power-bounded} if the set $\set{x^n}{n \geq 1}$ is bounded. The set of power-
bounded elements is denoted by $A^0$.

More generally, a subset $T$ of $A$ is called \emph{power-bounded} if $\bigcup_{n\in \NN}T(n)$ is bounded.
\end{defi}

\begin{rem}\label{PropBounded}
Let $A$ be a topological ring.
\begin{ali}
\item
Clearly, every finite subset of $A$ is bounded.
\item
Every subset of a (power-)bounded subset is (power-)bounded.
\item
Every finite union of (power-)bounded subsets is (power-)bounded. A finite union of bounded and topologically nilpotent subsets is again topologically nilpotent.
\item
A bounded and topologically nilpotent subset of $A$ is power-bounded.
\item
If $T_1$ is a power-bounded subset of $A$ and $T_2$ a topologically nilpotent subset of $A$. Then $\set{t_1t_2}{t_1 \in T_1, t_2 \in T_2}$ is topologically nilpotent.
\end{ali}
\end{rem}

\begin{example}
\begin{ali}
\item
Let $A = \CC$. Then a subset $B$ of $\CC$ is bounded if and only if there exists $C > 0$ such that $|z| < C$ for all 
$z \in \CC$. One has $A^o = \set{z \in \CC}{|z| \leq 1}$ and $A^{oo} = \set{z \in \CC}{|z|<1}$.
\item
Let $A$ be a ring and $\abs$ a non-trivial valuation on $A$.

Assume that $\abs$ takes values in $\RR^{>} \czero$ (i.e. $\abs$ is of height $1$). Endow $A$ with the structure of 
a topological ring such that the subgroups $\set{x \in A}{|x|<\eps}$ for $\eps \in \RR^{>0}$ form a fundamental 
system of neighborhoods of $0$ in $A$ (Example~\ref{FilteredRing}). Then again a subset $B$ of $A$ of set is 
bounded if and only if there exists $C > 0$ such that $|x| < C$ for all $x \in A$. One has $A^o = \set{x \in A}{|x| \leq 
1}$ and $A^{oo} = \set{x \in A}{|x|<1}$.

If $\abs$ has height $\geq 2$, then there exists a non-trivial convex subgroup $\Delta$ of $\Gamma_{\abs}$. If $x \in 
A$ is an element with $|x| \in \Delta$ then $x$ is never topologically nilpotent even if $|x| < 1$.
\end{ali}
\end{example}

\begin{prop}\label{GroupBounded}
Let $A$ be a non-archimedean topological ring.
\begin{ali}
\item
Let $T$ be a subset of $A$ and let $T'$ be the subgroup of $A$  generated by $T$. Then $T'$ is bounded (resp.~power-bounded, resp.~topologcially nilpotent) if and only if $T$ has this property.
\item
A subset $T$ of $A$ is power-bounded if and only if the subring of $A$ generated by $T$ is bounded.
\item
The union of all bounded subrings is $A^o$, and this is a subring of $A$.
\item
$A^o$ is integrally closed in $A$, and $A^{oo}$ is a radical ideal of $A^o$.
\end{ali}
\end{prop}

\begin{proof}
(1) follows from the definitions. To show (2), we may assume that $1 \in T$. Then the subring generated by $T$ is the subgroup of $A$ generated by $X := \bigcup_{n\in\NN} T(n)$. As $X$ is bounded, (2) follows from (1).

If $A_{1}$ and $A_{2}$ are bounded subrings, then by (1) the subring generated by $A_{1} \cup A_{2}$ is bounded. Thus the union $B$ of all bounded subrings is a subring. Now it follows from (2) that $B = A^o$.

It remains to show (4). Let $a \in A$ be integral over $A^o$. By (3), $a$ is integral over a bounded subring $B$ of $A$. Therefore there exists an integer $N$ such that $B[a] = B + Ba + \dots + Ba^{N-1}$. Thus $B[a]$ is bounded and hence $a$ is power-bounded. For $a \in A^{oo}$ the set $\{a\}$ is bounded and topologically nilpotent, hence it is power-bounded. Thus $A^{oo} \subseteq A^o$. By (1) it is a subgroup of $A^o$ and by Remark \ref{PropBounded} (5) it is an ideal of $A^o$. Let $a \in A$ such that $a^m \in A^{oo}$. We want to show that $a \in A^{oo}$. Let $U$ be a neighborhood of $0$ in $A$. Choose a neighborhood $V$ of $0$ such that $a^iv \subseteq U$ for all $i = 0,\dots,m-1$ and all $v \in V$. Choose $N \in\NN$ such that $(a^m)^n \in V$ for all $n \geq N$. Then $a^n \in U$ for all $n \geq mN$.
\end{proof}

%-----------------------------------------------------------------

\subsection{Completion of abelian topological groups and of topological rings}

\begin{defi}
\begin{ali}
\item
Let $X$ be a set. A \emph{filter basis on $X$} is a non-empty set $\Bcal$ of subsets of $X$ such that $\emptyset 
\notin \Bcal$ and such that for all $A_1, A_2 \in \Bcal$ there exists $A \in \Bcal$ with $A \subseteq A_1 \cap A_2$.

Example: Let $(x_n)_n$ be a sequence in $X$. Then $\set{\set{x_k}{k \geq n}}{n \in \NN}$ is a filter basis on $X$, 
called the \emph{filter basis of $(x_n)$}.
\item
Let $X$ be a topological space, $\Bcal$ a filter basis on $X$, and $x \in X$. Then \emph{$\Bcal$ converges to $x$} 
if for every open neighborhood $U$ of $x$ there exists $A \in \Bcal$ such that $A \subseteq U$.

Example: The filter basis of a sequence $(x_n)$ converges to $x$ if and only if $(x_n)$ converges to $x$.
\item
Let $G$ be an abelian topological group (written additively). A filter basis $\Bcal$ is called \emph{Cauchy} if for every open neighborhood $U$ of $0$ there exists $A \in \Bcal$ such that $x-y \in U$ for all $x,y \in A$.
\item
An abelian topological group $G$ is called \emph{complete} if $G$ is Hausdorff and every Cauchy filter basis converges.
\item
A topological ring $A$ is called \emph{complete} if $(A,+)$ is a complete topological group.
\end{ali}
\end{defi}

\Warning Here we deviate from Bourbaki's terminology, where complete topological groups are not necessarily Hausdorff.

\begin{propdef}\label{CompletionGroup}
Let $G$ be an abelian topological group (resp.~topological ring). Then there exists a complete abelian topological group $\Ghat$ (resp.~complete topological ring) and a continuous homomorphism $\iota\colon G \to \Ghat$ satisfying the following universal property. For every complete abelian topological group (resp.~complete topological ring) $H$ and every continuous homomorphism $\varphi\colon G \to H$ there exists a unique continuous homomorphism $\hat\varphi\colon \Ghat \to H$ such that $\hat\varphi \circ \iota = \varphi$. Clearly $(\Ghat,\iota)$ is unique up to unique isomorphism. Moreover:
\begin{ali}
\item
$\ker(\iota) = \overline{\{0\}} =: N$, $\im(\iota)$ is dense in $\Ghat$.
\item
If $\Ucal$ is a basis of neighborhoods of $0$ in $G$, then $\set{\text{$\overline{\iota(U)}$ closure in $\Ghat$}}{U \in 
\Ucal}$ is a basis of neighborhoods of $0$ in $\Ghat$.
\item
The formation of $\Ghat$ is functorial in $G$.
\end{ali}
\end{propdef}

\begin{proof}
\cite{Bou_TG} III, \S3.4.
\end{proof}

\begin{example}\label{CompletionSubgroup}
Let $G$ be an abelian topological group, where the topology is given by a family $(G_\gamma)$ of open subgroups of $G$ as in Example~\ref{TopBySubGroup}. Then $G/G_\gamma$ is discrete and the canonical homomorphism $G \to G/G_\gamma$ yields by the universal property of the completion (and of the projective limit) a continuous homomorphism $\Ghat \to \limproj_\gamma G/G_\gamma$. This is an isomorphism of topological groups and the canonical homomorphism $\iota\colon G \to \Ghat$ induces a homeomorphism $G/\overline{\{0\}} \iso \iota(G)$ (\cite{Bou_TG}~III, \S7.2). Moreover:
\begin{ali}
\item
The completion $\Ghat_\gamma$ of $G_\gamma$ is identified with the closure of $\iota(G_\gamma)$ (\cite{Bou_TG}
~II, \S3.9, Cor.~1 to Prop.~18). As $G_\gamma$ is closed in $G$, one has $G_\gamma = \iota^{-1}(\Ghat_\gamma) 
= \iota^{-1}(\Ghat_\gamma \cap \iota(G))$.
\item
The $\Ghat_\gamma$ form a fundamental system of neighborhoods of $0$ in $\Ghat$ and hence are open 
subgroups of $\Ghat$.
\item
As $\iota(G)$ is dense in $\Ghat$ and $\Ghat_\gamma$ is open, $\Ghat = \iota(G)\Ghat_\gamma$ for all $\gamma$. 
In particular one has $\bigcup_\gamma G_\gamma = G$ if and only if $\Ghat = \bigcup_\gamma \Ghat_\gamma$.
\end{ali}
\end{example}

\begin{proposition}\label{CharComplete}
Let $G$ be a filtered abelian group (written additively), $(G_n)_{n\in\ZZ}$ its filtration. Assume that $G$ is Hausdorff. 
Equivalent:
\begin{equivlist}
\item
$G$ is complete.
\item
Every Cauchy sequence in $G$ converges.
\item
For every sequence $(x_m)$ in $G$ converging to $0$ the series $\sum_m x_m$ converges.
\item
The canonical group homomorphism $G \to \limproj_n G/G_n$ is bijective.
\item
The canonical group homomorphism $G \to \limproj_n G/G_n$ is a homeomorphism if we endow $\limproj_n G/G_n$ 
with the projective limit topology.
\end{equivlist}
\end{proposition}

\begin{proof}
The equivalence of~(i) and~(ii) follows from the fact that $G$ is metrizable.

(ii) $\implies$ (iii): For every $n \geq 1$ there exists $M \geq 1$ such that $x_m \in G_n$ for all $m \geq M$. Hence 
$\sum_{m \in I} x_i \in G_n$ for all finite subsets $I$ of $\NN_{\geq M}$. This shows the convergence of the series.

(iii) $\implies$ (ii): Standard exercise in non-archimedean geometry.

The equivalence of (i),~(iv), and~(v) is \cite{Bou_TG} III, \S7.3.
\end{proof}

\begin{prop}\label{QuotComplete}
Let $G$ be a topological group whose topology is metrisable, and let $H \subseteq G$ be a closed subgroup. Then $G/H$ is metrisable. If $G$ is complete, then $G/H$ is complete.
\end{prop}

\begin{proof}
\cite{Bou_TG}~IX, \S3.1, Prop.~4.
\end{proof}

\begin{example}\label{CompletionIdeal}
Now assume that $A$ is a topological ring whose topology is defined by a family $(A_\gamma)$ of subgroups as in 
Example~\ref{FilteredRing} such that $\bigcup_\gamma A_\gamma = A$. Let $\iota\colon A \to \Ahat$ be the 
completion of $(A,+)$ as a topological group. Then there exists a unique multiplication on $\Ahat$ that makes $
\Ahat$ into a ring and such that $\iota$ is a homomorphism of rings (\cite{Bou_TG}~III, \S6.5).

For $\gamma,\gamma' \in \Gamma$ one has $\iota(A_\gamma)\iota(A_{\gamma}) \subseteq \iota(A_{\gamma+
\gamma'})$ and hence by continuity $\Ahat_\gamma\Ahat_{\gamma'} \subseteq \Ahat_{\gamma+\gamma'}$.
\end{example}

\begin{prop}\label{CompleteAdic}
Let $A$ be a ring, $I$ an ideal, let $\Ahat$ be its $I$-adic completion, and $\iota\colon A \to \Ahat$ the canonical 
homomorphism. For all $n > 0$ endow $I^n$ with the $I$-adic topology and consider their completion as ideals of $
\Ahat$. Assume that $I$ is finitely generated.
\begin{ali}
\item
For all $n > m \geq 0$ the canonical homomorphism $I^m/I^n \to \Ihat^m/\Ihat^n$ is an isomorphism.
\item
For all $n > 0$ one has $\widehat{I^n} = \Ihat^n = i(I^n)\Ahat$.
\item
The topology on $\Ahat$ is the $\Ihat$-adic topology.
\end{ali}
\end{prop} 

\begin{proof}
\cite{Bou_AC} III \S2.12
\end{proof}

\begin{proposition}\label{DefIdealMaximal}
Let $A$ be a complete non-archimedean ring. Let $\afr \subseteq A$ be an ideal such that every element of $\afr$ is topologically nilpotent. Then an element $x$ of $A$ is invertible if and only if its image in $A/\afr$ is invertible. In particular $\afr$ is contained in the Jacobson radical of $A$.
\end{proposition}

%----------------------------------------------

\subsection{The topology defined by a valuation}

\begin{propdef}
Let $A$ be a ring and let $v$ be a valution on $A$. For $\gamma \in \Gamma_v$ we set $A_\gamma := \set{a \in A}
{v(a) < \gamma}$. Then $(A_\gamma)_\gamma$ satisfies the conditions of Example~\ref{FilteredRing} and $
\bigcup_\gamma A_\gamma = A$. Thus $A$ becomes a topologcial ring with $(A_\gamma)$ as a fundamental 
system of neighborhoods of $0$. This topology is called the \emph{topology defined by $v$}. Moreover:
\begin{ali}
\item
If $\supp v = 0$, then $A$ is Hausdorff and totally disconnected.
\item
The group of units $A^{\times}$ endowed with the subspace topology is a topological group.
\end{ali}
\end{propdef}

\begin{proof}
One has $\bigcap_\gamma A_\gamma = \supp v$. This shows~(1). To show~(2) we have to prove that $A^{\times} 
\to A^{\times}$, $x \sends x^{-1}$ is continuous. Let $x_0 \in A^{\times}$. It suffices to show that for all $\gamma \in 
\Gamma_v$ and $x \in A^{\times}$ one has
\[
|x - x_0|_v < \min\{\gamma|x_0|_v^2,|x_0|_v\} \implies |x^{-1} - x_0^{-1}|_v < \gamma.
\]
But $x^{-1} - x_0^{-1} = x^{-1}(x_0 - x)x_0^{-1}$ and hence $|x^{-1} - x_0^{-1}|_v = |x - x_0|_v|x|_v^{-1}|x_0|_v^{-1}$. If 
$|x -x_0|_v < |x_0|_v$, one has $|x|_v = |x_0|_v$. If $|x -x_0|_v < \gamma|x_0|_v^2$ then
\[
|x^{-1} - x_0^{-1}|_v = |x - x_0|_v|x_0|_v^{-2} < \gamma.\qedhere
\]
\end{proof}

% \begin{rem}
% Let $A$ be a ring endowed with the topology defined by a valuation $v$. For all $\gamma \in \Gamma_v$, the sets $\set{a \in A}{v(a) \leq \gamma}$ and $\set{a \in A}{v(a) = \gamma}$ are open and closed in $A$.
% 
% Indeed, by definition $\set{a \in A}{v(a) < \gamma}$ is an open subgroup, hence it is closed. Thus its complement $\set{a \in A}{v(a) \geq \gamma}$ is open
% \end{rem}
% 
% 
\begin{cor}
Let $v$ be a valuation on a field $K$. Then the topology defined by $v$ on $K$ makes $K$ into a topological field.
\end{cor}

\begin{rem}
Let $v$ be a non-trivial valuation on a ring $A$ with value group $\Gamma$. Then $\set{a \in A}{v(a) \leq \gamma}$, 
$\set{a \in A}{v(a) \geq \gamma}$, $\set{a \in A}{v(a) > \gamma}$ are open and closed in the topology defined by $v
$.

Indeed, by definition $\set{a \in A}{v(a) < \gamma}$ is an open subgroup of the subgroup $\set{a \in A}{v(a) \leq 
\gamma}$. Thus $\set{a \in A}{v(a) \leq \gamma}$ is open and (being a subgroup) closed. The other sets are open 
and closed by taking complements.
\end{rem}

\begin{rem}\label{CompleteValField}
Let $K$ be a field and let $v$ be a valuation on $K$ with value group $\Gamma$. Endow $\Gamma$ with the 
discrete topology. Let $\Khat$ the completion ring. This is a topological field (\cite{Bou_AC}~VI, \S5.3, Prop.~5).

The homomorphism $v\colon K^{\times} \to \Gamma$ is continuous and thus can be extended uniquely to a 
continuous homomorphism $\vhat\colon \Khat^{\times} \to \Gamma$ and by continuity $\vhat$ defines a valuation 
on $\Khat$. For $\gamma \in \Gamma$ the closures of $\set{x \in K}{v(x) < \gamma}$ and $\set{x \in K}{v(x) \leq 
\gamma}$ in $\Khat$ are $\set{x \in \Khat}{\vhat(x) < \gamma}$ and $\set{x \in \Khat}{\vhat(x) \leq \gamma}$, 
respectively. In particular, the topology on $\Khat$ is defined by the valuation $\vhat$. Moreover, one obtains for its 
valuation ring $A(\vhat) = \widehat{A(v)}$ and for its maximal ideal $\mfr(\vhat) = \widehat{\mfr(v)}$.

Moreover $A(\vhat) = A + \mfr(\vhat)$ and in particular $\kappa(\vhat) = \kappa(v)$.
\end{rem}

\begin{defi}
Two valuations $v_1$ and $v_2$ on a field $K$ are called \emph{independent} if $K$ is the ring generated by 
$A(v_1)$ and $A(v_2)$. Otherwise they are called \emph{dependent}.
\end{defi}

The trivial valuation is independent of every other valuation. Clearly two equivalent non-trivial valuations are 
dependent. More generally, if $A(v_1) \subseteq A(v_2) \ne K$, then $v_1$ and $v_2$ are dependent.

\begin{prop}
Two valuations of height $1$ are equivalent if and only if they are dependent.
\end{prop}

\begin{proof}
\cite{Bou_AC}~VI, \S4.5, Prop.~6~(c).
\end{proof}

There are examples of dependent non-equivalent valuations which have both height $2$ (\cite{Bou_AC} VI, \S8, 
Exercise~1).

\begin{prop}
Two non-trivial valuations $v_1$ and $v_2$ on a field $K$ are dependent if and only if they define the same topology 
on $K$.
\end{prop}

\begin{proof}
\cite{Bou_AC} VI, \S7.2, Prop.~3
\end{proof}

\begin{remdef}\label{DefMicrobial}
Let $v$ be a valuation on a field $K$ and let $A$ its valuation ring. Then $v$ (or $A$) is called \emph{microbial} if 
the following equivalent assertions are satisfied.
\begin{eli}
\item
There exists a height $1$ valuation which is dependent to $v$.
\item
There exists exists a topologically nilpotent element $\ne 0$ in $K$.
\item
$A$ equipped with the valuation topology is non-discrete and adic.
\item
$A$ has a prime ideal of height $1$.
\item
There exists a convex subgroup $H$ of $\Gamma_v$ such that $\Gamma_v/H$ has height $1$.
\item
There exists a non-trivial homomorphisms of totally ordered groups $\Gamma_v \to \RR$.
\end{eli}
\end{remdef}

Every valuation of finite height $\ne 0$ is microbial.

%---------------------------------------------------------------

\subsection{Examples of non-archimedean rings}

For all examples in this section we fix the following notation. Let $A$ be a non-archimedean topological ring, let $I$ be an index set and fix a family $T = (T_i)_{i \in I}$ of subsets $T_i$ of $A$ such that for all $i \in I$, $m \in\NN$ and for every neighborhood $U$ of $0$ in $A$ the subgroup $T_i^mU$ is a neighborhood of $0$. This condition is automatic if all $T_i$ consist of units of $A$, e.g., if $T_i = \{1\}$ for all $i$, which is an important special case. For $\nu \in \NN_0^{(I)}$ we set
\[
T^{\nu} := \prod_{i\in I}T_i^{\nu(i)}.
\]
Then $T^{\nu}U$ is a neighborhood of $0$ for all $\nu$ and for all neighborhoods $U$ of $0$.

\begin{rem}
We claim that on the polynomial ring $A[X_i; i\in I]$ there is a (unique) topology that makes it into a topological ring and such that the following subgroups form a basis of neighborhoods of $0$.
\[
U_{[X]} := \set{\sum_{\nu}a_{\nu}X^{\nu}}{\text{$a_{\nu} \in T^{\nu}U$ for all $\nu$}},
\]
where $U$ runs through all open subgroups in $A$. We denote this non-archimedean topological ring by $A[X]_T$.

The inclusion $\iota\colon A \to A[X]_T$ is continuous and the set $\set{\iota(t)X_i}{I \in I, t \in T_i}$ is power-bounded. These properties characterize $A[X]_T$: Let $B$ be a non-archimedean topological ring, $f\colon A \to B$ a continuous ring homomorphism und sei $(x_i)_{i \in I}$ a family of elements $x_i \in B$ such that the set $\set{f(t)x_i}{I \in I, t \in T_i}$ is power-bounded in $B$. Then there exists a unique continuous ring homomorphism $g\colon A[X]_T \to B$ with $f = g \circ \iota$ and $g(X_i) = x_i$ for all $i \in I$.
\end{rem}

\begin{proof}[Proof of the claim]
We use Remark~\ref{DefineNATop}. For any open subgroup $V$ in $A$ there exists an open subgroup $U$ of $A$ such that $U\cdot U \subseteq V$. Then $U_{[X]}\cdot U_{[X]} \subseteq V_{[X]}$. Let $a \in A[X]$ and $U$ an open subgroup of $A$. It remains to show that there exists an open subgroup $V$ of $A$ such that $aV_{[X]} \subseteq U_{[X]}$. As $T^\nu U$ is a neighborhood of $0$ of $A$ for all $\nu$ and as $a_\nu = 0$ for almost all $\nu$, there exists a neighborhood $V$ of $0$ in $A$ such that $a_\nu V \subseteq T^\nu U$ for all $\nu$. Then $aV_{[X]} \subseteq U_{[X]}$.
\end{proof}

\begin{remdef}\label{ConvPower}
We consider now the ring of formal power series 
\[
A\dlbrack X \drbrack = A\dlbrack X_i ; i \in I \drbrack = \set{\sum_{\nu \in \NN_0^{(I)}}a_\nu X^\nu}{a_\nu \in A}
\]
We define a subring of $A\dlbrack X \drbrack$ as follows
\begin{equation}\label{DefConvPS}
A\langle X \rangle_T := \Bigl\{\sum_{\nu}a_\nu X^\nu \in A\dlbrack X \drbrack\ ;\
\begin{matrix}
\text{$a_\nu \in T^\nu\cdot U$ for all open subgroups}\\
\text{$U$ of $A$ and for almost all $\nu$}
\end{matrix}\Bigr\}.
\end{equation}
We endow $A\langle X \rangle_T$ with the (unique) structure of a topological ring such that the subgroups (for $U$ open subgroup in $A$)
\[
U_{\langle X\rangle} := \set{\sum_\nu a_\nu X^\nu \in A \langle X\rangle_T}{\text{$a_\nu \in T^\nu\cdot U$ for all $\nu \in \NN_0^{(I)}$}}
\]
form a fundamental system of neighborhoods in $A\langle X \rangle_T$. We also write simply $A\langle X \rangle$ if $T_i = \{1\}$ for all $i \in I$.
\end{remdef}

\begin{proof}
Easy (note that it is not entirely clear that $A\langle X \rangle_T$ is multiplicatively closed in $A\dlbrack X \drbrack$).
\end{proof}

\begin{prop}\label{PropConvPS}
\begin{ali}
\item
$A[X]_T$ is a dense subring of $A\langle X\rangle_T$, and the inclusion $A[X]_T \mono A\langle X\rangle_T$ is a topological embedding.
\item
If $A$ is Hausdorff and $T_i$ is bounded for all $i \in I$, then $A[X]_T$ and $A\langle X\rangle_T$ is Hausdorff.
\item
If $A$ complete and $T_i$ is bounded for all $i \in I$, then $A\langle X\rangle_T$ is the completion of $A[X]_T$.
\end{ali}
\end{prop}

\begin{proof}
Assertions~(1) is clear. Let us show~(2). The boundedness of the $T_i$ implies that $\bigcap_U T^\nu\cdot U \subseteq \bigcap_U U = \{0\}$, where $U$ runs through all open subgroups of $A$. This shows $\bigcap_U U_{\langle X\rangle} = \{0\}$.

To show~(3) it suffices to show that $A\langle X\rangle_T$ is complete. For a subset $B$ of $A\langle X\rangle_T$ let $B_\nu$ be the subset of $A$ consisting of the $\nu$-th coefficients of elements in $B$. Let $\Bcal$ be a Cauchy filter basis of $A\langle X\rangle_T$. Then for all $\nu \in \NN_0^{(I)}$ the set $\Bcal_\nu := \set{B_\nu}{B \in \Bcal}$ is a Cauchy filter basis of $A$. As $A$ is complete, $\Bcal_\nu$ converges to a unique element $a_{nu} \in A$. Let $a := \sum_{\nu \in \NN_0^{(I)}}a_{\nu}X^{\nu} \in A\dlbrack X \drbrack$. We want to show that $a \in A \langle X\rangle_T$ and that $\Bcal$ converges to $a$.

Let $U$ be an open subgroup of $A$. Choose $B \in \Bcal$ such that $b - c \in U_{\langle X \rangle}$ for all $b,c \in B$, i.e.~$b_\nu - c_\nu \in T^\nu\cdot U$ for all $\nu$ and all $b,c \in B$. Then
\begin{equation}
a_{\nu} - b_{\nu} \in T^{\nu}\cdot U \qquad\text{for all $b \in B$ and all $\nu$}.\tag{*}
\end{equation}
For any fixed $b \in B \subseteq A\langle X\rangle_T$, one has $b_{\nu} \in T^{\nu}\cdot U$ for almost all $\nu$. Hence (*) implies that $a_{\nu} \in T^{\nu}\cdot U$ for almost all $\nu$. This shows $a \in A \langle X\rangle_T$. Moreover, (*) also shows that $a - b \in U_{\langle X\rangle}$ for all $b \in B$. As $U$ was an arbitrary open subgroup of $A$, this shows that $\Bcal$ converges to $a$.
\end{proof}

\begin{cor}[Universal property of $A\langle X\rangle_T$]\label{UniversalConvPS}
Let $A$ be complete and let $T_i$ be bounded for all $i \in I$. Then the canonical homomorphism $\iota\colon A \to A\langle X \rangle_T$ is continuous and $\set{\iota(t)X_i}{i \in I, t \in T_i}$ is power-bounded in $A\langle X\rangle_T$. If $B$ is a complete ring, $f\colon A \to B$ a continuous homomorphism of rings and $(x_i)_i$ a family of elements $x_i \in B$ such that  $\set{f(t)x_i}{i \in I, t \in T_i}$ is power-bounded in $B$, then there exists a unique continuous ring homomorphism $g\colon A\langle X\rangle_T \to B$ with $f = g \circ \iota$ and $g(X_i) = x_i$ for all $i \in I$.
\end{cor}

The next class of examples are topological localizations of non-archimedean rings.

\begin{propdef}\label{TopLoc}
Let $S = (s_i)_{i\in I}$ be a family of elements of $A$ and let $R \subseteq A$ be the multiplicative subset generated by $\set{s_i}{i\in I}$. Then there exists on $R^{-1}A$ a non-archimedean ring topology making it into a topological ring
\[
A(\frac{T}{S}) = A(\frac{T_i}{s_i}\vert i \in I)
\]
such that $\set{\frac{t}{s_i}}{i \in I, t \in T_i}$ is power-bounded in $A(\frac{T}{S})$ and such that $A(\frac{T}{S})$ and the canonical  continuous homomorphism $\varphi\colon A \to A(\frac{T}{S})$ satisfy the following universal property. If $B$ is a non-archimedean topological ring and $f\colon A \to B$ is a continuous homomorphism such that $f(s_i)$ is invertible in $B$ for all $i \in I$ and such that the set $\set{f(t)f(s_i)^{-1}}{i \in I, t \in T_i}$ is power-bounded in $B$, then there exists a unique continuous ring homomorphism $g\colon A(\frac{T}{S}) \to B$ with $f = g \circ \varphi$.
\end{propdef}

The completion of $A(\frac{T}{S})$ is denoted by $A\langle\frac TS\rangle = A\langle\frac{T_i}{s_i}\vert i \in I\rangle$.

If $I$ consists of one element and $T = \{t_1,\dots,t_n\}$ is finite, we also write $A(\frac{t_1,\dots,t_n}{s})$ and $A\langle\frac{t_1,\dots,t_n}{s}\rangle$ instead of $A(\frac{\{t_1,\dots,t_n\}}{s})$ resp.\ $A\langle\frac{\{t_1,\dots,t_n\}}{s}\rangle$.

\begin{proof}
In the localization $R^{-1}A$ consider the subsets $E_i := \set{\frac{t}{s_i}}{t\in T_i}$ for $i \in I$ and $E := \bigcup_i E_i$. Let $D$ be the subring of $R^{-1}A$ generated by $E$. Endow $R^{-1}A$ with the group topology such that the subgroups $\set{D\cdot U}{\text{$U$ open subgroup of $A$}}$ is a fundamental system of neighborhoods of $0$. It is easy to check that this defines a ring topology on $R^{-1}A$ using as usual Remark~\ref{DefineNATop}. The canonical homomorphism $\varphi\colon A \to R^{-1}A$ is continuous as $1 \in D$. Moreover we have $D\cdot D = D$ which implies that $D$ is bounded in $R^{-1}A$. In particular, $E$ is power bounded in $R^{-1}A$.

We claim that $(R^{-1}A,\varphi)$ satisfies the desired universal property. Thus let $f\colon A \to B$ as in the proposition. Let $F$ be the subring generated by the power-bounded set $\set{f(t)f(s_i)^{-1}}{i \in I, t \in T_i}$. Then $F$ is bounded by Proposition~\ref{GroupBounded}~(2). Let $g\colon R^{-1}A \to B$ be the ring homomorphism such that $f = g \circ \varphi$. Then $g(D) \subseteq F$. Let $U$ be a neighborhood of $0$ in $B$. As $F$ is bounded, there exists a neighborhood $V$ of $0$ in $C$ with $F\cdot V \subseteq U$. As $f$ is continuous, there exists a neighborhood $W$ of $0$ in $A$ such that $f(W) \subseteq V$. Then $D\cdot W$ is a neighborhood of $0$ in $R^{-1}A$ and $g(D\cdot W) = g(D)\cdot f(W) \subseteq F\cdot V \subseteq U$. Therefore $g$ is continuous.
\end{proof}

\begin{rem}
We keep the notation of Proposition~\ref{TopLoc}.
\begin{ali}
\item
Let $J$ be the ideal of $A[X]_T$ generated by $\set{1-s_iX_i}{i \in I}$. We endow $A[X]_T/J$ with the quotient topology of $A[X]_T$. Then $A \to A[X]_T/J$ satisfies the universal property of $A \to A(\frac{T}{S})$ which shows that $A[X]_T/J = A(\frac{T}{S})$.
\item
One has $A(\frac{T_i}{s_i}\vert i \in I) = A(\frac{T_i \cup \{s_i\}}{s_i}\vert i \in I)$. Thus one may always assume that $s_i \in T_i$.
\end{ali}
\end{rem}

The last example class in this section is the following construction.

\begin{rem}\label{BoundedPS}
Consider the following subring of $A\dlbrack X \drbrack = A\dlbrack X_i \vert i \in I\drbrack$
\begin{align*}
A\dlangle X \drangle_T &:= A\dlangle X_i \vert i \in I \drangle_T \\
&:= \set{\sum_\nu a_\nu X^\nu \in A\dlbrack X \drbrack}{\text{$\exists\,K \subseteq A$ bounded: $a_\nu \in T^\nu\cdot K$ for all $\nu$}}.
\end{align*}
For every open subgroup $U$ of $A$ define a subgroup of $A\dlangle X \drangle_T$ by
\[
U_{\dlangle X \drangle} := \set{\sum_\nu a_\nu X^\nu \in A\dlangle X \drangle_T}{\text{$a_\nu \in T^\nu\cdot U$ for all $\nu$}}.
\]
Then the $U_{\dlangle X \drangle}$ form a fundamental system of neighborhoods of $0$ of a ring topology on $A\dlangle X \drangle_T$. If $T_i = \{1\}$ for all $i \in I$, we write $A\dlangle X \drangle = A\dlangle X_i \vert i\in I \drangle$ instead of $A\dlangle X \drangle_T$.
\end{rem}

%--------------------------------------------------------------

\subsection{Tate algebras}

In this subsection we let $k$ be a field which is complete with respect to a valuation $\abs\colon k \to \RR_{\geq 0}$. Its valuation ring is then Henselian (Remark~\ref{CompleteHensel}), and we endow any algebraic extension $k'$ of $k$ with the unique extension of $\abs$ to $k'$ (Proposition~\ref{ExtValHensel}). This extension is again called $\abs$.

\begin{example}\label{TateAlgebra}
Consider the special case $A = k$, $I = \{1,\dots,n\}$, and $T_i = \{1\}$ for all $i = 1,\dots,n$ of Definition~\ref{ConvPower}. Then
\[
T_n := T_{n,k} := k\langle X_1,\dots,X_n\rangle = \set{\sum_\nu a_\nu X^\nu}{\text{$a_\nu \to 0$ for $\sum_i \nu_i \to \infty$}}
\]
and the topology on $k\langle X_1,\dots,X_n\rangle$ is induced by the norm
\[
\dvert \sum a_\nu X^\nu \dvert := \max_\nu \vert a_\nu \vert
\]
which makes $T_n$ into a Banach algebra over $k$ (i.e. a complete normed algebra over $k$). It is called the \emph{Tate algebra of strictly convergent power series in $n$ variables} and the norm $\norm$ is called the \emph{Gau{\ss} norm}. It contains the polynomial algebra $k[X_1,\dots,X_n]$ as a dense $k$-subalgebra.
\end{example}

\begin{lem}\label{TateConvergent}
Let $\kbar$ be an algebraic closure of $k$ and let
\[
B^n(\kbar) := \set{(x_1,\dots,x_n) \in \kbar^n}{\text{$\vert x_i \vert \leq 1$ for all $i = 1,\dots,n$}}
\]
be the unit ball in $\kbar^n$. Then a formal power series $f = \sum_\nu a_\nu X^\nu \in k\dlbrack X_1,\dots,X_n\drbrack$ converges on $B^n(\kbar)$ if and only if $f \in k\langle X_1,\dots,X_n\rangle$.

More generally, for $i = 1,\dots,n$ let $T_i = \{t_i\}$ for some element $t_i \in k^\times$. Then a formal power series $f = \sum_\nu a_\nu X^\nu$ converges on $\set{x \in \kbar^n}{|x_i| \leq |t_i|^{-1}}$ if and only if $f \in k\langle X\rangle_{(t_1,\dots,t_n)}$.
\end{lem}

\begin{proof}
It suffices to show the special case. For $x \in B^n(\kbar)$ let $k \subseteq k' \subseteq \kbar$ be a finite subextension such that $x \in B^n(k')$. Then $|a_\nu x^\nu| \leq |a_\nu|$ so that the summands converge to zero, hence $\sum_\nu a_\nu x_\nu$ converges in $k'$ because $k'$ is complete. Conversely, if $f$ converges on $B^n(\kbar)$, then $f(1,...,1)$ converges so $(|a_\nu|)_\nu$ converges to zero.
\end{proof}

\begin{lem}\label{TnNoeth}
$T_n$ is noetherian.
\end{lem}

\begin{proof}
\cite{BGR}~6.1.1~Proposition~3.
\end{proof}

\begin{defi}\label{DefTateAlgebra}
A topological $k$-algebra $A$ is called \emph{topologically of finite type over $k$} if there exists a continuous open surjective $k$-algebra homomorphism $T_n \to A$.
\end{defi}

Thus every $k$-algebra $A$ topologically of finite type is of the form $T_n/\afr$, where $\afr$ is a closed ideal, and where the topology on $A$ is the quotient topology of the topology on $T_n$. Conversely:

\begin{prop}\label{DescribeTateAlgebra}
Let $\afr$ be an ideal of $T_n$. Let $\pi\colon T_n \to T_n/\afr$ be the canonical projection. Then $\afr$ is closed in $T_n$. Endow the $k$-algebra $A := T_n/\afr$ with the quotient topology. Then $A$ is a complete $k$-algebra and its topology is induced by the norm
\[
|a|_{\pi} := \inf\set{|f|}{\text{$f \in T_n$ with $\pi(f) = a$}}.
\]
The canonical projection $\pi$ is open and contractive.
\end{prop}

\begin{proof}
If $A$ is any commutative noetherian Banach $k$-algebra, then any finitely generated ideal is closed and $A/\afr$ is complete. Thus $\afr$ is closed in $T_n$ and $T_n/\afr$ is complete. Then clearly $\pi$ is continuous. It is not difficult to check that $|\cdot|_{\pi}$ is a norm, and it is immediate that $\pi$ is contractive. Further by Banach's theorem, any continuous $k$-linear surjective map of Banach spaces is open. In particular $\pi$ is open.
\end{proof}

Lemma~\ref{TnNoeth} implies:

\begin{prop}\label{TfNoetherian}
Every $k$-algebra topologically of finite type is noetherian.
\end{prop}

\begin{prop}\label{BasicTateAlgebra}
Every $k$-algebra homomorphism between $k$-algebras topologically of finite type is continuous.
\end{prop}

\begin{proof}
\cite{BGR}~6.1.3~Theorem~1.
\end{proof}

\begin{prop}\label{MaximalFinite}
Let $A$ be a $k$-algebra topologically of finite type. Then for every maximal ideal $\mfr$ the quotient $A/\mfr$ is a finite extension of $k$.
\end{prop}

%---------------------------------------------

%\subsection{Completed tensor products}
%
%Let $A$ be a non-archimedean topological ring and let $M$ and $N$ be topological $A$-modules such that there exists a fundamental system $(M_i)_i$ and $(N_j)_j$ of neighborhoods of zero of $M$ and $N$ consisting of subgroups of $(M,+)$ and $(N,+)$, respectively.
%

%=============================================

\section{$f$-adic rings and Tate rings}

\subsection{$f$-adic rings}
 
\begin{propdef}\label{CharFAdicRing}
A topological ring $A$ is called \emph{$f$-adic} if the following equivalent assertions are satisfied.
\begin{eli}
\item
There exists a subset $U$ of $A$ and a finite subset $T \subseteq U$ such that $\set{U^n}{n\geq 1}$ is a 
fundamental system of neighborhoods of $0$ in $A$ and such that $T\cdot U = U^2 \subseteq U$.
\item
$A$ contains on open subring $A_0$ such that the subspace topology on $A_0$ is $I$-adic, where $I$ is a finitely 
generated ideal of $A_0$.
\end{eli}
A subring $A_0$ (resp.~a pair $(A_0,I)$) as in~(ii) is called a \emph{ring of definition} (resp.~a \emph{pair of definition}).
\end{propdef}

Clearly any $f$-adic ring is non-archimedean.
 
\begin{proof}
We prove the equivalence together with the following lemma.
\end{proof}

\begin{lem}\label{DescribeRingDef}
Let $A$ be an $f$-adic ring. Then for a subring $A_0$ of $A$ (endowed with the subspace topology) the following 
assertions are equivalent.
\begin{dli}
\item
$A_0$ is a ring of definition.
\item
$A_0$ is open in $A$ and adic.
\item
$A_0$ is open in $A$ and bounded.
\end{dli}
\end{lem}
 
\begin{proof}
Let $A$ be a topological ring. We first remark that every open adic subring clearly is bounded. In particular~(b) implies~(c) for an arbitrary topological ring. The implication ``(a) $\implies$ (b)'' 
is trivial.
  
If $A$ satisfies~(ii) in Proposition~\ref{CharFAdicRing}, then we set $U := I$ and let $T$ be a finite system of 
generators of $I$ and see that $A$ satisfies~(i).
 
From now on let us assume that there exists a subset $U$ of $A$ and a finite subset $T \subseteq U$ such that $
\set{U^n}{n\geq 1}$ is a fundamental system of neighborhoods of $0$ in $A$ and such that $T\cdot U = U^2 
\subseteq U$. Let $W$ be the subgroup of $(A,+)$ generated by $U$. Then $W$ is open because $U$ is open. As 
$U^n \subseteq U$, $\ZZ\cdot1 + W$ is an open subring of $A$. It is also a bounded (Remark~\ref{GroupBounded}).
 
Thus we have shown all equivalences if we can prove that~(c) implies~(a). Let $A_0$ be an open and bounded 
subring of $A$. For every $n \in \NN$ put $T(n) = \set{t_1t_2\cdots t_n}{t_i \in T}$. Choose $k \in \NN$ such that 
$T(k) \subseteq A_0$ and let $I$ be the ideal of $A_0$ generated by $T(k)$. Let $l \in \NN$ such that $U^l 
\subseteq A_0$. Then we have for all $n \in \NN$
\[
I^n = T(nk)A_0 \supseteq T(nk)U^l = U^{l+nk}.
\]
Thus $I^n$ is an open neighborhood of $0$. Let $V$ be any neighborhood of $0$. As $A_0$ is bounded, there 
exists $m \in \NN$ such that $U^mA_0 \subseteq V$. Then $I^m \subseteq V$. Therefore $(I^n)_{n \in \NN}$ is a 
fundamental system of neighborhoods of $0$.
\end{proof}

\begin{rem}\label{FAdicMetrisable}
The topology of an $f$-adic ring $A$ is metrisable (Remark~\ref{TopFilteredGroup}). Thus Proposition~\ref{QuotComplete} shows that if $A$ is complete and $\afr \subseteq A$ is a closed ideal, then $A/\afr$ is complete.
\end{rem}

\begin{cor}\label{RingOfDef}
Let $A$ be an $f$-adic ring.
\begin{ali}
\item
If $A_0$ and $A_1$ are rings of definitions, then $A_0 \cap A_1$ and $A_0\cdot A_1$ are rings of definitions.
\item
Every open subring $C$ of $A$ is $f$-adic.
\item
$A^o$ is a subring of $A$ and it is the union of all rings of definitions of $A$.
\item
Let $B$ be a bounded subring of $A$ and let $C$ be an open subring of $A$ with $B \subseteq C$. Then there 
exists a ring of definition $A_0$ with $B \subseteq A_0 \subseteq C$.
\item
$A$ is adic if and only if it is bounded.
\end{ali}
\end{cor}

\begin{proof}
Assertions~(1) is clear, because if $A_0$ and $A_1$ are open and bounded, then $A_0 \cap A_1$ and $A_0\cdot A_1$ are open and bounded. Assertion~(2) is obvious, and~(3) follows from Proposition~\ref{GroupBounded}~(3). To show~(4) one may assume $A = C$ by~(2) and if $A_0$ is any ring of definition, then $A_0\cdot B$ is a ring of definition containing $B$. Finally Assertion~(5) follows immeditately from Lemma~\ref{DescribeRingDef}.
\end{proof}

\begin{example}\label{ExampleFAdic}
\begin{ali}
\item\label{ExampleFAdic1}
Every adic ring with a finitely generated ideal of definition is $f$-adic.
\item\label{ExampleFAdic2}
Let $k$ be a field endowed with the topolopgy defined by a valuation $v$. Then $k$ is $f$-adic if and only if $v$ is trivial or $v$ is microbial.
\item\label{ExampleFAdic3}
A field $k$ whose topology is defined by a microbial valuation $v$ is called a \emph{non-archimedean field}. Then the topology is also given by a nontrivial valuation $\abs$ with values in $\RR^{\geq 0}$. Then $k^o = \set{x \in k}{|x| \leq 1}$ and $k^{oo} =: \mfr$ is its maximal ideal. The subring $k^o$ is open in $k$ and it is an $I$-adic ring, where $I$ is an ideal generated by some element $\pi$ with $|\pi| < 1$.

If $\abs$ is a discrete valuation, then we may take for $I$ the maximal ideal. But in general $k^o$ does not carry the $\mfr$-adic topology (e.g., if $\Gamma_\abs \cong \QQ$, then $\mfr^n = \mfr$ for all $n \geq 1$).

The valuation ring of $v$ is contained in $k^o$.
\item\label{ExampleFAdic4}
Let $k$ be a non-archimedean field and let $T_n := k\langle X_1,\dots,X_n\rangle$ be the ring of convergent power series. Then $T_n$ is an $f$-adic ring. For the ring of power-bounded elements we have
\[
T_n^o = k^o\langle X_1,\dots,X_n\rangle,
\]
and $T_n^o$ is a ring of definition whose topology is the $I_{\langle X\rangle}$-adic topology (see Proposition~\ref{ConvPSTate} below).
\item\label{ExampleFAdic5}
Let $A$ be a ring, let $I$ be a finitely generated ideal of $A$, and define for $n \geq 0$ the following subgroup of the 
polynomial ring $A[T]$
\[
U_n := \set{\sum_k a_kT^k}{a_k \in I^{n+k}} 
\]
Then $U_nU_m \subseteq U_{n+m}$. Hence if we endow $A[T]$ with the topology such that the $(U_n)_n$ is a fundamental system of neighborhoods of $0$, $A[T]$ is a topological ring.
 
Then $A$ is an $f$-adic ring but it is not an adic ring if $I^m \ne I^{m+1}$ for all $m$.
\end{ali}
\end{example}

\begin{lem}\label{OpenIdeal}
Let $A$ be an $f$-adic ring and let $\nfr$ be the ideal of $A$ generated by the topologically nilpotent elements of $A$. Then an ideal $\afr$ of $A$ is open if and only if $\nfr \subseteq \rad(\afr)$.
\end{lem}

\begin{proof}
If $\afr$ is open, then for every topologically nilpotent element $f \in A$ there exists $n \in \NN$ such that $f^n \in \afr$, hence $\nfr \subseteq \rad(\afr)$. Conversely, assume that $\nfr \subseteq \rad(\afr)$. Let $A_0$ be a ring of definition of $A$ and let $I$ be a finitely generated ideal of definition of $A_0$. Then $I \subseteq \nfr \subseteq \rad(\afr)$ and thus there exists some power $I^m$ contained in $\afr$. But $I^m$ is open again and therefore $\afr$ is open.
\end{proof}

In particular one has for an $f$-adic ring $A$ and all ideals of definition $I$ of a ring of definition
\[
\set{\pfr \in \Spec A}{\text{$\pfr$ open in $A$}} = V(A^{oo}) = V(I),
\]
i.e., a prime ideal is open if and only if it contains every ideal of definition of every ring of definition of $A$.

\begin{remark}\label{QuotientFAdic}
Let $A$ be an $f$-adic ring, let $B$ be a topological ring and let $\pi\colon A \to B$ be a continuous surjective open homomorphism. Then $B$ is $f$-adic.
\end{remark}

\begin{rem}\label{FAdicCompletion}
Let $A$ be a $f$-adic ring, $B$ a ring of definition of $A$, and $I$ a finitely generated ideal of definition of $B$. We 
consider the completion $\Bhat$ as an open subring of $\Ahat$. By Proposition~\ref{CompleteAdic}, $\Bhat$ is adic, 
and $I\Bhat$ is an ideal of definition of $\Bhat$. Hence $\Ahat$ is an $f$-adic ring.
\end{rem}

\begin{prop}\label{PropCompletionFAdic}
Let $A$ be a $f$-adic ring, $B$ a ring of definition of $A$.
\begin{ali}
\item
The canonical ring homomorphism $A \otimes_B \Bhat \to \Ahat$ is an isomorphism.
\item
If $B$ is noetherian, then $A \to \Ahat$ is flat.
\item
If $A$ is a finitely generated $B$-algebra and $B$ is noetherian, then $\Ahat$ is noetherian.
\end{ali}
\end{prop}
 
\begin{proof}
Assertion~(1) is \cite{Hu_Cont}~Lemma~1.6. The remaining assertions follow immediately.
\end{proof}

%------------------------------------------------------------------

\subsection{Tate rings}

\begin{defi}\label{DefTate}
A topological ring is called \emph{Tate ring} if it is $f$-adic and has a topologically nilpotent unit.
\end{defi}

\begin{rem}\label{AdicTate}
An adic ring $A$ is a Tate ring if and only if its topology is the chaotic topology.
\end{rem}

\begin{example}\label{TateField}
Let $k$ be a field endowed with a non-trivial valuation $v$. Then $k$ with the valuation topology is a Tate ring if and only if $v$ is microbial (Definition~\ref{DefMicrobial}). Then $A(v)$ is a ring of definition.
\end{example}

\begin{example}\label{TateNormedAlgebra}
Let $k$ be a field endowed with a non-trivial valuation $\abs$. Assume that $\abs$ takes values in $\RR^{\geq 0}$. Every normed $k$-algebra $(A,\norm)$ is a Tate ring: $A_0 := \set{a \in A}{\dvert a\dvert \leq 1}$ is an open subring and if $r \in k^{\times}$ is any element with $v(r) < 1$, then $(r^nA_0)_{n\in\NN}$ is a fundamental system of neighborhoods of $0$ in $A_0$.

In particular every algebra topologically of finite type over $k$ (Example~\ref{TateAlgebra}) is a Tate ring.
\end{example}

\begin{prop}\label{DescribeTateRing}
Let $A$ be a Tate ring and let $B$ be a ring of definition. Then $B$ contains a topologically nilpotent unit $s$ of $A
$. For any such $s$ one has $A = B_s$ and $sB$ is an ideal of definition of $B$.
\end{prop}
 
\begin{proof}
As $B$ is an open neighborhood of $0$ there exists for every topologically nilpotent element $t$ an $n \in \NN$ 
such that $t^n \in B$. This shows the first assertion.
 
For every $a \in A$ there exists $n \in \NN$ such that $as^n \in B$, hence $A = B_s$. Multiplication with $s^n$ is a 
homeomorphism $A \to A$. This shows that $s^nB$ is open. Moreover, as $B$ is bounded, for every neighborhood 
$V$ of $0$ there exists $n \in \NN$ such that $s^nB \subseteq V$.
\end{proof}

\begin{cor}\label{DescribeTateAlg}
Let $k$ be a topological field whose topology is given by a valuation of height $1$. Let $\pi \in k^\times$ be a topologically nilpotent element and let $A$ be a topological $k$-algebra. Then $A$ is a Tate ring if and only if it contains a subring $A_0$ such that $(\pi^nA_0)_n$ is a fundamental system of open neighborhoods of $0$ of $A$.
\end{cor}

\begin{proof}
The condition is necessary by Proposition~\ref{DescribeTateRing}. Conversely, if $A_0$ is a subring such that $(\pi^nA_0)_n$ is a fundamental system of open neighborhoods of $0$ of $A$, then $A_0$ is an open subring and as linearly topologized ring it is bounded. Thus $A_0$ is a ring of definition.
\end{proof}

%-------------------------------------------------------------

\subsection{Banach's theorem for Tate rings}

Recall the following version of Banach's theorem:

\begin{thm}\label{Banach}
Let $A$ be a topological ring that has a sequence converging to $0$ consisting of units of $A$ (e.g., if $A$ is a Tate ring). Let $M$ and $N$ be Hausdorff topological $A$-modules that have countable fundamental systems of open neighborhoods of $0$. Assume that $M$ is complete. Let $u\colon M \to N$ be an $A$-linear map. Consider the following properties.
\begin{definitionlist}
\item
$N$ is complete.
\item
$u$ is surjective.
\item
$u$ is open.
\end{definitionlist}
Then any two of these properties imply the third.
\end{thm}

The proof will show that $u$ is already open and surjective if $N$ is complete and $u(M)$ is open in $N$.

\begin{proof}
Missing
%Note first, that $M$ and $N$ are both metrisable by Remark~\ref{MetrisableGroup}.
%
%Let us first assume that $u$ is open and surjective. Then $u$ is strict and induces a homeomorphism $M/\ker(u) \to N$. As $N$ is Hausdorff, $\ker(u)$ is closed in $M$ and hence $M/\ker(u)$ is complete by Proposition~\ref{QuotComplete}.
%
%Now assume that $N$ is complete and that $u(M)$ is open in $N$. As Baire's theorem (\cite{Bou_TG}~IX, \S5.3, Theorem~1) says that in a complete metrisable topological space no nonempty open subspace is meager (i.e. a countable union of subsets $A$ which are equal to their boundary), this implies that $u(M)$ is not meagre.
%
%Let $V$ be a neighborhood of zero in $M$. We claim that there $\overline{u(V)}$ is a neighborhood of zero in $N$.
%
%Let $W$ be a neighborhood of zero in $M$ such that $W = -W$ and $W + W \subseteq V$. Let $(a_n)_n$ a squence of units in $A$ converging to zero.
\end{proof}

\begin{prop}\label{SubmoduleClosed}
Let $A$ be a complete Tate ring, and let $M$ be a complete topological $A$-module that has a countable fundamental system of open neighborhoods of $0$. Then $M$ is noetherian if and only if every submodule of $M$ is closed. In particular $A$ is noetherian if and only if every ideal is closed.
\end{prop}

\begin{proof}
Missing
\end{proof}

\begin{prop}\label{LinearCont}
Let $A$ be a complete noetherian Tate ring.
\begin{assertionlist}
\item
Every finitely generated $A$-module has a unique $A$-module topology that is complete and that has a countable fundamental system of open neighborhoods of $0$.
\item
Let $f\colon M \to N$ be an $A$-linear map of finitely generated modules that are endowed with the topology from~(1). Then $f$ is continuous and the map $f\colon M \to f(M)$ is open.
\end{assertionlist}
\end{prop}

\begin{proof}
Missing
\end{proof}

\begin{rem}
Let $A$ be a complete noetherian Tate ring, $A_0$ a ring of definition and $s \in A_0$ a topologically nilpotent unit of $A$ (such that $A_0$ has the $sA_0$-adic topology). Let $M$ be a finitely generated $A$-module and choose a finitely generated $A_0$-submodule $M_0$ of $M$ such that $A\cdot M_0 = M$. Then $\set{s^nM_0}{n\in\NN}$ is a fundamental system of open neighborhoods of $0$ in $M$  for the topology defined in Proposition~\ref{LinearCont}.
\end{rem}

%-----------------------------------------------------

\subsection{Examples}

\begin{lem}\label{ZeroNeighbor}
Let $A$ be an $f$-adic ring and let $T$ be a subset such that ideal of $A$ generated by $T$ is open in $A$. Then $T^n\cdot U$ is a neighborhood of zero in $A$ for all $n \in \NN$ and for all open subgroups $U$ of $A$.
\end{lem}

\begin{proof}
Let $U$ and $n$ be given. As $T\cdot A$ is a neighborhood of zero, $T^n\cdot A = (T\cdot A)^n$ is a neighborhood of zero. Let $A_0$ be a ring of definition and $I$ an ideal of definition of $A_0$ such that $I \subseteq T^n\cdot A$. Let $L$ be a finite system of generators of $I$ and let $K$ be a finite subset of $A$ such that $L \subseteq T^n\cdot K$. Choose $m \in \NN$ with $K\cdot I^m \subseteq U$. Then $I^{m+1} = L\cdot I^m \subseteq (T^n\cdot K)\cdot I^m = T^n \cdot (K \cdot I^m) \subseteq T^n\cdot U$.
\end{proof}

\begin{prop}\label{ConvPSTate}
Let $A$ be an $f$-adic ring, let $B$ be a ring of definition and let $I$ a finitely generated ideal of definition of $B$. Let $\Lambda$ be a set, and let $T = (T_\lambda)_{\lambda\in \Lambda}$ be a family of subsets of $A$ such that the ideal of $A$ generated by $T_\lambda$ is open for all $\lambda$.
\begin{assertionlist}
\item
Then $A[X]_T$ is an $f$-adic ring, $B_{[X]}$ is a ring of definition, and $I_{[X]} = I\cdot B_{[X]}$ is a finitely generated ideal of definition. If $A$ is a Tate ring, then $A[X]_T$ is a Tate ring.
\item
Assume that $\Lambda$ is finite. Then $A\langle X\rangle_T$ is an $f$-adic ring, $B_{\langle X\rangle}$ is a ring of definition, and $I_{\langle X\rangle} = I\cdot B_{\langle X\rangle}$ is a finitely generated ideal of definition. If $A$ is a Tate ring, then $A\langle X\rangle_T$ is a Tate ring.
\item
Let $S = \set{s_\lambda}{\lambda \in \Lambda}$ be a family of elements of $A$. Then the topological localizations $A(\frac TS)$ and $A\langle\frac{T}{S}\rangle$ are $f$-adic rings. If $A$ is a Tate ring, then $A(\frac TS)$ and $A\langle\frac{T}{S}\rangle$ are also Tate rings.
\end{assertionlist}
\end{prop}

By Lemma~\ref{ZeroNeighbor} the topological rings $A[X]_T$, $A\langle X\rangle_T$, $A(\frac{T}{S})$, and $A\langle\frac{T}{S}\rangle$ are defined.

Assertion~(2) also holds for arbitrary index sets $\Lambda$ if one assumes in addition that $A$ is Hausdorff and that $T_\lambda$ is bounded for all $\lambda \in \Lambda$.

\begin{proof}
The inclusions $I_{[X]} = I\cdot B_{[X]}$ and $I\cdot B_{\langle X\rangle} \subseteq I_{\langle X\rangle}$ are clear. Let us show that $I_{\langle X\rangle} \subseteq I\cdot B_{\langle X\rangle}$. Let $i_1,\dots,i_m$ be generators of $I$. Then every $a \in T^\nu\cdot I^k$ ($\nu \in \NN_0^\Lambda$, $k \in \NN$) can be written as $a = i_1a_1 + \cdots + i_ma_m$ with $a_j \in T^\nu\cdot I^{k-1}$. As $\set{I^k}{k\in\NN}$ is a fundamentalsystem of neighborhoods of zero of $A$, every $u \in I_{\langle X\rangle}$ can be written as $u = i_1u_1+ \cdots + i_mu_m$ with $u_j \in B_{\langle X\rangle}$.

Clearly the rings $B_{[X]}$ and $B_{\langle X\rangle}$ are open subrings of $A_{[X]}$ and $A_{\langle X\rangle}$, respectively. Moreover, if $I$ is an ideal of definition of $B$, then $I^m$ is also an ideal of definition of $B$. Thus we have shown that $(I^m)_{[X]} = I^m\cdot B_{[X]}$ and $(I^m)_{\langle X\rangle} = I^m\cdot B_{\langle X\rangle}$. This shows that $B_{[X]}$ and $B_{\langle X\rangle}$ are adic rings. Thus we have seen that $A_{[X]}$ and $A_{\langle X\rangle}$ are $f$-adic rings.

Moreover, it is easy to see that $A(\frac{T}{S})$ is $f$-adic. Then $A\langle\frac{T}{S}\rangle$ is $f$-adic by Proposition~\ref{PropCompletionFAdic}.

Finally, there are continuous homomorphisms $A \to A[X]_T, A\langle X\rangle_T, A(\frac{T}{S}), A\langle\frac{T}{S}\rangle$. Thus if $A$ contains a nilpotent unit, all the other rings contain a nilpotent unit as well.
\end{proof}

\begin{rem}
Let $A$ be an adic ring. Then $A\langle X\rangle_T$ is in general not adic. In fact one can show that $A\langle X\rangle_T$ is adic if and only if $A\langle X\rangle_T$ is isomorphic to $A\langle X\rangle$. An analogous remark holds for $A[X]_T$.

Moreover, $A\langle\frac{T}{S}\rangle$ is adic if and only if $A\langle\frac{T}{S}\rangle$ is isomorphic to $A\langle\frac{\{1\}}{r} \vert r \in R\rangle$, where $R$ is some subset of $A$.
\end{rem}

%--------------------------------------------------------------------

\subsection{Adic homomorphisms}

\begin{defi}
Let $A$ and $B$ be $f$-adic rings. A ring homomorphism $\varphi\colon A \to B$ is called \emph{adic} if there exist 
rings of definitions $A_0$ of $A$ and $B_0$ of $B$ and an ideal of definition $I$ of $A_0$ such that $\varphi(A_0) 
\subseteq B_0$ and such that $\varphi(I)B_0$ is an ideal of definition of $B_0$.
\end{defi}
 
Any adic ring homomorphism is continuous. Conversely, for every continuous homomorphism $\varphi\colon A \to B$ there exist 
always rings of definitions $A_0$ of $A$ and $B_0$ of $B$ and finitely generated ideals of definition $I$ of $A_0$ and $J$ of $B_0$ such that $\varphi(A_0) \subseteq B_0$ and such that $\varphi(I) \subseteq J$. But in general $\varphi(I)B_0$ is not an ideal of definition of $B_0$.

\begin{example}\label{TrivExampleAdic}
Let $A$ be a discrete ring. Then $A$ is adic with ideal of definition $I = 0$. Any homomorphism $\varphi\colon A\to B$ to an $f$-adic ring $B$ is continuous. It is adic if and only if $B$ also carries the discrete topology.
\end{example}

\begin{prop}\label{HomomTate}
Let $\varphi\colon A \to B$ be a continuous ring homomorphism between $f$-adic rings. Assume that $A$ is a Tate ring. Then $B$ is a Tate ring, $\varphi$ is adic, and for every ring of definition $B_0$ of $B$ one has $\varphi(A)\cdot B_0 = B$.
\end{prop}
 
\begin{proof}
Let $A_0$ (resp.~$B_0$) be rings of definition of $A$ (resp.~$B$) such that $\varphi(A_0) \subseteq B_0$. Let $s \in A_0$ be a topologically nilpotent unit of $A$. Then $\varphi(s)$ is a topologically nilpotent unit of $B$ and hence $B$ is a Tate ring. Replacing $s$ be some power, we may assume that $s \in A_0$. Then $sA_0$ is an ideal of definition of $A_0$ and $\varphi(s)B_0$ is an ideal of definition of $B_0$ (Proposition~\ref{DescribeTateRing}). This shows that $\varphi$ is adic. Let $B_0 \subseteq B$ be an arbitrary ring of definition. Replacing $s$ by some power we may assuume that $\varphi(s) \in B_0$. Moreover one has $A = (A_0)_s$ and $B = (B_0)_{\varphi(s)}$ by Proposition~\ref{DescribeTateRing} and hence $\varphi(A)\cdot B_0 = B$.
\end{proof}

\begin{prop}\label{CompAdic}
Let $\varphi\colon A \to B$ and $\psi\colon B \to C$ be continuous ring homomorphisms of $f$-adic rings.
\begin{assertionlist}
\item
If $\varphi$ and $\psi$ are adic, then $\psi \circ \varphi$ is adic.
\item
If $\psi \circ \varphi$ is adic, then $\psi$ is adic.
\end{assertionlist}
\end{prop}

\begin{proof}
\cite{Hu_Cont}~Corollary~1.9
\end{proof}

\begin{remark}
Let $\varphi\colon A \to B$ be a continuous ring homomorphisms of $f$-adic rings. Let $A' \subseteq A$ and $B' \subseteq B$ be open subrings (then $A'$ and $B'$ are again $f$-adic rings by Corollary~\ref{RingOfDef}) such that $\varphi(A') \subseteq B'$. Then if $\varphi$ is adic, its restriction $A' \to B'$ is also adic.
\end{remark}

%------------------------------------------------------------------

\subsection{Homomorphisms topologically of finite type}\label{Sec:HomTFTfadic}

\begin{definition}
Let $A$ be an $f$-adic ring. A ring homomorphism $\varphi\colon A \to B$ to a complete $f$-adic ring $B$ is called \emph{strictly topologically of finite type} if there exist $n \in \NN_0$ and a surjective continuous open ring homomorphism
\[
\pi\colon \Ahat\langle X_1,\dots,X_n\rangle \epi B
\]
of $A$-algebras.
\end{definition}

\begin{propdef}\label{DefTFT}
Let $A$ be an $f$-adic ring. A ring homomorphism $\varphi\colon A \to B$ to a complete $f$-adic ring $B$ is called \emph{topologically of finite type} if the following equivalent conditions are satisfied.
\begin{equivlist}
\item
There exist $n \in \NN_0$, finite subsets $T_1,\dots,T_n$ of $A$ such that $T_i\cdot A$ is open in $A$ for all $i = 1,\dots,n$, and a surjective continuous open ring homomorphism
\[
\pi\colon \Ahat\langle X_1,\dots,X_n\rangle_{T_1,\dots,T_n} \epi B
\]
of $A$-algebras.
\item
The homomorphism $\varphi$ is adic, there exists a finite subset $M$ of $B$ such that $A[M]$ is dense in $B$, and there exist rings of definition $A_0 \subseteq A$ and $B_0 \subseteq B$ and a finite subset $N$ of $B_0$ such that $\varphi(A_0) \subseteq B_0$ and such that $A_0[N]$ is dense in $B_0$.
\item
There exist rings of definition $A_0 \subseteq A$ and $B_0 \subseteq B$ with $\varphi(A_0) \subseteq B_0$ such that $B_0$ is an $A_0$-algebra strictly topologically of finite type such and such that $B$ is finitely generated over $A\cdot B_0$.
\item
For every open subring $A_0$ of $A$ there exists an open subring $B_0$ of $B$ with $\varphi(A_0) \subseteq B_0$ such that $B_0$ is an $A_0$-algebra strictly topologically of finite type and such that $B$ is finitely generated over $A \cdot B_0$.
\end{equivlist}
\end{propdef}

\begin{remark}\label{RemTFT}
If $A_0$ is a ring of definition of $A$ in~(iv), then $B_0$ is automatically a ring of definition in $B$. Indeed, if $A_0$ is adic, then $A_0 \langle X_1,\dots,X_m\rangle$ is adic and hence the existence of a continuous surjective open homomorphism $A_0 \langle X_1,\dots,X_m\rangle \to B_0$ implies that $B_0$ is adic.
\end{remark}

\begin{proof}
\proofstep{``(i) $\implies$ (iii)''}
Let $(A_0,I)$ be a pair of definition of $A$. Then $\Ahat \langle X\rangle_T$ is the completion of $A[X]_T$ (Proposition~\ref{PropConvPS}). By Proposition~\ref{ConvPSTate}, $((A_0)_{[X]}, I(A_0)_{[X]}$ is a pair of definition of $A[X]_T$. Moreover $(A_0)_{[X]} = A_0[t_iX_i\,;\,i \in \{1,\dots,n\}, t_i \in T_i]$. Hence if we set $B_0 := \pi(\widehat{(A_0)_{[X]}})$, then there exists $\pi'$ as desired. Now by Proposition~\ref{PropCompletionFAdic} we have $\Ahat \langle X\rangle_T = A[X]_T \otimes_{(A_0)_{[X]}} \widehat{(A_0)_{[X]}}$. As $A[X]_T$ is of finite type over $A\cdot (A_0)_{[X]}$, this implies the claim.

\proofstep{``(iii) $\implies$ (iv)''}
\cite{Hu_Habil}~2.3.25

\proofstep{``(iv) $\implies$ (ii) $\implies$ (iii)''}
The fact that (iv) implies that $\varphi$ is adic follows from Proposition~\ref{ConvPSTate}. All other implications are clear.

\proofstep{``(iv) $\implies$ (i)''}
We may assume that $A$ is complete. By hypothesis $\varphi$ factors through a continuous open homomorphism $\sigma\colon A \langle X_1,\dots,X_n\rangle \to B$ such that $B$ is generated as $A \langle X_1,\dots,X_n\rangle$-algebra by finitely many elements $b_1,\dots,b_s$. Let $A_0$ be a ring of definition of $A$ and let $L \subseteq A_0$ a subset which generates an ideal of definition of $A_0$. Choose $k \in \NN$ such that $I^kb_j$ is power-bounded for all $j = 1,\dots,s$. We set $T := (T_1,\dots,T_{m+s})$ with $T_i := \{1\}$ for $i = 1,\dots,m$ and $T_i = I^k$ for $i = m+1,\dots,m+s$. Then the universal property of $A \langle X_1,\dots,X_{m+s}\rangle_T$ (Corollary~\ref{UniversalConvPS}) implies that there exists a continuous homomorphism of $A$-algebras $\rho\colon A \langle X_1,\dots,X_{m+s}\rangle_T \to B$ such that $\rho(X_i) = \sigma(X_i)$ for $i = 1,\dots,m$ and $\rho(X_{j+m}) = b_j$ for $j = 1,\dots,s$. Then $\rho$ is surjective. Moreover, for every open subgroup $U$of $A$ the image of $U_{\langle X_1,\dots,X_{m+s}\rangle} \subseteq A \langle X_1,\dots,X_{m+s}\rangle_T$ under $\rho$ contains the image of $U_{\langle X_1,\dots,X_m\rangle} \subseteq A \langle X_1,\dots,X_m\rangle$ under $\sigma$. Therefore $\rho$ is open because $\sigma$ is open.
\end{proof}

The same argument also shows the following strict variant.

\begin{prop}\label{PropStrictTFT}
Let $A$ be an $f$-adic ring. Let $\varphi\colon A \to B$ be a ring homomorphism to a complete $f$-adic ring $B$. Then the following conditions are equivalent.
\begin{equivlist}
\item
$\varphi$ is strictly topologically of finite type.
\item
There exist rings of definition $A_0 \subseteq A$ and $B_0 \subseteq B$ with $\varphi(A_0) \subseteq B_0$ such that $B_0$ is an $A_0$-algebra strictly topologically of finite type such and such that $B = A\cdot B_0$.
\item
For every open subring $A_0$ of $A$ there exists an open subring $B_0$ of $B$ with $\varphi(A_0) \subseteq B_0$ such that $B_0$ is an $A_0$-algebra strictly topologically of finite type and such that $B = A \cdot B_0$.
\end{equivlist}
\end{prop}

\begin{example}\label{LocalTFT1}
\begin{assertionlist}
\item
A ring homomorphism between discrete topological rings is topologically of finite type if and only if it is of finite type.
\item
Let $A$ be an $f$-adic ring, $s_1,\dots,s_n \in A$ and $T_i \subseteq A$ finite subsets such that $T_i\cdot A$ is open in $A$ for $i = 1,\dots,n$. Then the canonical homomorphism $A \to A \langle \frac{T_1}{s_1},\dots,\frac{T_n}{s_n}\rangle$ is topologically of finite type.

Indeed, this follows from Proposition~\ref{DefTFT}~(ii) by taking (with the notation there) $N = \set{\frac{t}{s_i}}{i \in \{1,\dots,n\}, t \in T_i}$ and $M = N \cup \{\frac{1}{s}\}$.
\end{assertionlist}
\end{example}

\begin{prop}\label{CompositionTFT}
Let $A$ be an $f$-adic ring, let $B$, $C$ be complete $f$-adic rings, and let $\varphi\colon A \to B$ and $\psi\colon B \to C$ be continuous ring homomorphisms.
\begin{assertionlist}
\item
If $\varphi$ and $\psi$ are (strictly) topologically of finite type, then $\psi \circ \varphi$ is (strictly) topologically of finite type.
\item
If $\psi \circ \varphi$ is topologically of finite type, then $\psi$ is topologically of finite type.
\end{assertionlist}
\end{prop}

\begin{proof}
The properties~(iv) in Proposition~\ref{DefTFT} and~(iii) in Proposition~\ref{PropStrictTFT} imply~(1). Assertion~(2) is clear.
\end{proof}

\begin{prop}\label{TateTFT}
Let $A$ be a Tate ring and let $B$ be a complete $f$-adic ring. Then a homomorphism $\varphi\colon A \to B$ is topologically of finite type if and only if it is strictly topologically of finite type.
\end{prop}

Note that if these equivalent conditions are satisfied, then $B$ is a Tate ring and $\varphi$ is adic.

\begin{proof}
Let $B_0$ be a ring of definition of $B$. Then Proposition~\ref{HomomTate} shows that $A\cdot B_0 = B$. This implies the claim by Proposition~\ref{DefTFT} and Proposition~\ref{PropStrictTFT}.
\end{proof}

\begin{proposition}\label{NoethTFT}
Let $A$ be an $f$-adic ring and let $B$ be an $f$-adic $A$-algebra topologically of finite type. If $A$ has a noetherian ring of definition, then $B$ has a noetherian ring of definition.
\end{proposition}

\begin{proof}
We may assume that $B = A \langle X\rangle_T$. Let $A_0$ be a noetherian ring of definition of $A$. Then $(A_0)_{[X]}$ is a ring of definition of $A[X]_T$ which is finitely generated of $A_0$ by $\set{t_iX_i}{i \in \{1,\dots,n\}, t_i \in T_i}$ and hence noetherian. Its completion is a ring of definition of $A \langle X\rangle_T$ and as the completion of a noetherian adic ring it is itself noetherian.
\end{proof}

%------------------------------------------------------------

\subsection{Strongly noetherian Tate rings}

\begin{propdef}
A Tate ring $A$ is called \emph{strongly noetherian} if the following equivalent conditions are satisfied.
\begin{equivlist}
\item
$\Ahat\langle X_1,\dots,X_n\rangle$ is noetherian for all $n \in \NN_0$.
\item
Every Tate ring topologically of finite type over $A$ is noetherian.
\end{equivlist}
\end{propdef}

\begin{proof}
The equivalence of both assertions follows immediately from \ref{TateTFT}.
\end{proof}

\begin{rem}
\begin{assertionlist}
\item
If $A$ is a strongly noetherian Tate ring, then every Tate ring topologically of finite type over $A$ is strongly noetherian by Proposition~\ref{CompositionTFT}.
\item
Every completely valued field $(k,v)$, where $v$ is of height $1$ (or, more generally, microbial), is strongly noetherian (\cite{BGR}~5.2.6~Theorem~1). Hence every Tate ring topologically of finite type over $k$ is strongly noetherian as well.
\item
Every Tate ring that has a noetherian ring of definition is strongly noetherian.
\end{assertionlist}
\end{rem}

\begin{example}\label{LocalTFT}
Let $A$ be an $f$-adic ring, let $s_1,\dots,s_n \in A$ and $T_1,\dots,T_n \subseteq A$ finite subsets such that $T_i\cdot A$ is open for all $i = 1,\dots,n$. Then the canonical ring homomorphism $\iota\colon A \to \Ahat\langle\frac{T_1}{s_1},\dots,\frac{T_n}{s_n}\rangle$ is topologically of finite type.

Indeed, $\iota$ is by definition adic, and if we set $M = \set{\frac{t}{s_i}}{i = 1,\dots,n; t \in T_i}$, then $A[M]$ is dense in $\Ahat\langle\frac{T_1}{s_1},\dots,\frac{T_n}{s_n}\rangle$ and $A_0[M]$ is dense in $\widehat{A_0}\langle\frac{T_1}{s_1},\dots,\frac{T_n}{s_n}\rangle$, if $A_0$ is a ring of definition of $A$.

If $A$ is a strongly noetherian Tate ring, we may represent $\Ahat\langle\frac{T_1}{s_1},\dots,\frac{T_n}{s_n}\rangle$ as a quotient of a ring of restricted power series. Set $C = \Ahat\langle X_{i,t} | i=1,\dots,n, t \in T_i\rangle$ and let $\afr$ be the ideal of $C$ generated by $\set{t - s_iX_{i,t}}{i = 1,\dots,n; t \in T_i}$. By hypothesis, $C$ is noetherian and hence $\afr$ is a closed ideal (Proposition~\ref{SubmoduleClosed}). As $A$ is a Tate ring, $T_i\cdot A = A$ by Remark~\ref{UnderstandRational}~(2). Hence the image of $s_i$ in $C/\afr$ is a unit. Now it is easily seen that the homomorphisms $A \to \Ahat\langle\frac{T_1}{s_1},\dots,\frac{T_n}{s_n}\rangle$ and $A \to C/\afr$ satisfy the same universal property.

In particular, $\Ahat\langle\frac{T_1}{s_1},\dots,\frac{T_n}{s_n}\rangle$ is again strongly noetherian.
\end{example}

\begin{example}\label{LaurentSeries}
Let $A$ be a Tate ring and let $\Ahat\langle X,X^{-1}\rangle$ be the ring of all formal series $\sum_{n \in \ZZ}a_nX^n$ such that $a_n \in \Ahat$ for all $n \in \ZZ$ and such that for every neighborhood $U$ of zero in $\Ahat$ there exist only finitely many $n \in \ZZ$ with $a_n \notin U$. Note that the product of two such series is well defined:
\[
(\sum_{n \in \ZZ}a_nX^n)(\sum_{n \in \ZZ}b_nX^n) = \sum_{n \in \ZZ}c_nX^n,
\]
where $c_n$ is the convergent series $\sum_{k+l=n}a_kb_l$. We endow $\Ahat\langle X,X^{-1}\rangle$ with the ring topology such that the sets \[
\sett{\sum_{n\in\ZZ}a_nX^n \in \Ahat\langle X,X^{-1}\rangle}{$a_n \in U$ for all $n \in \ZZ$}
\]
form a basis of neighborhoods of zero if $U$ runs through all open neighborhoods of zero of $A$.

Then $\Ahat\langle X,X^{-1}\rangle$ is an $A$-algebra topologically of finite type. More precisely, it is isomorphic to $\Ahat\langle X,Y\rangle/(XY - 1)$. In particular we see that, if $A$ is strongly noetherian, then $\Ahat\langle X,X^{-1}\rangle$ is a strongly noetherian Tate ring.

Indeed, it is not difficult to check that both algebras are initial objects in the category of homomorphisms from $A$ to complete Tate rings $B$ together with a distinguished unit $u \in B^\times$ such that $u$ and $u^{-1}$ are power-bounded.
\end{example}

%====================================================================

\section{The adic spectrum of an affinoid ring as topological space}

\subsection{The space $\Spv(A,I)$}
In this subsection we denote by $A$ a ring and by $I$ an ideal of $A$ such that there exists a finitely generated ideal $J$ of $A$ such that $\sqrt I = \sqrt J$.

Recall from Definition~\ref{DefCofinal} that an element $\gamma$ of a totally ordered group $\Gamma$ is called cofinal for a subgroup $H$ if for all $h \in H$ there exists $n \in \NN$ with $\gamma^n < h$.

\begin{lemma}\label{CofinalIdeal}
Let $v\colon A \to \Gamma_v \czero$ be a valuation on $A$. Let $H \subseteq \Gamma_v$ be a convex subgroup such that $c\Gamma_v \subsetneq H$. Then $\cfr := \sett{a \in A}{$v(a)$ is cofinal for $H$}$ is an ideal of $A$ and $\rad(\cfr) = \cfr$.
\end{lemma}

\begin{proof}
If $a \in \cfr$ and $b \in A$ with $v(b) < v(a)$, then $b \in \cfr$. Hence if $a,b \in \cfr$, then $v(a+b) \leq \max\{v(a),v(b)\}$ shows that $a + b \in \cfr$. Let $b \in A$, $a \in \cfr$. If $v(b) \leq 1$ then $v(ba) \leq v(a)$ and hence $ba \in \cfr$. If $v(b) > 1$, then $v(b) \in c\Gamma_v$ and hence $v(ab) = v(a)v(b)$ is cofinal for $H$ because $c\Gamma_v \subsetneq H$ (Remark~\ref{CofinalConvex}). The property $\rad(\cfr) = \cfr$ is clear.
\end{proof}

\begin{lem}\label{ExistConvexSubgroup}
Let $v\colon A \to \Gamma_v \czero$ be a valuation on $A$ with $v(I) \cap c\Gamma_v = \emptyset$. Then there exists a greatest convex subgroup $H$ of $\Gamma_v$ such that $v(a)$ is cofinal for $H$ for all $a \in I$. Moreover, if $v(I) \ne \{0\}$, then $H \supseteq c\Gamma_v$ and $H$ is the smallest convex subgroup of $\Gamma_v$ such that $v(I) \cap H \ne \emptyset$. 
\end{lem}

Note that $v(I) \cap c\Gamma_v = \emptyset$ is equivalent to $v(a) < c\Gamma_v$ for all $a \in I$ (if $v(a) \geq 1$ then $v(a) \in c\Gamma_v$ and if $1 \geq v(a) \geq \gamma$ for some $\gamma \in c\Gamma_v$, then $v(a) \in c\Gamma_v$ because $c\Gamma_v$ is by definition convex).

\begin{proof}
If $v(I) = \{0\}$, we may choose $H = \Gamma_v$. Thus we may assume that $v(I) \ne \{0\}$ and hence $c\Gamma_v \ne \Gamma_v$. 

Note that for any subgroup $H$ of $\Gamma_v$ one has $v(I) \cap H = \emptyset \iff v(\sqrt{I}) \cap H = \emptyset$. This and Lemma~\ref{CofinalIdeal} shows that we may assume that $I$ is finitely generated. Let $T$ be a finite set of generators of $I$ and let $H$ be the convex subgroup of $\Gamma_v$ generated by $h := \max \set{v(t)}{t \in T}$. As remarked above $v(h) < 1$ and hence $h$ is cofinal for $H$ and hence $v(t)$ is cofinal for $H$ for all $t \in T$. Moreover $h < c\Gamma_v < h^{-1}$ as remarked before the proof and hence $c\Gamma_v \subsetneq H$ because $H$ is convex. Thus Lemma~\ref{CofinalIdeal} shows that $v(a)$ is cofinal for $H$ for all $a \in I$. If $H'$ is any other convex subgroup of $\Gamma_v$ such that $v(a)$ is cofinal for $H'$ for all $a \in I$. Then in particular $h$ is cofinal for $H'$ which implies that $H' \subseteq H$.

Clearly one has $v(I) \cap H \ne \emptyset$. Let $H'$ be a convex subgroup of $\Gamma_v$ with $H' \cap v(I) \ne \emptyset$. We show that $H \subseteq H'$. As the set of totally convex subgroups is totally ordered and as $c\Gamma_v \cap v(I) = \emptyset$, we necessarily have $c\Gamma_v \subsetneq H'$. Choose $a \in I$ with $v(a) \in H'$ and write $a = \sum_{t\in T}a_tt$ with $a_t \in A$. Hence
\[
v(a) \leq \max\set{v(a_t)v(t)}{t \in T}.
\]
Let $t' \in T$ with $v(a) \leq v(a_{t'})v(t')$. If $v(a_{t'}) \leq 1$, then $v(a) \leq v(t') < 1$ which implies $v(t') \in H'$ (as $H'$ is convex). If $v(a_{t'}) > 1$, then $v(a_{t'}) \in c\Gamma_v \subseteq H'$ and hence $H' \ni v(a_{t'})^{-1}v(a) \leq v(t') < 1$ which again shows that $t' \in H$. By definition of $h$ and again by the convexity of $H'$ we get $h \in H'$. Thus $H \subseteq H'$.
\end{proof}

The lemma allows us to define for every valuation $v$ on $A$ the following subgroup of $\Gamma_v$.

\begin{definition}\label{DefGammaI}
Let $c\Gamma_v(I)$ be the group $c\Gamma_v$ if $v(I) \cap c\Gamma_v \ne \emptyset$. Otherwise let $c\Gamma_v(I)$ be the greatest convex subgroup $H$ of $\Gamma_v$ such that $v(a)$ is cofinal for $H$ for all $a \in I$.
\end{definition}

Then $c\Gamma_v(I)$ is a convex subgroup of $\Gamma_v$ which always contains $c\Gamma_v$.

\begin{lemma}\label{GammaIGamma}
The following conditions are equivalent.
\begin{equivlist}
\item
$c\Gamma_v(I) = \Gamma_v$.
\item
$v(a)$ is cofinal for $\Gamma_v$ for all $a \in I$, or $\Gamma_v = c\Gamma_v$.
\item
$v(a)$ is cofinal for $\Gamma_v$ for every element $a$ in a set of generators of an ideal $J$ such that $\sqrt{J} = \sqrt{I}$, or $\Gamma_v = c\Gamma_v$.
\end{equivlist}
Indeed, the equivalence of (i) and (ii) follows from the definition of $c\Gamma_v(I)$. For the equivalence of (ii) and (iii) follows from Lemma~\ref{CofinalIdeal}.
\end{lemma}

With this definition we set
\begin{equation}\label{DefSpvAI}
\Spv(A,I) := \set{v \in \Spv A}{c\Gamma_v(I) = \Gamma_v}
\end{equation}
and endow it with the subspace topology induced by $\Spv A$. We get the same subspace if we replace $I$ by any ideal $J$ with $\sqrt{I} = \sqrt{J}$. Moreover the map
\begin{equation}\label{DefRetraction}
r\colon \Spv A \to \Spv (A,I), \qquad v \sends v\rstr{c\Gamma_v(I)}
\end{equation}
is a retraction (i.e., $r(v) = v$ for all $v \in \Spv(A,I)$).

\begin{lem}\label{SpvAISpectral}
Let $A$ be an $f$-adic ring and let $I$ be an ideal of $A$ such that there exists a finitely generated ideal $J$ of $A$ with $\sqrt I = \sqrt J$.
\begin{ali}
\item\label{SpvAISpectral1}
The topological space $\Spv(A,I)$ is a spectral space and the set $\Rcal$ of the following sets form a basis of quasi-compact open subsets of the topology, which is stable under finite intersections.
\[
\Spv(A,I)(\frac{T}{s}) := \sett{v \in \Spv(A,I)}{$v(t) \leq v(s) \ne 0$ for all $t \in T$},
\]
where $s \in A$ and $T \subseteq A$ finite with $I \subseteq \sqrt{T\cdot A}$.
\item\label{SpvAISpectral2}
The retraction $r\colon \Spv A \to \Spv(A,I)$ is a continuous spectral map.
\item\label{SpvAISpectral3}
For $v \in \Spv A$ with $v(I) \ne 0$ one has $r(v)(I) \ne 0$.
\end{ali}
\end{lem}

\begin{proof}
We may assume that $I$ is finitely generated.

\proofstep{(i)}
As we endowed $\Spv(A,I)$ with the topology induced by $\Spv(A)$, it is clear that each set in $\Rcal$ is open in $\Spv(A,I)$. For $s$ and $T$ as above one has
\[
\Spv(A,I)(\frac{T}{s}) = R(\frac{T \cup \{s\}}{s}).
\]
Thus we may always assume that $s \in T$. If $s_1, s_2 \in A$ and let $T_1, T_2 \subseteq A$ finite subsets such that $I \subseteq \sqrt{T_i\cdot A}$ for $i = 1,2$. Setting $T := \set{t_1t_2}{t_i \in T_i}$, the ideal $\sqrt{T\cdot A}$ still contains $I$. Moreover, if we assume that $s_i \in T_i$, then
\[
\Spv(A,I)(\frac{T_1}{s_1}) \cap \Spv(A,I)(\frac{T_2}{s_2}) = \Spv(A,I)(\frac{T}{s_1s_2}).
\]
Therefore the intersection of two sets in $\Rcal$ is again in $\Rcal$.

\proofstep{(ii)}
We now show that $\Rcal$ is a basis of the topology of $\Spv(A,I)$. Let $v \in \Spv(A,I)$ and let $U$ be an open neighborhood of $v$ in $\Spv A$. Choose $g_0,\dots,g_n \in A$ such that
\[
v \in W := \sett{w \in \Spv A}{$w(g_i) \leq w(g_0) \ne 0$ for $i = 1,\dots,n$} \subseteq U,
\]
which is possible because the sets as $W$ form a basis of the topology of $\Spv A$.

Assume that $\Gamma_v = c\Gamma_v$. Then there exists $d \in A$ with $v(g_0d) \geq 1$. Hence
\[
v \in W' := \Spv(A,I)(\frac{g_1d,\dots,g_nd,1}{g_0d}) \subseteq W.
\]
Now assume that $\Gamma_v \ne c\Gamma_v$. Let $\{s_1,\dots,s_m\}$ be a set of generators of $I$. By Lemma~\ref{GammaIGamma} there exists $k \in \NN$ with $v(s_i)^k < v(g_0)$ for all $i = 1,\dots,m$. Then
\[
v \in \Spv(A,I)(\frac{g_1,\dots,g_n,s_1^k,\dots,s_m^k}{g_0}) \subseteq W.
\]

\proofstep{(iii)}
Let $s \in A$, $T \subseteq A$ finite with $I \subseteq \sqrt{T\cdot A}$. Set $U := \Spv(A,I)\frac{T}{s}$ and $W := \Spv(A)(\frac{T}{s})$. We claim that $W = r^{-1}(U)$.

Since $U \subseteq W$ and since every point $x \in r^{-1}(U)$ is a horizontal generization of $r(x) \in U$, one has $r^{-1}(U) \subseteq W$ because $W$ is open in $\Spv A$. Conversely, let $w \in W$ be given. It remains to show that $r(w) \in U$. If $w(I) = 0$, then $c\Gamma_v(I) = \Gamma_v$ and hence $r(w) = w$. Thus we may assume that $w(I) \ne 0$. As $r(w)$ is a horizontal specialization of $w$, one has $r(w)(t) \leq r(w)(s)$ for all $t \in T$. It remains to show that $r(w)(s) \ne 0$. But$r(w)(s) = 0$ would imply that $t \in \pfr := \supp(r(w))$ for all $t \in T$ and hence $I \subseteq \pfr$ as $I \subseteq \sqrt{T\cdot A}$. But by Lemma~\ref{ExistConvexSubgroup}, $w(I) \ne 0$ implies that $c\Gamma_w(I) \cap w(I) \ne \emptyset$, i.e., there exists $a \in I$ with $r(w)(a) \ne 0$. Contradiction.

We also have just seen that \eqref{SpvAISpectral3} holds.

\proofstep{(iv)}
Let $\hat\Rcal$ be the Boolean algebra of subsetes of $\Spv(A,I)$ generated by $\Rcal$ and let $X$ be the set $\Spv(A,I)$ equipped with the topology generated by $\hat\Rcal$. It follows from~(iii) that for every subset $C$ in $\hat\Rcal$ the preimage $r^{-1}(C)$ is constructible in $\Spv A$. Hence $r\colon (\Spv A)_{\rm cons} \to T$ is continuous. Since $(\Spv A)_{\rm cons}$ is quasi-compact (Proposition~\ref{PropConstTop}) and $r$ is surjective, $T$ is quasi-compact. By definition, every set of $\Rcal$ is open and closed in $T$. Moreover, $\Spv(A,I)$ is Kolmogorow as a subspace of the Kolmogorow space $\Spv A$. Hence we can apply Proposition~\ref{ConstructSpectral} to deduce that $\Spv(A,I)$ is a spectral space and that $\Rcal$ is a basis of open quasi-compact subsets. Moreover, (iii) then shows that $r\colon \Spv A \to \Spv(A,I)$ is spectral.
\end{proof}

\begin{remark}
The inclusion $\Spv(A,I) \mono \Spv(A)$ is not spectral in general
\end{remark}

%---------------------------------------------------------

\subsection{The spectrum of continuous valuations}
\begin{defi} Let $A$ be a topological ring, let $v$ be a valuation on $A$ and let $\Gamma$ be its valuation group. Then $v$ is called \emph{continuous} if it $\set{a \in A}{v(a) < \gamma}$ is open in $A$ for all $\gamma \in \Gamma$ (i.e., the topology of $A$ is finer than the topology defined by $v$).
 
We denote by $\Cont(A)$ the subspace of $\Spv(A)$ of continuous valuations on $A$
\end{defi}
 
\begin{remark}
Let $A$ be a topological ring.
\begin{assertionlist}
\item
A valuation $v$ on $A$ is continuous, if and only if the map $v\colon A \to \Gamma_v \czero$ is continuous, where $\Gamma_v \czero$ is endowed with the topology defined in Remark~\ref{TopTotOrdGroup}.
\item
If $A$ carries the discrete topology, then $\Cont(A) = \Spv(A)$.
\item
A valuation $v$ on $A$ is continuous if and only if for all $\gamma \in \Gamma_v$ the set $A_{\leq\gamma} := \set{a \in A}{v(a) \geq \gamma}$ is open in $A$.

Indeed, if $v$ is continuous, then the subgroup $A_{\leq\gamma}$ contains the open subgroup $A_{<\gamma}$. Conversely, $A_{<\gamma} = \bigcup_{\delta < \gamma}A_{\leq \delta}$.
\end{assertionlist}
\end{remark}

\begin{rem}
Clearly any continuous homomorphism $\varphi\colon A \to B$ of topological rings induces a continuous morphism $
\Cont(B) \to \Cont(B)$ by composition with $\varphi$.
\end{rem}

\begin{thm}\label{ContSpectral}
Let $A$ be an $f$-adic ring. Let $I$ be an ideal of definition of a ring of definition of $A$. Then
\[
\Cont(A) = \sett{v \in \Spv(A,I\cdot A)}{$v(a) < 1$ for all $a \in I$}.
\]
\end{thm}

\begin{proof}
Let $v \in \Cont(A)$. Every element $a \in I$ is topologically nilpotent and hence for every $\gamma \in \Gamma$ there exists $n \in \NN$ such that $a^n \in \set{f \in A}{v(f) < \gamma}$, i.e. $v(a)^n < \gamma$. Hence $v(a)$ is cofinal for $\Gamma_w$. This shows $v(a) < 1$ and that $c\Gamma_v(I\cdot A) = \Gamma_v$ by Lemma~\ref{GammaIGamma}~(iii). Hence $v \in \sett{v \in \Spv(A,I\cdot A)}{$v(a) < 1$ for all $a \in I$}$.

Conversely, let $v \in \Spv(A,I\cdot A)$ with $v(a) < 1$ for all $a \in I$. We claim that $v(a)$ is cofinal for $\Gamma_v$ for all $a \in I$. Indeed, if $c\Gamma_v \ne \Gamma_v$, then the claim follows from Lemma~\ref{GammaIGamma}. Thus we may assume that $c\Gamma_v = \Gamma_v$. Let $a \in I$ and $\gamma \in \Gamma_v$ be given. By our assumptions there exists $t \in A$ with $v(t) \ne 0$ and $v(t)^{-1} \leq \gamma$. Choose $n \in \NN$ such that $ta^n \in I$ (this exists because $\{t\}$ is bounded). Then $v(ta^n) < 1$ and hence $v(a)^n < \gamma$. This proves the claim.

Let $\gamma \in \Gamma$ be given. Let $T$ be a finite system of generators of $I$ and set $\delta := \max\set{v(t)}{t \in T}$. By the claim above there exists $n \in \NN$ such that $\delta^n < \gamma$. Then $v(a) < \gamma$ for all $a \in T^n\cdot I = I^{n+1}$. Hence $v$ is continuous.
\end{proof}

\begin{remark}\label{ContCofinal}
Let $A$ be an $f$-adic ring, $I$ an ideal of definition of a ring of definition of $A$.
\begin{assertionlist}
\item\label{ContCofinal1}
The proof of Theorem~\ref{ContSpectral} shows that a valuation $v$ on $A$ is continuous if and only if $v(a)$ is cofinal for $\Gamma_v$ for all $a \in I$.
\item\label{ContCofinal2}
Let $v$ be a continuous valuation on $A$ and let $H \subsetneq \Gamma_v$ be a proper convex subgroup. For $a \in I$ the image of $v(a)$ in $\Gamma_v/H$ is cofinal for $\Gamma_v/H$ (Corollary~\ref{CofinalQuot}). This shows that the vertical generalization $v_{/H}$ of $v$ is still continuous.
\end{assertionlist}
\end{remark}

\begin{corollary}
Let $A$ be an $f$-adic ring. Then $\Cont(A)$ is a spectral space.
\end{corollary}

\begin{proof}
Let $I$ be the ideal of $A$ generated by an ideal of definition of a ring of definition of $A$. Then $\Cont(A)$ is the complement of the open subset $\bigcup_{f \in I}\Spv(A,I)(\frac{1}{f})$ of $\Spv(A,I)$. Hence $\Cont(A)$ is closed in $\Spv (A,I)$ and hence a spectral space because $\Spv(A,I)$ is a spectral space by Lemma~\ref{SpvAISpectral}.
\end{proof}

\begin{example}
Let $K$ be a field, let $v$ be a valuation on $K$ and endow $K$ with the topology induced by $v$. Then $\Cont K = \set{w \in \Spv K}{\text{$w$ dependent on $v$}}$.
\end{example}

%-----------------------------------------------------------------

\subsection{Affinoid rings}

For a subring $B$ of a ring $A$ we denote by $B^{\rm int}$ the integral closure of $B$ in $A$.

\begin{defi}
\begin{ali}
\item
Let $A$ be an $f$-adic ring and let $A^o$ be the subring of power-bounded elements. A subring $B$ of $A$ is called \emph{ring of integral elements} if $B$ is open and integrally closed in $A$ and if $B \subseteq A^o$.
\item
An affinoid ring is a pair $(A,A^+)$, where $A$ is an $f$-adic ring and where $A^+$ is a ring of integral elements. Often we will simply write $A$ instead of $(A,A^+)$.
\item
An affinoid ring $(A,A^+)$ is called \emph{complete} (resp.\ \emph{adic}, resp.\ \emph{Tate}, resp.~...) if $A$ has this property.
\item
A morphism of affinoid rings $\affoid A \to \affoid B$ is a ring homomorphism $\varphi\colon A \to B$ such that $\varphi(A^+) \subseteq B^+$. It is called \emph{continuous} (resp.\ \emph{adic}) if $\varphi\colon A \to B$ is continuous (resp.\ adic).
\end{ali}
\end{defi}

\begin{rem}\label{LargestSmallestIntegral}
Let $A$ be an $f$-adic ring.
\begin{ali}
\item
Then $A^o$ is a ring of integral elements (Proposition~\ref{GroupBounded}). It is clearly the largest ring of integral elements.
\item
If $\Atilde$ is any subring of $A$ which is integrally closed in $A$, then $A^{oo} \subseteq \Atilde$ if and only if $\Atilde$ is open in $A$. Indeed, if $\Atilde$ is open in $A$, then for all $a \in A^{oo}$ there exists an $n \in \NN$ such that $a^n \in \Atilde$. But then $a \in \Atilde$ because $\Atilde$ is integrally closed in $A$. If $\Atilde$ contains $A^{oo}$, then for every ring of definition $A_0$ of $A$ and any ideal of definition $I$ of $A_0$ one has $I \subseteq A^{oo} \subseteq \Atilde$. Hence $\Atilde$ is open.

Thus if $A'$ is the integral closure of $\ZZ\cdot 1 + A^{oo}$ in $A$, then $A'$ is the smallest ring of integral elements of $A$.
\item
If $\Atilde$ is an subring of $A$ with $A^{oo} \subseteq \Atilde \subseteq A^o$, then its integral closure $A^+$ in $A$ is a ring of integral elements by~(2).
\end{ali}
\end{rem}

\begin{example}
Let $A$ be an adic ring which has a finitely generated ideal of definition. Then $(A,A)$ is an affinoid ring.
\end{example}

\begin{example}\label{TateAffinoidField}
Let $K$ be a field endowed with the topology given by a microbial valuation $v$ on $K$ and let $v_1$ be the unique equivalence class of a height $1$ valuation dependent on $v$. Then $(K,A(v))$ is a Tate affinoid ring and $K^o = A(v_1)$.
\end{example}

Let $A$ be an $f$-adic ring. For every subset $X \subseteq A$ we set
\[
S_X := \sett{v \in \Cont(A)}{$v(a) \leq 1$ for all $a \in X$}.
\]
As $S_{\{a\}} = \Spv(A,I)(\frac{1,a}{1}) \cap \Cont(A)$ ($I$ an ideal of definition of a ring of definition of $A$) is constructible for all $a \in A$, $S_X$ is pro-constructible for all subsets $X \subseteq A$.

\begin{prop}\label{ClassifyRingIntegral}
We set
\[
\Scal_A := \set{S_X}{X \subseteq A}
\]
and denote by $\Rcal_A$ the set of subrings of $A$ which are open and integrally closed in $A$.
\begin{ali}
\item
The maps
\begin{align*}
\sigma\colon \Rcal_A \to \Scal_A, &\qquad A' \sends \set{v \in \Cont A}{\text{$v(f) \leq 1$ for all $f \in A'$}},\\
\tau\colon \Scal_A \to \Rcal_A, &\qquad S \sends \set{f \in A}{\text{$v(f) \leq 1$ for all $v \in S$}},
\end{align*}
are mutually inverse bijections.
\item
Let $A^+$ be an element of $\Rcal_A$ with $A^+ \subseteq A^o$, i.e., $A^+$ is a ring of integral elements. Then every point of $\Cont(A)$ is a vertical specialization of a point in $\sigma(A^+)$ . In particular, $\sigma(A^+)$ is dense in $\Cont A$.
\item
If $A$ is a Tate ring and has a noetherian ring of definition, then also the converse of~(2) does hold: if $A' \in \Rcal_A$ such that $\sigma(A')$ is dense in $\Cont A$ then $A' \subseteq A^o$.
\end{ali}
\end{prop}

\begin{proof}
\cite{Hu_Cont}~Lemma~3.3
\end{proof}

\begin{prop}\label{ConvPSAffinoid}
Let $A = (A,A^+)$ be an affinoid ring and let $(T_i)_{i\in I}$ be a finite family of subsets of $A$ such that $T_i\cdot A$ is open for all $i$. Then $(A\langle X\rangle_T, (A^+\langle X\rangle_T)^{\rm int})$ is an affinoid ring which we denote by $A\langle X\rangle_T$.
\end{prop}

If $A$ is complete and $T_i$ is finite for all $i \in I$, then $A\langle X\rangle_T$ is complete by Proposition~\ref{PropConvPS}.

\begin{proof}
As $A^+$ is open in $A$, then $A^+_{\langle X\rangle}$ is open in $A\langle X\rangle_T$. Hence $(A^+_{\langle X\rangle})^{\rm int}$ is open. It remains to show that $(A^+_{\langle X\rangle})^{\rm int}$ is contained in $A\langle X\rangle_T^o$. As $A\langle X\rangle_T^o$ is integrally closed in $A\langle X\rangle_T$, it suffices to show that $A^+_{\langle X\rangle}$ is contained in $A\langle X\rangle_T^o$. This follows from the following lemma.
\end{proof}

\begin{lem}\label{AoPS}
Let $A$ be a $f$-adic ring and let $(T_i)_{i\in I}$ be a family of subsets of $A$ such that $T_i\cdot A$ is open in $A$ for all $i$. Then $(A^o)_{\langle X\rangle} \subseteq (A\langle X\rangle_T)^o$.
\end{lem}

It is not difficult to see that one has equality in Lemma~\ref{AoPS} if $I$ is finite and $T_i = \{1\}$ for all $i$.

\begin{proof}
Let $a \in (A^o)_{\langle X\rangle}$ and write $a = b + c$, where $b \in (A^o)_{[X]}$ and $c \in J_{\langle X\rangle}$, where $J$ is an ideal of definition of $A$. Then $b$ and $c$ are both power-bounded. Therefore $a$ is power-bounded.
\end{proof}

\begin{remark}\label{UniversalAffinoidPS}
With the notation of Proposition~\ref{ConvPSAffinoid} the canonical homomorphism $\iota \colon A \to A \langle X\rangle_T$ has the following universal property. Let $\varphi\colon A \to B$ a continuous homomorphism to a complete affinoid ring $B$ and let $b_i \in B$ for $i \in I$ such that $\varphi(t)b_i \in B^+$ for all $i \in I$ and $t \in T_i$. Then there exists a unique continuous homomorphism $\psi\colon A \langle X\rangle_T \to B$ of affinoid rings such that $\varphi = \psi \circ \iota$ and $\psi(X_i) = b_i$.
\end{remark}

\begin{remdef}\label{QuotientAffinoid}
Let $A = (A,A^+)$ be an affinoid ring and let $I$ be an ideal of $A$. Then $A/I := (A/I, (A^+/(A^+ \cap I))^{\rm int})$ is an affinoid ring, if $A/I$ is endowed with the quotient topology. If $A_0 \subseteq A$ is a ring of definition, its image in $A/I$ is a ring of definition in $A/I$. The canonical morphism $A \to A/I$ is an adic morphism of affinoid rings.
\end{remdef}

%----------------------------------------------------

\subsection{The adic spectrum as topological space}

\begin{defi}\label{DefAdicSpec}
Let $A = \affoid A$ be an affinoid ring. The \emph{adic spectrum of $A$} is the subspace
\[
\Spa A = \set{v \in \Cont(A)}{\text{$v(f) \leq 1$ for all $f \in A^+$}}
\]
of $\Cont(A)$.
\end{defi}

If $A$ is an $f$-adic ring and $A'$ is the integral closure of $\ZZ\cdot 1 + A^{oo}$ in $A$, then one has $\Spa (A,A') = \Cont A$.

\begin{rem}
Let $A$ be a ring, let $\Atilde \subseteq A$ be any subring, and let $A^+$ be its integral closure in $A$. Then for any $v \in \Spv A$ one has
\[
\text{$v(a) \leq 1$ for all $a \in A^+$} \iff \text{$v(a) \leq 1$ for all $a \in \Atilde$}
\]
Indeed if $a \in A^+$, then $a^n + b_{n-1}a^{n-1} + \cdots + b_0 = 0$ with $b_i \in \Atilde$. Assume $v(a) > 1$ but $v(b_i) \leq 1$ for all $i$. Then for all $i < n$ one has $v(a^n) = v(a)^n > v(a)^i \geq v(b_ia^i)$ and hence $v(a^n) > v(b_{n-1}a^{n-1} + \cdots + b_0)$; contradiction.

Thus if $A$ is $f$-adic and $\Atilde$ is a subring of $A$ such that its integral closure $A^+$ satisfies $A^{oo} \subseteq A^+ \subseteq A^o$ (i.e., $A^+$ is a ring of integral elements), then
\[
\Spa (A,A^+) = \sett{v \in \Cont(A)}{$v(a) \leq 1$ for all $a \in \Atilde$}.
\]
\end{rem}

\begin{remark}
Let $(A,A^+)$ be an affinoid ring. Then Proposition~\ref{ClassifyRingIntegral} shows that for every point $v$ of $\Cont(A)$ there exists a vertical generization $x$ which is in $\Spa(A,A^+)$. In fact, if $\supp v$ is open we can simply take for $x$ the trivial valuation with $\supp x = \supp v$. If $\supp v$ is not open, then it follows from Remark~\ref{ContCofinal}~\eqref{ContCofinal2} and Proposition~\ref{HeightOneSpa} below that there exists a vertical generization $x$ of $v$ of height $1$ (necessarily unique by Remark~\ref{DefVertical}) such that $x \in \Spa(A,A^+)$. 
\end{remark}

\begin{example}
Let $A$ be a discrete topological ring and let $A^+$ any subring of $A$ which is integrally closed in $A$ (e.g., $A^+ = A$). Then $(A,A^+)$ is an affinoid ring and $\Spa (A,A^+) = \set{v \in \Spv A}{\text{$v(a) \leq 1$ for all $a \in A^+$}}$.
\end{example}

\begin{rem}
Let $A$ be an affinoid ring. For $x \in \Spa A$ and $f \in A$ we sometimes write $|f(x)|$ instead of $x(f)$.
\end{rem}

\begin{remdef}
Any continuous morphism $\varphi\colon A \to B$ between affinoid rings defines via composition a continuous map $\Spa(\varphi)\colon\Spa B \to \Spa A$. Thus we obtain a contravariant functor from the category of affinoid rings to the category of topological spaces.
\end{remdef}

\begin{defi}\label{DefRational}
Let $A$ be an affinoid ring. The subsets of the form
\[
R(\frac{T}{s}) := \set{v \in \Spa A}{\forall\,t\in T\colon v(t) \leq v(s) \ne 0},
\]
where $s \in A$ and $T \subseteq A$ is a finite subset such that $T\cdot A$ is open in $A$, are called \emph{rational subsets}.
\end{defi}

By definition of the topology of $\Spa A$ as subspace topology of $\Spv A$, all rational subsets are open in $\Spa A$.

\begin{rem}\label{UnderstandRational}
Let $A$ be an affinoid ring.
\begin{ali}
\item\label{UnderstandRational1}
Let $I$ be the ideal of $A$ generated by the topologically nilpotent elements of $A$. Let $T$ be a finite subset of $A$. Then $T\cdot A$ is open in $A$ if and only if $I \subseteq \rad(T\cdot A)$ (Lemma~\ref{OpenIdeal}).
\item\label{UnderstandRational2}
Assume that $A$ is Tate. Then (1) implies immediately, that the ideal $T\cdot A$ of $A$ generated by a finite subset $T$ of $A$ is open if and only if $T\cdot A = A$. 
\item\label{UnderstandRational3}
For $s \in A$ and $T \subseteq A$ finite with $T\cdot A$ open in $A$ one has
\[
R(\frac{T}{s}) = R(\frac{T \cup \{s\}}{s}).
\]
Thus we may always assume that $s \in T$.
\item\label{UnderstandRational4}
For $s \in A^{\times}$ a unit and $T' \subseteq A$ an arbitrary finite subset, $R(\frac{T'}{s}) = R(\frac{T' \cup \{s\}}{s}$ is always a rational subset. In particular,
\[
R(\frac{f}{1}) = \set{x \in \Spa A}{|f(x)| \leq 1}
\]
is an open rational subset for all $f \in A$.
\item\label{UnderstandRational5}
Let $s_1, s_2 \in A$ and let $T_1, T_2 \subseteq A$ finite subsets such that $T_i\cdot A$ is open in $A$ for $i = 1,2$, i.e. $T_i\cdot A$ contains $I^{n_i}$, where $I$ is an ideal of definition of $A$ and $n_i \in \NN$. Setting $T := \set{t_1t_2}{t_i \in T_i}$ we see that $T\cdot A$ contains $I^{n_1+n_2}$ hence it is an open ideal.

Moreover, if we assume that $s_i \in T_i$, then
\[
R(\frac{T_1}{s_1}) \cap R(\frac{T_2}{s_2}) = R(\frac{T}{s_1s_2}).
\]
Therefore the intersection of two rational subsets is again a rational subset.
\end{ali}
\end{rem}

\begin{lem}\label{NonInfinitesimalZero}
Let $A = (A,A^+)$ be an affinoid ring, $X \subset \Spa A$ a quasi-compact subset and $f \in A$ such that $|f(x)| \ne 0$ for all $x \in X$. Then there exists a neighborhood of zero $I$ in $A$ such that $|a(x)| < |f(x)|$ for all $x \in X$ and $a \in I$.
\end{lem}

\begin{proof}
Let $T$ be a finite subset of $A^{oo}$ such that $T\cdot A^{oo}$ is open (e.g., one may choose $T$ as a system of generators of an ideal of definition). For every $n \in \NN$ set
\[
X_n := \sett{x \in \Spa A}{$|t(x)| \leq |f(x)| \ne 0$ for all $t \in T^n$}.
\]
Then $X_n$ is open in $\Spa A$, and $X \subseteq \bigcup_n X_n$. Hence $X \subseteq X_m$ for some $m \in \NN$. Then $I := T^m\cdot A^{oo}$ has the desired properties.
\end{proof}

\begin{cor}\label{TateZeroNeighbor}
Let $A = (A,A^+)$ be a Tate affinoid ring, $Y \subset \Spa A$ a quasi-compact subset and $s \in A$ such that $|s(y)| \ne 0$ for all $y \in Y$. Then there exists a unit $\pi \in A^{\times}$ such that $|\pi(y)| < |s(y)|$ for all $y \in Y$.
\end{cor}

\begin{proof}
Let $I \subseteq A$ be a neighborhood of zero $I$ such that $|a(y)| < |s(y)|$ for all $y \in Y$ and $a \in I$. As $A$ is Tate, there exists a unit $\pi$ of $A$ in $I$ and thus we have $|\pi(y)| \leq |s(y)|$ for all $y \in Y$. 
\end{proof}

\begin{rem}\label{RationalBoundedBelow}
Let $A$ be Tate, $T \subseteq A$ finite with $T\cdot A = A$, $s \in A$. Let $x \in \Spa A$ with $|t(x)| \leq |s(x)|$ for all $t \in T$. Then there exists a unit $\pi$ of $A$ such that $|\pi(x)| \leq |s(x)|$. In particular
\[
R(\frac{T}{s}) = \set{x \in \Spa A}{\forall\,t \in T: |t(x)| \leq |s(x)|}.
\]
Indeed, let $T = \{f_1,\dots,f_n\}$ and $g_1,\dots,g_n \in A$ with $\sum_i g_if_i = 1$. As $A^+$ is open in $A$, there exists a topologically nilpotent unit $\pi \in A$ such that $\pi g_i \in A^+$ for all $i$. Then we have
\[
|\pi(x)| = |\sum_i(\pi g_i)(x)f_i(x)| \leq \max_i|(\pi g_i)(x)||f_i(x)| \leq |s(x)|.
\]
\end{rem}

\begin{proposition}\label{RationalOpenElements}
Let $A$ be a complete affinoid ring, $s \in A$, $T = \{t_1,\dots,t_n\} \subseteq A$ finite with $T\cdot A$ open in $A$. Then there exists a neighborhood $J$ of $0$ in $A$ such that for all $s' \in s + J$, $t'_i \in t_i + J$ the ideal of $A$ generated by $t'_1,\dots,t'_n$ is open in $A$ and such that
\[
R(\langle \frac{t_1,\dots,t_n}{s}\rangle) = R(\langle \frac{t'_1,\dots,t'_n}{s'}\rangle.
\]
\end{proposition}

\begin{proof}
\cite{Hu_Cont}~3.10
\end{proof}

\begin{thm}\label{AdicSpecSpectral}
Let $A$ be an affinoid ring.
\begin{assertionlist}
\item
$\Spa A$ is a spectral space.
\item
The rational subsets of $\Spa A$ form a basis of quasi-compact open subsets of $\Spa A$ which is stable under finite intersection.
%\item
%If $A$ is Tate, then the rational subsets are of the form $R(\frac Ts)$ with $s \in A$ and $T$ a finite subset of $A$ such that $T\cdot A = A$.
\end{assertionlist}
\end{thm}

\begin{proof}
\proofstep{(1)}
Let $I$ be the ideal of $A$ generated by $(A)^{oo}$. Then $\Cont(A)$ is a closed subspace of $\Spv(A,I)$. For all $a \in A$ the subset $\Spv(A,I)(\frac{a}{1}) = \Spv(A,I)(\frac{a,1}{1})$ is an open constructible subset of $\Spv(A,I)$ by Lemma~\ref{SpvAISpectral}. Thus we see that
\[
\Spa A = \Cont(A) \cap \bigcap_{a \in A^+}\Spv(A,I)(\frac{a}{1})
\]
is a pro-constructible subset of the spectral space $\Spv(A,I)$. Thus it is pro-constructible in $\Spv(A,I)$. In particular it is a spectral space.

\proofstep{(2)}
By Remark~\ref{UnderstandRational}~(1), one has $R(\frac{T}{s}) = \Spa A \cap \Spv(A,I)(\frac{T}{s})$. As the sets of the form $\Spv(A,I)(\frac{T}{s})$ form a basis of open quasi-compact subsets of $\Spv(A,I)$ and as $\Spa A$ is pro-constructible in $\Spv(A,I)$, the rational subsets form a basis of quasi-compact open subsets of $\Spa A$. Moreover, this basis is stable under finite intersection by Remark~\ref{UnderstandRational}~(4).
%
%\proofstep{(3)}
%This follows immediately from Remark~\ref{UnderstandRational}~(2).
\end{proof}

\begin{cor}
Let $A$ be an affinoid ring. Then a subset of $\Spa A$ is constructible if and only if it is in the Boolean algebra generated by the rational subsets.
\end{cor}

\begin{example}\label{AdicSpecField}
Let $K$ be a field, endowed with the topology induced by a microbial valuation $v$, and let $A(v)$ be its valuation ring. Then
\[
\Spa(K,A(v)) = \set{x \in \Spv K}{A(v) \subseteq A(x) \ne K},
\]
and the open subsets are the sets $U_w := \set{x \in \Spv K}{A(w) \subseteq A(x)}$, where $w \in \Spa(K,A(v))$.

In particular, $\Spa(K,A(v))$ consists of $h \in \NN$ points if and only if $v$ is of height $h$.
\end{example}

The next proposition is easy to check.

\begin{prop}\label{ClosedSub}
Let $A = (A,A^+)$ be an affinoid ring and let $\afr \subset A$ be an ideal. Let $A/\afr$ be the quotient affinoid ring (Definition~\ref{QuotientAffinoid}) and let $\pi\colon A \to A/\afr$ be the canonical homomorphism. Then $\Spa(\pi)\colon \Spa A/\afr \to \Spa A$ is a homeomorphism of $\Spa A/\afr$ onto the closed subset of points $x \in \Spa A$ with $\supp(x) \supseteq \afr$.
\end{prop}

Not all closed subspaces of $\Spa A$ are of this form.

%---------------------------------------------------

\subsection{Analytic Points}

\begin{defi}\label{DefAnalytic}
Let $A$ be a topological ring. A point $x \in \Cont A$ is called \emph{analytic} if $\supp x$ is not open in $A$.

If $A = \affoid A$ is an affinoid ring, then the subset of analytic points in $\Spa A$ is denoted by $(\Spa A)\an$, its complement in $\Spa A$ is denoted by $(\Spa A)\na$.
\end{defi}

Any valuation with open support is automatically continuous.

\begin{rem}\label{RemAnalytic}
Let $A = (A,A^+)$ be an affinoid ring.
\begin{ali}
\item\label{RemAnalytic1}
For $x \in \Cont(A)$ the following assertions are equivalent.
\begin{eli}
\item
$x$ has non-open support.
\item
There exists $a \in (A)^{oo}$ such that $x(a) \ne 0$.
\end{eli}
\item\label{RemAnalytic2}
$(\Spa A)\an$ is an open quasi-compact subset of $\Spa A$.
\item\label{RemAnalytic3}
If $A$ is Tate, then $\Cont A = (\Cont A)\an$ and $\Spa A = (\Spa A)\an$.
\item\label{RemAnalytic4}
In $(\Cont A)\an$ (and hence in $(\Spa A)\an$) there are no proper horizontal specializations. Thus all specializations are vertical.
\item\label{RemAnalytic5}
Every valuation $x$ in $(\Cont A)\an$ has height $\geq 1$ and is microbial.
\item\label{RemAnalytic6}
Let $x \in (\Cont A)\an$. Then $x$ has height $= 1$ if and only if $x$ is a maximal point of $(\Spa A)\an$ (i.e., the closure of $\{x\}$ in $(\Spa A)\an$ is an irreducible component).
\end{ali}
\end{rem}

\begin{proof}
Assertion~(1) follows from Lemma~\ref{OpenIdeal}, which in particular shows
\[
\set{\pfr \in \Spec A}{\text{$\pfr$ is open}} = V((A)^{oo}).
\]
Let $T \subseteq (A)^{oo}$ be a finite subset such that $T\cdot A$ is open in $A$ (e.g., if $T$ is a system of generators of an ideal of definition of a ring of definition of $A$). Then
\[
(\Spa A)\an = \bigcup_{t\in T} R(\frac{T}{t}).
\]
This shows~(2). Assertion~(3) follows from~(1).

Let $v$ be a continuous valuation on $A$ with value group $\Gamma$, and let $v\rstr{H}$ be a horizontal specialization, where $H$ is a proper convex subgroup of $\Gamma$ (containing $c\Gamma_v$). Then $\supp(v\rstr{H})$ is the union of $\set{a \in A}{v(a) < \gamma}$, where $\gamma$ runs through all elements of $\Gamma \setminus H$. Thus $v\rstr{H}$ is non-analytic. This shows~(4).

A continuous valuation $x$ of height $0$ is a trivial valuation whose support is $\set{a\in A}{x(a) < 1}$ and hence it is open. Assume that $x$ is analytic, i.e., $A/(\supp x)$ endowed with the quotient topology is not discrete. Then there exists by \eqref{RemAnalytic1} a topologically nilpotent element $\ne 0$ in $K(x)$. Thus the valuation $x$ on $K(x)$ is microbial. This shows~(5).

Let us prove~(6). By~(4) it remains to show that an analytic continuous valuation $x$ has height $1$ if there exists no proper vertical generization of $x$. But the vertical generizations are those within $\Spv K(x)$ (Remark~\ref{DefVertical}).
\end{proof}

\begin{proposition}\label{HeightOneSpa}
Let $A$ be an $f$-adic ring and let $x \in \Cont(A)\an$ be of height $1$. Then $x(a) \leq 1$ for all $a \in A^o$. In particular $x \in \Spa(A,A^+)\an$ for every ring of integral elements $A^+$.
\end{proposition}

\begin{proof}
Let $a \in A^o$ and assume that $x(a) > 1$. Choose $b \in A^{oo}$ with $x(b) \ne 0$ (possible by Remark~\ref{RemAnalytic}~\eqref{RemAnalytic1}). As $\Gamma_x$ has height $1$, it is archimedean (Proposition~\ref{Height1Groups}). Hence there exists $n \in \NN$ with $x(a^n) > x(b)^{-1}$, i.e., $x(a^nb) > 1$. But as $a$ is power-bounded, $a^nb \in A^{oo}$ and thus the continuity of $x$ implies that $x(a^nb) < 1$. Contradiction.
\end{proof}

\begin{rem}
Let $A = (A,A^+)$ be an affinoid ring and let $x \in \Spa A$.
\begin{assertionlist}
\item
Let $H \subsetneq \Gamma_x$ be a proper convex subgroup. Then the vertical generization $x_{/H}$ is a point of $\Spa A$.

Indeed, $x_{/H}$ is continuous by Remark~\ref{ContCofinal}~\eqref{ContCofinal2}, and if $x(a) \leq 1$ for all $a \in A^+$ then clearly $x_{/H}(a) \leq 1$ for all $a \in A^+$.
\item
If $x$ is analytic, it is microbial (Remark~\ref{RemAnalytic}~\eqref{RemAnalytic5}) and hence there exists by (1) a vertical generization of height $1$.
\item
If $x$ is non-analytic, then the vertical generization $x_{/\Gamma_x}$ is the trivial valuation with support $\supp(x)$. As $\supp(x)$ is an open prime ideal, $x_{/\Gamma_x}$ is continuous and hence a point of $\Spa A$.
\end{assertionlist}
\end{rem}

\begin{example}
Let $A$ be an affinoid ring and let $x \in \Spa A$. Endow $K(x) = \Frac A/(\supp x)$ with the topology induced by $x$. As $x$ is continuous, the canonical ring homomorphism $\gamma\colon A \to K(x)$ is continuous.

Assume that $x$ is analytic. By Remark~\ref{RemAnalytic}~\eqref{RemAnalytic5} the valuation $x$ on $K(x)$ is microbial. Thus $(K(x),A(x))$ is an affinoid ring and the condition $|f(x)| \leq 1$ for all $f \in A^+$ then implies, that $\gamma$ is a morphism of affinoid rings.
\end{example}

\begin{lem}\label{IsoOutsideOpen}
Let $A$ be an $f$-adic ring and let $B$ be an open subring of $A$. Denote by $f\colon \Spec A \to \Spec B$ and $g\colon \Spv A \to \Spv B$ be the morphisms induced by the inclusion $B \mono A$.
\begin{assertionlist}
\item\label{IsoOutsideOpen1}
Set $T := \set{\qfr \in \Spec B}{\text{$\qfr$ open}}$. Then $f^{-1}(T) = \set{\pfr \in \Spec A}{\text{$\pfr$ open}}$ and the restriction $f\colon (\Spec A)\setminus f^{-1}(T) \to (\Spec B) \setminus T$ is a homeomorphism.
\item\label{IsoOutsideOpen2}
$\Cont(A) = g^{-1}(\Cont(B))$.
\item\label{IsoOutsideOpen3}
$g$ yields by restriction a homeomorphism $\Cont(A)\an \iso \Cont(B)\an$ and one has for every $v \in \Cont(A)\an$ that $\Gamma_v = \Gamma_{g(v)}$.
\end{assertionlist}
\end{lem}

\begin{proof}
\proofstep{(1)}
As $B$ is open in $A$, a prime ideal $\pfr$ of $A$ is open if and only if $\pfr \cap B$ is open in $B$. This shows $f^{-1}(T) = \set{\pfr \in \Spec A}{\text{$\pfr$ open}}$.

Let $\qfr \subseteq B$ be a non-open prime ideal. Then there exists $s \in B^{oo}$ with $s \notin \qfr$. As $B$ is open, there exists for all $a \in A$ an integer $n \geq 1$ such that $s^na \in B$. This shows that the injective homomorphism $B_s \to A_s$ is also surjective. This shows the claim because $(\Spec B) \setminus T = \bigcup_{s \in B^{oo}} \Spec B_s$ and $f^{-1}(\Spec B_s) = \Spec A_s$.

\proofstep{(2)}
Clearly one has $\Cont(A) \subseteq g^{-1}(\Cont(B))$. Let $v$ be a valuation on $A$ such that $w := v\rstr{B}$ is continuous. If $\supp(v)$ is open in $A$, then $v$ is continuous. Assume that $\pfr := \supp(v)$ is not open in $A$ and set $\qfr := B \cap \pfr$. Then (1) shows that $B_{\qfr} \to A_{\pfr}$ is an isomorphism and in particular $\Gamma_v = \Gamma_w$. Hence the subgroup $\set{a \in A}{v(a) < \gamma}$ contains the open subgroup $\set{b \in B}{w(b) < \gamma}$ for all $\gamma \in \Gamma_v = \Gamma_w$. Therefore $v$ is continuous.

%$\pfr$ is the unique prime ideal of $A$ such that $\qfr := B \cap \pfr = \supp w$. Hence the image of $\pfr$ in $A/\qfr A$ is the nil radical. This shows that $\qfr A$ and $\pfr$ have the same radical.
%
%Now $B + \pfr$ is open in $A$ because it contains the open subgrup $B$. Hence we may consider $B/\qfr$ as an open subring of $A/\pfr$. To show the claim we may therefore replace $B$ by $B/\qfr$ and $A$ by $A/\pfr$ (with the induced topologies which are not discrete by assumption). Hence $\kappa(\qfr)

\proofstep{(3)} This follows from~(1) and the proof of~(2).
\end{proof}

%\begin{lemma}\label{ValHeight1Dominate}
%Let $A$ be a ring, let $\pfr \subsetneq \qfr$ be a prime ideal of $A$. Then there exists a valuation ring $R$ of height $1$ with $\Frac R = K := \Frac A/\pfr$ dominating $(A/\pfr)_{\qfr/\pfr}$. If $A$ is noetherian, then we may choose $R$ to be a discrete valuation ring.
%\end{lemma}
%
%\begin{proof}
%We may replace $A$ by $(A/\pfr)_{\qfr/\pfr}$ and can thus assume that $A$ is a local integral domain with maximal ideal $\qfr$ and $\pfr = 0$. Let $Z$ be the blow-up of $\Spec A$ in the closed subscheme $\{\pfr\}$. Then $Z$ is again integral, the canonical morphism $Z \to \Spec A$ is birational, and $D := \pi^{-1}(\{x\})$ is a divisor. Let $\eta$ be a maximal point of $D$ and let $B = \Oscr_{Z,\eta}$. Then $B$ is a local ring of dimension $1$ with $\Frac B = K$ dominating $A$. By \cite{Bou_AC}~Chap.~6, \S4, Exercise~6 there exists a valuation ring $R$ of $K$ of height $1$ with $B \subseteq R$. As we cannot have $\mfr_R \cap B = \{0\}$ (because then the inclusion $B \mono R$ would induce an inclusion $K \mono R$), $R$ dominates $B$.
%\end{proof}

\begin{lemma}\label{ExistAnalyticPoint}
Let $A$ be a complete affinoid ring. Let $\pfr$ be a non-open prime ideal of $A$. Then there exists an analytic point $x \in \Spa A$ of height $1$ such that $\supp x \supseteq \pfr$.

If $A$ has a noetherian ring of definition, we may assume in addition that $x$ is a discrete valuation and that $\supp x = \pfr$.
\end{lemma}

\begin{proof}
Let $A_0$ be a ring of definition of $A$ and let $I \subseteq A_0$ be a finitely generated ideal of definition. Set $\pfr_0 := \pfr \cap A_0$. As $\pfr_0$ is not open, one has $I \not\subseteq \pfr_0$.

\proofstep{Construction of $x$ in the general case}
Let $\mfr$ be a maximal ideal of $A_0$ containing $\pfr_0$. As $A_0$ is complete for the $I$-adic topology, one has $I \subseteq \mfr$ (Proposition~\ref{DefIdealMaximal}). In particular $\pfr_0 \subsetneq \mfr$. Let $u_0$ be a valuation of $\Frac A_0/\pfr_0$ such that its valuation ring dominates $(A_0/\pfr_0)_{\mfr/\pfr_0}$. Let $u$ be the corresponding valuation of $A_0$ such that $\supp u = \pfr_0$. Then $u(a) < 1$ for all $a \in \mfr$.

Let $r\colon \Spv(B_0) \to \Spv(B_0,I)$ be the retraction from \eqref{DefRetraction}. As $u(a) < 1$ for all $a \in I$ one also has $r(u) < 1$ for all $a \in I$ by definition of $r(u)$. Hence $r(u) \in \Cont(B_0)$ by Theorem~\ref{ContSpectral}. Moreover by Lemma~\ref{SpvAISpectral}~\eqref{SpvAISpectral3} one has $I \not\subseteq \supp(r(u))$. Hence $r(u)$ is not analytic. By Lemma~\ref{IsoOutsideOpen}~\eqref{IsoOutsideOpen3} there exists a unique continuous analytic valuation $v$ on $A$ such that $v\rstr{A_0} = r(u)$. Then $\supp(v) \cap A_0 = \supp r(u) \supseteq \supp u = \pfr_0$ because $r(u)$ is a (horizontal) specialization of $u$.

By Remark~\ref{RemAnalytic}~\eqref{RemAnalytic5} $v$ is microbial and thus there exists a vertical generalization $x$ of $v$ of height $1$ within $\Spv A$ (Remark~\ref{DefVertical}). Then $x$ is continuous Remark~\ref{ContCofinal}~\eqref{ContCofinal2}. As $x$ is of height $1$, Proposition~\ref{HeightOneSpa} shows that $x \in \Spa A$.

\proofstep{Construction of $x$ if $A_0$ is noetherian}
Let $A_0$ be noetherian. Let $J$ be the image of $I$ in $A_0/\pfr_0$. Let $\pi\colon Z \to X := \Spec (A_0/\pfr_0)$ be the blow-up of $X$ in $Y := V(J)$. As $J$ contains a regular element, $\pi$ is projective and birational, and $\pi^{-1}(Y)$ is a non-empty divisor of $Z$. As $A_0$ is noetherian, $Z$ is noetherian. Thus if $z$ is a maximal point of $\pi^{-1}(Y)$, then $\Oscr_{Z,z}$ is a local noetherian ring of dimension $1$. By the theorem of Krull-Akizuki its normalization is a discrete valuation ring $R$. We have $\Frac R = \Frac A_0/\pfr_0$ because $\pi$ is birational. Let $w$ be the valuation on $A_0$ with $\supp w = \pfr_0$ corresponding to $R$. Then $w(a) < 1$ for all $a \in I$. As $\Gamma_w$ is of height $1$ (in fact isomorphic to the totally ordered group $\ZZ$) this implies that $w(a)$ is cofinal for $\Gamma_w$ for all $a \in I$. Hence $w$ is continuous (Remark~\ref{ContCofinal}). As $I \not\subseteq \pfr_0$, $w$ is analytic. Therefore Lemma~\ref{IsoOutsideOpen}~\eqref{IsoOutsideOpen3} shows that there exists a unique continuous valuation $x \in \Cont(A)$ with $x\rstr{A_0} = w$. Moreover $\Gamma_w = \Gamma_x$ and hence $x$ is a discrete valuation. Again Proposition~\ref{HeightOneSpa} shows that $x \in \Spa A$.
\end{proof}

\begin{lemma}\label{AdicSpectral}
Let $\varphi\colon A \to B$ be a continuous homomorphism between affinoid rings and let $f = \Spa(\varphi)\colon X := \Spa B \to Y := \Spa A$ be the attached continuous map.
\begin{assertionlist}
\item\label{AdicSpectral1}
$f(X\na) \subseteq Y\na$. If $\varphi$ is adic, then $f(X\an) \subseteq Y\an$.
\item\label{AdicSpectral2}
If $B$ is complete and $f(X\an) \subseteq Y\an$, then $\varphi$ is adic.
\item\label{AdicSpectral3}
If $\varphi$ is adic, then for every rational subset $V$ of $Y$ the preimage $f^{-1}(V)$ is rational. In particular $f$ is spectral.
\end{assertionlist}
\end{lemma}

\begin{proof}
\proofstep{(1)}
The first assertion of~(1) is clear and the second follows from Remark~\ref{RemAnalytic}~\eqref{RemAnalytic1}.

\proofstep{(2)}
Assume that $B$ is complete and that $\varphi$ is not adic. We want to show that $f(X\an) \not\subseteq Y\an$. Let $(A_0,I)$ and $(B_0,J)$ pairs of definition of $A$ and $B$, respectively, such that $\varphi(A_0) \subseteq B_0$ and $\varphi(I) \subseteq J$. Since $\varphi$ is not adic, $\varphi(I)B_0$ and $J$ do not have the same radical. Therefore there exists a prime ideal $\pfr$ of $B_0$ such that $\varphi(I) \subseteq \pfr$ and $J \not\subseteq \pfr$ (i.e., $\pfr$ is non-open). By Lemma~\ref{ExistAnalyticPoint} there exists a point $x \in X\an$ such that $(\supp x) \cap B_0 \supseteq \pfr$. But $\supp f(x) = \varphi^{-1}(\supp x) \supseteq I$ and hence $f(x) \notin Y\an$.

%Since $B_0$ is complete in the $J$-adic topology, there exists a prime ideal $\qfr$ with $\pfr \subseteq \qfr$ and $J \subseteq \qfr$ (\ref{???}). Let $u$ be a continuous valuation of $B$ such that such that $\supp u = \pfr$ and such that the valuation ring $B_0(u)$ of $u$ dominates the local ring $(B_0/\pfr)_{\qfr/\pfr}$. $J \not\subseteq \supp u_0$. Let $r\colon \Spv B_0 \to \Spv(B_0,J)$ be the corresponding retraction~\eqref{DefRetraction}. Then $r(u)$ is a continuous valuation of $B_0$ with $J \subsetneq \supp(r(u))$. By Lemma~\ref{IsoOutsideOpen} there exists a continuous valuation $v$ on $B$ with $r(u) = v\rstr{B}$.
%
%We first claim that if $\pfr \subseteq \qfr$ are prime ideals of $B$ and $I \subseteq \qfr$, then $I \subseteq \pfr$. Assume that $I \subsetneq \pfr$. Let $u$ be a valuation of $B$ such that $\supp u = \pfr$ and such that the valuation ring $B(u)$ of $u$ dominates the local ring $(B/\pfr)_{\qfr/\pfr}$. Let $r\colon \Spv A \to \Spv (A,I)$ be the retraction~\eqref{DefRetraction}. Then $r(u)$ is a continuous valuation with $I \subsetneq \supp(r(u))$. By Lemma~\ref{IsoOutsideOpen} there exists a continuous valuation $v$ on $A$ with $r(u) = v\rstr{B}$. The vertical generization $w$ of $v$ of height $1$ is then an element of $X\an$. Contradiction.

\proofstep{(3)}
Let $s  \in A$ and $T \subseteq A$ a finite subset such that $T\cdot A$ is open in $A$. If $\varphi$ is adic, then $\varphi(T)\cdot B$ is open in $B$. Hence $f^{-1}(R(\frac{T}{s}))$ is the rational subset $R(\frac{\varphi(T)}{s})$ of $\Spa B$.
\end{proof}

%---------------------------------------------------

\subsection{First properties of $\Spa A$}

\begin{lemma}\label{SubgroupCompletionFAdic}
Let $A$ be an $f$-adic ring, let $\iota\colon A \to \Ahat$ be its completion. Then under the bijection (Example~\ref{CompletionSubgroup})
\begin{align*}
\{\text{open subgroups $G$ of $A$}\} &\bijective \{\text{open subgroups $G'$ of $\Ahat$}\}, \\
G \sends \Ghat = \widehat{\iota(G)}, &\qquad \iota^{-1}(G') \lsends G'
\end{align*}
the following subgroups correspond to each other.
\begin{assertionlist}
\item
$A^o$ and $(\Ahat)^o$.
\item
$A^{oo}$ and $(\Ahat)^{oo}$.
\item
Rings of definition of $A$ and rings of definition of $\Ahat$.
\item
Rings of integral elements of $A$ and rings of integral elements of $\Ahat$.
\end{assertionlist}
\end{lemma}

\begin{proof}
\cite{Hu_Habil}~2.4.3.
\end{proof}

If $A$ is an affinoid ring, then $\widehat{A^+}$ is a ring of integral elements in $\widehat{A}$ and we set $\Ahat := (\widehat{A},\widehat{A^+})$. This is an affinoid ring by Lemma~\ref{SubgroupCompletionFAdic}.

\begin{prop}\label{SpaCompletion}
Let $A$ be an affinoid ring. The canonical map $\Spa \Ahat \to \Spa A$ is a homeomorphism which maps rational subsets to rational subsets.
\end{prop}

\begin{proof}
\cite{Hu_Cont}~Prop.~3.9
\end{proof}

\begin{prop}\label{SpaEmpty}
Let $A = (A,A^+)$ be an affinoid ring.
\begin{ali}
\item
The following assertions are equivalent.
\begin{equivlist}
\item
$\Spa A = \emptyset$.
\item
$\Cont(A) = \emptyset$.
\item
$A/\cl{0} = 0$.
\end{equivlist}
\item
The following assertions are equivalent.
\begin{equivlist}
\item
$(\Spa A)\an = \emptyset$.
\item
$(\Cont A)\an = \emptyset$.
\item
The Hausdorff ring $A/\cl{0}$ attached to $A$ has the discrete topology.
\end{equivlist}
\end{ali}
\end{prop}

The equivalent conditions in~(1) (resp.~in~(2)) imply that $\Ahat = 0$ (resp.~that $\Ahat$ is discrete).

\begin{proof}
First note that in (1) and (2) the implications ``(iii) $\implies$ (ii) $\implies$ (i)'' are clear. We set $X := \Spa A$. Let $B$ be a ring of definition of $A$ and let $I$ be a finitely generated ideal of definition of $B$. 

\proofstep{(2)}
Assume that $X\an = \emptyset$. 

We first claim that if $\pfr \subseteq \qfr$ are prime ideals of $B$ and $I \subseteq \qfr$, then $I \subseteq \pfr$. The proof is similar as the proof of the general case of Lemma~\ref{ExistAnalyticPoint}. Assume that $I \not\subseteq \pfr$. Let $u$ be a valuation of $B$ such that $\supp u = \pfr$ and such that the valuation ring $B(u)$ of $u$ dominates the local ring $(B/\pfr)_{\qfr/\pfr}$. Let $r\colon \Spv B \to \Spv (B,I)$ be the retraction~\eqref{DefRetraction}. Then $r(u)$ is a continuous valuation with $I \subsetneq \supp(r(u))$. By Lemma~\ref{IsoOutsideOpen} there exists a continuous valuation $v$ on $A$ with $r(u) = v\rstr{B}$. The vertical generization $w$ of $v$ of height $1$ is then an element of $X\an$. Contradiction.

We now prove~(2) using the claim. Set $S := 1 + I$ and let $\varphi\colon B \to C := S^{-1}B$ the localization. As in $\Spec C$ every prime ideal specializes to a prime ideal containing $\varphi(I)$, the claim implies that $\varphi(I)$ is contained in every prime ideal of $C$. Hence there exists $n \in \NN$ with $\varphi(I^n) = 0$ because $I$ is finitely generated. This means that there exists $i \in I$ with $(1+i)I^n = 0$ in $B$. But then $I^k = I^n$ for all $k \geq n$. This shows that $A/\cl{0}$ is discrete.

\proofstep{(1)}
Assume that $X = \emptyset$. Then by~(2) the ideal $\cl{0}$ is open in $A$. Then every trivial valuation whose support contains $\cl{0}$ is an element of $X$. This shows that $A$ has no prime ideal containing $\cl{0}$, i.e. $A = \cl 0$.
\end{proof}

\begin{remark}
Let $A$ be a complete affinoid ring such that $(\Spa A)\an = \emptyset$. Then $A$ is discrete by Proposition~\ref{SpaEmpty}. This can also be seen as follows. Let $(B,I)$ be a pair of definition of $A$. By Lemma~\ref{ExistAnalyticPoint} all prime ideals of $B$ are open. Thus $I$ is contained in all prime ideals of $B$ and hence in the nil radical of $B$. As $I$ is finitely generated, this shows $I^n = 0$.
\end{remark}

\begin{prop}\label{MaxComplete}
Let $A$ be a complete affinoid ring and let $\mfr \subset A$ be a maximal ideal. Then $\mfr$ is closed and there exists $v \in \Spa A$ with $\supp v = \mfr$.
\end{prop}

\begin{proof}
The set of $(A)^{oo}$ of topologically nilpotent elements of an $f$-adic ring is open (as the union of all definition ideals of all definition rings). Hence $1 + A^{oo}$ is open in $A$. As $A$ is complete, $1 + A^{oo}$ is a subgroup of the group of units of $A$. This shows that $A^\times$ is open in $A$. If $\mfr$ was not closed, then its closure (being an ideal) would be $A$, which contradicts that $A \setminus A^{\times}$ is closed. Hence $\mfr$ is closed.

Thus $A/\mfr$ is a Hausdorff and $\set{v \in \Spa A}{\supp(v) = \mfr} =\Spa A/\mfr \ne \emptyset$ by Proposition~\ref{SpaEmpty}.
\end{proof}

\begin{prop}\label{FunctionInvertible}
Let $A = \affoid A$ be an affinoid ring, and $f \in A$.
\begin{ali}
\item
Then $|f(x)| \leq 1$ for all $x \in \Spa A$ if and only if $f \in A^+$.
\item
Assume that $A$ is complete. Then $f$ is a unit if and only if $|f(x)| \ne 0$ for all $x \in \Spa A$.
\end{ali}
\end{prop}

\begin{proof}
Assertion~(1) is simply a reformulation of Proposition~\ref{ClassifyRingIntegral}~(1). Assertion~(2) is a reformulation of Proposition~\ref{MaxComplete}.
\end{proof}

\begin{cor}\label{OpenCovering}
Let $A$ be a complete affinoid ring and let $T \subseteq A$ be a finite subset. Then the ideal generated by $T$ equals $A$ if and only if for all $x \in \Spa A$ there exists $t \in T$ with $|t(x)| \ne 0$. If these equivalent conditions are satisfied, $(R(\frac{T}{t}))_{t \in T}$ is an open covering of $\Spa A$.
\end{cor}

\begin{proof}
If $T\cdot A$ is contained in a maximal ideal $\mfr$ then there exists $x \in \Spa A$ with $|t(x)| = 0$ for all $t \in T$ by Proposition~\ref{MaxComplete}. Conversely, if $\sum_i a_it_i = 1$ with $t_i \in T$ and $a_i \in A$ and there exists $x \in \Spa A$ with $|t_i(x)| = 0$ for all $i$. Then $1 \leq \max\{|t_i(x)||a_i(x)|\} = 0$. Contradiction.

Assume now that for all $x \in X$ there exists $t \in T$ with $|t(x)| \ne 0$. Let $t_x \in T$ such that $|t_x(x)| = \max\set{|t(x)|}{t \in T}$. Then $x \in R(\frac{T}{t_x})$.
\end{proof}

\begin{lem}\label{RationalRefinement}
Let $A$ be a complete affinoid ring and let $(V_j)_{j\in J}$ be an open covering of $\Spa A$. Then there exist $f_0,\dots,f_n \in A$ generating $A$ as an ideal such that for all $i \in \{0,\dots,n\}$ the rational subset $R(\frac{f_0,\dots,f_n}{f_i})$ is contained in some $V_j$.
\end{lem}

By Corollary~\ref{OpenCovering} the rational subsets $R(\frac{f_0,\dots,f_n}{f_i})$ form an open covering of $\Spa A$.

\begin{proof}
\cite{Hu_Gen}~Lemma~2.6.
\end{proof}

\begin{rem}\label{DevissageRational}
Let $A = (A,A^+)$ be a Tate affinoid ring, $s \in A$ and $T = \{t_1,\dots,t_n\} \subset A$ a finite subset such that $T\cdot A$ is open in $A$. Let $U = \sett{x \in \Spa A}{$x(t_i) \leq x(s) \ne 0$ for $i = 1,\dots,n$}$ be the corresponding rational subset. Since $U$ is quasi-compact, there exists by Corollary~\ref{TateZeroNeighbor} a unit $u \in A^{\times}$ such that $|u(x)| < |s(x)|$ for all $x \in U$. We set
\[
X_0 := \set{x \in \Spa A}{1 \leq x(\frac{s}{u})}
\]
(a rational subset of $\Spa A$). Then $x(s) \ne 0$ for all $x \in X_0$, thus $s\rstr{X_0}$ is a unit in $A \langle \frac{su^{-1}}{1}\rangle$ (Proposition~\ref{FunctionInvertible}). Define now inductively rational subsets $X_1,\dots,X_n$ of $\Spa A$ by
\[
X_i := \set{x \in X_{i-1}}{x(\frac{t_i}{s}) \leq 1}, \qquad i =1,\dots,n.
\]
Thus we obtain a chain of rational subsets $\Spa A \supseteq X_0 \supseteq X_1 \supseteq \dots \supseteq X_n = U$.
\end{rem}

%-------------------------------------------------------------------

\subsection{Examples of adic spectra}\label{ExampleAdicSpec}

\begin{defi}
A \emph{non-archimedean field} is a topological field $k$ whose topology is induced by a valuation of height $1$.

An \emph{affinoid $k$-algebra topologically of finite type} (\emph{tft} for short) is an affinoid ring $A = \affoid A$ such that $A$ is a $k$-algebra topologically of finite type and such that $A^+ = A^o$.
\end{defi}

Every non-archimedean field has then a continuous valuation $\abs\colon k \to \RR^{\geq 0}$ which is unique up to passing from $\abs$ to $\abs^t$ for some $t \in \RR^{>0}$. Moreover, $k^o = \set{c \in k}{|c| \leq 1}$ and $k^{oo} = \set{c \in k}{|c|<1}$.

\begin{example}
Let $k$ be a non-archimedean field and fix the norm $\abs\colon k \to \RR^{\geq 0}$. Let $|k^\times|$ be its value group. For $x \in k$ and $r \in \RR^{\geq 0}$ we set $D(x,r) := \set{y \in k}{|y-x| \leq r}$. Such subsets are called discs. Note that for any $x' \in D(x,r)$ one has $D(x',r) = D(x,r)$. We also set $D^o(x,r) := \set{y \in k}{|x-y| < r}$. Then $k^{oo} = D^0(0,1)$ is the maximal ideal of $k^o$ and we denote by $\kappa = k^o/k^{oo}$ the residue field. We assume that $k$ is algebraically closed and complete. Then $\kappa$ is also algebraically closed.

Let $A = k\langle t\rangle$ be the ring of convergent power series in one variable and $A^+ = A^o = k^o \langle t\rangle$. In $X = \Spa(A,A^o)$ there are $5$ different kind of points. We visualize $X$ as a tree.
\begin{assertionlist}
\item[(1)]
The classical points (end points): Let $x \in k^o$, i.e., $x \in D(0,1)$. Then for any $f \in k\langle t\rangle$ we can evaluate $f$ at $x$ to get a map $k\langle t\rangle \to k$, $f \sends f(x)$. Composing with the norm on $k$, one gets a valuation $f \sends |f(x)|$ on $k\langle t\rangle$, which is continuous and $\leq 1$ for all $f \in R^+$. These classical points correspond to the maximal ideals of $A$: If $\mfr \subset A$ is a maximal ideal, then $A/\mfr = k$ (\cite{BGR}~6.1.2~Corollary~3; because $k$ is algebraically closed) and $\mfr$ is of the form $(t - x)$ for a unique $x \in D(0,1)$.

These are the end points of the branches.
\item[(2),(3)]
Points on the limbs: Let $0 \leq r \leq 1$ be some real number and $x \in k^o$. Then
\[
x_r\colon f = \sum_n a_n(t-x)^n \sends \sup\limits_n |a_n|r^n = \sup_{\substack{y \in k^o,\\ |y-x|\leq r}} |f(y)|
\]
is a point of $X$. It depends only on $D(x,r)$. For $r = 0$ it agrees with the classical point defined by $x$, for $r = 1$ the disc $D(x,1)$ is independent of $x$ and we obtain the Gau{\ss} norm $\sum a_nt^n \sends \sup_n |a_n|$ as the ``root'' of the tree.

If $r \in |k^\times|$, then the point $x_r$ corresponding to $D(x,r)$ is a branching point. These are the points of type~(2), the other points are of type~(3).
\item[(4)]
Dead ends: Let $D_1 \supset D_2 \supset \cdots$ be a sequence of closed discs such that $\bigcap_i D_i = \emptyset$ (such sequences exist if $k$ is not spherically complete, e.g. if $k$ is $\CC_p = \widehat{\overline{\QQ_p}}$). Then
\[
f \sends \inf_i \sup_{y \in D_i}|f(y)|
\]
is a point of $X$.
\item[(5)]
Finally there are some valuations of height $2$. Let $x \in k^o$ and fix a real number $r$ with $0 < r \leq 1$ and let $\Gamma_{<r}$ be the abelian group $\RR^{>0} \times \gamma^{\ZZ}$ endowed with the unique total order such that $r' < \gamma < r$ for all $r' < r$. Then
\[
x_{< r}\colon f = \sum_n a_n(t-x)^n \sends \max_n |a_n|\gamma^n \in \Gamma_{<r}
\]
is a point of $X$ which depends only on $D^0(x,r)$. For $0 < r < 1$ let $\Gamma_{>r}$ be the abelian group $\RR^{>0} \times \gamma^{\ZZ}$ endowed with the unique total order such that $r' > \gamma > r$ for all $r' > r$. Then
\[
x_{> r}\colon f = \sum_n a_n(t-x)^n \sends \max_n |a_n|\gamma^n \in \Gamma_{>r}
\]
is a point of $X$ which depends only on $D(x,r)$ (and hence $x_{>r} = x'_{>r}$ if and only if $x_r = x'_r$).

If $r \notin |k^{\times}|$, then $x_{<r} = x_{>r} = x_r$. But for each point $x_r$ of type~(2) this gives exactly one additional point for each ray starting from $x_r$. The points $x_{<r}$ correspond to rays towards the classical points, and the point $x_{>r}$ corresponds to the ray towards the Gau{\ss} point.
\end{assertionlist}
All points except the points of type~(2) are closed. The closure of such a point $x_r$ consists of $x_r$ and all points of type~(5) around it. The closure of the Gau{\ss} point is homeomorphic to the scheme $\AA^1_\kappa$ with the Gauss point as generic point. The closure of all other points of type~(2) is homeomorphic to the scheme $\PP^1_\kappa$.

The subspace, consisting only of points of type (1) -- (4), is the Berkovich space $X^{\rm an}$ attached to $k\langle t \rangle$. Note that $X^{\rm an}$ is Hausdorff.

In case~(5) one also could have defined a valuation $x_{>1}$ (note that if $1 < \gamma < r'$ for all $r' > 1$, then $\lim_n |a_n|\gamma^n = 0$ if $a_n \to 0$, thus $\max_n |a_n|\gamma^n$ is still defined). This is a continuous valuation on $k\ang t$ but $x_{>1}(t) > 1$, thus it is not a point in $\Spa(A,A^o)$. In fact, one can show that $\Spa(A,A^o) \cup \{x_{>1}\} = \Spa(A,A^+)$, where $A^+$ is the smallest ring of integral elements of A which contains $k^o$, i.e. the integral closure of $A^{oo} + k^o$.
\end{example}

\begin{example}\label{AdditionalPoint}
Let $k$ be a complete non-archimedean field and let $\abs\colon k \to \RR^{\geq 0}$ be a height $1$ valuation. Let $A$ be a $k$-algebra topologically of finite type over $k$. Every maximal ideal element $\mfr$ in
\[
\Max(A) := \set{\mfr \subset A}{\text{$\mfr$ maximal ideal}}
\]
defines a point of $X := \Spa(A,A^o)$ (again called ``classical point''): For $x \in \Max(A)$ the quotient $A/x$ is a finite extension of $k$ by Proposition~\ref{MaximalFinite}. Thus there exists a unique extension of $\abs$ to $A/x$. Its composition with $A \to A/\mfr$ is denoted by $\abs_x$ and it is a point of $X$. Via $x \sends \abs_x$ we consider $\Max(A)$ as a subset of $X$. We also set
\[
L_A := \set{v\in \Spv A}{\text{$v(a) \leq 1$ for all $a \in A^o$ and $v(a) < 1$ for all $a \in A^{oo}$}}.
\]
Then
\[
\Max(A) \subseteq \Spa(A,A^o) \subseteq L_A \subseteq \Spv A.
\]
One can show that $L_A$ is the closure of $\Max A$ in the constructible topology of $\Spv(A)$.

Instead of $A^o$ we can also consider the following ring of integral elements $A^+$ which is defined to be the integral closure of $k^o + A^{oo}$. Then $A^+$ is the smallest ring of integral elements of $A$ which contains $k^o$. Then $X_c := \Spa(A,A^+)$ contains $X$ as an open dense subset.
\end{example}

\begin{example}[cf.~\cite{Hu_Habil}~(3.7.1)]\label{AdicExample}
Let $A$ be an adic ring with a finitely generated ideal of definition $I$. We consider
\[
X = \Spa (A,A) = \set{x \in \Cont(A)}{\text{$x(f) \leq 1$ for all $f \in A$}}.
\]
Let $Y = X_{\rm triv}$ be the subspace of trivial valuations. Such a trivial valuation is continuous if and only if its support is an open prime ideal. Thus we can identify $Y = \Spf A := \set{\pfr \in \Spec A}{\text{$\pfr$ open in $A$}}$. It follows from Remark~\ref{SpvSpec} that $Y$ carries also the subspace topology of $\Spec A$. 

The subspace $Y$ is pro-constructible, in particular $Y$ is spectral and the inclusion $Y \mono X$ is spectral. There is a spectral retraction $r\colon X \to Y$, given by $x \sends x\rstr{c\Gamma_x}$.

Once we have defined the structure sheaf $\Oscr_X$ on $X$, we will see that attaching to a complete noetherian adic ring $A$ the adic space $\Spa (A,A)$ yields an embedding $t$ of the category of noetherian formal affine schemes into the category of adic spaces. The inverse functor (on the essential image of $t$) is given by attaching to $(X,\Oscr_X)$ the topologically ringed space $(X_{\rm triv},r_*\Oscr_X)$.
\end{example}

%----------------------------------------------------------

\section{Adic spaces}

\subsection{The presheaf $\Oscr_X$}

We recall the following special case of Definition~\ref{TopLoc}. Let $A$ be an $f$-adic ring, $s \in A$ and let $T = \{t_1,\dots,t_n\} \subseteq A$ be a finite subset such that $T\cdot A$ is open in $A$. Then there exists on $A_s$ a non-archimedean ring topology making it into a topological ring
\[
A(\frac{T}{s}) = A(\frac{t_1,\dots,t_n}{s})
\]
such that $\set{\frac{t}{s}}{t \in T}$ is power-bounded in $A(\frac{T}{s})$ and such that $A(\frac{T}{s})$ and the canonical homomorphism $\varphi\colon A \to A(\frac{T}{s})$ satisfy the following universal property. If $B$ is a non-archimedean topological ring and $f\colon A \to B$ is a continuous homomorphism such that $f(s)$ is invertible in $B$ and such that the set $\set{f(t)f(s)^{-1}}{t \in T}$ is power-bounded in $B$, then there exists a unique continuous ring homomorphism $g\colon A(\frac{T}{s}) \to B$ with $f = g \circ \varphi$.

The completion of $A(\frac{T}{s})$ is denoted by $A\langle\frac Ts\rangle = A\langle\frac{t_1,\dots,t_n}{s}\rangle$.

The ring $A(\frac{T}{s})$ is constructed as follows. Let $(A_0,I)$ be a pair of definition. In the localization $A_s$ consider the subring $D$ generated by $A_0$ and $E := \{\frac{t_1}{s},\dots,\frac{t_n}{s}\}$, i.e.,
\[
D = A_0[\frac{t_1}{s},\dots,\frac{t_n}{s}].
\]
Then $(I^n\cdot D)_n$ is a fundamental system of neighborhoods of $0$ in $A(\frac{T}{s})$. 

The construction also shows that $A(\frac{T}{s})$ is again an $f$-adic ring. The pair $(D,I\cdot D)$ is a pair of definition of $A(\frac{T}{s})$.

Now let $A^+$ be a ring of integral elements of $A$, i.e. $A = (A,A^+)$ is an affinoid ring. Let $C$ be the integral closure of $A^+[\frac{t_1}{s},\dots,\frac{t_n}{s}]$ in $A_s$. Then $C$ is a ring of integral elements in $A(\frac{T}{s})$. We denote the affinoid ring $(A(\frac{T}{s}),C)$ simply by $A(\frac{T}{s})$ and its completion by $A\langle\frac{T}{s}\rangle$.

The canonical ring homomorphism $A \to A_s$ is a continuous homomorphism of affinoid rings $h\colon A \to A(\frac{T}{s})$, one has $h(s) \in A(\frac{T}{s})^\times$, and $\frac{h(t_i)}{h(s)}) \in C = A(\frac{T}{s})^+$. It is universal for these properties.

In a similar way, the canonical continuous homomorphism of affinoid rings $\rho\colon A \to A\langle\frac{T}{s}\rangle$ is universal for continuous homomorphism of affinoid rings $\varphi\colon A \to B$, with $B$ complete, $\varphi(s) \in B^\times$ and $\frac{\varphi(t_i)}{\varphi(s)} \in B^+$. But in fact $\rho$ has also a more geometric universal property.

\begin{lem}\label{RationalUniversal}
Let $U = R\langle\frac{T}{s}\rangle \subseteq \Spa A$ be the corresponding rational subset. Then $\Spa\rho\colon \Spa A\langle\frac{T}{s}\rangle \to \Spa A$ factors through $U$ and whenever $\varphi\colon A \to B$ is a continuous homomorphism from $A$ to a complete affinoid ring $B$ such that $\Spa(\varphi)$ factors through $U$, then there exists a unique continuous ring homomorphism $\psi\colon A\langle\frac{T}{s}\rangle \to B$ such that $\psi \circ \rho = \varphi$.
\end{lem}

\begin{proof}
The definition of $A\langle\frac{T}{s}\rangle$ shows that for $x \in \Spa A\langle\frac{T}{s}\rangle$ one has $|\varphi(t)|_x \leq |\varphi(s)|_x$ for all $t \in T$. This means $\Spa(\rho)(x) \in U$.

As $\Spa(\varphi)$ factors through $U$, we have $|\varphi(t)|_w \leq |\varphi(s)|_w \ne 0$ for all $w \in \Spa B$ and for all $t \in T$. This implies $\varphi(s) \in B^\times$ by Proposition~\ref{FunctionInvertible}.
%Then there exists a maximal ideal $\mfr \subset B$ containing $\varphi(s)$. By Proposition~\ref{MaxComplete} there exists $w \in \Spa B$ with $\supp w = \mfr$ which contradicts $\varphi(s)|_w \ne 0$. Thus $\varphi(s) \in B^\times$. 
Moreover, for all $w \in \Spa B$ we have $|\frac{\varphi(t)}{\varphi(s)}|_w \leq 1$. This implies $\frac{\varphi(t)}{\varphi(s)} \in B^+$ by Proposition~\ref{FunctionInvertible}. Thus the claim follows from the universal property of $A \to A\langle\frac{T}{s}\rangle$.
\end{proof}

\begin{prop}\label{RationalFunction}
Let $A = (A,A^+)$ be an affinoid ring. Let $s,s' \in A$ and $T, T' \subseteq A$ finite subsets such that $T\cdot A$ and $T'\cdot A$ are open. Let $U = R(\frac{T}{s})$ and $U' = R(\frac{T'}{s'})$ be the corresponding rational subsets. Let $\rho\colon A \to A\ang{\frac{T}{s}}$ and $\rho'\colon A \to A\ang{\frac{T'}{s'}}$ be the canonical continuous homomorphisms of affinoid rings.
\begin{assertionlist}
\item
If $U' \subseteq U$, then there exists a unique continuous homomorphism $\sigma\colon A\ang{\frac{T}{s}} \to A\ang{\frac{T'}{s'}}$ such that $\sigma \circ \rho = \rho'$.
\item
The map $\Spa(\rho)\colon \Spa A\langle\frac{T}{s}\rangle \to \Spa A$ is a homeomorphism of $\Spa A\langle\frac{T}{s}\rangle$ onto $R(\frac{T}{s})$, and it induces a bijection between rational subsets of $\Spa A\langle\frac{T}{s}\rangle$ and rational subsets of $\Spa A$ contained in $R(\frac{T}{s})$.
\end{assertionlist}
\end{prop}

\begin{proof}
(1)~follows immediately from Lemma~\ref{RationalUniversal}. Let us prove~(2). We set $j := \Spa(\rho)$. We factorize $\rho$ into
\[
A \ltoover{\rho'} A(\frac{T}{s}) \ltoover{\iota} A\langle\frac{T}{s}\rangle.
\]
As $\rho'$ is adic by definition and $\iota$ is clearly adic, the composition $\rho$ is adic by Proposition~\ref{CompAdic}. Thus for every rational subset $V$ of $\Spa A$ its inverse image $j^{-1}(V)$ is a rational subset of $\Spa A(\frac{T}{s})$ (Lemma~\ref{AdicSpectral}). Moreover $\Spa (\iota)$ is a homeomorphism mapping rational subsets to rational subsets by Proposition~\ref{SpaCompletion}. Thus it remains to show that $j' := \Spa(\rho')$ is a homeomorphism from $\Spa A(\frac{T}{s})$ onto $U$ mapping rational subsets to rational subsets.

Write $T = \{t_1,\dots,t_n\}$. A valuation $v$ on $A$ extends (necessarily uniquely) to $A_s$ if and only if $v(s) \ne 0$. Moreover, if $v$ is continuous, such an extension is continuous with respect to the topology defined on $A_s = A(\frac{T}{s})$ if and only if $v(t_i) \leq v(s)$ for all $i = \{1,\dots,n\}$. Finally, it satisfies $v(f) \leq 1$ for all $f \in A(\frac{T}{s})^+ = A^+[\frac{t_1}{s},\dots,\frac{t_n}{s}]^{\rm int}$ if and only $v(f) \leq 1$ for all $f \in A^+$ and $v(t_i) \leq v(s)$ for all $i$. This shows that $j'$ is injective and that the image of $j'$ is $R(\frac{T}{s})$. Thus it remains to show that $j'$ maps rational subsets to rational subsets.

Let $V = R(\frac{g_1,\dots,g_m}{r})$ be a rational subset of $\Spa A(\frac{T}{s})$ for $r,g_1,\dots,g_m \in A_s = A(\frac{T}{s})$. Multiplying $r,g_1,\dots,g_m$ with a suitable power of $s$ we may assume that all these elements lie in the image of $\rho'$, say $r = \rho'(q)$ and $\{g_1,\dots,g_m\} = \rho'(H)$ for some $q \in A$ and some finite subset $H$ of $A$. As $V$ is quasi-compact, $j'(V)$ is quasi-compact. Now every $x \in j'(V)$ is of the form $v \circ \rho'$ for some $v \in V$. Thus by definition one has $x(q) \ne 0$ for all $x \in j'(V)$. Thus Lemma~\ref{NonInfinitesimalZero} shows that there exists a neighborhood $E$ of $0$ in $A$ with $v(\rho'(e)) \leq v(q)$ for all $v \in V$ and $e \in E$. Let $D \subseteq E$ be a finite subset such that $D\cdot A$ is open in $A$ (e.g. a finite system of generators of an ideal of definition contgained in $E$) and put $W := R(\frac{H \cup D}{q})$. Then $j'(V) = U \cap W$, hence $j'(V)$ is a rational subset of $\Spa A$.
\end{proof}

Let $A = (A,A^+)$ be an affinoid ring. Proposition~\ref{RationalFunction} now allows us to define a presheaf on the basis of rational subsets of $X := \Spa A$ by setting for $s \in A$ and $T \subseteq A$ finite with $T\cdot A$ open:
\begin{equation}\label{EqPreSheaf}
\Oscr_X(R(\frac{T}{s})) := A\langle\frac{T}{s}\rangle.
\end{equation}
This is a presheaf with values in the category of complete topological rings and continuous ring homomorphisms. For an arbitrary open subset $V$ of $\Spa A$ we set
\[
\Oscr_X(V) := \limproj_U \Oscr_X(U),
\]
where $U$ runs through rational subsets of $\Spa A$ contained in $V$. We equip $\Oscr_X(V)$ with the projective limit topology and obtain a presheaf on $\Spa A$ with values in the category of complete topological rings.

\begin{rem}\label{GlobalSections}
One has $X = R(\frac{1}{1})$ and hence $\Oscr_X(X) = \Ahat$.
\end{rem}

\begin{rem}\label{PresheafOpenRational}
Let $R(\frac{T}{s}) \subseteq \Spa A$ be a rational subset, let $B$ be the affinoid ring $A\langle\frac{T}{s}\rangle$, and let $j \colon U := \Spa B \to X := \Spa A$ be the canonical continuous map with image $R(\frac{T}{s})$. Let $V$ be a rational subset of $\Spa A$ with $V \subseteq R(\frac{T}{s})$. Then the unique continuous ring homomorphism $\sigma\colon \Oscr_X(V) \iso \Oscr_U(j^{-1}(V))$ making the diagram
\[\xymatrix{
A \ar[d] \ar[r] & A\langle\frac{T}{s}\rangle \ar[d] \\
\Oscr_X(V) \ar[r]^\sigma & \Oscr_U(j^{-1}(V))
}\]
commutative, is an isomorphism of complete topological rings.
\end{rem}

\begin{remdef}\label{DefStalk}
Let $A$ be an affinoid ring, $X = \Spa A$. For $x \in X$ let
\[
\Oscr_{X,x} := \limind_{\text{$x \in U$ open}} \Oscr_X(U) = \limind_{\text{$x \in U$ rational}} \Oscr_X(U)
\]
be the stalk. Here we take the inductive limit in the category of rings. Thus $\Oscr_{X,x}$ is not endowed with a topology.

For every rational subset $U$ of $X$ with $x \in U$ the valuation $x\colon A \to \Gamma_x \czero$ extends uniquely to a valuation $v_U\colon \Oscr_X(U) \to \Gamma_x \czero$ (Proposition~\ref{RationalFunction}~(2)). By passing to the inductive limit one obtains a valuation
\begin{equation}\label{ValStalk}
v_x\colon \Oscr_{X,x} \to \Gamma_x \czero.
\end{equation}
For every rational subset $U$ with $x \in U$ we obtain the homomorphism
\begin{equation}\label{EqCanStalk}
\rho_x\colon A \to \Oscr_X(U) \to \Oscr_{X,x}
\end{equation}
which is independent of the choice of $U$. By definition one has
\begin{equation}\label{ValuesAtStalk}
v_x(\rho_x(f)) = x(f) \in \Gamma_x \czero
\end{equation}
for $f \in A$ and $x \in X$. We usually omit $\rho_x$ from the notation and simply write $v_x(f)$.
\end{remdef}

\begin{prop}\label{PropStalk}
For every $x \in X = \Spa A$, the stalk $\Oscr_{X,x}$ is a local ring and the maximal ideal of $\Oscr_{X,x}$ is the support of $v_x$.
\end{prop}

\begin{proof}
Let $U \subseteq \Spa A$ be an open subset with $x \in U$ and let $f \in \Oscr_X(U)$ with $v_x(f) \ne 0$. We have to show that the image of $f$ in $\Oscr_{X,x}$ is a unit. Let $x \in W \subseteq U$ be a rational subset. Then $v_x$ defines by restriction a continuous valuation $v_W$ on $\Oscr_X(W)$. As $v_W(f) \ne 0$ there exists a finite subset $T$ of $B := \Oscr_X(W)$ such that $T\cdot B$ is open in $B$ and $v_W(t) \leq v_W(f)$ for all $t \in T$. Then we have in $Y := \Spa B$ the rational subset $V:= R(\frac{T}{f})$ with $v_W \in V$ and $f \in \Oscr_Y(V)^\times$. Hence there exists a rational subset $S$ of $\Spa A$ such that $x \in S \subseteq W$ and such that $f$ is a unit in $\Oscr_X(S)$.
\end{proof}

We denote by $k(x)$ the residue field of the local ring $\Oscr_{X,x}$. The valuation $v_x$ induces a valuation on $k(x)$ which is again denoted by $v_x$. Its valuation ring is denoted by $k(x)^+$.

\begin{def}\label{DefVPre}
We denote by $\Vscr^{\rm pre}$ the category of tuples $X = (X,\Oscr_X, (v_x)_{x \in X})$, where
\begin{definitionlist}
\item
$X$ is a topological space,
\item
$\Oscr_X$ is a presheaf of complete topological rings on $X$ such that the stalk $\Oscr_{X,x}$ of $\Oscr_X$ (considered as a presheaf of rings) is a local ring,
\item
$v_x$ is an equivalence class of valuations on the stalk $\Oscr_{X,x}$ such that $\supp(v_x)$ is the maximal ideal of $\Oscr_{X,x}$.
\end{definitionlist}
The morphisms $f\colon X \to Y$ are pairs $(f,f^{\flat})$, where $f$ is a continuous map of topological spaces $X \to Y$ and $f^{\flat}\colon \Oscr_Y \to f_*\Oscr_X$ is a morphisms of pre-sheaves of topological rings (i.e., for all $V \subseteq Y$ open, $\varphi_V\colon \Oscr_Y(V) \to \Oscr_X(f^{-1}(V))$ is a continuous ring homomorphism) such that for all $x \in X$ the induced ring homomorphism $f^{\flat}_x\colon \Oscr_{Y,f(x)} \to \Oscr_{X,x}$ is compatible with the valuation $v_x$ and $v_{f(x)}$, i.e., $v_{f(x)} = v_x \circ f^{\flat}_x$. (Note that this implies that $f^{\flat}_x$ is a local homomorphism and that $\Gamma_{v_{f(x)}} \subseteq \Gamma_{v_x}$.)
\end{def}

We have seen above that for an affinoid ring $A$ its adic spectrum $\Spa A$ is an object in $\Vscr^{\rm pre}$.

\begin{remdef}
Let $X$ be a an object in $\Vscr^{\rm pre}$ and let $U$ be an open subset of the underlying topological space of $X$. Then $(U,\Oscr_X\rstr{U},(v_x)_{x\in U})$ is again an object of $\Vscr^{\rm pre}$.

A morphism $j\colon Y \to X$ in $\Vscr^{\rm pre}$ is called an \emph{open immersion} if $j$ is a homeomorphism of $Y$ onto an open subspace $U$ of $X$ which induces an isomorphism
\[
(Y,\Oscr_Y,(v_y)_{y \in Y}) \iso (U,\Oscr_X\rstr{U},(v_x)_{x\in U})
\]
in $\Vscr^{\rm pre}$.
\end{remdef}

\begin{rem}\label{RationalOpenImmersion}
Let $A$ be an affinoid ring, $X = \Spa$ and let $U = R(\frac{T}{s})$ be a rational subset. Then $\Spa A\langle\frac{T}{s}\rangle \to \Spa A$ is an open immersion with image $U$.
\end{rem}

\begin{remdef}
Let $X$ be a topological space and let $\Bcal$ be a basis of the topology. Let $\Fscr$ be a presheaf on $X$ with values in a category, where projective limits exist. We call $\Fscr$ \emph{adapted to $\Bcal$}if for every open subset $V$ of $X$ the restriction maps $\Fscr(V) \to \Fscr(U)$ with $U \in \Bcal$ and $U \subseteq V$ yield an isomorphism
\[
\Fscr(V) \iso \limproj_{U} \Fscr(U).
\]
In this case, $\Fscr$ is a sheaf if and only if $\Fscr$ is a sheaf on $\Bcal$.
\end{remdef}

\begin{remdef}\label{DefPreAdic}
An \emph{affinoid pre-adic space} is an object of $\Vscr^{\rm pre}$ which is isomorphic to $\Spa A$ for an affinoid ring $A$.

Let $X$ be an object in $\Vscr^{\rm pre}$ such that there exists an open covering $(U_i)_i$ of $X$ such that $(U,\Oscr_X\rstr{U},(v_x)_{x\in U})$ is an affinoid pre-adic space. An open subset $U$ of $X$ is called \emph{open affinoid subspace} if $(U,\Oscr_X\rstr{U},(v_x)_{x\in U})$ is an affinoid pre-adic space. The open affinoid subspaces form a basis of the topology of the underlying topological space of $X$.

We call $X$ a \emph{pre-adic space} if in addition the sheaf of topological rings $\Oscr_X$ is adapted to the basis of open affinoid subspaces. The full subcategory of $\Vscr^{\rm pre}$ of pre-adic spaces is called the \emph{category of pre-adic spaces} and it is denoted by $\PreAd$.

For every affinoid ring $A$, $\Spa A$ is a pre-adic space by definition.
\end{remdef}

\begin{remark}
Let $X$ be a pre-adic space.
\begin{assertionlist}
\item
For every open subset $U \subseteq X$ the object $(U,\Oscr_X\rstr{U},(v_x)_{x \in U})$ in $\Vscr$ is again a pre-adic space.
\item
Let $U$ and $V$ be open affinoid adic subspaces of $X$. Then for all $x \in U \cap V$ there exists an open neighborhood $W \subseteq U \cap V$ of $x$ such that $W$ is a rational subset of $U$ and of $V$.
\end{assertionlist}
\end{remark}

\begin{rem}\label{PrincipalOpen}
Let $U \subseteq X = \Spa A$ be open and $f,g \in \Oscr_X(U)$. Then $V := \set{x \in U}{v_x(f) \leq v_x(g) \ne 0}$ is open in $X$.

Indeed, we can assume that $U = R(\frac{T}{s})$ is rational. Replacing $A$ by $A\langle\frac{T}{s}\rangle$ we may assume that $U = X$ and $f,g \in A$ (Proposition~\ref{RationalFunction}). But then $V$ is open by definition of the topology on $\Spa A$ as subspace of $\Spv A$.
\end{rem}

\begin{remdef}\label{DefOPlus}
Let $X$ be a pre-adic space. For every open subset $U \subseteq X$ we set
\begin{equation}\label{EqDefOPlus}
\Oscr_X^+(U) := \sett{f \in \Oscr_X(U)}{$v_x(f) \leq 1$ for all $x \in U$},
\end{equation}
endowed with the topology induced by $\Oscr_X(U)$. This is a sub-presheaf $\Oscr^+_X$ of topological rings of $\Oscr_X$. For every $x \in X$ let $\Oscr_{X,x}^+$ denote the stalk of $\Oscr_X^+$ at $x$. It follows from Remark~\ref{PrincipalOpen} that
\begin{equation}\label{StalkPlus}
\Oscr_{X,x}^+ = \set{f \in \Oscr_{X,x}}{v_x(f) \leq 1}.
\end{equation}
Proposition~\ref{PropStalk} then shows that $\Oscr_{X,x}^+$ is the inverse image of the valuation ring of $(k(x),v_x)$ under the canonical homomorphism $\Oscr_{X,x} \to k(x)$. In particular we see that $\Oscr_{X,x}^+$ is a local ring with maximal ideal $\set{f \in \Oscr_{X,x}}{v_x(f) < 1}$.
\end{remdef}

\begin{lemma}\label{MorphismAdic}
Let $X$ and $Y$ be pre-adic spaces. Let $(f,f^{\flat})\colon (X,\Oscr_X) \to (Y,\Oscr_Y)$ be a pair such that $f\colon X \to Y$ is continuous and such that $f^{\flat}\colon \Oscr_Y \to f_*\Oscr_X$ is a local morphism of presheaves of topological rings. Then $(f,f^{\flat})$ is a morphism of pre-adic spaces if and only if the following two conditions are satisfied.
\begin{definitionlist}
\item
$f^{\flat}(\Oscr^+_Y) \subseteq f_*\Oscr^+_X$.
\item
The induced morphism $\Oscr^+_Y \to f_*\Oscr^+_X$ is a local morphism of presheaves of rings.
\end{definitionlist}
\end{lemma}

\begin{proof}
It is clear that the condition is necessary. Assume that condition~(a) and ~(b) are satisfied. We have to show that for all $x \in X$ the local homomorphism $f^{\flat}_x\colon \Oscr_{Y,f(x)} \to \Oscr_{X,x}$ is compatible with valuations, i.e. $v_x \circ f^{\flat}_x = v_{f(x)}$. Now $\mfr_{f(x)} = \supp v_{f(x)}$ and $\mfr_{x} = \supp v_{x}$ (Proposition~\ref{PropStalk}). As $f^{\flat}_x$ is local, this shows that $v_x \circ f^{\flat}_x$ and $v_{f(x)}$ both have the support $\mfr_{f(x)}$. Let $K := \Oscr_{Y,f(x)}/\mfr_{f(x)}$ and let $A_x$ and $A_{f(x)}$ in $K$ the valuation rings of $f^{\flat}_x \circ v_x$ and $v_{f(x)}$, respectively. We have to show that $A_{f(x)} = A_x$. Condition~(a) and \eqref{StalkPlus} show that $A_{f(x)} \subseteq A_x$, and Condition~(b) shows that $A_x$ dominates $A_{f(x)}$. Hence $A_x = A_{f(x)}$.
\end{proof}

\begin{proposition}
Let $(A,A^+)$ be an affinoid Tate ring and let $\pi \in A$ be a topologically nilpotent unit. Then $\pi \in A^+$ and the $\pi$-adic completion of $\Oscr^+_{X,x}$ is equal to the $\pi$-adic completion of $k(x)^+$.
\end{proposition}

\begin{proposition}\label{SectionsRational}
For every rational subset $U = R(\frac{T}{s})$ one has
\[
(\Oscr_X(U),\Oscr_X^+(U)) = (A\langle\frac{T}{s}\rangle, A\langle\frac{T}{s}\rangle^+).
\]
\end{proposition}

\begin{proof}
One has $\Oscr_X(U) = A\langle\frac{T}{s}\rangle$ by definition and $\Spa A\langle\frac{T}{s}\rangle \to \Spa A$ is an open immersion with image $U$. Hence
\begin{align*}
\Oscr_X^+(U) &= \sett{f \in A\langle\frac{T}{s}\rangle}{$v(f) \leq 1$ for all $v \in \Spa A\langle\frac{T}{s}\rangle$}\\
&= A\langle\frac{T}{s}\rangle^+,
\end{align*}
where the second equation holds by Proposition~\ref{FunctionInvertible}~(1).
\end{proof}

\begin{rem}\label{NoSheaf}
In general $\Oscr_X$ is not a sheaf of rings even if $A$ is a Tate ring (see the end of section~3.2 in~\cite{Hu_Habil}).
\end{rem}

Let $\varphi\colon A \to B$ be a continuous morphism of affinoid rings. Set $X := \Spa B$, $Y := \Spa A$ considered as objects in $\Vscr^{\rm pre}$. Let $f\colon X \to Y$ be the continuous map of the underlying topological spaces attached to $\varphi$.

Let $V$ be a rational subset of $Y$ and $U = R(\frac{T}{s})$ be a rational subset of $X$ with $f(U) \subseteq V$. Then $U = \Spa B\langle\frac{T}{s}\rangle$ and as $f\rstr{U}$ factors through $V$ there exists a unique continuous ring homomorphism $\varphi_{V,U}\colon\Oscr_Y(V) \to B\langle\frac{T}{s}\rangle = \Oscr_X(U)$ such that the following diagram commutes
\[\xymatrix{
A \ar[r]^{\varphi} \ar[d] & B \ar[d] \\
\Oscr_Y(V) \ar[r]^{\varphi_{V,U}} & \Oscr_X(U).
}\]
These $\varphi_{V,U}$ yield a morphism of presheaves of complete topological rings $\Oscr_Y \to f_*\Oscr_X$ on the basis of rational subsets on $X$ resp. $Y$ and hence a morphism or presheaves  $f^{\flat}\colon \Oscr_Y \to f_*\Oscr_X$ of complete topological rings.

For every $x \in X$ the induced ring homomorphism $f_x\colon \Oscr_{Y,f(x)} \to \Oscr_{X,x}$ is compatible with the valuations $v_{f(x)}$ and $v_x$. Thus we obtain a morphism
\[
f := {}^{a}\varphi := \Spa(\varphi)\colon X \to Y
\]
in $\Vscr^{\rm pre}$.

The formation $\varphi \sends {}^a\varphi$ is functorial. If we denote by $\Affd$ the category of affinoid rings with continuous morphisms of affinoid rings, we obtain a contravariant functor
\[
\Affd \to \PreAd, \qquad A \sends \Spa A.
\]

\begin{rem}\label{IsomCompletion}
Let $A$ be an affinoid ring and $\iota\colon A \to \Ahat$ the canonical morphism into the completion. Then ${}^a\iota\colon \Spa \Ahat \to \Spa A$ is an isomorphism in $\Vscr^{\rm pre}$.
\end{rem}

\begin{prop}\label{AdicFullyFaithful}
Let $A$ and $B$ affinoid rings. If $B$ is complete, the map
\[
\Hom_{\Affd}(A,B) \to \Hom_{\Vscr^{\rm pre}}(\Spa B, \Spa A), \qquad \varphi \sends {}^a\varphi
\]
is bijective.
\end{prop}

\begin{proof}
The inverse map is given by
\[
\Hom_{\Vscr^{\rm pre}}(\Spa B, \Spa A) \to \Hom_{\Affd}(\Ahat,\Bhat) = \Hom_{\Affd}(A,B),\qquad f \sends f^{\flat}_Y,
\]
where the equality $\Hom_{\Affd}(\Ahat,\Bhat) = \Hom_{\Affd}(A,B)$ is due to the completeness of $B$. For details we refer to~\cite{Hu_Cont}.
\end{proof}

%---------------------------------------------------------

\subsection{Adic spaces}

\begin{rem}\label{SheafTopRing}
Let $X$ be a topological space, and let $\Oscr_X$ be a presheaf of topological rings (i.e. $\Oscr_X(U)$ is a topological ring for $U \subseteq X$ open and the restriction maps are continuous). Then $\Oscr_X$ is a sheaf of topological rings (i.e., for every topological ring $T$, $U \sends \Hom(T,\Oscr_X(U))$ is a sheaf of sets on $X$) if and only if $\Oscr_X$ is a sheaf of rings and for every open subset $U$ of $X$ and for every open covering $(U_i)_{i\in I}$ of $U$ the canonical map $\Oscr_X(U) \to \prod_{i \in I}\Oscr_X(U_i)$ is a topological embedding (where the product is endowed with the product topology).
\end{rem}

Let $\Vscr$ be the full subcategory of $\Vscr^{\rm pre}$ consisting of those objects $(X,\Oscr_X, (v_x)_{x \in X})$ in $\Vscr^{\rm pre}$ such that $\Oscr_X$ is a sheaf of topological rings.

Let $A$ be an affinoid ring, set $X := \Spa A$, and for $x \in X$ let $v_x$ be the valuation defined on $\Oscr_{X,x}$ in Definition~\ref{DefStalk}. If $\Oscr_X$ is a sheaf of topological rings, then $(X,\Oscr_X,(v_x)_x)$ is an object of $\Vscr$ which we again denote by $\Spa A$.

\begin{defi}\label{DefAffinoidSpace}
An \emph{affinoid adic space} is an object of $\Vscr$ which is isomorphic to $X := \Spa A$ for some affinoid ring $A$ such that $\Oscr_X$ is a sheaf of topological rings.
\end{defi}

\begin{defi}\label{DefAdicSpace}
An \emph{adic space} is an object $X$ of $\Vscr$ such that there exists an open covering $(U_i)_{i\in I}$ of $X$ such that $(U_i,\Oscr\rstr{U_i},(v_x)_{x\in U_i})$ is an affinoid adic space for all $i \in I$. The full subcategory of $\Vscr$ of adic spaces is denoted by $\Adic$.
\end{defi}

\begin{rem}
Let $X$ be an adic space.
\begin{assertionlist}
\item
The open affinoid adic subspaces form a basis of the topology of $X$.
\item
For every open subset $U \subseteq X$ the object $(U,\Oscr_X\rstr{U},(v_x)_{x \in U})$ in $\Vscr$ is again an adic space.
\item
Let $U$ and $V$ be open affinoid adic subspaces of $X$. Then for all $x \in U \cap V$ there exists an open neighborhood $W \subseteq U \cap V$ of $x$ such that $W$ is a rational subset of $U$ and of $V$.
\end{assertionlist}
\end{rem}

\begin{rem}\label{HomSheaf}
The morphisms of adic spaces form a sheaf. More precisely, if $X$ and $Y$ are adic spaces, then the presheaf on $X$ of sets
\[
U \sends \Hom_{\Adic}(U,Y),
\]
with the obvious restriction maps, is a sheaf.
\end{rem}

For two pairs $(A,A')$, $(B,B')$ consisting of a topological ring $A$ (resp.~$B$) and a subring $A'$ of $A$ (resp.~$B'$ of $B$) we denote by $\Hom((A,A'),(B,B'))$ the set of continuous ring homomorphisms $\varphi\colon A \to B$ such that $\varphi(A') \subseteq B'$.

\begin{prop}\label{UniversalAffinoid}
Let $X$ be an adic space, let $A$ be an affinoid ring and set $Y = \Spa A$. Then the map
\begin{align*}
\Hom_{\Adic}(X,Y) &\to \Hom((\Ahat,\Ahat^+),(\Oscr_X(X),\Oscr_X^+(X))) = \Hom((\Ahat,\Ahat^+),(\Oscr_X(X),\Oscr_X^+(X))),\\
f &\sends f^{\flat}_Y
\end{align*}
is a bijection.
\end{prop}

\begin{proof}
This follows from Proposition~\ref{AdicFullyFaithful} and Remark~\ref{HomSheaf}.
\end{proof}

\begin{definition}\label{DefSheafy}
We call an $f$-adic ring $A$ \emph{sheafy}, if for every ring $A^+$ of integral elements of $\Ahat$ the presheaf $\Oscr_{\Spa (\Ahat,A^+)}$ is a sheaf of topological rings.

We call $A$ \emph{stably sheafy} if $B$ is sheafy for every $\Ahat$-algebra topologically of finite type.
\end{definition}

\begin{remark}\label{RemSheafy}
Let $X$ be a pre-adic space such that one of the following conditions hold.
\begin{definitionlist}
\item
For every open affinoid subspace $U = \Spa A$ of $X$ the $f$-adic ring $A$ is sheafy.
\item
There exists a covering by open affinoid subspaces of the form $\Spa A$ such that the $f$-adic ring $A$ is stably sheafy.
\end{definitionlist}
Then $X$ is an adic space.
\end{remark}

\begin{theorem}\label{AdicSheaf}
Let $A = (A,A^+)$ be an affinoid ring and $X = \Spa A$. Assume that $A$ satisfies one of the following conditions.
\begin{definitionlist}
\item
The completion $\Ahat$ has a noetherian ring of definition.
\item
$A$ is a strongly noetherian Tate ring.
\item
$\Ahat$ has the discrete topology.
\end{definitionlist}
Then $\Oscr_X$ is a sheaf of complete topological rings. Moreover, one has $H^q(U,\Oscr_X) = 0$ for all $q \geq 1$ and all rational subsets $U$ of $X$.
\end{theorem}

We will give the proof only in the cases~(b) and~(c). In case~(c) the proof is esay:

\begin{proof}[Proof in the discrete case]
We may and do assume that $A$ is complete. Assume that $A$ carries the discrete topology. Then $A^+$ can be any subring of $A$ which is integrally closed in $A$, and $\Spa A = \sett{v \in \Spv A}{$v(f) \leq 1$ for all $f \in A^+$}$.

If $U = R(\frac{T}{s})$ is a rational subset of $\Spa A$, then $\Oscr_X(U) = A_s$, endowed with the discrete topology, and the $\Oscr_X$-acyclicity of rational covers follows from the analogue statement for the structure sheaf of an affine scheme and open coverings by principally open subsets $D(s) = \Spec A_s$. In fact, in this case the structure sheaf on $\Spa A$ is simply the pullback of the structure sheaf on $\Spec A$ under the continuous and surjective map $\Spa A \to \Spec A$, $x \sends \supp(x)$.
\end{proof}

For the proof if $A$ is a strongly noetherian Tate algebra we need some preparations. We will use some general results on \v{C}ech cohomology, recalled in Appendix~\ref{AppCech}.

\begin{rem}\label{ModuleConvPS}
Let $A$ be a complete noetherian Tate ring and let $M$ be a finitely generated $A$-module endowed with its canonical topology (Proposition~\ref{LinearCont}~(1)). 
Denote by $M\langle X\rangle$ the $A\langle X\rangle$-module of elements $\sum_{\nu \geq 0} m_{\nu}X^{\nu}$ with $m_\nu \in M$ for all $\nu$ and such that for every neighborhood of zero $U$ in $M$ one has $m_{\nu} \in U$ for almost all $\nu$.

We claim that the homomorphism of $A\langle X\rangle$-modules
\[
\mu_M \colon M \otimes_A A\langle X\rangle \to M\langle X\rangle, \qquad m \otimes a \sends ma
\]
is bijective. Indeed, this is clear if $M$ is a finitely generated free $A$-module. In general, we find a presentation
\[
A^n \ltoover{u} A^m \ltoover{p} M \to 0
\]
(because $A$ is noetherian). Proposition~\ref{LinearCont}~(2) shows that $u$ and $p$ are continuous and open onto its image. Thus we obtain an exact sequence $A^m\langle X\rangle \to A^n\langle X\rangle \to M\langle X\rangle \to 0$. Then the 5-lemma implies our claim.
\end{rem}

\begin{prop}\label{RationalFlat}
Let $A = (A,A^+)$ be a strongly noetherian Tate affinoid ring, and let $U \subseteq V \subseteq X := \Spa A$ be two rational subsets. Then the restriction homomorphism $\Oscr_X(V) \to \Oscr_X(U)$ is flat.
\end{prop}

\begin{proof}
By Example~\ref{LocalTFT} we know that $\Oscr_X(V)$ is again a strongly noetherian Tate ring. Thus we may assume that $X = V$ and that $A$ is complete. By Remark~\ref{DevissageRational} we may moreover assume that $U$ is either of the form $U_1 = R(\frac{f}{1}) = \set{x \in X}{x(f) \leq 1}$ or of the form $U_2 = R(\frac{1}{f}) = \set{x \in X}{x(f) \geq 1}$ for some $f \in A$.

In Example~\ref{LocalTFT} we have seen that $\Oscr_X(U_1) = \Ahat\langle X\rangle/(f - X)$ and $\Oscr_X(U_2) = \Ahat\langle X\rangle/(1 - fX)$. Thus it suffices to show the following lemma.
\end{proof}

\begin{lem}\label{FlatLocalization}
Let $A$ be a noetherian complete Tate ring.
\begin{assertionlist}
\item
The ring $A\langle X\rangle$ is faithfully flat over $A$.
\item
For all $f \in A$ the rings $A\langle X\rangle/(f - X)$ and $A\langle X\rangle/(1 - fX)$ are flat over $A$.
\end{assertionlist}
\end{lem}

\begin{proof}
\proofstep{(1)}
Let $i\colon N \mono M$ be an injective homomorphism of finitely generated $A$-modules. Then Remark~\ref{ModuleConvPS} implies that $i \otimes \id_{A\langle X\rangle}\colon M \otimes_A A\langle X\rangle \to N \otimes_A A\langle X\rangle$ is again injective. This shows that $A\langle X\rangle$ is flat over $A$.

If $\pfr$ is a prime ideal of $A$, then the set of $\sum_\nu a_\nu X^{\nu} \in A\langle X\rangle$ such that $a_0 \in \pfr$ is a prime ideal $\qfr$ of $A\langle X\rangle$ with $\qfr \cap A = \pfr$. Thus $A\langle X\rangle$ is faithfully flat over $A$.

\proofstep{(2)}
We first show the following claim. Let $g \in A\langle X\rangle$ and assume that for every finitely generated $A$-module $M$ the multiplication $w_g\colon M\langle X\rangle \to M\langle X\rangle$ is injective. Then $B := B_g := A\langle X\rangle/(g)$ is flat over $A$.

To show the claim consider the sequence of $A$-modules
\[
0 \to A\langle X\rangle \ltoover{v} A\langle X\rangle \ltoover{p} B \to 0,
\]
where $v = v_g$ is the multiplication with $g$. It is exact by our assumption applied to $M = A$. If we want to show that $\Tor^A_1(M,B) = 0$ for every finitely generated $A$-module $M$, then it suffices to show that $w_g := \id_M \otimes_A v_g$ is injective (use that we have already seen in~(1) that $A\langle X\rangle$ is flat over $A$ and the long exact Tor-sequence). But using Remark~\ref{ModuleConvPS} we see that $w_g$ is the homomorphism $M\langle X\rangle \to M\langle X\rangle$ given by multiplication with $g$. This proves the claim.

If $g = 1 - fX$, then $w_g$ is easily checked to be injective.

Let $g = f - X$, $u = \sum_{\nu}m_{\nu}X^{\nu} \in M\langle X\rangle$ with $(f - X)u = 0$. Then
\[
fm_0 = 0, \qquad fm_\nu = m_{\nu-1} \quad\text{for all $\nu \geq 1$}.\tag{+}
\]
Let $M'$ be the submodule of $M$ generated by $\set{m_\nu}{\nu \in \NN_0}$. As $M$ is noetherian, $M'$ is generated by finitely many elements, say $m_0,\dots,m_l$. Then $(+)$ implies that $M'$ is generated by $m_l$ and $f^{l+1}m_l = 0$. Write $m_{2l+1} = am_l$ for some $a \in A$. Then $m_l = f^{l+1}m_{2l+1} = af^{l+1}m_l = 0$. Thus $M' = 0$ and hence $u = 0$.
\end{proof}

\begin{cor}\label{ExactInZero}
Let $A$ be a strongly noetherian Tate affinoid ring, $X = \Spa A$, and let $(U_i)_{1\leq i\leq n}$ a finite covering of $X$ be rational subsets. Then the homomorphism
\[
\Oscr_X(X) \to \prod_{i=1}^n \Oscr_X(U_i), \qquad f \sends (f\rstr{U_i})_{1\leq i\leq n}
\]
is faithfully flat (and in particular injective).
\end{cor}

\begin{lem}\label{TrivialChech}
Let $A = (A,A^+)$ be a strongly noetherian Tate affinoid ring, $X = \Spa A$. Let $f \in A$ and set $U_1 = \set{x \in X}{x(f) \leq 1}$ and $U_2 = \set{x \in X}{x(f) \geq 1}$. Then the augmented \v Cech complex (with alternating cochains) for $\Oscr_X$ and the open covering $\{U_1,U_2\}$ of $X$
\[
0 \to \Oscr_X(X) \ltoover{\eps} \Oscr_X(U_1) \times \Oscr_X(U_2) \ltoover{\delta} \Oscr_X(U_1 \cap U_2) \to 0
\]
is exact.
\end{lem}

\begin{proof}
We may assume that $A$ is complete (to simplify the notation). We have already seen that $\eps$ is injective (Corollary~\ref{ExactInZero}). Moreover, by Examples~\ref{LocalTFT} and~\ref{LaurentSeries} we have
\begin{equation}\label{EqDescribeSection}
\begin{aligned}
\Oscr_X(U_1) &= A\langle \zeta\rangle/(f - \zeta),\\
\Oscr_X(U_2) &= A\langle \eta\rangle/(1 - f\eta),\\
\Oscr_X(U_1 \cap U_2) &= A\langle \zeta,\eta\rangle/(f-\zeta,1 - f\eta) = A\langle \zeta,\eta\rangle/(f - \zeta,1 - \zeta\eta) \\
&= A\langle \zeta,\zeta^{-1}\rangle/(f - \zeta).
\end{aligned}
\end{equation}
Consider the following commutative diagram
\[\xymatrix{
& & 0 \ar[d] & 0 \ar[d] \\
& & (f - \zeta)A\langle\zeta\rangle \times (1- f\eta)A\langle\eta\rangle
\ar[d] \ar[r]^{\lambda'} & (f - \zeta)A\langle\zeta,\zeta^{-1}\rangle \ar[d] \ar[r] & 0 \\
0 \ar[r] & A \ar[r]^{\iota} \ar@{=}[d] & A\langle\zeta\rangle \times A\langle\eta\rangle \ar[d] \ar[r]^{\lambda} & A\langle\zeta,\zeta^{-1}\rangle \ar[d] \ar[r] & 0 \\
0 \ar[r] & A \ar[r]^{\eps} & \Oscr_X(U_1) \times \Oscr_X(U_2) \ar[d] \ar[r]^{\delta} & \Oscr_X(U_1 \cap U_2) \ar[d] \ar[r] & 0 \\
& & 0 & 0
}\]
Here $\iota$ is the canonical injection, $\lambda$ is the map $g((\zeta), h(\eta)) \sends g(\zeta) - h(\zeta^{-1})$, and $\lambda'$ is induced by $\lambda$. The columns are exact by~\eqref{EqDescribeSection}. A diagram chase shows that if we the first and second row are exact, then the third row is exact (note that we know already the injectivity of $\eps$).

The equations
\begin{align*}
A\langle\zeta,\zeta^{-1}\rangle &= A\langle\zeta\rangle + \zeta^{-1}A\langle\zeta^{-1}\rangle,\\
(f - \zeta)A\langle\zeta,\zeta^{-1}\rangle &= (f - \zeta)A\langle\zeta\rangle + (1 - f\zeta^{-1})A\langle\zeta^{-1}\rangle
\end{align*}
show the surjectivity of $\lambda$ and $\lambda'$ (and in particular the exactness of the first row). Finally, the equality
\[
0 = \lambda(\sum_{k\geq0}a_k\zeta^k,\sum_{k\geq0}b_k\eta^k) = \sum_{k\geq0}a_k\zeta^k - \sum_{k\geq0}b_k\zeta^{-k}
\]
is equivalent to $a_k = b_k = 0$ for $k > 0$ and $a_0 = b_0$. Thus $\im(\iota) = \ker(\lambda)$.
\end{proof}

We can now prove that $\Oscr_X$ is a sheaf if $A$ is a strongly noetherian Tate algebra.

\begin{proof}[Proof of Theorem~\ref{AdicSheaf} if $A$ is a strongly noetherian Tate algebra]
By Proposition~\ref{CechZero} it suffices to show that every open covering by rational subsets is $\Oscr_X$-acyclic. We may assume that $A$ is complete. Every open covering of $X$ has a refinement $\Ucal = (U_t)_{t \in T}$ of the form $U_t := R(\frac{T}{t})$ with $T \subseteq A$ generating $A$ as an ideal (Lemma~\ref{RationalRefinement}). Let us call such a cover \emph{the rational cover generated by $T$}. If $U$ is any rational subset, then $\Ucal\rstr{U}$ is the rational cover of $U = \Spa \Oscr_X(U)$ generated by the set of images of $t$ in $\Oscr_X(U)$ for $t \in T$. Thus by Proposition~\ref{AcyclicProduct}~(2) it suffices to show the following lemma.
\end{proof}

\begin{lem}
Let $A$ be a complete strongly noetherian Tate ring and $\Ucal$ be a rational cover generated by some finite subset $T \subseteq A$ with $T\cdot A = A$. Then $\Ucal$ is $\Oscr_X$-acyclic.
\end{lem}

\begin{proof}
\proofstep{(i)}
For $f \in A$ let $\Ucal_f$ be the open covering of $X = \Spa A$ consisting of $R(\frac{f}{1})$ and $R(\frac{1}{f})$. Then $\Ucal_f$ is $\Oscr_X$-acyclic by Lemma~\ref{TrivialChech}.

Moreover, if $U = R(\frac{T}{s})$ is any rational subset then $\Ucal_f\rstr{U} = \Ucal_{f\rstr{U}}$, where $f\rstr{U}$ is the image of $f$ under the homomorphism $A \to A(\frac{T}{s})$. Thus $\Ucal_f\rstr{U}$ is $(\Oscr_X)$-acyclic.

Using Proposition~\ref{AcyclicProduct}~(3) it follows by induction that all open covers of the form $\Vcal := \Ucal_{f_1} \times \dots \times \Ucal_{f_r}$ are $\Oscr_X$-acyclic. Such a cover is called a \emph{Laurent cover generated by $f_1,\dots,f_r$}. It is the rational cover generated by $T = \set{\prod_{j \in J}f_j}{J \subseteq \{1,\dots,r\}}$.

If $U$ is any rational subset of $X$, then $\Vcal\rstr{U}$ is the Laurent cover generated by $f_1\rstr{U},\dots,f_r\rstr{U}$. Thus we have seen that for every Laurent cover $\Vcal$ of $X$ and every open rational subset $U$ the restriction $\Vcal\rstr{U}$ is $\Oscr_X$-acyclic (more precisely, $\Oscr_X\rstr{U}$-acyclic).

\proofstep{(ii)}
We show the following claim. Let $T = (f_0,\dots,f_n) \subseteq A$ be finite such that $T$ generates $A$ as ideal and let $\Ucal$ be the rational cover of $X$ generated by $T$. Then there exists a Laurent cover $(V_j)_{j\in J}$ of $X$ such that $\Ucal\rstr{V_j}$ is a rational cover generated by units of $\Oscr_X(V_j)$ for all $j \in J$.

Indeed, for all $x \in X$ there exists $f_i$ such that $x(f_i) \ne 0$. Thus by Corollary~\ref{TateZeroNeighbor} there exists a unit $s \in A^{\times}$ such that for all $x \in X$ an $i \in \{0,\dots,n\}$ exists with $x(s) < x(f_i)$. Then the Laurent cover generated by $s^{-1}f_1,\dots,s^{-1}f_r$ satisfies the claim.

\proofstep{(iii)}
Every rational cover $\Ucal$ of $X$ which is generated by units $f_0,\dots,f_n$ of $A$ has a refinement by a Laurent cover.

Indeed, the Laurent cover generated by $\set{f_if_j^{-1}}{0 \leq i,j \leq n}$ is a refinement of $\Ucal$.

\proofstep{(iv)}
As the restriction of Laurent covers to arbitrary rational subsets are $\Oscr_X$-acyclic by~(i), it follows from~(iii) that all restrictions to rational subsets of all rational covers generated by units are $\Oscr_X$-acyclic (Proposition~\ref{AcyclicProduct}~(2)). 

Now let $\Ucal$ be a rational cover generated by some finite subset $T \subseteq A$ with $T\cdot A = A$ and let $\Vcal$ be a Laurent cover such that $\Ucal\rstr{V}$ is a rational cover generated by a finite set of units for all $V$ in $\Vcal$ (which exists by~(ii)). Then we have just seen that $\Ucal\rstr{V}$ is $\Oscr_X$-ayclic. Moreover, by~(i) $\Vcal\rstr{U}$ is $\Oscr_X$-acyclic for every $U$ in $\Uscr$. Thus the $\Oscr_X$-acyclicity of $\Vcal$ implies the $\Oscr_X$-acyclicity of $\Ucal$ by Proposition~\ref{AcyclicProduct}~(1).
\end{proof}

\begin{corollary}\label{StablySheafy}
Let $A$ be an $f$-adic ring satisfying one of the following properties.
\begin{definitionlist}
\item
The completion $\Ahat$ has a noetherian ring of definition.
\item
$A$ is a strongly noetherian Tate ring.
\item
$\Ahat$ has the discrete topology.
\end{definitionlist}
Then $A$ is stably sheafy.
\end{corollary}

\begin{proof}
If $A$ has one of these properties, then every $A$-algebra topologically of finite type over $A$ has the same property.
\end{proof}

%------------------------------------------------

\subsection{Analytic Points of adic spaces}

\begin{propdef}
Let $X$ be an adic space. A point $x \in X$ is called \emph{analytic} if the following equivalent conditions are satisfied.
\begin{equivlist}
\item
There exists an open neighborhood $U$ of $x$ such that $\Oscr_X(U)$ contains a topologically nilpotent unit.
\item
For every open affinoid neighborhood $U = \Spa A$ of $x$, the point $\supp x \subset A$ is not open in $A$ (i.e., $x \in \Spa A$ is analytic in the sense Definition~\ref{DefAnalytic}.)
\end{equivlist}
We set $X\an := \sett{x \in X}{$x$ is analytic}$ and $X\na := X \setminus X\an$.
\end{propdef}

\begin{proof}
We may assume that $X = \Spa A$, $A$ complete affinoid ring. 

\proofstep{``(ii) $\implies$ (i)''}
Let $x \in \Spa A$ such that $\supp x$ is not open in $A$. By Remark~\ref{RemAnalytic}~\eqref{RemAnalytic1} there exists a topologically nilpotent element $s$ of $A$ with $x(s) \ne 0$. Then $U := \set{y \in \Spa A}{y(s) \ne 0}$ is an open neighborhood of $x$ in $\Spa A$. As the restriction $A = \Oscr_X(X) \to \Oscr_X(U)$ is a continuous ring homomorphism, the image of $s$ in $\Oscr_X(U)$ is again a topologically nilpotent unit.

\proofstep{``(i) $\implies$ (ii)''}
Let $x \in \Spa A$ such that $\supp x$ is open $A$, and let $U$ be an open neighborhood of $x$. We have to show that $\Oscr_X(U)$ has no topologically nilpotent unit. Let $V$ be a rational subset of $\Spa A$ with $x \in V \subseteq U$ and set $\pfr := \set{f \in \Oscr_X(V)}{v_x(f) = 0}$. Then $\pfr$ is a prime ideal of $\Oscr_X(V)$ with $\pfr \cap A = \supp x$. As $\supp x$ contains an ideal of definition of (a ring of definition of) $A$, $\pfr$ contains an ideal of definition of $\Oscr_X(V)$ by definition of $\Oscr_X(V)$. Thus $\pfr$ is an open prime ideal of $\Oscr_X(V)$ and contains therefore all topologically nilpotent elements of $\Oscr_X(V)$. As $\pfr$ contains no units, $\Oscr_X(V)$ contains no topologically nilpotent units. Hence $\Oscr_X(U)$ does not contain a topologically nilpotent unit.
\end{proof}

\begin{remark}\label{AnalyticConstructible}
Let $X$ be an adic space.
\begin{assertionlist}
\item\label{AnalyticConstructible1}
Remark~\ref{RemAnalytic}~\eqref{RemAnalytic2} implies that $X\an$ is an open constructible subset of $X$.
\item\label{AnalyticConstructible2}
For every open subspace $U$ of $X$ one has $U\an = X\an \cap U$ and $U\na = X\na \cap U$.
\end{assertionlist}
\end{remark}

%------------------------------------------------

\subsection{Adic morphisms of adic spaces}

\begin{definition}\label{DefAdicMorphism}
A morphism $f\colon X \to Y$ of adic spaces is called \emph{adic} if for every $x \in X$ there exist open affinoid neighborhoods $U$ of $x$ in $X$ and $V$ of $f(x)$ in $Y$ with $f(U) \subseteq V$ such that the ring homomorphism of $f$-adic rings $\Oscr_Y(V) \to \Oscr_X(U)$ induced by $f$ is adic.
\end{definition}

\begin{proposition}\label{PropertyAdic}
Let $f\colon X \to Y$ be a morphism of adic spaces.
\begin{assertionlist}
\item\label{PropertyAdic1}
$f$ is a adic if and only if $f(X\an) \subseteq Y\an$.
\item\label{PropertyAdic2}
$f(X\na) \subseteq Y\na$.
\end{assertionlist}
\end{proposition}

\begin{proof}
It follows from Remark~\ref{AnalyticConstructible}~\eqref{AnalyticConstructible2} that we may assume that $X$ and $Y$ are affinoid. But in this case we have already shown the all results in Lemma~\ref{AdicSpectral}.
\end{proof}

\begin{corollary}\label{AffinoidAdic}
Let $f\colon X \to Y$ be an adic morphism of adic spaces. Then for all open affinoid subspaces $U \subseteq X$ and $V \subseteq Y$ with $f(U) \subseteq V$ the ring homomorphism $\Oscr_Y(V) \to \Oscr_X(U)$ induced by $f$ is adic.
\end{corollary}

\begin{proof}
It follows from Proposition~\ref{PropertyAdic}~\eqref{PropertyAdic1} that $f(U\an) \subseteq V\an$. Thus Lemma~\ref{AdicSpectral}~\eqref{AdicSpectral2} implies the claim.
\end{proof}

Let $X$ be an adic space, $x \in X$. Recall that we denote by $k(x)$ the residue field of $\Oscr_{X,x}$. If $x$ is non-analytic, we endow $k(x)$ with the discrete topology. If $x$ is analytic, then we endow $k(x)$ with the topology induced by $v_x$. Then $k(x)$ is an $f$-adic field (Example~\ref{ExampleFAdic}~\eqref{ExampleFAdic3}). In fact, this is the unique structure of an $f$-adic field on $k(x)$ such that for every open affinoid neighborhood $U$ of $x$ the canonical homomorphism $\Oscr_X(U) \to k(x)$ is adic (\cite{Hu_Habil}~3.8.10). We also have defined $k(x)^+ := \set{a \in k(x)}{v_x(a) \leq 1}$ and we set
\[
\kappa(x) := (k(x),k(x)^+)
\]
which is an affinoid field. Moreover, the underlsing topological space of $\Spa \kappa(x)$ is totally ordered by specialization. The closed point is $v_x$. If $x$ is analytic, the generic point is the unique height $1$ valuation on $k(x)$ that is dependent on $v_x$. If $x$ is non-analytic, the generic point of $\Spa \kappa(x)$ is the trivial valuation on $k(x)$.

The canonical morphism $i_x\colon \Spa \kappa(x) \to X$ of adic spaces is adic. It maps the closed point of $\Spa \kappa(x)$ to $x \in X$ and it yields a homeomorphism of $\Spa \kappa(x)$ onto the subspace consisting of the vertical generizations of $x$.

%------------------------------------------------

\subsection{Morphisms of finite type}

Let $A$ be an $f$-adic ring, let $M_1,\dots,M_n$ be finite subsets of $A$ such that $M_i\cdot A$ is open in $A$ for $i = 1,\dots,n$. Then $M_i^mU$ is a neighborhood of $0$ for all neighborhoods $U$ of $0$ and for all $m \geq 1$ (Lemma~\ref{ZeroNeighbor}) and we defined in Definition~\ref{ConvPower} the ring
\[
A \langle X\rangle_M  = A \langle X_1,\dots,X_n\rangle_{M_1,\dots,M_n}.
\]
We saw in Proposition~\ref{PropConvPS} that $A \langle X\rangle_M$ is complete if $A$ is complete. We also described in Corollary~\ref{UniversalConvPS} the universal property of $A \langle X\rangle_M$.

Recall also the definition of homomorphisms topologically of finite type and some of their properties between $f$-adic rings from Section~\ref{Sec:HomTFTfadic}. We now define the notion of homomorphisms topologically of finite type between affinoid rings.

\begin{remdef}
Let $A = (A,A^+)$ be an affinoid ring and let $M_1,\dots,M_n$ be finite subsets of $A$ such that $M_i\cdot A$ is open in $A$ for $i = 1,\dots,n$. Recall that we set $M^{\nu} := M_1^{\nu_1}\cdots M_n^{\nu_n}$ for all $\nu \in \NN_0^n$. Then
\[
B := \sett{\sum_{\nu \in \NN_0^n}a_{\nu}X^{\nu} \in A \langle X\rangle_M}{$a_{\nu} \in M^{\nu}\cdot \widehat(A^+)$ for all $\nu \in \NN_0^n$}
\]
is a subring of $A \langle X\rangle_M$. Its integral closure $C$ in $A \langle X\rangle_M$ is a ring of integral elements of $A \langle X\rangle_M$. The affinoid ring $(A \langle X\rangle_M,C)$ is simply denoted by $A \langle X\rangle_M$. We set $A \langle X\rangle := A \langle X\rangle_{\{1\},\dots,\{1\}}$. 
\end{remdef}

\begin{definition}
A homomorphisms $\pi\colon (C,C^+) \to (B,B^+)$ between affinoid rings is called a \emph{quotient mapping} if $\pi\colon C \to B$ is surjective continuous and open and if $B^+$ is the integral closure of $\pi(C^+)$ in $B$.

A homomorphism $\varphi\colon (A,A^+) \to (B,B^+)$ of affinoid rings with $B$ complete is called \emph{topologically of finite type} if there exist $n \in \NN_0$ and finite subsets $T_1,\dots,T_n$ of $A$ with $T_i \cdot A$ open in $A$ and a quotient mapping $\pi\colon A \langle X_1,\dots,X_n\rangle_{T_1,\dots,T_n} \to B$ such that $\varphi = \pi \circ \iota$, where $\iota\colon A \to A \langle X\rangle_T$ is the canonical homomorphism.
\end{definition}

\begin{proposition}\label{CharAffinoidTFT}
Let $A$ and $B$ affinoid rings with $B$ complete. Let $\varphi\colon A \to B$ be a homomorphism of affinoid rings. Then $\varphi$ is topologically of finite type if and only if the homomorphism $A \to B$ of $f$-adic rings is topologically of finite type and there exists an open subring $C$ of $B^+$ such that $B^+$ is integral over $C$, $\varphi(A^+) \subseteq C$, and $\varphi\colon A^+ \to C$ is topologically of finite type.
\end{proposition}

\begin{proof}
\cite{Hu_Gen}~Lemma~3.5.
\end{proof}

\begin{corollary}\label{LocalAffinoidTFT}
Let $A$ be an affinoid ring, $s_1,\dots,s_n \in A$, $T_i \subseteq A$ finite subsets such that $T_i\cdot A$ for all $i = 1,\dots,n$. Then the canonical homomorphism $A \to A \langle X_1,\dots,X_n\rangle_T$ is topologically of finite type.
\end{corollary}

The following proposition follows from the corresponding result for homomorphism of $f$-adic rings (Proposition~\ref{TateTFT}).

\begin{proposition}\label{TateAffinoidTFTRem}
Let $A$ be a Tate affinoid ring, then a homomorphism $\varphi\colon A \to B$ to a complete affinoid ring $B$ is topologically of finite type if and only if $\varphi$ factors through a quotient mapping $A \langle X_1,\dots,X_n\rangle \to B$ for some $n \in \NN_0$.
\end{proposition}

\begin{proposition}\label{CompositionAffinoidTFT}
Let $A$, $B$, and $C$ be affinoid rings with $B$ and $C$ complete and let $\varphi\colon A \to B$ and $\psi\colon B \to C$ be continuous homomorphism of affinoid rings.
\begin{assertionlist}
\item
If $\varphi$ and $\psi$ are topologically o finite type, then $\psi \circ \varphi$ is topologically of finite type.
\item
If $\psi \circ \varphi$ is topologically of finite type, then $\psi$ is topologically of finite type.
\end{assertionlist}
\end{proposition}

\begin{proof}
\cite{Hu_Gen}~Lemma~3.5.
\end{proof}

\begin{proposition}\label{TateAffinoidTFT}
Let $A$ and $B$ complete $f$-adic rings and let $\varphi\colon A \to B$ be a continuous ring homomorphism which is topologically of finite type. Assume that $A$ is a Tate ring and has a noetherian ring of definition. Then there is a unique ring of integral elements $B^+$ of $B$ such that $\varphi\colon (A,A^o) \to (B,B^+)$ is a homomorphism of affinoid rings which is topologically of finite type, namely $B^+ = B^o$.
\end{proposition}

\begin{proof}
\cite{Hu_Habil}~2.4.17
\end{proof}

\begin{definition}\label{DefFinType}
Let $f\colon X \to Y$ be a morphism of adic spaces.
\begin{assertionlist}
\item
Then $f$ is called \emph{locally of finite type} (resp.~\emph{locally of weakly finite type}), if for every $x \in X$ there exists an open affinoid neighborhood $U = \Spa B$ of $x$ in $X$ and an open affinoid subspace $V = \Spa A$ of $Y$ with $f(U) \subseteq V$ such that the induced homomorphism of affinoid rings $(A,A^+) \to (B,B^+)$ is topologically of finite type (resp.~such that the induced homomorphism of $f$-adic rings $A \to B$ is topologically of finite type).
\item
The morphism $f$ is called \emph{of (weakly) finite type}, if $f$ is locally of (weakly) finite type and quasi-compact.
\end{assertionlist}
\end{definition}

\begin{remark}\label{RemLocallyFT}
Let $X$ be an adic space.
\begin{assertionlist}
\item
Let $U \subseteq X$ be an open subspace. Then the inclusion $U \mono X$ is locally of finite type.
\item
Every morphism locally of weakly finite type is adic (Proposition~\ref{DefTFT}).
\end{assertionlist}
\end{remark}

\begin{proposition}\label{RestrictFinType}
Let $f\colon X \to Y$ be a morphism of adic spaces locally of (weakly) finite type. Let $U' \subseteq X$ and $V' \subseteq Y$ be open subspaces with $f(U') \subseteq V'$. Then the morphism of adic spaces $U' \to V'$ obtained by restriction from $f$ is locally of (weakly) finite type.
\end{proposition}

\begin{proof}
We have to show that for all $x \in U'$ there exist open affinoid neighborhoods $U = \Spa B \subseteq U'$ of $x$ and $V = \Spa A \subseteq V'$ of $f(x)$ with $f(U) \subseteq V$ such that the induced homomorphism of affinoid rings $A \to B$ is topologically of finite type.

Indeed, by hypothesis there exist open affinoid neighborhoods $\Utilde = \Spa \Btilde$ of $x$ and $\Vtilde = \Spa \Atilde$ with $f(\Utilde) \subseteq \Vtilde$ such that $\Atilde \to \Btilde$ is topologically of finite type. Choose open rational subsets $R(\frac{T}{s})$ of $\Spa \Atilde$ contained in $V'$ and $R(\frac{M}{r})$ of $\Spa \Btilde$ contained in $U'$ with $f(R(\frac{M}{r})) \subseteq R(\frac{T}{s})$. Then the composition $\Atilde \to \Btilde \to \Btilde \langle \frac{M}{r}\rangle$ is topologically of finite type (Corollary~\ref{LocalAffinoidTFT} and Proposition~\ref{CompositionAffinoidTFT}~(1)) and thus Proposition~\ref{CompositionAffinoidTFT}~(2) implies that the induced homomorphism $\Atilde \langle \frac{T}{s}\rangle \to \Btilde \langle \frac{M}{r}\rangle$ is topologically of finite type.

The same proof shows the claim in the ``locally of weakly finite type'' case (using Proposition~\ref{CompositionTFT} instead of Proposition~\ref{CompositionAffinoidTFT} and Example~\ref{LocalTFT} instead of Corollary~\ref{LocalAffinoidTFT}).
\end{proof}

\begin{proposition}\label{GlobalFinType}
Let $f\colon X \to Y$ be a morphism of adic spaces locally of finite type. Let $U = \Spa B \subseteq X$ and $V = \Spa A \subseteq Y$ be open affinoid subspaces with $f(U) \subseteq V$. Then the induced homomorphism of affinoid rings $(A,A^+) \to (B,B^+)$ is topologically of finite type.
\end{proposition}

\begin{proof}
Replacing $f$ by its restriction $U \to V$ (again locally of finite type by Proposition~\ref{RestrictFinType}), we may assume that $X$ and $Y$ are affinoid. Then the proposition is shown in~\cite{Hu_Habil}~Prop.~3.8.15.
\end{proof}

\begin{remark}\label{NotGlobalWeaklyFinType}
Let $X = \Spa B$ amd $Y = \Spa A$ be affinoid adic spaces, and let $f\colon X \to Y$ be weakly of finite type. Then the induced homomorphism $A \to B$ of $f$-adic rings is not necessarily topologically of finite type (\cite{Hu_Habil}~3.8.18).
\end{remark}

\begin{corollary}\label{OpenAffinoidFinType}
Let $X = \Spa A$ be an affinoid adic space, $U = \Spa B$ an open affinoid subspace of $X$. Then the induced homomorphism $(A,A^+) \to (B,B^+)$ is topologically of finite type.
\end{corollary}

%------------------------------------------------

\subsection{Fiber products of adic spaces}

\begin{definition}\label{DefStableSpace}
An adic space $X$ is called \emph{stable} if there exists an open  covering by affinoid subspaces of the form $\Spa (A,A^+)$ , where the $f$-adic ring is stably sheafy.
\end{definition}

\begin{rem}\label{RemStableSpace}
Let $X$ be a stable adic space.
\begin{assertionlist}
\item
The underlying topological space of $X$ has a basis of open affinoid subspaces $\Spa (A,A^+)$, where $A$ is stably sheafy.

Indeed, if $U = \Spa (A,A^+)$ with $A$ stably sheafy, then for every open rational subspace $V = \Spa (B,B^+)$ of $U$, the $f$-adic ring $B$ is stably sheafy because $B$ is topologically of finite type over $A$.
\item
If $f\colon Y \to X$ is a morphism locally of weakly topologically space, then $Y$ is stable adic space.
\end{assertionlist}
\end{rem}

\begin{theorem}\label{FiberProduct}
Let $f\colon X \to S$ and $g\colon Y \to S$ be morphisms of adic spaces. Assume $Y$ is stable and that one of the following consitions is satisfied.
\begin{definitionlist}
\item
$f$ is locally of finite type.
\item
$f$ is locally of weakly finite type and $g$ is adic.
\end{definitionlist}
Then the fiber product $X \times_S Y$ in the category of adic spaces exists. Moreover, $X \times_S Y$ is a stable adic space.
\end{theorem}

\begin{proof}
\proofstep{(i)}
As usual (using the fact that morphisms between adic spaces can be glued) the existence of the fibre product for $f$ and $g$ follows if there exist an open coverings $S = \bigcup_i S_i$ and open coverings $f^{-1}(S_i) = \bigcup_j X_{ij}$ and $g^{-1}(S_i) = \bigcup_k Y_{ik}$ such that all fiber products $X_{ij} \times_{S_i} Y_{ik}$ exist. Thus we may assume that $X = \Spa (A,A^+)$, $Y = \Spa (B,B^+)$ and $S = \Spa (R,R^+)$ are affinoid with $R$, $A$, and $B$ complete and $B$ stably sheafy. Let $\varphi\colon R \to A$ and $\psi\colon R \to B$ be the homomorphisms of $f$-adic rings induced by $f$ and $g$, respectively.

\proofstep{(ii)}
We first assume that (b) is satisfied. Then $\psi$ is adic (Corollary~\ref{AffinoidAdic}). Moreover, we may assume that the homomorphism $\varphi\colon R \to A$ of $f$-adic rings is topologically of finite type. Let $R_0 \subseteq R$, $A_0 \subseteq A$, and $B_0 \subseteq B$ be rings of definitions with $\varphi(R_0) \subseteq A_0$ and $\psi(R_0) \subseteq B_0$. Let $I$ be an ideal of definition of $R_0$.

Set $C' := A \otimes_R B$, let $C'_0$ be the image of $A_0 \otimes_{R_0} B_0$ in $C'$, and let $C^{\prime+}$ be the integral closure of the image of $A^+ \otimes_{R^+} B^+$ in $C'$. Endow $C'$ with the topology such that $\{I^n\cdot C\}_{n\in\NN}$ is a fundamental system of neighborhoods of zero. Let $C$ be the completion of $C'$. As $A$ is topologically of finite type over $R$, $C$ is topologically of finite type over $B$. Thus $C$ is an $f$-adic ring by Proposition~\ref{ConvPSTate} and Remark~\ref{QuotientFAdic}. Moreover $C$ is stably sheafy because $B$ is stably sheafy. The completion $C^+$ of $C^{\prime+}$ is a ring of integral elements in $C$.

By Proposition~\ref{UniversalAffinoid} it suffices to show that the natural diagram
\[\xymatrix{
(C,C^+) & (B,B^+) \ar[l] \\
(A,A^+) \ar[u] & (R,R^+) \ar[l]_{\varphi} \ar[u]_{\psi} 
}\]
is cocartesian in the category of complete affinoid rings. This is proved in \cite{Hu_Habil}~2.4.18.

\proofstep{(iii)}
For the proof under condition~(a) we refer to~\cite{Hu_Etale}~1.2.2.
\end{proof}

The construction of the fiber product in loc.~cit.\ under condition~(a) does not show that in this case the fiber product of affinoids adic spaces is again affinoid.

\begin{remark}\label{PropFiberProd}
Let
\[\xymatrix{
Z \ar[r]^q \ar[d]_p & Y \ar[d]^g \\
X \ar[r]^f & S
}\]
be a cartesian diagram of adic spaces such that $f$ is locally of weakly finite type and $g$ is adic.
\begin{assertionlist}
\item
The proof of Theorem~\ref{FiberProduct} shows that $q$ is locally of weakly finite type and $p$ is adic.
\item
If $f$ (resp.~$g$) is quasi-compact, then $q$ (resp.~$p$) is quasi-compact.
\end{assertionlist}
\end{remark}

\begin{lemma}\label{PointsFiberProduct}
Notation and hypotheses of Remark~\ref{PropFiberProd}. Then for all $x \in X$ and $y \in Y$ with $f(x) = g(y)$ there exists $z \in Z$ with $p(z) = x$ and $q(z) = y$.
\end{lemma}

\begin{proof}
This shown as in \cite{Hu_Gen}~Lemma~3.9~(i).
\end{proof}

\begin{definition}\label{DefFiber}
Let $f\colon X \to Y$ be a morphism locally of weakly finite type, let $y \in Y$, and let $\Spa \kappa(y) \to Y$ be the canonical adic morphism. The fiber product $X \times_Y \Spa \kappa(y)$ is denoted by $X_y$ or by $f_a^{-1}(y)$ and is called the \emph{adic fiber of $f$ in $y$}.
\end{definition}

\begin{proposition}\label{DescriptionFiber}
Let $f\colon X \to Y$ be a morphism of adic spaces locally of weakly finite type, let $y \in Y$, and let $S$ be the set of vertical generizations of $y$. Then the projection $f_a^{-1}(y) \to X$ induces a homeomorphism of $f^{-1}_a(y)$ onto the topological space $f^{-1}(S)$.
\end{proposition}

\begin{proof}
\cite{Hu_Habil}~3.10.4
\end{proof}

%------------------------------------------------

\subsection{Fiber products of adic spaces and schemes}

We denote by $(X,\Oscr_X,(v_x)_{x\in X}) \sends \Xline := (X,\Oscr_X)$ the forgetful functor from the category of adic spaces to the category of locally ringed spaces.

\begin{propdef}
Let $X$ and $Y$ be schemes, $S$ a stable adic space and let $f\colon X \to Y$ and $g\colon \Sline \to Y$ be morphisms of locally ringed spaces. Assume that $f$ is scheme morphism locally of finite type. Then there exists an adic space $R$, a morphism $p\colon R \to S$ of adic spaces, and a morphism $q\colon \Rline \to X$ of locally ringed spaces such that
\begin{equation}\label{FiberSchemeDiagram}
\begin{aligned}\xymatrix{
\Rline \ar[r]^q \ar[d]_{\pline} & X \ar[d]^f \\
\Sline \ar[r]^g & Y 
}\end{aligned}
\end{equation}
commutes and such that the following universal property is satisfied. For every adic space $U$, every morphism $u\colon U \to S$ of adic spaces, and every morphism $v\colon \Uline \to X$ of locally ringed spaces with $g \circ \uline = f \circ v$ there exists a unique morphism $w\colon U \to R$ of adic spaces with $p \circ w = u$ and $q \circ \wline = v$.

Moreover, the morphism $p$ is locally of finite type and $R$ is a stable adic space.

The adic space $R$ is denoted by $X \times_Y S$ and is called the \emph{fiber product of $X$ and $S$ over $Y$}.
\end{propdef}

\begin{proof}
We may assume that $Y = \Spec B$ and $X = \Spec B[X_1,\dots,X_n]/I$ are affine and that $S = \Spa A$ is affinoid such that the affinoid ring $A$ is stably sheafy. Then $g$ is induced by a ring homomorphism $\lambda\colon B \to A$. Let $E$ be a finite set of topologically nilpotent elements of $A$ such that $E\cdot A$ is open in $A$.

For every $k \in \NN$ let $A(k)$ be the affnoid ring $A \langle X_1,\dots,X_n\rangle_{E^k,\dots,E^k}$ (where $E^k = \sett{e_1\cdot \dots \cdot e_k}{$e_i \in E$}$). Let $\lambda_k\colon B[X_1,\dots,X_n] \to A(k)$ be the extension of $\lambda$ with $\lambda_k(X_i) = X_i$ for $i = 1,\dots,n$. Let
\[
\pi_k\colon A(k) \to A(k)/\lambda_k(I)\cdot A(k) := A_k
\]
be the canonical homomorphism of $f$-adic rings (i.e. $A_k$ is endowed with the $f$-adic topology such that $\pi_k$ is open and continuous). Then $A_k$ is topologically of finite type over $A$ and hence stably sheafy. Let $A_k^+$ be the integral closure of $\pi_k(A(k)^+)$ in $A_k$. We obtain an affinoid ring $(A_k,A_k^+)$. Let $R_k = \Spa (A_k,A_k^+)$. 

For $k \leq l$ let $\varphi_{lk}\colon R_k \to R_l$ be the morphism induced by the continous $A$-homomorphism $A_l \to A_k$ with $\pi_l(X_i) \sends \pi_k(X_i)$ for $i = 1,\dots,n$. Then $\varphi_{lk}$ is an isomorphism of $R_k$ onto the rational subset
\[
\sett{x \in R_l}{$v_x(e\pi_l(X_i)) \leq 1$ for all $e \in E(k)$, $i = 1,\dots,n$}.
\]
of $R_l$. Hence there exists in the category of adic spaces over $S$ the inductive limit $R$ of the system $(R_k, \varphi_{lk})$. Let $p\colon R \to S$ be the structure morphism.

The ring homomorphisms $\pi_k \circ \lambda_k\colon B[X_1,\dots,X_n] \to A_k$ induce morphisms of locally ringed spaces $\underline{R_k} \to X$ of locally ringed spaces which glue together to a morphism of locally ringed spaces $q\colon \Rline \to X$. We have $g \circ \pline = f \circ q$. It is not difficult to check that $R$, $p$, and $q$ satisfy the universal property.
\end{proof}

\begin{remark}
Let $Y$ be a scheme, let $X_1$ and $X_2$ two $Y$-schemes locally of finite type, and let $h\colon X_1 \to X_2$ be a morphism of $Y$-schemes. Let $S$ be an adic space and let $g\colon \Sline \to Y$ be morphisms of locally ringed spaces. Then $h$ induces by the universal property of the fiber product a canonical morphism of adic spaces $X_1 \times_Y S \to X_2 \times_Y S$ which is denoted by $h_{S}$.
\end{remark}

\begin{proposition}\label{PointsFiberScheme}
With the notations of the commutative diagram~\eqref{FiberSchemeDiagram} let $x \in X$ and $s \in S$ with $f(x) = g(s)$. Then there exists $r \in R$ with $q(r) = x$ and $p(r) = s$.
\end{proposition}

\begin{proof}
We may assume $Y = \Spec B$, $X = \Spec C$, $S = \Spa A$. Then $f$ and $g$ are given by ring homomorphisms $\varphi\colon B \to C$ and $\psi\colon B \to A$. Then $\pfr := \pfr_x \in \Spec C$, $v := v_s \in \Spa A$ with $\varphi^{-1}(\pfr_x) = \psi^{-1}(\supp(v_s))$. Set $D := C \otimes_B A$. Choose $\qfr \in \Spec D$ whose image in $\Spec A$ is $\supp(v)$ and whose image in $\Spec C$ is $x$. Let $r \in \Spv(D)$ be a valuation satisfying the following conditions.
\begin{definitionlist}
\item
$\supp(r) = \qfr$.
\item
$r$ extends the valuation $v_s$.
\item
$\Gamma_v$ and $\Gamma_r$ have the same divisible hull (in particular, then convex subgroup generated by $\Gamma_v$ in $\Gamma_r$ is all of $\Gamma_r$).
\end{definitionlist}
We claim that $r \in R$.

Let $(A_0,I)$ be a pair of definition of $A$ and let $C' \subseteq C$ be a finite set of generators of the $B$-algebra $C$. As the convex hull of $\Gamma_v$ is $\Gamma_r$ we can find $k \in \NN$ such that $r(c \otimes i) \leq 1$ for all $c \in C'$ and $i \in I^k$.

Set $D_0 := A_0[c \otimes i\,;\,c\in C', i\in I^k] \subseteq D$. Let $D^+$ be the integral closure of $A^+[c \otimes i\,;\,c\in C', i\in I^k]$ in $D$. We endow $D$ with the grup topology such that $\set{I^n\cdot D_0}{n\in \NN}$ is a fundamental system of neighborhoods of $0$. Then $D$ is an $f$-adic ring and $D^+$ is a ring of integral elements of $D$. One has $r \in \Cont(D)$ and $r(d) \leq 1$ for all $d \in D^+$ by definition. Hence $r \in \Spa (D,D^+)$. As $\Spa (D,D^+)$ is an open subspace of $R$ the claim is proved.
\end{proof}

\begin{defi}\label{SchemeToAdicSpace}
Let $k$ be a non-archimedean field and let $X$ be a $k$-scheme locally of finite type. The identity $k \to k$ corresponds to a morphism of locally ringed spaces $\Spa (k,k^o) \to \Spec k$. Thus we obtain the adic space
\[
X^{\rm ad} := X \times_{\Spec k} \Spa(k,k^o),
\]
which is locally of finite type over $\Spa(k,k^o)$. It is called the \emph{adic spaces associated to $X$}.
\end{defi}

\begin{example}\label{SchemeTrivialAdicSpace}
Let $R$ be a ring endowed with the discrete topology and let $X$ be a $k$-scheme locally of finite type. Let $R^+ \subseteq R$ be a ring of integral elements. The identity $R \to R$ corresponds to a morphism of locally ringed spaces $S := \Spa (R,R^+) \to \Spec R$ and we obtain the adic space $X \times_{\Spec R} \Spa(R,R)$. It can be described as follows.

Set $T := \sett{(x,v) \in Spv X}{$v(a) \leq 1$ for all $a \in R^+$}$. Let $f\colon T \to X$ be the restriction of $\supp\colon \Spv X \to X$ which sends $(x,v)$ to $x$ and define $\Ascr := f^{-1}(\Oscr_X)$ (a sheaf of rings). We denote again by $\Ascr$ its topologicalization (see Remark~\ref{PseudoDiscrete} below). For each $t = (x,v) \in T$ let $v_t$ be the valuation on the stalk $\Ascr_{(x,v)} = \Oscr_{X,x}$ whose support is $\mfr_x$ and which induces $v$. Then $R := (T,\Ascr, (v_t)_{t\in T})$ is the fiber product $X \times_{\Spec R} S$.

If $R = R^+ = k$ is a field, then $S = \Spa (k,k)$ consists only of the trivial valuation and the underlaing topological space of $X \times_{\Spec k} \Spa(k,k)$ is the subspace of $\Spv X$ consisting of the points $(x,v)$ such that the restriction to $k$ is trivial.\end{example}

%================================================

\section{Formal schemes as adic spaces}

\subsection{Formal (adic) schemes}

\begin{definition}
A topological ring $A$ is called \emph{pre-admissible} if it satisfies the following properties.
\begin{definitionlist}
\item
$A$ is linearly topologized (i.e., there exists a fundamental system of neighborhood of $0$ consisting of ideals of $A$).
\item
There exists an open ideal $I \subseteq A$ such that $(I^n)_{n}$ tends to zero (i.e., for every neighborhood of zero $V$ there exists $n \geq 1$ with $I^n \subseteq V$).
\end{definitionlist}
An ideal as in~(b) is called ideal of definition. $A$ is called \emph{admissible}, if $A$ is in addition complete.
\end{definition}

\begin{remark}
If $I$ is an ideal of definition in a pre-admissible ring and if $J$ is an open ideal, then $I \cap J$ is again an ideal of definition. This shows that the ideals of definition in $A$ form a fundamental system of open ideals.

Note that $I^n$ is not necessarily open.
\end{remark}

An adic ring is a topological ring $A$, such that an ideal $I \subseteq A$ exists such that $(I^n)_{n\in\NN}$ is a fundamental system of neighborhoods of zero. Such an ideal $I$ is then an ideal of definition of $A$.

Any adic ring is pre-admissible.

Recall the following properties (Remark~\ref{TopNilpotentLinTop}).

\begin{remark}\label{PropAdicRing}
Let $A$ be a pre-admissible ring.
\begin{assertionlist}
\item
The set of topologically nilpotent elements $A^{oo}$ is an open ideal of $A$, and $A^{oo}$ is the union of all ideals of definition.
\item
For every ideal of definition $I$ one has
\[
\Spec(A/I) = \Spec(A/A^{oo}) = \sett{\pfr \in \Spec A}{$\pfr$ open in $A$}.
\]
In particular, the left hand side is independent of $I$.
\end{assertionlist}
\end{remark}

\begin{remark}\label{PseudoDiscrete}
Let $X$ be a topological space. Let $\Fscr$ be a sheaf of sets on $X$. If one endows $\Fscr(U)$ with the discrete topology for all $U \subseteq X$ open, then $\Fscr$ is a presheaf of topological spaces. This is not a sheaf of topological spaces in general: If $U$ is an arbitrary open subset and $(U_i)_i$ is an open covering, then the discrete topology on $\Fscr(U)$ is in general not the coarsest topology making all restriction maps $\Fscr(U) \to \Fscr(U_i)$ continuous.

Now assume that the set $\Bcal$ of open quasi-compact subsets of $X$ is a basis of the topology of $X$. Then it is easy to check that the restriction of $\Fscr$ to $\Bcal$ is a sheaf of topological spaces. Thus there exists a sheaf $\Fscr'$ of topological spaces on $X$ (unique up to unique isomorphism) extending $\Fscr\rstr{\Bcal}$. For every open subset $U$ the underlying set of $\Fscr'(U)$ is equal to $\Fscr(U)$, but the topology is not necessarily the discrete topology if $U$ is not quasi-compact.

We call $\Fscr'$ the \emph{topologicalization of $\Fscr$}.

If $u\colon \Fscr \to \Gscr$ is a morphism of sheaves of sets, then $u\colon \Fscr' \to \Gscr'$ is a morphism of sheaves of topological spaces.

The same considerations hold for sheafs of groups or of rings.
\end{remark}

Now let $A$ be an admissible ring. We attach a topologically ringed space $(\Xcal,\Oscr_{\Xcal})$ to $A$ as follows. Let $(I_{\lambda})_{\lambda}$ be a fundamental system of ideals of definition in $A$. The underlying topological space $\Xcal$ is the subspace $\Spec(A/I_{\lambda})$ of $\Spec(A)$ of open prime ideals of $A$. Let $\Oscr_\lambda$ be the topologicalization of $\Oscr_{\Spec(A/I_{\lambda})}$ on $\Xcal$ and set $\Oscr_{\Xcal} := \limproj_\lambda \Oscr_\lambda$ (projective limit in the category of sheaves of topological rings).

\begin{definition}\label{DefFormalSpec}
The topologically ringed space $(\Xcal,\Oscr_{\Xcal})$ is called the \emph{formal spectrum of $A$}. It is denoted by $\Spf A$.
\end{definition}

\begin{remark}
Let $A$ be a linearly topologized ring, $(I_{\lambda})_{\lambda}$ be a fundamental system of ideals and let $S \subseteq A$ be a multipicative set. Let $S_\lambda$ be the image of $S$ in $A_\lambda := A/I_\lambda$. Then the $S_{\lambda}^{-1}A_\lambda$ form a projective system with surjective transition homomorphisms. We denote by $A\{S^{-1}\}$ its projective limit.

If $S = \set{f^n}{n\in\NN_0}$ for $f \in A$, we write $A_{\{f\}}$.
\begin{assertionlist}
\item
The canonical homomorphism $S^{-1}A \to S_{\lambda}^{-1}A_{\lambda}$ is surjective and has as its kernel $S^{-1}I_{\lambda}$. This shows that $A\{S^{-1}\}$ is the completion of $S^{-1}A$ for the topology given by the fundamental system $S^{-1}I_{\lambda}$.

In particular $A_{\{f\}} = A \langle \frac{1}{f}\rangle$.
\item
Let $\afr \subseteq A$ be an open ideal. Then $S^{-1}\afr$ is an open ideal in $S^{-1}A$. Its completion is denoted by $\afr\{S^{-1}\}$.
\item
If $A$ is pre-admissible, then $A\{S^{-1}\}$ is an admissible ring. If $I$ is an ideal of definition, then $I\{S^{-1}\}$ is an ideal of definition of $A\{S^{-1}\}$.
\item
If $A$ is an adic with an ideal of definition $I$ such that $I/I^2$ is a finitely generated $A/I$-module (e.g., if $I$ is a finitely generated $A$-module), then $A\{S^{-1}\}$ is an adic ring and $I' := I\{S^{-1}\}$ is an ideal of definition. Moreover $(I')^n = I^n\{S^{-1}\}$ (\cite{EGA} $0_I$~(7.6.11)).
% Here $\Ihat_f = \limproj I_f/(I^n)_f$ is the completion of $I_f$ (equal to the closure of $I_f$ in $A_{\{f\}}$).
\item
If $A$ is adic and noetherian, then $A\{S^{-1}\}$ is adic and noetherian.
\end{assertionlist}
\end{remark}

\begin{remark}\label{PropAffineFormal}
Let $A$ be an admissible ring, $\Xcal = \Spf(A)$.
\begin{assertionlist}
\item
$\Xcal$ is a spectral space (because $\Xcal \subseteq \Spec A$ closed). For $f \in A$ set $\Dcal(f) := D(f) \cap \Xcal$. Then $(\Dcal(f))_{f\in A}$ is a basis of quasi-compact open subsets, stable under finite intersections.
\item
For all $f \in A$ one has
\[
\Gamma(\Dcal(f),\Oscr_{\Xcal}) = A_{\{f\}} := A \langle \frac{1}{f}\rangle
\]
as topological rings. In particular $\Gamma(\Xcal,\Oscr_{\Xcal}) = A$.
\item
If $A$ is in adic with a finitely generated ideal of definition, then $A_{\{f\}}$ is again a complete $f$-adic adic ring.
\item
For all $U \subseteq \Xcal$ the topologically ring $\Gamma(U,\Oscr_{\Xcal})$ is complete.
\item
For $x \in \Xcal$ let
\[
\Oscr_{\Xcal,x} := \limind_{f \notin \pfr_x} A_{\{f\}}
\]
be the stalk of $\Oscr_{\Xcal}$ as sheaf of rings (without topology). Then $\Oscr_{\Xcal,x}$ is a local ring with residue field $k(x) = A_{\pfr_x}/\pfr_xA_{\pfr_x}$. If $A$ is noetherian, then $\Oscr_{\Xcal,x}$ is noetherian.
\end{assertionlist}
\end{remark}

\begin{example}
Let $A$ be a ring endowed with the discrete topology. Then $\Spf A = \Spec A$ (as locally ringed spaces).
\end{example}

\begin{definition}
Let $(\Xcal,\Oscr_{\Xcal})$ be a topologically ringed space.
\begin{assertionlist}
\item
It is called \emph{affine formal scheme} if it is isomorphic to $\Spf(A)$ for an admissible ring $A$. An affine formal scheme is called \emph{adic} (resp.~$f$-adic, resp.~\emph{noetherian}) if it is isomorphic to $\Spf A$, where $A$ is an adic (resp.~$f$-adic and adic, resp.~noetherian and adic) ring.
\item
An open subset $\Ucal \subseteq \Xcal$ is called \emph{(adic, resp.~$f$-adic, resp.~noetherian) affine formal open}, if $(\Ucal,\Oscr_{\Xcal}\rstr{\Ucal})$ is an (adic, resp.~$f$-adic, resp.~noetherian) affine formal scheme.
\item
$(\Xcal,\Oscr_{\Xcal})$ is called an \emph{(adic resp.~$f$-adic resp.~locally noetherian) formal scheme} if every point $x \in \Xcal$ has an open neighborhood that is an (adic resp.~$f$-adic resp.~noetherian) affine formal open. 
\end{assertionlist} 
\end{definition}

It follows from Remark~\ref{PropAffineFormal} that every formal scheme is a locally ringed space.

\begin{remark}\label{BasisFormal}
By Remark~\ref{PropAffineFormal} the $f$-adic (resp.~noetherian) affine formal open subsets of an $f$-adic (resp.~locally noetherian) formal scheme $\Xcal$ form a basis of the topology of $\Xcal$. In particular if $\Ucal \subseteq \Xcal$ is an open subset, then $(\Ucal,\Oscr_{\Xcal}\rstr{\Ucal})$ is again an $f$-adic (resp.~locally noetherian) formal scheme.
\end{remark}

\begin{definition}
A morphism of formal schemes is a morphism of locally and topologically ringed spaces.
\end{definition}

Let $A$, $B$ be complete adic rings, $\Xcal = \Spf B$, $\Ycal = \Spf A$, and let $\varphi\colon A \to B$ be a continuous homomorphism. Then ${}^a\varphi := \Spec(\varphi)$ maps $\Spf B$ into $\Spf A$. Moreover, for every $f \in A$ one has $\Dcal(\varphi(f)) = {}^a\varphi^{-1}(\Dcal(f))$, and $\varphi$ induces a continuous homomorphism $A_{\{f\}} \to B_{\{\varphi(f)\}}$ compatible with restrictions corresponding by passage of $f$ to $fg$ for some $g \in A$. Thus one obtains a homomorphism of sheaves of topological rings $\Oscr_{\Ycal} \to {}^a\varphi_*\Oscr_{\Xcal}$ which is easily seen to induce a local homomorphism on stalks. This construction is compatible with composition of continuous ring homomorphisms and we obtain a functor $\Spf$ from the category of complete adic rings into the category of toplogically and locally ringed spaces.

\begin{proposition}
The functor $\Spf$ is fully faithful.
\end{proposition}

\begin{remark}
Let $A$ be a complete $f$-adic adic ring, $\Xcal = \Spf A$, $f \in A$, and let $\varphi\colon A \to A_{\{f\}}$ be the canonical homomorphism. Then $\Spf(\varphi)$ yields an isomorphism $\Spf A_{\{f\}} \iso (\Dcal(f),\Oscr_{\Xcal}\rstr{\Dcal(f)})$.
\end{remark}

\begin{example}\label{CompletionAlongSubscheme}
Let $X$ be a scheme and let $Y \subseteq X$ be a closed subschemes, defined by a quasi-coherent ideal $\Iscr$ of finite type. Let $\Oscr_{X/Y}$ be the sheaf of topological rings which is obtained by restriction to $Y$ of the sheaf $\limproj_n \Oscr_X/\Iscr^{n+1}$. Then $(Y,\Oscr_{X/Y}$ is an $f$-adic formal scheme (\cite{EGA}~I~(10.8.3)). It is called the \emph{completion of $X$ along $Y$} and it is denoted by $X_{/Y}$.

If $X = \Spec A$ is affine, $Y = V(I)$, then $X_{/Y} = \Spf(\Ahat)$, where $\Ahat$ is the completion of $A$ for the $I$-adic topology.
\end{example}

Let $A$ be an admissible ring, $\Xcal := \Spf(A)$. Let $I \subseteq A$ an open ideal and let $\{I_{\lambda}\}_{\lambda}$ be the partially ordered set of ideals of definition contained in $I$. Set $X_{\lambda} = \Spec A/I_{\lambda}$ and let $\widetilde{I/I_{\lambda}}$ be the quasi-coherent ideal of $\Oscr_{X_{\lambda}}$ corresponding to $I/I_{\lambda}$. We denote by $I^{\Delta}$ the projective limit of the sheaves induced by $\widetilde{I/I_{\lambda}}$ on $\Xcal$. This is an ideal of $\Oscr_{\Xcal}$.

For all $f \in A$ one has
\[
I^{\Delta}(\Dcal(f)) = \limproj_{\lambda}S_f^{-1}I/S^{-1}_fI_{\lambda} = I_{\{f\}} \subseteq A_{\{f\}}.
\]
In particular we have $I^{\Delta}(\Xcal) = I$. As the $\Dcal(f)$ form a basis of the topology, one has
\begin{equation}\label{IDeltaDf}
I^{\Delta}\rstr{\Dcal(f)} = (I_{\{f\}})^{\Delta}.
\end{equation}

\begin{definition}
Let $\Xcal$ be a formal scheme. An ideal $\Iscr$ of $\Oscr_{\Xcal}$ is called \emph{ideal of definition} if for every point there exists an open affine neighborhood $\Ucal = \Spf(A)$ and an ideal of definition $I$ of $A$ such that $\Iscr\rstr{\Ucal} = I^{\Delta}$.
\end{definition}

\begin{proposition}
Let $\Xcal$ be a formal locally noetherian scheme. Then there exists a largest ideal of definition $\Tscr$ of $\Xcal$. This is the unique idealof definition $\Iscr$ of $\Oscr_{\Xcal}$ such that the scheme $(\Xcal,\Oscr_{\Xcal}/\Iscr)$ is reduced.
\end{proposition}

\begin{definition}
The scheme $(\Xcal,\Oscr_{\Xcal}/\Tscr)$ is denoted by $\Xcal_{\rm red}$. It is called the \emph{underlying reduced subscheme of the formal scheme $\Xcal$}.
\end{definition}

\begin{remark}\label{FormalInductiveSystem}
Let $\Xcal$ and $\Ycal$ be formal schemes, let $\Iscr$ (resp.~$\Kscr$) be an ideal of definition of $\Oscr_{\Xcal}$ (resp.~of $\Ycal$). Let $f\colon \Xcal \to \Ycal$ be a morphism of formal schemes such that $f^*(\Kscr)\Oscr_{\Xcal} \subseteq \Iscr$. Then one has $f^*(\Kscr^n)\Oscr_{\Xcal} = (f^*(\Kscr)\Oscr_{\Xcal})^n \subseteq \Iscr^n$ for all $n \in \NN$. Thus $f$ induces for all $n \geq 0$ morphisms of schemes $f_n\colon X_n \to Y_n$, where $X_n := (\Xcal,\Oscr_{\Xcal}/\Iscr^{n+1})$ and $Y_n := (\Ycal,\Oscr_{\Ycal}/\Kscr^{n+1})$. For all $m \leq n$ the following diagram is commutative
\begin{equation}\label{LevelFormal}
\begin{aligned}\xymatrix{
X_m \ar[r]^{f_m} \ar[d] & Y_m \ar[d] \\
X_n \ar[r]^{f_n} & Y_n 
}\end{aligned}
\end{equation}
\end{remark}

\begin{proposition}\label{MorphismInductive}
Let $\Xcal$ and $\Ycal$ be adic formal schemes. Then the construction in Remark~\ref{FormalInductiveSystem} yields a bijection between the set of morphisms of formal schemes $\Xcal \to \Ycal$ such that $f^*(\Kscr)\Oscr_{\Xcal} \subseteq \Iscr$ and of the set of sequences $(f_n)$ of morphisms $f_n\colon X_n \to Y_n$ making~\eqref{LevelFormal} commutative.
\end{proposition}

\begin{remark}
If $\Xcal$ is locally noetherian with largest ideal of definition $\Tscr$. Then $f^*(\Kscr)\Oscr_{\Xcal} \subseteq \Tscr$ for every ideal of definition $\Kscr$ of $\Ycal$.
\end{remark}

Recall (\cite{EGA} $0_I$, 7.7) that for a linearly topologized ring $R$ and linearly topologiczed topological $R$-modules $M$ one defines a topolpogy on $M \otimes_R N$ by declaring that the $\im(U \otimes_R N) + \im(M \times_R V)$ form a fundamental system of neighborhoods of zero, if $U \subseteq M$ and $V \subseteq N$ run through the set of open $R$-submodules. We denote by $M \hat\otimes_R N$ its completion.

If $M = A$ and $N = B$ are $R$-algebras, then $A \otimes_R B$ and $A \hat\otimes_R B$ are topological rings. If $A$ and $B$ are pre-admissible $R$-algebras, then $A \otimes_R B$ is pre-admissible and $A \hat\otimes_R B$ is admissible.

\begin{proposition}
Let $\Xcal = \Spf(A)$, $\Ycal = \Spf(B)$ affine formal schemes over an affine formal scheme $\Scal = \Spf(R)$. Setze $\Zcal := \Spf(A \hat\otimes_R B)$ and let $p\colon \Zcal \to \Xcal$ and $q\colon \Zcal \to \Ycal$ be the $\Scal$-morphisms corresponding to the canonical continuous $R$-algebra homomorphisms $A \to A \hat\otimes_R B$ and $B \to A \hat\otimes_R B$. Then $(\Zcal,p,q)$ is a fiber product of $\Xcal$ and $\Ycal$ over $\Scal$ in the category of formal schemes.
\end{proposition}

\begin{proposition}
In the category of formal schemes arbitrary fiber products exist.
\end{proposition}

\begin{definition}
Let $\Xcal$, $\Scal$ be locally noetherian schemes. A morphism $f\colon \Xcal \to \Scal$ is called \emph{adic} if there exists an ideal of definition $\Iscr$ of $\Oscr_{\Scal}$ such that $f^*(\Iscr)\Oscr_{\Xcal}$ is an ideal of definition. One also then says that $\Xcal$ is an \emph{adic $\Scal$-scheme}.
\end{definition}

\begin{remark}
Let $\Xcal$, $\Scal$ be locally noetherian formal schemes and let $f\colon \Xcal \to \Scal$ be an adic morphism. Then for every ideal of definition $\Iscr$ of $\Oscr_{\Scal}$ one has that $f^*(\Iscr)\Oscr_{\Xcal}$ is an ideal of definition.

This implies that for two adic $\Scal$-schemes $\Xcal$ and $\Ycal$ any morphism $f\colon \Xcal \to \Ycal$ is adic.
\end{remark}

Let $\Scal$ be a locally noetherian formal $\Scal$-scheme and let $\Iscr$ be an ideal of definition of $\Scal$. We set $S_n = (\Scal,\Oscr_{\Scal}/\Iscr^{n+1})$ (a locally noetherian scheme). An inductive system $(f_n\colon X_n \to S_n)_n$ of locally noetherian $S_n$-schemes is called an \emph{adic inductive $(S_n)$-system} if for all $m \leq n$ the diagram
\[\xymatrix{
X_m \ar[r] \ar[d]_{f_m} & X_n \ar[d]^{f_n} \\
S_m \ar[r] & S_n
}\]
is cartesian.

If $f\colon \Xcal \to \Scal$ is a formal locally noetherian adic $\Scal$-scheme, then $\Kscr := f^*(\Iscr)\Oscr_{\Xcal}$ is an ideal of definition and the inductive system $(X_n)$ with $X_n := (\Xcal,\Oscr_{\Xcal}/\Kscr^{n+1})$ is an adic inductive $(S_n)$-system.

\begin{theorem}\label{FormalAdicSystem}
The above construction yields an equivalence of the category of locally noetherian formal adic $\Scal$-schemes and the category of adic inductive $(S_n)_n$-systems.
\end{theorem}

In particular for all (locally noetherian) adic $\Scal$-schemes $\Xcal$ and $\Ycal$ the following canonical map is bijective
\[
\Hom_{\Scal}(\Xcal,\Ycal) \isom \limproj_n \Hom_{S_n}(X_n,Y_n).
\]

\begin{proof}
\cite{EGA} I, (10.12.3).
\end{proof}

\begin{example}
Let $R$ be a noetherian ring, $I \subseteq A$ an ideal, let $\Rhat$ be the $I$-adic completion of $R$ and set $S := \Spec R$, $S_n := \Spec R/I^{n+1}$. Let $f\colon X \to S$ be a locally noetherian $S$-scheme and set $X_n := S_n \times_S X$. Then $(X_n)_n$ is an adic inductive $(S_n)$-system corresponding under the equivalence of Theorem~\ref{FormalAdicSystem} to an adic $\Scal$-scheme $\Xcal$, where $\Scal = \Spf(\Rhat)$. Moreover, $\Xcal$ is the formal completion of $X$ along $f^{-1}(V(I))$.
\end{example}

\begin{definition}
Let $\Ycal$ be a locally noetherian formal scheme, $\Kscr$ an ideal of definition of $\Oscr_{\Ycal}$, let $f\colon \Xcal \to \Ycal$ be an adic morphism and set $\Iscr := f^*(\Kscr)\Oscr_{\Xcal}$ (an ideal of definition of $\Xcal$). Let $f_0\colon (\Xcal,\Oscr_{\Xcal}/\Iscr) \to (\Ycal,\Oscr_{\Ycal}/\Kscr)$ be the induced morphism of schemes.

Then $f$ is called \emph{(locally) of finite type}, if $f_0$ is (locally) of finite type.
\end{definition}

\begin{proposition}\label{DescribeFormalSchemeFinType}
Let $\Ycal$ be a locally noetherian formal scheme. Then a morphism of formal schemes $f\colon \Xcal \to \Ycal$ is locally of finite type (resp.~of finite type) if and only if every point of $\Ycal$ has an open affine neighborhood $V = \Spf A$ such that $f^{-1}(V)$ is the union of a family (resp.~of a finite family) of open affine formal subschemes $U_i = \Spf B_i$ such that $B_i$ is strictly topologically of finite type (i.e., there exists a surjective continuous open ring homomorphism $A \langle X_1,\dots,X_n\rangle \to B_i$).
\end{proposition}

\begin{proof}
\cite{EGA}~(I, 10.13.1).
\end{proof}

\begin{definition}\label{DefLocallyFormallyFiniteType}
Let $\Ycal$ be a locally noetherian formal scheme. Then a morphism of formal schemes $f\colon \Xcal \to \Ycal$ is called \emph{locally of formally finite type} if every point of $\Ycal$ has an open affine neighborhood $V = \Spf A$ such that $f^{-1}(V)$ is the union of a family of open affine formal subschemes $U_i = \Spf B_i$ such that there exists a surjective continuous open ring homomorphism $A \dlbrack T_1,\dots,T_n\drbrack \langle X_1,\dots,X_n\rangle \to B_i$, where $A \dlbrack T_1,\dots,T_n\drbrack$ is the noetherian adic ring whose ideal of definition is generated by an ideal of definition of $A$ and the $T_i$.
\end{definition}

\begin{example}
Let $V$ be a discrete valuation ring with uniformizing element $\pi$, let $X$ be a scheme locally of finite type over $V$, and let $Y$ be a closed subscheme of $X$ such that $\pi$ is locally nilpotent on $Y$ (e.g., if $Y$ is contained in the special fiber of $X$). Let $\Xcal$ be the formal completion of $X$ along $Y$ (which is the same as the formal completion of $X$ along the special fiber of $Y$). Then $\Xcal$ is locally of formally finite type over $\Spf V$. It is locally of finite type over $\Spf V$ if and only if the underlying topological space of $Y$ is equal to a connected component of the special fiber of $X$.
\end{example}

%---------------------------------------------------------

\subsection{The adic space attached to a formal scheme}

\begin{definition}
A complete affinoid ring $(A,A^+)$ is called \emph{sheafy} if $\Oscr_{\Spa(A,A^+)}$ is a sheaf of topological rings.

An adic ring $A$ is called \emph{sheafy}, if $A$ is complete, $f$-adic (i.e. it has a finitely generated ideal of definition) and if the affinoid ring $(A,A)$ is sheafy.

A formal scheme $\Xcal$ is called \emph{sheafy} if for every affine formal open subset $\Ucal = \Spf(A)$ the adic ring $A$ is sheafy.
\end{definition}

\begin{example}
A noetherian adic ring is sheafy. A ring with the discrete topology is sheafy.
\end{example}

If $Y$ is an adic space, the subsheaf $\Oscr^+_Y$ of $\Oscr_Y$ is considered as a sheaf of topological rings by endowing $\Oscr^+_Y(U)$ with the subspace topology of $\Oscr_Y(U)$ for $U \subseteq Y$ open. Then $(Y,\Oscr^+_Y)$ is a locally and topologically ringed space. Every morphism $f\colon Z \to Y$ of adic spaces induces a morphism $f^+\colon (Y,\Oscr^+_Y) \to (Z,\Oscr^+_Z)$ of locally and topologically ringed spaces (Proposition~\ref{MorphismAdic}).

\begin{proposition}\label{AdicFormalAffine}
Let $A$ be a sheafy adic ring and set $\Xcal := \Spf A$. Then for any adic space $Y$ the map
\begin{align*}
\Hom_{\Ltr}((Y,\Oscr_{Y}^+),(\Xcal,\Oscr_{\Xcal})) &\to \Hom(A,\Oscr_{Y}^+(Y)),\\
(f,f^{\flat}) &\sends f^{\flat}_{\Xcal}
\end{align*}
is bijective. Here the set on the left hand side denotes morphism of locally and topologically ringed spaces, and the set of the right hand side denotes continuous ring homomorphisms.
\end{proposition}

\begin{proof}
As morphisms of adic spaces can be glued, we may assume that $Y = \Spa(B,B^+)$, where $(B,B^+)$ is a complete affinoid ring. Let $\varphi\colon A \to \Oscr_Y^+(Y)$ be a continuous ring homomorphism. For every $y \in Y$ we set $g(y) := \set{a \in A}{v_y(\varphi(a)) < 1}$. Then $g(y)$ is an open prime ideal of $A$. This defines a map $g\colon Y \to \Xcal$. This map is continuous because for all $s \in A$ one has
\[
g^{-1}(\Dcal(s)) = g^1(\set{x \in \Xcal}{s \notin \pfr_x}) = \set{y \in Y}{v_y(\varphi(s)) \geq 1} = Y(\frac{1}{\varphi(s)}).
\]
Let $s \in A = \Oscr_{\Xcal}(\Xcal)$. Then $\Oscr_{\Xcal}(\Dcal(s)) = A \langle \frac{1}{s}\rangle$ and $\Oscr^+_Y(g^{-1}\Dcal(s)) = B^+ \langle \frac{1}{\varphi(s)}\rangle$. Thus by universal properties of the topological localization there exists a unique continuous ring homomorphism $\psi_U$ making the following diagram commutative.
\[\xymatrix{
\Oscr_{\Xcal}(\Xcal) \ar[r]^{\varphi} \ar[d] & \Oscr^+_Y(Y) \ar[d] \\
\Oscr_{\Xcal}(\Dcal(s)) \ar[r]^{\varphi_U} & \Oscr^+_Y(g^{-1}\Dcal(s))
}\]
The $\psi_U$ define a homomorphism of sheaves of topological rings $\psi\colon \Oscr_{\Xcal} \to g_*\Oscr^+_Y$. It remains to show that for all $y \in Y$ the induced homomorphism on stalks $\psi_y\colon \Oscr_{\Xcal,g(y)} \to \Oscr^+_{Y,y}$ is local. Let $\mfr_y \subset \Oscr^+_{Y,y}$ and $\mfr_{g(y)} \subset \Oscr_{\Xcal,g(y)}$ be the maximal ideals and let $\iota\colon A \to \Oscr_{\Xcal,g(y)}$ be the caonical ring homomorphism. Then by definition of $g$ one has $\iota^{-1}(\mfr_{g(y)}) = \iota^{-1}(\psi_y^{-1}(\mfr_y))$ which implies $\mfr_{g(y)} = \psi_y^{-1}(\mfr_y)$.
\end{proof}

\begin{theorem}
Let $\Xcal$ be a sheafy formal scheme. Then there exists an adic space $t(\Xcal)$ and a morphism of locally and topologically ringed spaces
\[
\pi = \pi_{\Xcal}\colon (t(\Xcal),\Oscr^+_{t(\Xcal)}) \to (\Xcal,\Oscr_{\Xcal})
\]
satisfying the following universal property. For every adic space $Z$ and for every morphism $f\colon (Z,\Oscr^+_Z) \to (\Xcal,\Oscr_{\Xcal})$ of locally and topologically ringed spaces there exists a unique morphism of adic spaces $g\colon Z \to t(\Xcal)$ making the following diagram commutative
\[\xymatrix{
(Z,\Oscr^+_Z) \ar[r]^f \ar@{.>}[rd]_{g^+} &  (\Xcal,\Oscr_{\Xcal}) \\
& (t(\Xcal),\Oscr^+_{t(\Xcal)}) \ar[u]^\pi
}\]
\end{theorem}

\begin{proof}
As we can glue morphisms of adic spaces, we may assume that $\Xcal = \Spf(A)$ is affine. We set $t(\Xcal) := \Spa(A,A)$. By Proposition~\ref{AdicFormalAffine}, the identity $A \to A$ corresponds to a morphism of locally and topologically ringed spaces $\pi\colon (t(\Xcal),\Oscr^+_{t(\Xcal)}) \to (\Xcal,\Oscr_{\Xcal})$.

For every adic space $Z$ we then have identifications
\begin{align*}
\Hom_{\Adic}(Z,t(\Xcal)) &\eqann{\ref{UniversalAffinoid}} \Hom_{\Affd}((A,A),(\Oscr_Z(Z),\Oscr^+_Z(Z))) \\
&\eqann{\phantom{\ref{UniversalAffinoid}}} \Hom_{(\textup{TopRing})}(A,\Oscr^+_Z(Z)) \\
&\eqann{\ref{AdicFormalAffine}} \Hom_{\Ltr}((Z,\Oscr^+_Z),(t(\Xcal),\Oscr^+_{t(\Xcal)})),
\end{align*}
and the composition is given by $f \sends \pi \circ f^+$.
\end{proof}

\begin{remark}
Let $f\colon \Xcal \to \Ycal$ be a morphism of sheafy formal schemes.  The universal property, applied to the composition
\[
(t(\Xcal),\Oscr^+_{t(\Xcal)}) \ltoover{\pi_{\Xcal}} (\Xcal,\Oscr_{\Xcal}) \ltoover{f} (\Ycal,\Oscr_{\Ycal}),
\]
yields a unique morphism $t(f)\colon t(\Xcal) \to t(\Ycal)$ of adic spaces such that the diagram
\[\xymatrix{
(t(\Xcal),\Oscr^+_{t(\Xcal)}) \ar[r]^{t(f)^+} \ar[d]_{\pi_{\Xcal}} & (t(\Ycal),\Oscr^+_{t(\Ycal)}) \ar[d]_{\pi_{\Ycal}} \\
(\Xcal,\Oscr_{\Xcal}) \ar[r]^f & (\Ycal,\Oscr_{\Ycal})
}\]
commutes. The uniqueness of $t(f)$ implies that the formation of $t(\cdot)$ is compatible with composition. Hence we obtain a functor
\[
t\colon \bigl(\text{sheafy formal schemes}\bigr) \to \Adic.
\]
\end{remark}

\begin{example}
Let $A$ be a complete valuation ring of height $1$ endowed with its valuation topology. Let $\mfr$ be its maximal ideal, $k := A/\mfr$, $K := \Frac A$.

Then $A$ is a $\gamma$-adic ring for all $0 \ne \gamma \in \mfr$. Set $\Xcal := \Spf(A) = \{\mfr\}$ and $\Oscr_{\Xcal}(\Xcal) = A$.

For $x \in X := \Spa (A,A)$ with $\supp v_x = \mfr$ one has $v_x(\abar) = 1$ for all $\abar \in k^{\times}$ (because $v(a) \leq 1$ for all $a \in A$). Thus the only point $x \in \Spa(A,A)$ with $\supp(v_x) = \mfr$ is the trivial valuation with support $\mfr$.

For $y \in \Spa(A,A)$ with $\supp v_y = (0)$ one obtains an induced valuation $v_y$ on $K = \Frac A$. As $v_y$ is continuous (and $(0)$ is not open in $A$), $v_y$ is not trivial. As $v_y(a) \leq 1$ for all $a \in A$, $A$ is contained in the valuation ring $A_y$ of $v_y$. As $A$ is of height $1$, it is maximal among all valuation rings $\ne K$ of $K$ and hence $A = A_y$. Thus the only point $y \in \Spa(A,A)$ with $\supp(v_y) = (0)$ is the valuation of $A$.

Clearly, $y$ is a horizontal generization of $x$. One has $\Oscr_X(X) = \Oscr^+_X(X) = A$, $\Oscr_X(\{y\}) = \Frac A$, and $\Oscr^+_X(\{y\}) = A$.
\end{example}

\begin{remark}\label{AdicAdicSpace}
Let $X$ be an adic space.
\begin{assertionlist}
\item\label{AdicAdicSpace1}
Set
\[
X_{\rm triv} := \sett{x \in X}{$v_x$ is trivial valuation}.
\]
The value group of $v_x$ is the same as the value group of the corresponding valuation on $A$, where $U = \Spa (A,A^+) \subseteq X$ is an open affinoid neighborhood of $x$. Hence
\[
X_{\rm triv} \cap U = U_{\rm triv} = \sett{v \in \Spa A}{$v$ trivial}.
\]
\item\label{AdicAdicSpace2}
Let $X = \Spa A$ be affinoid. A trivial valuation on $A$ is continuous if and only if its support is open. Hence the map $v \sends \supp v$ yields a map $X_{\rm triv} \sends \sett{\pfr \in \Spec A}{$\pfr$ open in $A$}$ which is a homeomorphism (Remark~\ref{SpvSpec}).
\item\label{AdicAdicSpace3}
Now assume that $X = \Spa (A,A)$, where $A$ is an adic ring and let $I$ be a finitely generated ideal of definition. Then supp yields a homeomorphism
\[
X_{\rm triv} \iso \Spf(A),
\]
and $\Spf A$ is a closed constructible subset of $\Spec A$ (because $\Spf A = V(I)$ and $I$ is finitely generated). The inclusion $X_{\rm triv} \mono X$ is spectral as the intersection of a rational subset $R(\frac{T}{s}) \subseteq X$ with $X_{\rm triv}$ is $\Dcal(s)$.
\item\label{AdicAdicSpace4}
For $X = \Spa (A,A)$, $A$ adic ring, one has $v_x(f) \leq 1$ for all $x \in X$. In particular one has for all $s \in A$ and $T \subseteq A$ finite with $T\cdot A \subseteq A$ open that
\[
R(\frac{1}{s}) = \set{v \in \Spa (A,A)}{v(s) = 1} \subseteq R(\frac{T}{s}).
\]
Moreover
\begin{equation}\label{EqOPlusR1}
\begin{aligned}
\Oscr^+_X(R(\frac{1}{s})) &= \sett{f \in A_{\{s\}}}{$v_x(f) \leq 1$ for all $x \in X$ with $1 \leq v_x(s)$} \\
&= A_{\{s\}} = \Oscr_X(R(\frac{1}{s})).
\end{aligned}
\end{equation}
\item\label{AdicAdicSpace5}
Let $A$ be a sheafy adic ring, $\Xcal := \Spf A$ and let $i\colon X_{\rm triv} \mono X = \Spa(A,A)$ be the inclusion. The restriction of the canonical morphism of locally and topologically ringed spaces $\pi\colon (X,\Oscr^+_X) \to (\Xcal,\Oscr_{\Xcal})$ to $X_{\rm triv}$ defines an isomorphism of locally ringed spaces
\begin{equation}\label{RetractAdicFormal}
\sigma_X\colon (X_{\rm triv}, i^{-1}\Oscr^+_X) \iso (\Xcal, \Oscr_{\Xcal}).
\end{equation}
Indeed, it is a homeomorphism by~\eqref{AdicAdicSpace2}. For all $x \in X_{\rm triv}$ the open rational neighborhoods of $x$ of the form $R(\frac{1}{s})$ for $s \notin \supp v_x$ are cofinal in the set of all open rational neighborhoods of $x$ by \eqref{AdicAdicSpace4}. Thus we have for the stalk in $x$:
\begin{align*}
(i^{-1}\Oscr^+_X)_x = \Oscr^+_{X,x} &\eqann{\phantom{\eqref{EqOPlusR1}}} \limind_{s \notin \supp v_x}\Oscr^+_X(R(\frac{1}{s})) \\
&\eqann{\eqref{EqOPlusR1}} \limind_{s \notin \supp v_x} A_{\{s\}} \\
&\eqann{\phantom{\eqref{EqOPlusR1}}} \Oscr_{\Xcal,\sigma_X(x)}.
\end{align*}
\item\label{AdicAdicSpace6}
The homomorphisms induced by $\sigma^{\flat}_X$ on rings of sections are bijective continuous ring homomorphisms. The morphism $\sigma_X$ is an isomorphism of locally and topologically ringed spaces if $A$ is noetherian by the following lemma applied to the isomorphisms of $A_{\{s\}}$-modules $\sigma^{\flat}_{\Dcal(s)}\colon \Oscr_{\Xcal}(\Dcal(s)) \iso \sigma_*(i^{-1}\Oscr^+_X)(\Dcal(s))$ for all $s \in A$.
\end{assertionlist}
\end{remark}

\begin{lemma}\label{NoetherianStrict}
Let $A$ be an adic noetherian ring, $I$ an ideal of definition. Let $M$, $N$ be $A$-modules and let $N$ be finitely generated. Then every $A$-linear map $M \to N$ is strict for the $I$-adic topologies.
\end{lemma}

We denote by $\lnSchf$ the category of locally noetherian formal schemes.

\begin{proposition}\label{FormalFullyFaithful}
The functor $t\colon \lnSchf \to \Adic$ is fully faithful.
\end{proposition}

\begin{proof}
For a locally noetherian formal scheme $\Zcal$ set $Z_t := \sett{z \in t(\Zcal)}{$v_z$ is trivial}$ and let $Z_0$ be the locally and topologically ringed space $(Z_t,\Oscr^+_{t(\Zcal)}\rstr{Z_t})$. By Remark~\ref{AdicAdicSpace}~\eqref{AdicAdicSpace6} the morphism $\sigma_Z\colon Z_0 \to \Zcal$ is an isomorphism of locally and topologically ringed spaces.

Let $\Xcal$ and $\Ycal$ be locally noetherian formal schemes and let $f\colon t(\Xcal) \to t(\Ycal)$ be a morphism of adic spaces. As $f$ is compatible with valuations, ona has $f(X_t) \subseteq Y_t$ and thus $f$ induces a morphism of locally and topologically ringed spaces $f_0\colon X_0 \to Y_0$ which we consider via $\sigma_X$ and $\sigma_Y$ as morphism $g\colon \Xcal \to \Ycal$ of formal schemes. It is easy to see that $t(g) = f$.

Let $g_1,g_2\colon \Xcal \to \Ycal$ be morphisms of sheafy formal schemes with $t(g_1) = t(g_2) =: f$. Then $g_i \circ \sigma_X = \sigma_Y \circ f_0$ for $i = 1,2$ and hence $g_1 = g_2$. This shows that $t$ is fully faithful.
\end{proof}

To determine the essential image of the functor $t$ we make the following definition.

\begin{definition}\label{DefSaturatedSubCat}
A full subcategory $\Ccal$ of the category $\Adic$ of adic spaces is called \emph{saturated} if it satisfies the following properties.
\begin{definitionlist}
\item
If an adic space $X$ is isomorphic to $Y \in \Ccal$, then $X \in \Ccal$.
\item
If $X \in \Ccal$ and $U \subseteq X$ is an open subspace, then $u \in \Ccal$.
\item
If $X$ is an adic space that has an open covering $(U_i)_{i\in I}$ with $U_i \in \Ccal$ for all $i \in I$, then $X \in \Ccal$.
\end{definitionlist}
\end{definition}

Clearly the intersection saturated subcategories is again a saturated subcategory.

\begin{proposition}
Let $\Ccal$ be the smallest saturated subcategory of $\Adic$ such that $t(\lnSchf) \subseteq \Ccal$. Then the objects of $\Ccal$ are those adic spaces $X$ such that every $x \in X$ has an open affinoid neighborhood $U$ such that the following conditions hold.
\begin{definitionlist}
\item
$\Oscr_X(U)$ has a noetherian ring of definition $A$.
\item
$\Oscr_X(U)$ is a finitely generated $A$-algebra.
\item
$\Oscr^+_X(U)$ is the integral closure of $A$ in $\Oscr_X(U)$.
\end{definitionlist}
\end{proposition}

\begin{proof}
\cite{Hu_Gen}~Prop.~4.2.
\end{proof}

\begin{proposition}
Let $f\colon \Xcal \to \Ycal$ be a morphism of locally noetherian formal schemes.
\begin{assertionlist}
\item
$f$ is adic if and only if $t(f)$ is adic.
\item
$f$ is locally of finite type if and only if $t(f)$ is locally of finite type.
\end{assertionlist}
\end{proposition}

\begin{proof}
\cite{Hu_Gen}~4.2.
\end{proof}

%================================================

\section{Rigid analytic spaces as adic spaces}

\begin{remark}
Let $k$ be a non-archimedean field, let $X$ be a scheme locally of finite type over $k$ and let $X^{\rm rig}$ be the attached rigid analyit variety. Then the adic space attached to $X^{\rm rig}$ is isomorphic to the adic space $X^{\rm ad}$ associated to the scheme $X$.
\end{remark}

\begin{remark}
Let $V$ be a complete discrete valuation ring (of height $1$) and let $k$ be the field of fractions of $V$. Let $\Fcal_f$ be the category of formal schemes locally of finite type over $\Spv V$. Let $\Rcal_k$ be the category of rigid analytic spaces over $k$. Finally let $\Ecal$ be the category of locally noetherian formal schemes whose morphisms are the adic morphisms between formal schemes.

Let $t_a\colon \Ecal \to \Adic$ be the functor $X \sends t(X)\an$ (note that because all morphisms are adic, the respect to open subset of analytic points). Let ${\rm rig}\colon \Fcal_f \to \Rcal_k$ be the generic fiber functor defined by Raynaud and let $u \colon \Fcal_f \to \Ecal$ be the inclusion functor. Then the diagram
\[\xymatrix{
\Fcal_f \ar[r]^{\rm rig} \ar[d]_u & \Rcal_k \ar[d]_{r_k} \\
\Ecal \ar[r]^{t_a} & \Adic
}\]
is 2-commutative.

The functor $t_a$ can be extended and described differently as follows (cf.~\cite{Hu_Etale}~1.2.2). The adic space $S := \Spa (V,V)$ attached to $\Spf V$ consists of an open and a closed point, and the canonical morphism $S^0 = \Spa (k,V) \to S$ is an open immersion onto the open point (in particular it is a morphism of adic spaces of finite type).

Let $\Fcal_{ff}$ be the category of formal schemes locally of formally finite type over $V$. Then Berthelot has extended Raynaud's generic fiber functor to a functor ${\rm rig}\colon \Fcal_{ff} \to \Rcal_k$. If $\Xcal \to \Spf V$ is an object in $\Fcal_{ff}$ we obtain an associated morphism of adic spaces $t(\Xcal) \to S$ and we can form the fiber product $t(\Xcal) \times_S S^0$. Then Berthelot's generic fiber functor and the functor $\Xcal \sends t(\Xcal) \times_S S^0$ are isomorphic.
\end{remark}

%================================================

\section{Perfectoid spaces}

In this section we explain a special class of adic spaces, namely the perfectoid spaces introduced by P.~Scholze. They allow to ``tilt'' from characteristic zero to characteristic $p$ and vice versa in presence of high ramification.

%================================================

\appendix

\section{\v{C}ech cohomology}\label{AppCech}

Let $X$ be a topological space. For a family of open subsets $\Ucal = (U_i)_{i\in I}$ we set
\[
U_{i_0\dots i_q} := U_{i_0} \cap \dots \cap U_{i_q}
\]
for all $(i_0,\dots,i_q) \in I^{q+1}$. If $Y$ is a subspace of $X$ we set $\Ucal\rstr{Y} := (U_i \cap Y)_{i\in I}$. If $\Vcal = (V_j)_{j\in J}$ is a second family of open subsets of $X$, we set
\[
\Ucal \times \Vcal := (U_i \cap V_j)_{i\in I,j\in J}.
\]
If $\Ucal$ and $\Vcal$ are open coverings of $X$, then $\Ucal \times \Vcal$ is an open covering of $X$.

An open covering $\Vcal = (V_j)_{j\in J}$ of $X$ is called \emph{refinement} of an open covering $\Ucal = (U_i)_{i\in I}$ of $X$ if there exists a map $\tau\colon J \to I$ such that $V_j \subseteq U_{\tau(j)}$ for all $j \in J$.

Now let $\Fscr$ be a presheaf of abelian groups on $X$. For an open covering $\Ucal$ of $X$ we define the abelian group of $q$-cochains
\[
\Cvee^q(\Ucal,\Fscr) := \prod_{(i_0,\dots,i_q) \in I^{q+1}}\Fscr(U_{i_0\dots i_q}).
\]
If $\Fscr$ is a presheaf of rings, then $\Cvee^q(\Ucal,\Fscr)$ is an $\Fscr(X)$-module. For $f \in \Cvee^q(\Ucal,\Fscr)$ we denote by $f_{i_0,\dots i_q}$ its $(i_0,\dots,i_q)$-component. We call $f$ \emph{alternating} if for all permutations $\pi$ of $\{0,1,\dots,q\}$ one has
\[
f_{i_{\pi(0)},\dots,i_{\pi(q)}} = (\sgn \pi)f_{i_0,\dots i_q}
\]
and if $f_{i_0,\dots i_q} = 0$ whenever the indices $i_0,\dots,i_q$ are not pairwise distinct. The alternating $q$-cohains form a submodule $\Cvee^q_a(\Ucal,\Fscr)$ of $\Cvee^q(\Ucal,\Fscr)$.

For $q \geq 0$ we define $d^q\colon \Cvee^q(\Ucal,\Fscr) \to \Cvee^{q+1}(\Ucal,\Fscr)$ by
\[
d^q(f)_{i_0\dots i_{q+1}} := \sum_{j=0}^{q+1}(-1)^jf_{i_0\dots\hat{i}_j\dots i_{q+1}}.
\]
Then $d^q \circ d^{q-1} = 0$ for all $q \geq 1$ and we obtain the \emph{\v{C}ech complex of cochains on $\Ucal$ with values in $\Fscr$}. The alternating cochains form a subcomplex $\Cvee^{\bullet}_a(\Ucal,\Fscr)$.

Recall that the inclusion $\Cvee_a^{\bullet}(\Ucal,\Fscr) \mono \Cvee^{\bullet}(\Ucal,\Fscr)$ is a quasi-isomorphism, i.e., it yields for all $q \geq 0$ an isomorphism
\[
\Hvee^q(\Ucal,\Fscr) := H^q(\Cvee_a^{\bullet}(\Ucal,\Fscr)) = H^q(\Cvee^{\bullet}(\Ucal,\Fscr)).
\]
One has a homomorphism
\[
\eps\colon \Fscr(X) \to \Cvee_a^0(\Ucal,\Fscr) = \Cvee^0(\Ucal,\Fscr), \qquad s \sends (s\rstr{U})_{U \in \Ucal}
\]
called the augmentation homomorphism.

Let $\Vcal = (V_j)_{j\in J}$ be a refinement of $\Ucal$ and let $\tau\colon J \to I$ be a map such that $V_j \subseteq U_{\tau(j)}$ for all $j \in J$. Then $\tau$ induces homomorphisms of complexes $\tau^{\bullet}\colon \Cvee^{\bullet}(\Ucal,\Fscr) \to \Cvee^{\bullet}(\Vcal,\Fscr)$ and $\tau^{\bullet}\colon \Cvee_a^{\bullet}(\Ucal,\Fscr) \to \Cvee_a^{\bullet}(\Vcal,\Fscr)$, and the induced map
\[
H^q(\tau^{\bullet})\colon \Hvee^q(\Ucal,\Fscr) \to \Hvee^q(\Vcal,\Fscr)
\]
is independent of the choice of $\tau$. Thus one may define the \v{C}ech cohomology on $X$ with values in $\Fscr$ by
\[
\Hvee^q(X,\Fscr) := \limind_{\Ucal} \Hvee^q(\Ucal,\Fscr),
\]
where $\Ucal$ runs through the set of open coverings of $X$, preordered by refinement.

\begin{defi}
Let $\Ucal$ be an open covering of $X$. Then the covering $\Ucal$ is called \emph{$\Fscr$-acyclic} if the the augmented \v Cech complex
\[
0 \to \Fscr(U) \to \Cvee^0(\Ucal,\Fscr) \to \Cvee^1(\Ucal,\Fscr) \to \Cvee^2(\Ucal,\Fscr) \to \dots
\]
is exact. In other words, the augmentation morphism $\eps$ yields an isomorpism $\Fscr(X) \iso \Hvee^0(\Ucal,\Fscr)$ and $\Hvee^q(\Ucal,\Fscr) = 0$ for $q \geq 1$.
\end{defi}

The trivial open covering $\Ucal_0$ consisting of $X$ only is always $\Fscr$-acyclic. An arbitrary open covering $\Ucal$ is always a refinement of $\Ucal_0$ (in a unique way) and $\Ucal$ is $\Fscr$-acyclic if and only if $\Cvee^{\bullet}(\Uscr_0,\Fscr) \to \Cvee^{\bullet}(\Uscr,\Fscr)$ is a quasi-isomorphism. This implies in particular the following remark.

\begin{rem}\label{AcyclicMutuallyRefinement}
Let $\Ucal$ and $\Vcal$ be open coverings that are refinements of each other. Then $\Ucal$ is $\Fscr$-acyclic if and only if $\Vcal$ is $\Fscr$-acyclic.
\end{rem}

The preceeding remark show that every open covering $(U_i)_{i}$ of $X$ such that there exists $i_0 \in I$ with $U_{i_0} = X$ is $\Fscr$-acyclic.

\begin{prop}\label{AcyclicProduct}
Let $\Fscr$ be a presheaf of abelian groups and let $\Ucal = (U_i)_{i\in I}$ and $\Vcal = (V_j)_{j\in J}$ be open coverings of $X$ such that $\Vcal\rstr{U_{i_0\dots i_q}}$ is $\Fscr$-acyclic (or, more precisely, $\Fscr\rstr{U_{i_0\dots i_q}}$-acyclic) for all $(i_0,\dots,i_q) \in I^{q+1}$ and for all $q \geq 0$.
\begin{assertionlist}
\item
Assume that also $\Ucal\rstr{V_{j_0\dots j_q}}$ is $\Fscr$-acyclic for all $(j_0,\dots,j_q) \in J^{q+1}$ and all $q \geq 0$. Then $\Ucal$ is $\Fscr$-acycic if and only if $\Vcal$ is $\Fscr$-acyclic.
\item
If $\Vcal$ is a refinement of $\Ucal$, then $\Ucal$ is $\Fscr$-acyclic if and only if $\Vcal$ is $\Fscr$-acyclic.
\item
Then $\Ucal \times \Vcal$ is $\Fscr$-acyclic if and only if $\Ucal$ is $\Fscr$-acyclic.
\end{assertionlist}
\end{prop}

\begin{proof}
Under the hypotheses in Assertion~(1) one knows more generally that $\Hvee^q(\Ucal,\Fscr) = \Hvee^q(\Vcal,\Fscr)$ for all $q \geq 0$. In particular~(1) holds.

If $\Vcal$ is a refinement of $\Ucal$, then $\Ucal\rstr{V_{j_0\dots j_q}}$ and the trivial covering of $V_{j_0\dots j_q}$ are refinements of each other. This shows~(2).

Let us show~(3). The coverings $\Vcal\rstr{U_{i_0\dots i_q}}$ and $(\Ucal \times \Vcal)\rstr{U_{i_0\dots i_q}}$ are refinements of each other. Thus $(\Ucal \times \Vcal)\rstr{U_{i_0\dots i_q}}$ is $\Fscr$-acyclic. Thus we may apply~(2) to $\Ucal$ and its refinement $\Ucal \times \Vcal$.
\end{proof}

\begin{prop}\label{CechZero}
Let $\Bscr$ be a basis of the topology of $X$ that is stable under finite intersections. Let $\Fscr'$ be a presheaf of abelian groups on $\Bscr$ and define a presheaf $\Fscr$ on $X$ by setting for $V \subseteq X$ open
\[
\Fscr(V) = \limproj_{U \subseteq V, U \in \Bscr} \Fscr(U).
\]
Assume that for every $U \in \Bscr$ and for every open covering $\Ucal = (U_i)_{i\in I}$ of $U$ by open subsets $U_i \in \Bscr$ the presheaf $\Fscr$ is $\Ucal$-acyclic, i.e., the augmented \v Cech complex
\[
0 \to \Fscr(U) \to \Cvee^0(\Ucal,\Fscr) \to \Cvee^1(\Ucal,\Fscr) \to \Cvee^2(\Ucal,\Fscr) \to \dots
\]
is exact. Then $\Fscr$ is a sheaf on $X$ and for all $q \geq 0$ and for every open subset $U$ of $X$ the canonical homomorphisms $\Hvee^q(U,\Fscr) \to H^q(U,\Fscr)$ is an isomorphism. In particular $H^q(U,\Fscr) = 0$ for all $U \in \Bscr$ and $q \geq 1$.
\end{prop}

\begin{proof}
The exactness of $0 \to \Fscr(U) \to \Cvee^0(\Ucal,\Fscr) \to \Cvee^1(\Ucal,\Fscr)$ shows that $\Fscr'$ is a sheaf on $\Bscr$ and thus $\Fscr$ is a sheaf on $X$. The exactness of the augmented \v Cech complex then shows that $\Hvee^q(\Uscr,\Fscr) = 0$ for all $q \geq 1$ and for every open covering $\Uscr$ of $U \in \Bscr$ by open subsets in $\Uscr$. This shows $\Hvee^q(U,\Fscr) = 0$ for all $q \geq 1$ and all $U \in \Bscr$ and thus a result of Cartan (e.g.~\cite{Godement}~II~5.9.2) shows that $\Hvee^q(U,\Fscr) \to H^q(U,\Fscr)$ is an isomorphism for every open subset $U \subseteq X$ and all $q \geq 0$.
\end{proof}

%==================================================================


\begin{thebibliography}{BouLie}
\setlength{\itemsep}{\itemsepamount}
\bibitem[BGR]{BGR} S. Bosch, U. G\"untzer, R. Remmert, \emph{Non-Archimedean Analysis}, Springer (1984).
\bibitem[BouA]{Bou_Alg} N. Bourbaki, {\em Algebra}, chap. I-III, Springer (1989);
chap. IV-VII, Springer (1988); chap. VIII, Hermann (1958); chap. IX,
Hermann (?); chap. X, Masson (1980).
\bibitem[BouTG]{Bou_TG} N. Bourbaki, {\em Topologie g\'en\'erale}, chap. I - IV, Springer (2007); chap. V - X, 
Springer (2007).
\bibitem[BouAC]{Bou_AC} N. Bourbaki, {\em Alg\`ebre commutative}, chap. I - VII,
Springer (1989); chap. VIII-IX, Masson (1983); chap. X, Masson (1998).
\bibitem[EGA]{EGA} A.~Grothendieck, J.~Dieudonn\'e: {\em El\'ements de G\'eom\'etrie Alg\'ebrique}, I Grundlehren 
der Mathematik {\bf 166} (1971) Springer,
II-IV Publ.~Math.~IHES {\bf 8} (1961), {\bf 11} (1961), {\bf 17} (1963), {\bf 20} (1964), {\bf 24} (1965), {\bf 28} (1966), 
{\bf 32} (1967).
\bibitem[FuSa]{FuchsSalce} L.~Fuchs, L.~Salce, {\em Modules over non-Notherian domains}, Mathematical 
Surveys and Monographs {\bf84}, AMS (2001).
\bibitem[God]{Godement} R.~Godement, {\em Topologie alg\'ebrique et th\'eorie des faisceaux}, Hermann (1958).
\bibitem[Ho]{Ho_PrimeSpec} M.~Hochster, {\em Prime ideal structure in commutative rings}, Trans. Amer. Math. 
Soc. {\bf142} (1969), 43--60.
\bibitem[Hu1]{Hu_Habil} R.~Huber, {\em Bewertungsspektrum und rigide Geometrie}, Regensburger Math. Schriften 
{\bf23} (1993).
\bibitem[Hu2]{Hu_Cont} R.~Huber, {\em Continuous valuations}, Math. Zeitschrift {\bf212} (1993), 445--477.
\bibitem[Hu3]{Hu_Gen} R.~Huber, {\em A generalization of formal schemes and rigid analytic varieties}, Math. 
Zeitschrift {\bf217} (1994), 513--551.
\bibitem[Hu4]{Hu_Etale} R.~Huber, {\em \'Etale Cohomology of Rigid Analytic Varieties and Adic Spaces}, Aspects 
of Mathematics {\bf E30}, Vieweg (1996).
\bibitem[HuKn]{HuKne} R.~Huber, M.~Knebusch, {\em On valuation spectra}, Contemporary mathematics {\bf155} 
(1994), 167--206.
\bibitem[Sch]{Scholze_Perfect} P.~Scholze, \emph{Perfectoid Spaces}, arXiv:1111.4914v1.
\bibitem[Str]{Str_Adic} M.~Strauch, {\em Introductory lectures on adic spaces}, Lecture script (2005).
\end{thebibliography}
\end{document}